\documentclass[11pt,letterpaper]{amsart}
\usepackage{amsmath,amssymb,amsthm}
\usepackage{graphicx}
\usepackage{placeins}
\usepackage[percent]{overpic}
\newcommand{\R}{\mathbb{R}}
\newcommand{\Z}{\mathbb{Z}}
\newcommand{\Q}{\mathbb{Q}}
\newcommand{\C}{\mathbb{C}}

\newcommand{\ZN}[1]{\Z_{#1}}
\newcommand{\ex}[1]{e^{2\pi i #1}}

\newcommand{\cyc}[1]{\Phi_{#1}}
\newcommand{\zroot}[1]{\zeta_{#1}}
\newcommand{\muS}{\mu_S}
\newcommand{\mask}[1]{#1(X)}

\newcommand{\nset}{[n]}
\newcommand{\one}{\mathbf{1}}
\newcommand{\slice}[1]{\mathcal{T}_{#1}}
\newcommand{\Rat}{R(\ell)}
\newcommand{\Den}{D(\ell)}
\newcommand{\kzero}{k_0}
\newcommand{\Kn}{K(n)}
\newcommand{\Krel}{K_\ell}
\newcommand{\dimQ}{\dim_{\Q}}
\newcommand{\clsgraph}{G}
\newcommand{\nloop}{n_{\mathrm{loop}}}
\newcommand{\nfwd}{n_{\mathrm{fwd}}}
\newcommand{\nback}{n_{\mathrm{back}}}
\newcommand{\Out}[1]{\mathrm{Out}_{#1}}
\newcommand{\In}[1]{\mathrm{In}_{#1}}
\newcommand{\ev}[1]{e_{#1}}
\newcommand{\Dop}{\Delta}
\newcommand{\zeroset}[1]{Z(#1)}
\newcommand{\Tsub}[1]{T_{#1}}
\newcommand{\FL}[1]{(\mathrm{F}_{#1})}
\newcommand{\UC}[1]{(\mathrm{U}_{#1})}
\newcommand{\CP}[1]{(\mathrm{CP}_{#1})}
\newcommand{\cycset}[1]{E(#1)}
\newcommand{\Ssev}{S'_{72}}
\newcommand{\OmegaEight}{\Omega_8}
\newcommand{\Rp}[1]{R_{#1}}

\makeatletter
\def\subsection{\@startsection{subsection}{2}\z@{.7\linespacing\@plus.7\linespacing}{.3\linespacing}{\normalfont\bfseries}}
\makeatother
\theoremstyle{plain}
\newtheorem{theorem}{Theorem}[section]
\newtheorem{proposition}[theorem]{Proposition}
\newtheorem{lemma}[theorem]{Lemma}
\newtheorem{corollary}[theorem]{Corollary}
\newcommand{\mainletter}{}
\newtheorem*{maintheoreminner}{Theorem~\mainletter}
\makeatletter
\newenvironment{maintheorem}[1]{\renewcommand{\mainletter}{#1}\begin{maintheoreminner}\phantomsection\def\@currentlabel{#1}}{\end{maintheoreminner}}
\newtheorem*{maincorollaryinner}{Corollary~\mainletter}
\newenvironment{maincorollary}[1]{\renewcommand{\mainletter}{#1}\begin{maincorollaryinner}\phantomsection\def\@currentlabel{#1}}{\end{maincorollaryinner}}
\makeatother
\theoremstyle{definition}
\newtheorem{definition}[theorem]{Definition}
\theoremstyle{remark}
\newtheorem{remark}[theorem]{Remark}
\usepackage[hidelinks,pdftitle={Spectral unions of intervals with few rational relations tile the real line},pdfauthor={Alexander Okezue Bell}]{hyperref}
\title[Spectral unions of intervals with few rational relations]{Spectral unions of intervals with few rational relations tile the real line}
\author{Alexander Okezue Bell}
\thanks{Department of Mathematics, Stanford University, 450 Jane Stanford Way, Building 380, Stanford, CA 94305-2125, USA. Corresponding author: Alexander Okezue Bell. Email: \href{mailto:okezue@stanford.edu}{okezue@stanford.edu}.}
\begin{document}
\begin{abstract}
We prove that every spectral finite union of pairwise disjoint bounded intervals whose lengths satisfy at most four linearly independent rational relations tiles the real line, and that each of its spectra has at most five cosets modulo its least period. This includes every spectral union of at most five bounded intervals. The proof combines specialization preserving subset-sum equalities with rigidity of exchanges between finite spectral fibers. We also obtain common tiling complements for all spectra of primitive sets of at most seven elements in finite cyclic groups, and establish the sharp eleven-point range in which consecutive pairs force periodicity of the dual set.
\end{abstract}
\maketitle
\section{Introduction}
\label{sec:intro}

Fuglede proved that a bounded measurable subset of Euclidean space tiles by a lattice iff it is spectral with the dual lattice as a spectrum, and he asked whether spectrality and tiling by translations remain equivalent without the lattice restriction \cite{fuglede1974}. On the line the question is open in general. It is settled for a union of two intervals by \L{}aba \cite{laba2001}, for a union of three intervals by Bose, Anil Kumar, Krishnan and Madan in the direction from tiling to spectrality \cite{bose-etal2010} and by Shi in the direction from spectrality to tiling \cite[Theorem~1.1]{shi2026}, and in the direction from spectrality to tiling for a union of four intervals by Bell \cite{bell2026four}, whose methods are used throughout this paper.

\begin{theorem}[{\cite[Theorem~1.1]{bell2026four}}]
\label{thm:four}
Let \(\Omega \subset \R\) be a bounded measurable set which, up to a null set, is a union of four bounded intervals. If \(\Omega\) is spectral, then \(\Omega\) tiles \(\R\) by translations.
\end{theorem}

A set is spectral when the normalized exponentials at the frequencies of some spectrum form an orthonormal basis of \(L^2(\Omega)\), and it tiles by translations when countably many translates partition the line up to null sets; the definitions are recalled in Section~\ref{sec:fibers}.

The extension to arbitrarily many intervals rests on the arithmetic of their lengths. Equalities between subset sums of the lengths constrain the possible periods of a spectrum. With at most \(d\) independent rational relations among the lengths, the length vector admits a nonzero, nonnegative specialization with at most \(d+1\) positive coordinates that preserves every original subset-sum equality. Taking \(d=4\) thus brings the arithmetic of five lengths into the problem, regardless of the number of intervals. Together with the geometric arguments below, this gives the following theorem.

\begin{maintheorem}{A}
\label{thm:A}
Let \(\Omega \subset \R\) be a bounded measurable set which, up to a null set, is a finite union of pairwise disjoint bounded intervals with lengths \(\ell_1, \dots, \ell_n\) satisfying at most four linearly independent rational relations. If \(\Omega\) is spectral, then \(\Omega\) tiles \(\R\) by translations. Moreover, every spectrum of \(\Omega\) is a union of at most five cosets of its minimal period.
\end{maintheorem}

\begin{maincorollary}{B}
\label{cor:B}
Let \(\Omega \subset \R\) be a bounded measurable set which, up to a null set, is a union of at most five bounded intervals. If \(\Omega\) is spectral, then \(\Omega\) tiles \(\R\) by translations.
\end{maincorollary}

For the length vector \(\ell=(\ell_1,\dots,\ell_n)\), write
\[
\Krel = \Bigl\{ r \in \Q^n : \sum_{j=1}^{n} r_j \ell_j = 0 \Bigr\}
\]
for the space of rational linear relations. Theorem~A's hypothesis is therefore \(\dimQ\Krel\le4\). This dimension is \(n-\rho\), where \(\rho\) is the dimension of the rational span of the lengths; it is \(n-1\) when the lengths are commensurable and \(0\) when they are linearly independent over \(\Q\). In particular, \(\dimQ\Krel\le n-1\), so Corollary~B is the case \(n\le5\) of Theorem~A. Both are proved in Section~\ref{sec:assembly}. The bound of five cosets is sharp: it is attained for six, seven and eight intervals (Proposition~\ref{prop:sharp}). Three further results come out of the same method and are stated at the end of this introduction: the fiber lemmas for six and seven cosets (Theorem~C), the structure of finite spectral pairs whose spectrum contains two consecutive integers, with its failure at twelve points (Theorem~D), and the reduction of the general case to statements about finite spectral pairs of finitely many possible cardinalities (Theorem~E).

\begin{figure}[ht]
\centering
\input{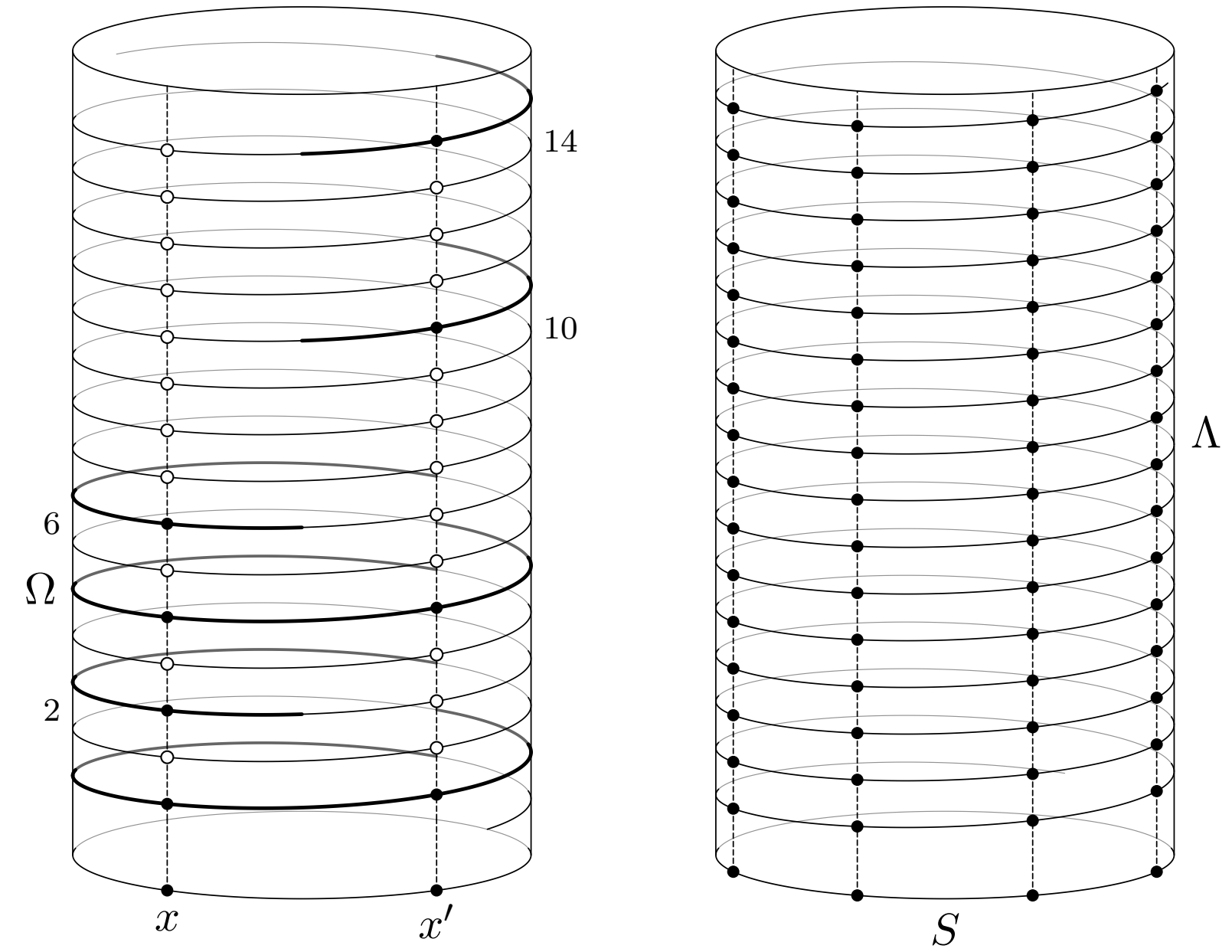}
\caption{The union \(\Omega = [0,1) \cup [2, 2.6) \cup [4,5) \cup [6, 6.6) \cup [10.6, 11) \cup [14.6, 15)\) of six intervals, drawn heavy on the line wound on a cylinder of circumference \(h = 1\), with the parts behind the cylinder in gray. The dashed lines over \(x\) and \(x'\) meet the winding at the points \(x + m\) and \(x' + m\), \(m \in \Z\), drawn as dots, filled when the point lies in \(\Omega\); the filled dots form the fibers \(A_x = \{0, 2, 4, 6\}\) and \(A_{x'} = \{0, 4, 10, 14\}\). The two points are separated by the endpoint class of size two at \(0.6\): the fiber loses the heights \(2\) and \(6\), labeled at the left edge of the cylinder, and gains \(10\) and \(14\), labeled at its right edge, and the character sums of the two pairs agree on \(S - S\). On the right, the spectrum \(\Lambda = \tfrac18\{0,1,2,3\} + \Z\) is wound on a cylinder of circumference \(h^{-1}\); it lies on \(k = 4\) vertical lines over the points of \(S\).}
\label{fig:fibers}
\end{figure}

We outline the proof of Theorem~A. A spectrum \(\Lambda\) of a bounded measurable set on the line is periodic, by the theorem of Iosevich and Kolountzakis \cite{ik2013}. Let \(h^{-1}\) be its least period, so that \(\Lambda = \Gamma + h^{-1}\Z\) for a finite set \(\Gamma\) of \(k\) frequencies, and let \(S = h\Gamma\), a set of \(k\) points of the circle \(\R/\Z\). Cutting the line along the period turns the continuous problem into finitely many finite ones: for each point \(x\), the integers \(m\) with \(x + mh \in \Omega\) form the fiber \(A_x\), and \(\Lambda\) is a spectrum of \(\Omega\) exactly when almost every fiber is a spectrum of the one fixed set \(S\), in the sense that \(\lvert A_x\rvert = k\) and the character sums of \(A_x\) vanish on \(S - S\) away from \(0\). Then \(\lvert\Omega\rvert = kh\), so \(\Omega\) covers almost every point of the line exactly \(k\) times by the translates \(\Omega + mh\). For a union of intervals the fiber is a union of at most \(n\) runs of consecutive integers, one for each interval, and it changes only when \(x\) crosses an endpoint of an interval modulo \(h\): at such a class of endpoints, the fiber loses the last point of each run that ends there and gains the first point of each run that begins there, and the two sets of heights exchanged have equal character sums on \(S - S\). This conservation law, which is proved in Section~\ref{sec:fibers} together with its converse, is the whole content of spectrality once the period is fixed. Figure~\ref{fig:fibers} shows a union of six intervals, its spectrum, and two fibers on either side of a class where two intervals change at once.

The lengths of the intervals control \(k\) through the arithmetic of their subset sums. Shi's boundary unitary and its secular determinant \cite{shi2026}, which are available for any number of intervals, give an exponential polynomial in the frequency whose exponents are the subset sums \(\ell_T\), \(T \subset \nset\), and whose zero set is the spectrum. A cofactor identity from \cite[Theorem~4.1]{bell2026four}, whose proof applies to any number of intervals, says that at each value \(\alpha\) of a subset sum the sum of the subset indicators \(1_T\), weighted by their coefficients over \(\ell_T = \alpha\), is proportional to the all-ones vector. An interior level \(0<\alpha<\lvert\Omega\rvert\) can therefore carry a nonzero coefficient only when it is effective, meaning that the indicators of the subsets at that level span the all-ones vector over \(\Q\), and then the ratio \(\alpha/\lvert\Omega\rvert\) is a rational number determined by those subsets alone, the pinned ratio of the level. The periodicity of the spectrum forces every exponent that survives into a lattice, and the upshot is the verdict rule of Section~\ref{sec:levels}: \(k\) divides the integer \(\kzero(\ell)\), the least common multiple of the denominators of the pinned ratios. For four lengths the effective levels are classified in \cite[Theorem~4.4]{bell2026four}. For five lengths Section~\ref{sec:fivelengths} does the same, and finds that every pinned ratio has denominator in \(\{2, 3, 4, 5, 7, 8, 9\}\) and that \(\kzero(\ell) \le 10\).

Specialization extends the five-length bound to every \(n\) with \(\dimQ\Krel\le4\). A coincidence \(\ell_T = \ell_{T'}\) is a rational relation, and the relations form the space \(\Krel\); moving the length vector inside the linear space of vectors satisfying all the original relations keeps every coincidence, and an extreme ray of the cone of nonnegative such vectors has at most \(\dimQ \Krel + 1\) positive coordinates. The specialization lemma of Section~\ref{sec:fivelengths} shows that every pinned ratio of \(\ell\) is a pinned ratio of the specialized vector, so \(\kzero(\ell) \le 10\).

The second ingredient is the consecutive pair. When a fiber has more points than there are intervals, two of its points are consecutive integers, since the fiber is a union of at most \(n\) runs. A finite spectral pair \((A, S)\) with \(A \subset \Z\) has rational differences in \(S\), by the theorem of Fu and Song \cite{fu-song2026}, so \(S\) is, after translation and scaling, a set of residues modulo some \(N\), and the consecutive pair in \(A\) says that the \(k\) corresponding roots of unity sum to zero. Section~\ref{sec:consecutive} proves that for \(k \le 10\) this forces \(S\) to have a translation symmetry \(S + 1/m = S\), through the classification of vanishing sums of at most ten roots of unity and the orthogonality of the other elements of \(A\), and Section~\ref{sec:twelve} extends this to \(k = 11\). But \(S = h\Gamma\) has no translation symmetry at all, because \(h^{-1}\) is the least period of \(\Lambda\). So \(k \le n\) whenever \(k \le 11\), and in particular whenever \(\kzero(\ell) \le 10\).

The third ingredient is merging. A pair of endpoints \(b_i \equiv a_j \pmod{h}\) is rigid when the slide \(d = (b_i - a_j)/h\) satisfies \(d(S - S) \subset \Z\). Sliding the interval \([a_j, b_j)\) by \(dh\) then preserves the spectrum, and the choice of the slide makes the moved interval land exactly against \([a_i, b_i)\), while a dimension count excludes overlap with every other interval. Thus the two intervals merge into one. A class of endpoints of size one is always rigid. Iterating, every spectral union of intervals merges, without changing \(\Lambda\) or \(k\), either to a single interval, where \(k = 1\), or to a fully mixing configuration, in which no pair is rigid and every endpoint class has size at least two. The number of classes is then at most \(n/2\), and since the lengths lie in the rational span of \(h\) and the class differences \(c_i-c_1\), \(\dimQ \Krel \ge n/2\). The bound \(\dimQ\Krel\le4\) therefore leaves at most eight intervals in a fully mixing configuration, and its class pattern is one of a short list (Section~\ref{sec:classgraph}).

For the resulting fully mixing configurations of six, seven or eight intervals, the bounds above give \(k \le n - 1 \le 7\), and what is special about \(n \ge 6\) is that several intervals may change their runs at the same class, so that the fiber changes by an exchange of two or more points at once. Section~\ref{sec:mechanisms} analyzes these exchanges on the fiber side, in the space \(\C^S\) of functions on \(S\), where a spectrum of \(S\) is an orthogonal basis of exponentials. A loop at a class of size two exchanges two exponentials for two others in a two-dimensional space and forces a consecutive pair; a class moving fewer than \(k\) points cannot change a spectrum of \(S\) when the spectra of \(S\) are cosets, which is the case for three and five points; for four points the two-coset structure of the spectra propagates along the classes and again needs a class of size four; and around a closed walk of classes the exchange spaces return to themselves twisted by a power of the diagonal unitary of \(S\), which bounds \(k\) by the number of backward steps. These mechanisms give \(k \le 2, 3, 4\) at \(n = 8, 7, 6\). In every case \(k \le 5\), and for \(k \le 5\) \cite[Theorem~5.8]{bell2026four} gives one translation set, depending only on \(S\), that tiles with every fiber at once, so \(\Omega\) tiles \(\R\).

The proof of Theorem~A therefore uses the consecutive-pair statement only up to ten points and the fiber lemma only up to five cosets. The results stated next, the consecutive-pair statement at eleven points and its failure at twelve, the description of the failures by gluing, and the fiber lemmas for six and seven cosets, do not enter it; they are proved after it, in Sections~\ref{sec:twelve} to~\ref{sec:general}.

We turn to the three further results. The fiber lemma is the statement that a bounded measurable set covering the line \(k\) times by the translates \(F + mh\) and having a spectrum made of \(k\) cosets of \(h^{-1}\Z\) tiles \(\R\) by translations. For \(k\le5\), it was proved by Dutkay and Haussermann \cite{dutkay-haussermann2015}, whose argument constructs a common complement for the finite fibers; we use the formulation in \cite[Theorem~5.8]{bell2026four}. In finite terms it is a common tiling complement for all spectra of a \(k\)-point set of residues, and Section~\ref{sec:fiberlemmas} proves it through a six-point dichotomy and seven-point rigidity.

\begin{maintheorem}{C}
\label{thm:C}
Let \(1 \le k \le 7\) be an integer, let \(h > 0\), and let \(F \subset \R\) be a bounded measurable set with \(\sum_{m \in \Z} 1_F(x + mh) = k\) for almost every \(x\) and with a spectrum \(\Gamma + h^{-1}\Z\), where \(\lvert\Gamma\rvert = k\) and the elements of \(\Gamma\) are distinct modulo \(h^{-1}\Z\). Then \(F\) tiles \(\R\) by translations, by a translation set \(hB\) with \(B \subset \Z\) depending only on \(h\Gamma\). Moreover, for \(k = 6\), every spectrum \(A\) of a primitive six-point set \(S' \subset \ZN{N}\) satisfies \(A + N/2 = A\) or \(A + N/3 = A\); and for \(k = 7\), every primitive seven-point spectral set \(S' \subset \ZN{N}\) is a complete system of residues modulo \(7\), and each of its spectra is a coset of the subgroup of order \(7\).
\end{maintheorem}

The consecutive-pair statement is what makes the count \(k \le n\) work, and it has a precise range of validity.

\begin{maintheorem}{D}
\label{thm:D}
\leavevmode
\begin{enumerate}
\item[\rm (i)] Let \((A, S)\) be a spectral pair with \(A \subset \Z\), \(S \subset \R/\Z\) and \(\lvert A\rvert = \lvert S\rvert = k \le 11\), and suppose that \(A\) contains two consecutive integers. Then \(S + 1/m = S\) for some integer \(m \ge 2\). For \(k \in \{7, 11\}\), \(S\) is a coset of \((1/k)\Z/\Z\).
\item[\rm (ii)] The statement fails for \(k = 12\). The set
\[
\Ssev = \{5, 41, 11, 47, 35, 71\} \cup \{6, 30, 54\} \cup \{10, 34, 58\} \subset \ZN{72}
\]
has no period, and \(A' = \{0, 1, 8, 9, 17, 28\} + \{0, 36\}\) is a spectrum of it.
\item[\rm (iii)] Every spectral pair \((A, S')\) in \(\ZN{N}\), with \(\lvert A\rvert=\lvert S'\rvert=k\), whose spectrum \(A\) is periodic, with stabilizer of order \(m\), arises by gluing \(m\) spectral pairs of size \(k/m\) along a common spectrum, the pieces of \(S'\) being its intersections with the cosets of \(m\ZN{N}\). The glued set \(S'\) is periodic exactly when, for some prime \(p \mid N\), translation by \(N/p\) carries every piece onto a piece: onto itself when \(p \mid N/m\), so that the pieces share the period \(N/p\), and onto another piece, in orbits of length \(p\), when \(p \nmid N/m\).
\item[\rm (iv)] Taking the displayed representatives of \(A'\) in \(\{0,\dots,71\}\), the set \(\OmegaEight = A' + [0, 1)\) is a spectral union of eight intervals whose spectrum has \(k = 12\) cosets, and \(\OmegaEight\) tiles \(\R\) by translations.
\end{enumerate}
\end{maintheorem}

Part (i) is proved for \(k \le 10\) in Section~\ref{sec:consecutive} and for \(k = 11\) in Section~\ref{sec:twelve}, which also proves parts (ii) and (iii); part (iv) is proved in Section~\ref{sec:general}. The example in (iv) shows that the inequality \(k \le n\), which holds for \(n \le 5\) and is the hinge of the proof of Theorem~A, fails for eight intervals, and Corollary~\ref{cor:alllargern} extends this failure to every number \(n \ge 8\) of intervals.

For arbitrary lengths, the method still reduces the problem to statements about finite spectral pairs of bounded cardinality. The ambient cyclic-group order is not bounded by this reduction. Section~\ref{sec:levels} proves that the denominator of a pinned ratio of \(n\) lengths divides a nonzero determinant of a \(0\)-\(1\) matrix of order at most \(n + 1\), which Hadamard's inequality bounds, and Section~\ref{sec:general} draws the consequence. Write \(\FL{k}\) for the fiber lemma for \(k\) cosets, as in Theorem~C.

\begin{maintheorem}{E}
\label{thm:E}
Let \(\Omega\) be a spectral union of \(n\) pairwise disjoint non-adjacent bounded intervals with lengths \(\ell = (\ell_1, \dots, \ell_n)\), with spectrum \(\Lambda\) of minimal period \(k/\lvert\Omega\rvert\). Then \(k\) divides the integer \(\kzero(\ell)\) determined by the subset sums of the lengths, and \(\kzero(\ell)\) divides
\[
\Kn = \operatorname{lcm}(1, 2, \dots, d(n)), \qquad d(n) = \bigl\lfloor (n+2)^{(n+2)/2} / 2^{n+1} \bigr\rfloor;
\]
if \(\dimQ \Krel \le 4\), then \(\kzero(\ell) \le 10\). If the fiber lemma \(\FL{k}\) holds for every \(k \le \Kn\), then every spectral union of \(n\) intervals tiles \(\R\) by translations.
\end{maintheorem}

Relative to the work of Shi \cite{shi2026} and of Bell \cite{bell2026four}, the new elements are the following. The level theory of the subset sums (effective levels, pinned ratios and the verdict rule) is set up for every \(n\), together with the bound \(\Kn\). The specialization lemma reduces the arithmetic of any number of lengths with at most four rational relations to that of five lengths, which is classified in full. The fiber-side mechanisms, the loop exchange, the propagation of cosets and lineages along classes, and the monodromy of the exchange spaces around a closed walk, replace the case analysis of \cite{bell2026four} by arguments on the space \(\C^S\). The six-point dichotomy and seven-point rigidity give the fiber lemmas up to seven cosets. And the gluing construction explains every failure of the consecutive-pair statement with a periodic spectrum, starting at twelve points.

Section~\ref{sec:general} records the three places where the method stops. At twelve points the consecutive-pair statement is false, so the count \(k \le n\) cannot be proved in general by this route, and indeed eight intervals can carry twelve cosets. For \(\dimQ \Krel > 4\), the number of cosets is bounded only by \(\Kn\), and the fiber lemmas proved here cover only \(k \le 7\). And the fiber lemma for a given \(k\) contains the discrete direction of Fuglede's conjecture for \(k\)-element sets of integers, so the general problem on the line is at least as hard as that discrete problem.

Section~\ref{sec:fibers} fixes the definitions, restates the results of \cite{bell2026four} that are used, and proves the conservation law, the dimension count and the merging lemma for \(n\) intervals. Section~\ref{sec:levels} develops the effective levels and proves the verdict rule and the bound \(\Kn\). Section~\ref{sec:fivelengths} classifies the effective levels of five lengths and proves the specialization lemma. Section~\ref{sec:consecutive} proves the consecutive-pair statement up to ten points. Section~\ref{sec:classgraph} introduces the class graph and reduces Theorem~A to fully mixing configurations of at most eight intervals. Section~\ref{sec:mechanisms} contains the fiber-side mechanisms and the bounds on \(k\) for six, seven and eight intervals. Section~\ref{sec:assembly} assembles the proof of Theorem~A and Corollary~B and gives the sharpness example. Section~\ref{sec:twelve} proves Theorem~D(i) to (iii): eleven points, the failure at twelve, and the gluing. Section~\ref{sec:fiberlemmas} proves Theorem~C. Section~\ref{sec:general} proves Theorem~E and Theorem~D(iv), and closes with the limits of the method and the open questions. Appendix~\ref{app:cases} contains the case analyses behind the classification of five lengths and the list of class graphs at eight intervals.

\section{Fibers, classes and merging for \texorpdfstring{$n$}{n} intervals}
\label{sec:fibers}

We recall the periodicity and fiber descriptions, then establish the conservation law and merging arguments used in Sections~\ref{sec:classgraph} to~\ref{sec:assembly}. The results imported from \cite{bell2026four} are stated with their original references.

\begin{definition}[standing reduction]
\label{def:standing}
Let \(\Omega \subset \R\) be a bounded measurable set which, up to a null set, is a finite union of bounded intervals. Spectrality and tiling by translations are unchanged when \(\Omega\) is altered on a null set, so the intervals may be taken half-open, and two of them that overlap or share an endpoint may be replaced by their union. We therefore write
\[
\Omega = \bigsqcup_{j=1}^{n} [a_j, b_j), \qquad \ell_j = b_j - a_j > 0, \qquad L = \lvert\Omega\rvert = \sum_{j=1}^{n} \ell_j,
\]
with the \(n\) components pairwise disjoint and non-adjacent, that is, \(b_i \ne a_j\) \(\forall i, j\). This is the standing reduction, and \(n\) denotes the number of components after it. We write \(\nset = \{1, \dots, n\}\); for \(T \subset \nset\) the subset sum is \(\ell_T = \sum_{j \in T} \ell_j\), the indicator of \(T\) is \(1_T \in \{0,1\}^n\), and \(\one = 1_{\nset}\). The space of rational linear relations among the lengths is
\[
\Krel = \Bigl\{ r \in \Q^n : \sum_{j=1}^{n} r_j \ell_j = 0 \Bigr\},
\]
and \(\rho = \dimQ \operatorname{span}_{\Q}(\ell_1, \dots, \ell_n)\), so that \(\dimQ \Krel = n - \rho\), since \(\Krel\) is the kernel of the \(\Q\)-linear map \(r \mapsto \sum_j r_j \ell_j\) from \(\Q^n\) onto that span.
\end{definition}

We first fix the spectral and fiber notation, which applies to every bounded measurable set of positive measure.

\begin{definition}
\label{def:spectrum}
For \(\lambda \in \R\) the normalized exponential on \(\Omega\) is \(e_\lambda(x) = L^{-1/2} \ex{\lambda x}\). A set \(\Lambda \subset \R\) is a \emph{spectrum} of \(\Omega\), and \(\Omega\) is \emph{spectral}, when \(\{e_\lambda : \lambda \in \Lambda\}\) is an orthonormal basis of \(L^2(\Omega)\). Since \(\langle e_\lambda, e_{\lambda'} \rangle = L^{-1} \widehat{1_\Omega}(\lambda' - \lambda)\) with \(\widehat{1_\Omega}(\xi) = \int_\Omega \ex{-\xi t}\,dt\), orthogonality is the vanishing of \(\widehat{1_\Omega}\) on \((\Lambda - \Lambda) \setminus \{0\}\). The set \(\Omega\) \emph{tiles \(\R\) by translations} when there is a countable set \(T \subset \R\) such that the translates \(\Omega + t\), \(t \in T\), are pairwise disjoint up to null sets and cover \(\R\) up to a null set.
\end{definition}

\begin{definition}
\label{def:fiber}
Let \(h > 0\) and let \(\Gamma \subset \R\) be a finite set whose elements are pairwise distinct modulo \(h^{-1}\Z\). Write \(k = \lvert\Gamma\rvert\), \(\Lambda = \Gamma + h^{-1}\Z\) and \(S = h\Gamma \subset \R/\Z\), so that \(\lvert S\rvert = k\). For \(x \in \R\) the \emph{fiber} of \(\Omega\) at \(x\) is the finite set
\[
A_x = \{ m \in \Z : x + mh \in \Omega \},
\]
and its character sum is \(\widehat{A_x}(\sigma) = \sum_{m \in A_x} \ex{\sigma m}\) for \(\sigma \in \R/\Z\). The functions \(x \mapsto \lvert A_x\rvert\) and \(x \mapsto \widehat{A_x}(\sigma)\) are bounded and measurable on \([0, h)\), and \(A_{x+h} = A_x - 1\), so that \(\widehat{A_{x+h}}(\sigma) = \ex{-\sigma}\,\widehat{A_x}(\sigma)\).
\end{definition}

\begin{definition}
\label{def:spectralpair}
Let \(S \subset \R/\Z\) be finite, of cardinality \(k\), and let \(A \subset \Z\) be finite. Then \(A\) is a \emph{spectrum} of \(S\), and \((A, S)\) is a \emph{spectral pair}, when \(\lvert A\rvert = k\) and \(\widehat{A}(\sigma) = 0\) \(\forall \sigma \in (S - S) \setminus \{0\}\). Equivalently, \(H = (\ex{sa})_{s \in S,\, a \in A}\) is a complex Hadamard matrix of order \(k\), that is, \(HH^* = kI\); then also \(H^*H = kI\), so that \(\widehat{\muS}(a - a') = \sum_{s \in S} \ex{s(a - a')} = 0\) for \(a \ne a'\) in \(A\), where \(\muS = \sum_{s \in S} \delta_s\).
\end{definition}

The periodicity theorem of Iosevich and Kolountzakis \cite[Theorem~1.5]{ik2013}, in the form of \cite[Theorem~2.4]{bell2026four}, and the fiber reformulation turn the problem into a finite one.

\begin{theorem}[{\cite[Theorem~2.4]{bell2026four}}]
\label{thm:periodic}
Let \(\Lambda\) be a spectrum of a bounded measurable set \(\Omega \subset \R\) of measure \(L\). Every positive period of \(\Lambda\) lies in \(L^{-1}\Z_{>0}\), and the set of all periods of \(\Lambda\) is \(h^{-1}\Z\) for a unique \(h > 0\). With \(\Gamma = \Lambda \cap [0, h^{-1})\), \(k = \lvert\Gamma\rvert\) and \(S = h\Gamma \subset \R/\Z\), one has \(\Lambda = \Gamma + h^{-1}\Z\), \(\lvert S\rvert = k\), and \(S\) has no nontrivial translation symmetry: \(S + t = S\) with \(t \in \R/\Z\) forces \(t = 0\).
\end{theorem}

\begin{theorem}[{\cite[Theorem~2.5]{bell2026four}}]
\label{thm:fiber}
Let \(\Omega \subset \R\) be bounded and measurable, of measure \(L\), and let \(h\), \(\Gamma\), \(k\), \(\Lambda\) and \(S\) be as in Definition~\ref{def:fiber}. Then \(\Lambda\) is a spectrum of \(\Omega\) iff, for almost every \(x \in [0, h)\), the fiber \(A_x\) is a spectrum of \(S\). When these hold, \(L = kh\).
\end{theorem}

Throughout the paper, when \(\Lambda\) is a spectrum of \(\Omega\), the letters \(h\), \(\Gamma\), \(k\) and \(S\) refer to the minimal period of Theorem~\ref{thm:periodic}: \(h^{-1}\) is the least positive period of \(\Lambda\), \(L = kh\), and \(k\) is the number of cosets of \(h^{-1}\Z\) in \(\Lambda\).

\subsection{The fibers of a union of intervals}

\begin{proposition}
\label{prop:intervalfibers}
Let \(\Omega = \bigsqcup_{j=1}^{n} [a_j, b_j)\) be as in Definition~\ref{def:standing}, and let \(h > 0\). For \(x \in \R\) and \(j \in \nset\) set
\[
I_j(x) = \Z \cap \bigl[ \lceil (a_j - x)/h \rceil,\ \lceil (b_j - x)/h \rceil \bigr).
\]
Then \(A_x = \bigsqcup_{j=1}^{n} I_j(x)\), each \(I_j(x)\) is a set of consecutive integers, and \(\lvert I_j(x)\rvert\) equals \(\lfloor \ell_j/h \rfloor\) or \(\lceil \ell_j/h \rceil\); in particular every fiber is a union of at most \(n\) runs. Let
\[
E = \bigcup_{j=1}^{n} \bigl( (a_j + h\Z) \cup (b_j + h\Z) \bigr).
\]
The map \(x \mapsto A_x\) is right-continuous and constant on each interval between consecutive points of \(E\). For \(c \in E\), let \(\In{c}\) be the set of starting heights \((a_j - c)/h\) over the \(j\) with \(a_j \in c + h\Z\), and \(\Out{c}\) the set of ending heights \((b_i - c)/h\) over the \(i\) with \(b_i \in c + h\Z\); these are sets of integers whose cardinalities are the numbers of left and of right endpoints in \(c + h\Z\). If \(A_{c-}\) denotes the constant value of the fiber on a left neighborhood of \(c\), then \(\Out{c} \subset A_{c-}\), \(\In{c}\) is disjoint from \(A_{c-} \setminus \Out{c}\), and
\begin{equation}
\label{eq:swap}
A_c = \bigl(A_{c-} \setminus \Out{c}\bigr) \sqcup \In{c}.
\end{equation}
Consequently \(\widehat{A_c}(\sigma) - \widehat{A_{c-}}(\sigma) = \widehat{\In{c}}(\sigma) - \widehat{\Out{c}}(\sigma)\) \(\forall \sigma \in \R/\Z\). Replacing \(c\) by \(c + h\) shifts \(\In{c}\), \(\Out{c}\), \(A_{c-}\) and \(A_c\) by \(-1\). Finally, \(\int_0^h \lvert A_x\rvert\,dx = L\); the cardinality \(\lvert A_x\rvert\) is constant in \(x\) iff \(\lvert\In{c}\rvert = \lvert\Out{c}\rvert\) \(\forall c \in E\), and the constant is then \(L/h\).
\end{proposition}

\begin{proof}
An integer \(m\) lies in \(A_x\) iff \(x + mh\) lies in exactly one component \([a_j, b_j)\), the components being disjoint, and \(a_j \le x + mh < b_j\) is the condition \(\lceil (a_j - x)/h \rceil \le m < \lceil (b_j - x)/h \rceil\) because \(m\) is an integer. This is the decomposition of \(A_x\). For the length, write \(\ell_j/h = N + \theta\) with \(N = \lfloor \ell_j/h \rfloor\) and \(\theta \in [0, 1)\), and set \(u = (a_j - x)/h\). Then
\[
\lceil u + N + \theta \rceil - \lceil u \rceil = N + \bigl(\lceil u + \theta \rceil - \lceil u \rceil\bigr),
\]
and \(\lceil u + \theta \rceil - \lceil u \rceil \in \{0, 1\}\) because \(u \le u + \theta < u + 1\); when \(\theta = 0\) the difference is \(0\). Thus \(\lvert I_j(x)\rvert\) is \(\lfloor \ell_j/h \rfloor\) or \(\lceil \ell_j/h \rceil\).

For fixed \(e \in \R\) and \(x \in [e - mh,\ e - (m-1)h)\), the number \((e - x)/h\) lies in \((m - 1,\ m]\), so its ceiling equals \(m\). These intervals partition \(\R\) as \(m\) runs through \(\Z\). The function \(x \mapsto \lceil (e - x)/h \rceil\) is therefore integer valued, non-increasing and right-continuous, constant between consecutive points of \(e + h\Z\), and it drops by exactly \(1\) at each such point. Each of the \(2n\) ceilings defining the \(I_j(x)\) is constant on every interval between consecutive points of \(E\), and so is \(A_x\). The set \(E\) is a finite union of translates of \(h\Z\), hence closed and discrete.

Now fix \(c \in E\) and let \(\varepsilon > 0\) be smaller than the distance from \(c\) to \(E \setminus \{c\}\). If \(a_j = c + hp\) with \(p \in \Z\), the previous paragraph gives \(\lceil (a_j - x)/h \rceil = p + 1\) on \([c - h,\ c)\) and \(= p\) on \([c,\ c + h)\). The same holds for a right endpoint \(b_i = c + hq\), with the values \(q + 1\) on the left and \(q\) on the right. An endpoint outside \(c + h\Z\) contributes no change on \((c - \varepsilon,\ c + \varepsilon)\).

Fix a component \(j\), and write \(\alpha\) and \(\beta\) for the two ceilings of \(I_j\) on \((c - \varepsilon,\ c)\), so that \(I_j(c-) = \Z \cap [\alpha, \beta)\). If \(a_j \in c + h\Z\) and \(b_j \notin c + h\Z\), then \(I_j(c) = \Z \cap [\alpha - 1,\ \beta)\), and \(I_j\) gains the point \(p = \alpha - 1\). If \(b_j \in c + h\Z\) and \(a_j \notin c + h\Z\), then \(\beta = q + 1\), and \(q \in I_j(c-)\) because \(c - \varepsilon' + qh = b_j - \varepsilon'\) lies in \([a_j, b_j)\) \(\forall \varepsilon' \in (0, \ell_j)\); thus \(I_j(c) = \Z \cap [\alpha,\ \beta - 1)\), and \(I_j\) loses the point \(q\). If both endpoints lie in \(c + h\Z\), then \(I_j(c) = \Z \cap [\alpha - 1,\ \beta - 1)\), and \(1_{I_j(c)} - 1_{I_j(c-)} = 1_{\{\alpha - 1\}} - 1_{\{\beta - 1\}}\) whether or not \(\alpha = \beta\). If neither does, \(I_j\) is unchanged. Summing over \(j\) and using the disjointness of the \(I_j\) at each \(x\),
\[
1_{A_c} - 1_{A_{c-}} = \sum_{j :\, a_j \in c + h\Z} 1_{\{p_j\}} - \sum_{i :\, b_i \in c + h\Z} 1_{\{q_i\}}.
\]
The heights \(q_i\) in the second sum are the maximal points of the pairwise disjoint integer intervals \(I_i(c-)\), hence distinct, and each lies in \(A_{c-}\); the heights \(p_j\) in the first sum are the minimal points of the pairwise disjoint intervals \(I_j(c)\), hence distinct. So the two sums are \(1_{\In{c}}\) and \(1_{\Out{c}}\), the cardinalities are as stated, and \(\Out{c} \subset A_{c-}\).

If \(p \in \In{c} \cap A_{c-}\), evaluation of the identity at \(p\) gives \(1_{A_c}(p) - 1 = 1 - 1_{\Out{c}}(p)\), so \(1_{A_c}(p) = 2 - 1_{\Out{c}}(p) \le 1\), which forces \(p \in \Out{c}\). This is the disjointness of \(\In{c}\) and \(A_{c-} \setminus \Out{c}\), and with it the identity of indicators is \eqref{eq:swap}. Multiplying by \(\ex{\sigma m}\) and summing over \(m \in \Z\) gives the formula for the character sums, and \(\sigma = 0\) gives \(\lvert A_c\rvert - \lvert A_{c-}\rvert = \lvert\In{c}\rvert - \lvert\Out{c}\rvert\). The behavior under \(c \mapsto c + h\) is the shift \(A_{x+h} = A_x - 1\) of Definition~\ref{def:fiber}.

Finally, \(\int_0^h \lvert A_x\rvert\,dx = \sum_{m \in \Z} \int_0^h 1_\Omega(x + mh)\,dx = \int_\R 1_\Omega = L\). The function \(x \mapsto \lvert A_x\rvert\) is a right-continuous step function whose only possible jumps are at the points of \(E\), where the jump is \(\lvert\In{c}\rvert - \lvert\Out{c}\rvert\); it is constant on \(\R\) iff every jump vanishes, and the constant is then \(L/h\) by the integral just computed.
\end{proof}

\begin{lemma}
\label{lem:everyfiber}
Let \(\Omega\) be as in Definition~\ref{def:standing}, let \(h\), \(\Gamma\), \(k\), \(\Lambda\) and \(S\) be as in Definition~\ref{def:fiber}, and suppose that \(\Lambda\) is a spectrum of \(\Omega\). Then \(L = kh\), and every fiber \(A_x\), \(x \in \R\), is a spectrum of \(S\). In particular \(\lvert A_x\rvert = k\) \(\forall x\), so that \(\sum_{m \in \Z} 1_\Omega(x + mh) = k\) \(\forall x \in \R\): the set \(\Omega\) is a \(k\)-fold multitile of \(\R\) by \(h\Z\).
\end{lemma}

\begin{proof}
Theorem~\ref{thm:fiber} gives \(L = kh\) and a set of full measure in \([0, h)\) on which \(A_x\) is a spectrum of \(S\), that is, \(\lvert A_x\rvert = k\) and \(\widehat{A_x}(\sigma) = 0\) for \(\sigma \in (S - S) \setminus \{0\}\). By Proposition~\ref{prop:intervalfibers}, \(x \mapsto A_x\) is a step function on \([0, h)\) with finitely many pieces, each an interval of positive length. Every value of \(A_x\) on \([0, h)\) is therefore attained on a set of positive measure, which meets the full-measure set; so every fiber over \([0, h)\) is a spectrum of \(S\). The shift \(A_{x+h} = A_x - 1\) preserves cardinality and the vanishing of the character sums on \(S - S\), and extends the conclusion to every \(x \in \R\).
\end{proof}

\subsection{Endpoint classes and the conservation law}

\begin{definition}
\label{def:classes}
Let \(\Omega\) be as in Definition~\ref{def:standing}, let \(\Lambda\) be a spectrum of \(\Omega\), and let \(h\) be as in Theorem~\ref{thm:periodic}. The \emph{endpoint classes} are the residues modulo \(h\Z\) of the \(2n\) endpoints \(a_j, b_j\); their representatives in \([0, h)\) are written \(c_1 < \dots < c_r\), and a class is identified with its representative. For \(j \in \nset\) let \(u(j), v(j) \in [r]\) be the classes of \(a_j\) and \(b_j\), and write
\[
a_j = c_{u(j)} + h p_j, \qquad b_j = c_{v(j)} + h q_j, \qquad d_j = q_j - p_j,
\]
with \(p_j, q_j \in \Z\); then \(\ell_j = c_{v(j)} - c_{u(j)} + h d_j\), so \(d_j\) is \(\lfloor \ell_j/h \rfloor\) or \(\lceil \ell_j/h \rceil\). At a class \(c\), the set of \emph{leaving heights} is \(\Out{c} = \{q_j : c_{v(j)} = c\}\) and the set of \emph{entering heights} is \(\In{c} = \{p_j : c_{u(j)} = c\}\); these are the sets of Proposition~\ref{prop:intervalfibers} at the representative \(c\). By Lemma~\ref{lem:everyfiber} the cardinality \(\lvert A_x\rvert\) is constant, so Proposition~\ref{prop:intervalfibers} gives the \emph{counting law} \(\lvert\In{c}\rvert = \lvert\Out{c}\rvert\) at every class: the multiset of residues of the left endpoints equals the multiset of residues of the right endpoints. The common value \(m_c = \lvert\In{c}\rvert = \lvert\Out{c}\rvert \ge 1\) is the \emph{size} of the class, the number of intervals that start, equivalently end, at \(c\). Since every interval starts at exactly one class, \(\sum_c m_c = n\) and \(r \le n\). The multiset of sizes \((m_c)\) is the \emph{pattern} of the configuration, a partition of \(n\).
\end{definition}

On \([c_1,c_1+h)\), the fiber is constant on each of the \(r\) intervals \([c_i,c_{i+1})\), \(i<r\), and \([c_r,c_1+h)\). Consecutive values are related by the swap \eqref{eq:swap} of \(m_c\) heights for \(m_c\) heights; after one full turn the first value is shifted according to \(A_{x+h}=A_x-1\).

\begin{theorem}[conservation law]
\label{thm:conservation}
Let \(\Omega = \bigsqcup_{j=1}^{n} [a_j, b_j)\) be as in Definition~\ref{def:standing}, of measure \(L\), and let \(h\), \(\Gamma\), \(k\), \(\Lambda\) and \(S\) be as in Definition~\ref{def:fiber}, with the classes, \(\In{c}\) and \(\Out{c}\) computed from this \(h\). Say that the conservation law holds at the class \(c\) when
\begin{equation}
\label{eq:conservation}
\sum_{q \in \Out{c}} \ex{\sigma q} = \sum_{p \in \In{c}} \ex{\sigma p} \qquad \forall \sigma \in S - S,
\end{equation}
including \(\sigma = 0\), where it reads \(\lvert\Out{c}\rvert = \lvert\In{c}\rvert\). If \(\Lambda\) is a spectrum of \(\Omega\), then \eqref{eq:conservation} holds at every class. Conversely, if \(L = kh\) and \eqref{eq:conservation} holds at every class, then every fiber is a spectrum of \(S\) and \(\Lambda\) is a spectrum of \(\Omega\). The same conclusion holds if a single fiber \(A_{x_0}\) is a spectrum of \(S\) and \eqref{eq:conservation} holds at every class.
\end{theorem}

\begin{proof}
Both sides of \eqref{eq:conservation} are multiplied by \(\ex{-\sigma}\) when the representative \(c\) is replaced by \(c + h\), by the last assertion of Proposition~\ref{prop:intervalfibers}, so the condition depends only on the class. Suppose that \(\Lambda\) is a spectrum of \(\Omega\). By Lemma~\ref{lem:everyfiber}, both \(A_c\) and \(A_{c-}\) are spectra of \(S\), so for \(\sigma \in S - S\) their character sums agree: both vanish when \(\sigma \ne 0\), and both equal \(k\) when \(\sigma = 0\). The difference formula of Proposition~\ref{prop:intervalfibers} then gives \(\widehat{\In{c}}(\sigma) = \widehat{\Out{c}}(\sigma)\), which is \eqref{eq:conservation}.

Conversely, fix \(\sigma \in S - S\) and set \(\psi_\sigma(x) = \widehat{A_x}(\sigma)\). By Proposition~\ref{prop:intervalfibers} this is a right-continuous step function on \(\R\) whose only possible jumps are at the points of \(E\), and the jump at \(c\) equals \(\widehat{\In{c}}(\sigma) - \widehat{\Out{c}}(\sigma)\), which vanishes by \eqref{eq:conservation}. Hence \(\psi_\sigma\) is constant on \(\R\). On the other hand \(\psi_\sigma(x + h) = \ex{-\sigma}\,\psi_\sigma(x)\) by Definition~\ref{def:fiber}. If \(\sigma \ne 0\) in \(\R/\Z\), then \(\ex{-\sigma} \ne 1\), and the constant is \(0\). For \(\sigma = 0\), the identity \eqref{eq:conservation} says \(\lvert\In{c}\rvert = \lvert\Out{c}\rvert\) at every class, so Proposition~\ref{prop:intervalfibers} makes \(\lvert A_x\rvert\) constant and equal to \(L/h = k\). Every fiber is therefore a spectrum of \(S\), and Theorem~\ref{thm:fiber} says that \(\Lambda\) is a spectrum of \(\Omega\).

If instead a single fiber \(A_{x_0}\) is a spectrum of \(S\), the constancy of each \(\psi_\sigma\) gives \(\psi_\sigma(x) = \psi_\sigma(x_0)\) \(\forall x\), which is \(0\) for \(\sigma \ne 0\) and \(k\) for \(\sigma = 0\). Every fiber is a spectrum of \(S\), so \(\lvert A_x\rvert = k\) and \(L = \int_0^h \lvert A_x\rvert\,dx = kh\), and Theorem~\ref{thm:fiber} applies again.
\end{proof}

With \(L=kh\), the positions \(c_1, \dots, c_r\) of the classes do not enter \eqref{eq:conservation}: spectrality with the period \(h^{-1}\) is determined by the integer heights \((\In{c}, \Out{c})_c\) and by \(S\).

\begin{lemma}[dimension count]
\label{lem:dimension}
Let \(S \subset \R/\Z\) be finite, of cardinality \(k \ge 1\). Let \(B\) be a nonempty finite multiset in \(\Z\), with multiplicity function \(\mu_B \colon \Z \to \Z_{\ge 0}\), support \(\operatorname{supp} B = \{ b : \mu_B(b) > 0 \}\), cardinality \(\lvert B\rvert = \sum_b \mu_B(b)\) and character sum \(\widehat{B}(\sigma) = \sum_b \mu_B(b) \ex{\sigma b}\). If \(\widehat{B}(\sigma) = 0\) \(\forall \sigma \in (S - S) \setminus \{0\}\), then
\[
k \le \lvert\operatorname{supp} B\rvert \le \lvert B\rvert.
\]
A multiset of cardinality \(k\) with this vanishing is a set, and is a spectrum of \(S\). If \(A\) is a spectrum of \(S\) and \(B\) is a nonempty multiset with \(\mu_B \le 1_A\) and with the same vanishing, then \(B = A\).
\end{lemma}

\begin{proof}
For \(s \in S\) define \(u_s \in \C^{\operatorname{supp} B}\) by \((u_s)_b = \mu_B(b)^{1/2} \ex{sb}\); the exponential is well defined because \(b \in \Z\). Then
\[
\langle u_s, u_{s'} \rangle = \sum_b \mu_B(b) \ex{(s - s') b} = \widehat{B}(s - s'),
\]
which vanishes for \(s \ne s'\), while \(\lVert u_s\rVert^2 = \lvert B\rvert > 0\). The \(k\) vectors \(u_s\) are pairwise orthogonal and nonzero in \(\C^{\operatorname{supp} B}\), so \(k \le \lvert\operatorname{supp} B\rvert\). The comparison of support and cardinality is immediate.

If \(\lvert B\rvert = k\), the two inequalities force \(\lvert\operatorname{supp} B\rvert = \lvert B\rvert\), hence \(\mu_B \le 1\), and the vanishing hypothesis with the cardinality are the conditions of Definition~\ref{def:spectralpair}. If \(\mu_B \le 1_A\) for a spectrum \(A\) of \(S\), the inequalities give \(k \le \lvert B\rvert \le \lvert A\rvert = k\), so \(\lvert B\rvert = \lvert A\rvert\) and \(\mu_B = 1_A\).
\end{proof}

\subsection{Rigid pairs and merging}

In this subsection \(\Omega\) is as in Definition~\ref{def:standing} with \(n \ge 2\), \(\Lambda\) is a spectrum of \(\Omega\), and \(h\), \(k\), \(\Gamma\), \(S\) are as in Theorem~\ref{thm:periodic}; the classes, heights and sizes are those of Definition~\ref{def:classes}.

\begin{definition}
\label{def:rigid}
Let \(i, j \in \nset\) satisfy \(b_i \equiv a_j \pmod{h\Z}\), and set \(d = (b_i - a_j)/h \in \Z\). The pair \((i, j)\) is \emph{rigid} when \(d(S - S) \subset \Z\), that is, when \(\ex{\sigma d} = 1\) \(\forall \sigma \in S - S\); the integer \(d\) is its \emph{slide}. The endpoints \(b_i\) and \(a_j\) then lie in one class \(c\), with \(q_i \in \Out{c}\), \(p_j \in \In{c}\) and \(d = q_i - p_j\). A class is \emph{rigid} when it contains a rigid pair, and \emph{mixing} otherwise. The configuration is \emph{fully mixing} when no pair is rigid.
\end{definition}

\begin{lemma}
\label{lem:rigidpair}
Let \(i, j \in \nset\) satisfy \(b_i \equiv a_j \pmod{h\Z}\), and set \(d = (b_i - a_j)/h\).
\begin{enumerate}
\item[\rm (a)] \(d \ne 0\). If \((i, j)\) is rigid, then \(S\) lies in a single coset of \((1/\lvert d\rvert)\Z/\Z\), so that \(k \le \lvert d\rvert\) and \(L \le \lvert b_i - a_j\rvert\).
\item[\rm (b)] No pair \((i, i)\) is rigid: if \(\ell_i \in h\Z\), then \((\ell_i/h)(S - S) \not\subset \Z\).
\item[\rm (c)] A class of size one is rigid, and the two indices of its rigid pair are distinct.
\item[\rm (d)] Let \(c\) be a class of size two, with \(\Out{c} = \{q_1, q_2\}\) and \(\In{c} = \{p_1, p_2\}\). Then \(c\) is rigid iff one of the two matchings of \(\Out{c}\) with \(\In{c}\) satisfies \eqref{eq:conservation} term by term for every \(\sigma \in S - S\), and in that case both pairs of that matching are rigid. For each single \(\sigma \in S - S\), either one of the two matchings satisfies \eqref{eq:conservation} term by term at \(\sigma\), or \(\ex{\sigma(q_2 - q_1)} = -1 = \ex{\sigma(p_2 - p_1)}\), so that both sides of \eqref{eq:conservation} vanish at \(\sigma\).
\end{enumerate}
\end{lemma}

\begin{proof}
(a) The equality \(d = 0\) would mean \(b_i = a_j\), which the standing reduction excludes. Suppose \(d(S - S) \subset \Z\). For \(s, s' \in S\) the element \(s - s'\) of \(\R/\Z\) is killed by multiplication by \(d\), hence lies in the subgroup \((1/\lvert d\rvert)\Z/\Z\), which has \(\lvert d\rvert\) elements. Thus \(S\) lies in one coset of that subgroup, and \(k = \lvert S\rvert \le \lvert d\rvert\). Multiplying by \(h\) and using \(L = kh\) gives \(L \le \lvert d\rvert h = \lvert b_i - a_j\rvert\).

(b) Suppose that \(\ell_i = dh\) with \(d \ge 1\) and \(d(S - S) \subset \Z\). For every \(x\) the window \([(a_i - x)/h,\ (a_i - x)/h + d)\) of length \(d\) contains exactly \(d\) integers, so \(I_i(x)\) is a block of \(d\) consecutive integers \(\{m_0, \dots, m_0 + d - 1\}\). For \(\sigma \in (S - S) \setminus \{0\}\) one has \(\ex{\sigma} \ne 1\) and \(\ex{\sigma d} = 1\), so the character sum of the block is the geometric sum
\[
\widehat{I_i(x)}(\sigma) = \ex{\sigma m_0}\,\frac{\ex{\sigma d} - 1}{\ex{\sigma} - 1} = 0.
\]
By Lemma~\ref{lem:everyfiber}, \(A_x\) is a spectrum of \(S\), and Lemma~\ref{lem:dimension}, applied to the nonempty sub-multiset \(I_i(x)\) of \(A_x\), gives \(I_i(x) = A_x\) \(\forall x\). Therefore
\[
\lvert\Omega \setminus [a_i, b_i)\rvert = \int_0^h \lvert A_x \setminus I_i(x)\rvert\,dx = 0,
\]
which contradicts the presence of a second component of positive length.

(c) Let \(c\) be a class with \(\Out{c} = \{q_i\}\) and \(\In{c} = \{p_j\}\). The conservation law \eqref{eq:conservation} at \(c\) reads \(\ex{\sigma q_i} = \ex{\sigma p_j}\) \(\forall \sigma \in S - S\), so \(\ex{\sigma(q_i - p_j)} = 1\) and \((i, j)\) is rigid with slide \(q_i - p_j = (b_i - a_j)/h\). The indices are distinct by (b).

(d) If, say, the pair with heights \((q_1, p_1)\) is rigid, then \(\ex{\sigma q_1} = \ex{\sigma p_1}\) \(\forall \sigma \in S - S\), and subtracting this from \eqref{eq:conservation} gives \(\ex{\sigma q_2} = \ex{\sigma p_2}\) for every such \(\sigma\): the matching \(q_1 \leftrightarrow p_1\), \(q_2 \leftrightarrow p_2\) holds term by term and both its pairs are rigid. Conversely, a termwise matching consists of two rigid pairs. For the last assertion, fix \(\sigma\) and apply the following observation to \(u_r = \ex{\sigma q_r}\), \(v_r = \ex{\sigma p_r}\): if \(u_1, u_2, v_1, v_2\) are complex numbers of modulus one with \(u_1 + u_2 = v_1 + v_2 = w\), then either \(\{u_1, u_2\} = \{v_1, v_2\}\) as multisets, or \(w = 0\). Indeed, if \(w \ne 0\), then \(\overline{u_r} = 1/u_r\) gives \(\overline{w} = 1/u_1 + 1/u_2 = w/(u_1 u_2)\), hence \(u_1 u_2 = w/\overline{w}\), and likewise \(v_1 v_2 = w/\overline{w}\); both pairs are the multiset of roots of \(z^2 - wz + w/\overline{w}\). Equal multisets means that one of the two matchings holds at \(\sigma\); \(w = 0\) means \(u_2 = -u_1\) and \(v_2 = -v_1\).
\end{proof}

\begin{lemma}[the slide]
\label{lem:slide}
Let \((i, j)\) be a rigid pair with slide \(d\), set \(J = [a_j, b_j)\) and \(\Omega' = (\Omega \setminus J) \cup (J + dh)\). Then \(J + dh = [b_i, b_i + \ell_j)\) is disjoint from \(\Omega \setminus J\). The set \(\Omega'\) is the disjoint union of the \(n - 1\) half-open intervals \([a_m, b_m)\), \(m \notin \{i, j\}\), and \([a_i, b_i + \ell_j)\); it has measure \(L\) and, after the standing reduction, at most \(n - 1\) components. The set \(\Lambda\) is a spectrum of \(\Omega'\), and \(h\), \(k\), \(\Gamma\) and \(S\) are the same for \(\Omega'\) as for \(\Omega\).
\end{lemma}

\begin{proof}
By Lemma~\ref{lem:rigidpair}(b), \(i \ne j\). The identity \(a_j + dh = b_i\) gives \(J + dh = [b_i, b_i + \ell_j)\). Consider the function \(g = 1_{\Omega \setminus J} + 1_{J + dh}\), with values in \(\{0, 1, 2\}\). Its fiber multiset at \(x\) is
\[
G_x = \bigl(A_x \setminus I_j(x)\bigr) \sqcup \bigl(I_j(x) + d\bigr),
\]
because \(x + mh \in J + dh\) iff \(m - d \in I_j(x)\). For every \(\sigma \in S - S\),
\[
\widehat{G_x}(\sigma) = \widehat{A_x}(\sigma) - \widehat{I_j(x)}(\sigma) + \ex{\sigma d}\,\widehat{I_j(x)}(\sigma) = \widehat{A_x}(\sigma),
\]
since \(\ex{\sigma d} = 1\) by rigidity. This is the only use of rigidity. By Lemma~\ref{lem:everyfiber}, \(A_x\) is a spectrum of \(S\) \(\forall x\), so \(G_x\) is a multiset of cardinality \(k\) whose character sums vanish on \((S - S) \setminus \{0\}\), and Lemma~\ref{lem:dimension} says that \(G_x\) is a set. In particular \((I_j(x) + d) \cap (A_x \setminus I_j(x)) = \emptyset\) \(\forall x\).

For a bounded measurable set \(F \subset \R\), the decomposition of Lebesgue measure along the cosets of \(h\Z\) gives \(\lvert F\rvert = \int_0^h \#\{m : x + mh \in F\}\,dx\). Applied to \(F = (J + dh) \cap (\Omega \setminus J)\), the integrand is the cardinality just shown to vanish, so \(F\) is a null set. The intersection of \(J + dh\) with any one component \([a_m, b_m)\), \(m \ne j\), is an intersection of two half-open intervals, hence of the form \([u, v)\), and a null set of this form is empty. Thus \(J + dh\) is disjoint from \(\Omega \setminus J\), the function \(g\) is the indicator of \(\Omega'\), and
\[
[a_i, b_i) \cup [b_i, b_i + \ell_j) = [a_i, b_i + \ell_j).
\]
The \(n - 1\) resulting intervals are pairwise disjoint, so \(\lvert\Omega'\rvert = L\). Among the intervals \([a_m, b_m)\), \(m \notin \{i, j\}\), no two are adjacent, and none is adjacent to \([a_i, b_i + \ell_j)\) at its left end \(a_i\), by the convention on \(\Omega\). The only possible new adjacency is \(b_i + \ell_j = a_m\) for some \(m\), and coalescing adjacent intervals does not increase their number, so \(\Omega'\) has at most \(n - 1\) components after the standing reduction. For every \(x\) the fiber of \(\Omega'\) is the set \(G_x\), a spectrum of \(S\), so Theorem~\ref{thm:fiber} says that \(\Lambda\) is a spectrum of \(\Omega'\). The minimal period of \(\Lambda\), and with it \(h\), \(\Gamma\), \(k\) and \(S\), is determined by \(\Lambda\) alone (Theorem~\ref{thm:periodic}).
\end{proof}

\begin{figure}[ht]
\centering
\input{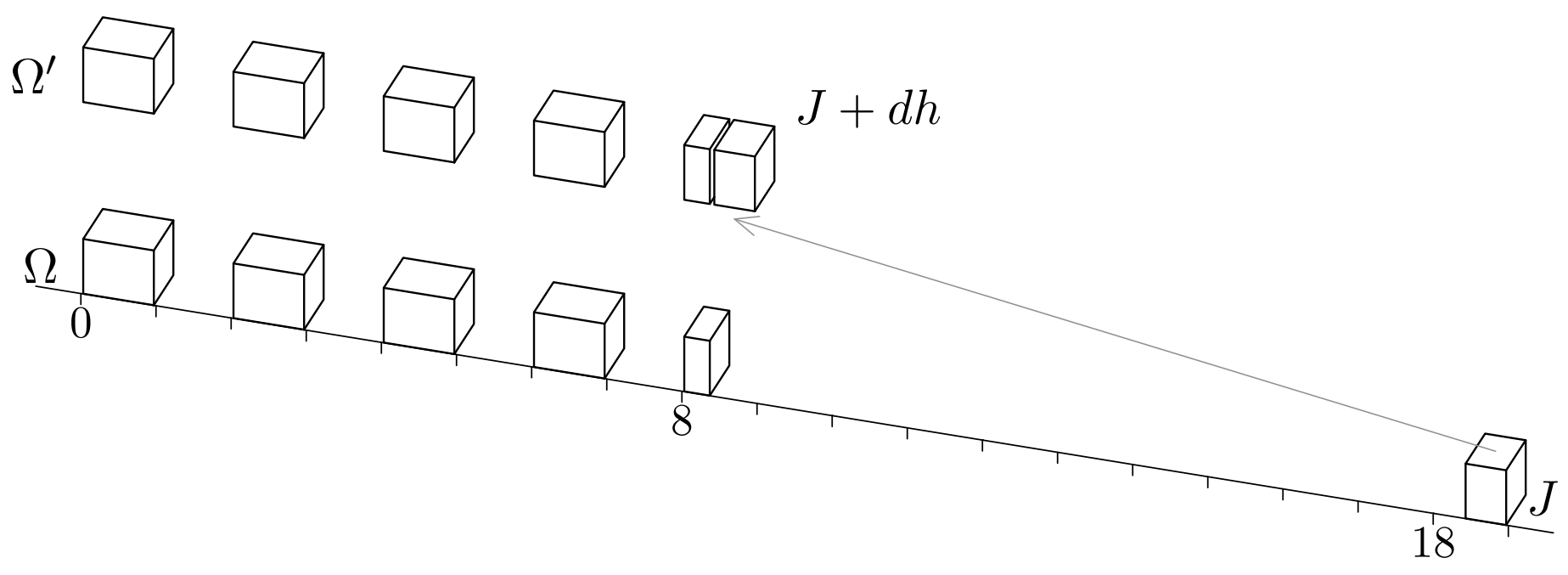}
\caption{The slide of Lemma~\ref{lem:slide} for the set \(\Omega\) of Proposition~\ref{prop:sharp}, drawn schematically with \(a \approx 0.4\), for which \(h = 1\). Below, \(\Omega\), in which \([8, 8 + a)\) and \(J = [18 + a, 19)\) form a rigid pair with slide \(d = -10\). Above, the set \(\Omega'\) obtained by the slide: \(J\) moves to \(J + dh = [8 + a, 9)\), which completes \([8, 8 + a)\) to \([8, 9)\), so that \(\Omega' = \{0, 2, 4, 6, 8\} + [0, 1)\). The gray arrow carries \(J\) to \(J + dh\). In every fiber that contains \(18\), the point \(18\) is replaced by \(8\), which leaves the character sums on \(S - S\) unchanged since \(10(S - S) \subset \Z\).}
\label{fig:slide}
\end{figure}

Figure~\ref{fig:slide} shows a slide in the example of Proposition~\ref{prop:sharp}.

Iterating the slide gives the reduction with which the proof of Theorem~\ref{thm:A} starts.

\begin{proposition}[iterated merging]
\label{prop:mergechain}
Let \(\Omega\) be as in Definition~\ref{def:standing} with \(n \ge 1\) components, let \(\Lambda\) be a spectrum of \(\Omega\), and let \(h\), \(k\), \(\Gamma\), \(S\) be as in Theorem~\ref{thm:periodic}. There is a finite sequence \(\Omega = \Omega_0, \Omega_1, \dots, \Omega_t\) of sets, each of the form of Definition~\ref{def:standing}, of measure \(L\), with \(\Lambda\) a spectrum of each and with the same \(h\), \(k\), \(\Gamma\) and \(S\), in which \(\Omega_{s+1}\) is obtained from \(\Omega_s\) by the slide of Lemma~\ref{lem:slide} along a rigid pair followed by the standing reduction, the numbers of components strictly decrease, and the last set \(\Omega_t\), with \(n_r\) components, satisfies one of the following.
\begin{enumerate}
\item[\rm (i)] \(n_r = 1\). Then \(k = 1\), \(h = L\), and \(\Omega\) tiles \(\R\) by the lattice \(L\Z\).
\item[\rm (ii)] \(n_r \ge 2\) and \(\Omega_t\) is fully mixing. Then every class of \(\Omega_t\) has size at least two, and \(\Omega_t\) has at most \(n_r/2\) classes.
\end{enumerate}
\end{proposition}

\begin{proof}
Start with \(\Omega_0 = \Omega\). If the current set \(\Omega_s\) has at least two components and a rigid pair, Lemma~\ref{lem:slide} followed by the standing reduction produces \(\Omega_{s+1}\) with fewer components, of measure \(L\), with the spectrum \(\Lambda\) and the same \(h\), \(k\), \(\Gamma\), \(S\); rigidity for \(\Omega_{s+1}\) is defined with respect to these same \(h\) and \(S\), so the step can be repeated. Since the number of components decreases at each step, the process stops, at a set \(\Omega_t\) with \(n_r = 1\) components or with \(n_r \ge 2\) components and no rigid pair.

In case (i), \(\Omega_t = [a, b)\) is one interval of length \(L = kh\), and by Proposition~\ref{prop:intervalfibers} every fiber of \(\Omega_t\) is a block of \(\lfloor L/h\rfloor = \lceil L/h \rceil = k\) consecutive integers \(\{m_0, \dots, m_0 + k - 1\}\). By Lemma~\ref{lem:everyfiber} this block is a spectrum of \(S\). If \(k \ge 2\), take \(\sigma \in (S - S) \setminus \{0\}\); then \(\ex{\sigma} \ne 1\), and the vanishing of the geometric sum \(\ex{\sigma m_0}(\ex{\sigma k} - 1)/(\ex{\sigma} - 1)\) forces \(\ex{\sigma k} = 1\). Thus \(k(S - S) \subset \Z\), so \(S\) lies in one coset of \((1/k)\Z/\Z\), and since \(\lvert S\rvert = k\) it is that whole coset, which satisfies \(S + 1/k = S\). Theorem~\ref{thm:periodic} excludes a nontrivial translation symmetry, so \(k = 1\). Then \(h = L\), and Lemma~\ref{lem:everyfiber} for the original \(\Omega\) says \(\sum_{m} 1_\Omega(x + mL) = 1\) \(\forall x \in \R\): the translates \(\Omega + mL\), \(m \in \Z\), partition \(\R\).

In case (ii), a class of size one would be rigid by Lemma~\ref{lem:rigidpair}(c), so every class of \(\Omega_t\) has \(m_c \ge 2\), and \(\sum_c m_c = n_r\) gives at most \(n_r/2\) classes.
\end{proof}

\subsection{Three statements about \texorpdfstring{$k$}{k} points, and the reduction modulo \texorpdfstring{$N$}{N}}

\begin{definition}
\label{def:namedstatements}
Let \(k \ge 1\) be an integer. The following three statements are referred to throughout the paper.
\begin{enumerate}
\item[\(\FL{k}\)] (the fiber lemma) Let \(h > 0\) and let \(F \subset \R\) be a bounded measurable set with \(\sum_{m \in \Z} 1_F(x + mh) = k\) for almost every \(x \in \R\). If \(\Gamma + h^{-1}\Z\) is a spectrum of \(F\) for a set \(\Gamma\) of \(k\) real numbers, pairwise distinct modulo \(h^{-1}\Z\), then \(F\) tiles \(\R\) by translations.
\item[\(\UC{k}\)] (the common complement) Let \(N \ge 1\) and let \(S' \subset \ZN{N}\) be a primitive set of \(k\) elements that has a spectrum in \(\ZN{N}\). Then there is a set \(B \subset \ZN{N}\) with \(A \oplus B = \ZN{N}\) for every spectrum \(A\) of \(S'\).
\item[\(\CP{k}\)] (the consecutive pair) Let \((A, S)\) be a spectral pair with \(A \subset \Z\), \(S \subset \R/\Z\) and \(\lvert A\rvert = \lvert S\rvert = k\). If \(A\) contains two consecutive integers, then \(S + 1/m = S\) for some integer \(m \ge 2\).
\end{enumerate}
The terms spectrum and primitive in \(\UC{k}\) are those of Definition~\ref{def:zeroset} below. In \(\FL{k}\) the multitile hypothesis is automatic: by Theorem~\ref{thm:fiber} a spectrum of the stated form forces \(\lvert F\rvert = kh\) and \(k\) points in almost every fiber.
\end{definition}

\begin{theorem}[{Fu and Song \cite[Theorem~1.2]{fu-song2026}; \cite[Theorem~6.1]{bell2026four}}]
\label{thm:rational}
Let \((A, S)\) be a spectral pair with \(A \subset \Z\) and \(S \subset \R/\Z\). Then \(S - S \subset \Q/\Z\).
\end{theorem}

\begin{lemma}[{\cite[Lemma~5.2]{bell2026four}}]
\label{lem:reduction}
Let \(S \subset \R/\Z\) be a finite set with at least two elements and with \(S - S \subset \Q/\Z\). Let \(N\) be the least common multiple of the denominators in lowest terms of the nonzero elements of \(S - S\), and fix \(s_0 \in S\). Then \(S - S \subset (1/N)\Z/\Z\), and \(S' = N(S - s_0) \subset \ZN{N}\) has \(\lvert S\rvert\) elements, contains \(0\), and satisfies \(\gcd(N, S') = 1\). Reduction modulo \(N\) is injective on every spectrum \(A \subset \Z\) of \(S\), and the image is a spectrum of \(S'\) in \(\ZN{N}\). Conversely, a subset of \(\Z\) with \(\lvert S\rvert\) elements whose reduction modulo \(N\) is a spectrum of \(S'\) is a spectrum of \(S\).
\end{lemma}

\begin{definition}
\label{def:zeroset}
Let \(N \ge 1\) and \(S' \subset \ZN{N}\) be a finite set. Its character sum is \(\widehat{\mu_{S'}}(t) = \sum_{s \in S'} \zroot{N}^{st}\) for \(t \in \ZN{N}\), and its \emph{zero set} is
\[
\zeroset{S'} = \{ t \in \ZN{N} \setminus \{0\} : \widehat{\mu_{S'}}(t) = 0 \}.
\]
A \emph{spectrum} of \(S'\) in \(\ZN{N}\) is a set \(A \subset \ZN{N}\) with \(\lvert A\rvert = \lvert S'\rvert\) and \((A - A) \setminus \{0\} \subset \zeroset{S'}\); equivalently, \((\zroot{N}^{sa})_{s \in S',\, a \in A}\) is a complex Hadamard matrix, since for a square matrix the orthogonality of the columns is equivalent to that of the rows. The set \(S'\) is \emph{periodic} when \(S' + N/p = S'\) for some prime \(p\) dividing \(N\); equivalently, \(S' + t = S'\) for some \(t \ne 0\), because the stabilizer of \(S'\) is then a nontrivial subgroup of \(\ZN{N}\) and contains a subgroup of prime order. The set \(S'\) is \emph{primitive} when \(S' - S'\) generates \(\ZN{N}\). The set \(S'\) of Lemma~\ref{lem:reduction} contains \(0\) and has \(\gcd(N, S') = 1\), so it is primitive; and since every translation symmetry \(t\) of \(S\) lies in \(S - S \subset (1/N)\Z/\Z\), the set \(S\) has a nontrivial translation symmetry iff \(S'\) is periodic.
\end{definition}

The fiber lemma through five cosets, including a common complement for the fibers, is due to Dutkay and Haussermann \cite{dutkay-haussermann2015}. We use the following formulation from \cite[Theorem~5.8]{bell2026four}.

\begin{theorem}
\label{thm:fiberfew}
The fiber lemma \(\FL{k}\) holds for \(k \le 5\). More precisely, let \(h > 0\), \(k \le 5\), and let \(F \subset \R\) be a bounded measurable set with a spectrum \(\Gamma + h^{-1}\Z\), where \(\lvert\Gamma\rvert = k\) and the elements of \(\Gamma\) are distinct modulo \(h^{-1}\Z\); put \(S = h\Gamma\). Then there is a set \(B \subset \Z\), depending only on \(S\), such that \(\sum_{b \in B} 1_F(t - hb) = 1\) for almost every \(t \in \R\). Thus \(F\) tiles \(\R\) by the translation set \(hB\).
\end{theorem}

A common complement for the finite spectra gives the fiber lemma by the following pullback argument; see also \cite{dutkay-haussermann2015} and \cite[Theorem~5.8]{bell2026four}.

\begin{lemma}
\label{lem:ucfl}
Let \(k \ge 2\). If \(\UC{k}\) holds, then \(\FL{k}\) holds, with a translation set \(hB\) in which \(B \subset \Z\) depends only on \(S = h\Gamma\).
\end{lemma}

\begin{proof}
Let \(F\), \(h\), \(\Gamma\) and \(S\) be as in \(\FL{k}\), and let \(A_x = \{m \in \Z : x + mh \in F\}\). By Theorem~\ref{thm:fiber} there is a set \(X_0 \subset [0, h)\) of full measure such that \(A_x\) is a spectrum of \(S\) \(\forall x \in X_0\). Fix one \(x_1 \in X_0\). Theorem~\ref{thm:rational} applied to \((A_{x_1}, S)\) gives \(S - S \subset \Q/\Z\), and Lemma~\ref{lem:reduction} produces \(N\) and the primitive set \(S' = N(S - s_0) \subset \ZN{N}\), which has the spectrum \(A_{x_1} \bmod N\). By \(\UC{k}\) there is \(B' \subset \ZN{N}\) with \(A \oplus B' = \ZN{N}\) for every spectrum \(A\) of \(S'\). Let \(B = \{ b \in \Z : b \bmod N \in B' \}\), which depends only on \(S\).

Let \(x \in X_0\) and \(m \in \Z\). The integers \(b \in B\) with \(m - b \in A_x\) correspond to the elements \(a \in A_x\) with \((m - a) \bmod N \in B'\). Since reduction modulo \(N\) is injective on \(A_x\) by Lemma~\ref{lem:reduction}, and \((A_x \bmod N) \oplus B' = \ZN{N}\), there is exactly one such \(a\). Hence \(A_x \oplus B = \Z\) \(\forall x \in X_0\). Now let \(t = x + hm\) with \(x \in [0, h)\) and \(m \in \Z\). For \(x \in X_0\) and \(b \in B\), the point \(t - hb\) lies in \(F\) iff \(m - b \in A_x\), so \(\sum_{b \in B} 1_F(t - hb)\) is the number of representations of \(m\) in \(A_x \oplus B\), which is \(1\). The exceptional set \(([0, h) \setminus X_0) + h\Z\) is null.
\end{proof}

\section{Effective levels and the verdict rule}
\label{sec:levels}

Throughout this section \(\Omega = \bigsqcup_{j=1}^{n}[a_j,b_j)\) is a spectral union of \(n\) pairwise disjoint non-adjacent bounded intervals in the normal form of Definition~\ref{def:standing}, with lengths \(\ell_j = b_j - a_j > 0\) and measure \(L = \sum_j \ell_j\), and \(\Lambda\) is a spectrum of \(\Omega\) with minimal period \(h^{-1}\), \(\Lambda = \Gamma + h^{-1}\Z\), \(k = \lvert\Gamma\rvert\) and \(S = h\Gamma \subset \R/\Z\) as in Theorem~\ref{thm:periodic}. By Theorem~\ref{thm:fiber}, \(L = kh\). The periodicity theorem gives no bound on \(k\). The bound comes from the lengths, through the boundary unitary of Shi, and this section explains the mechanism for any number of intervals. The lengths enter only through the coincidences among their subset sums \(\ell_T = \sum_{j \in T}\ell_j\), \(T \subset \nset\), and the outcome is an integer \(\kzero(\ell)\), determined by these coincidences alone, which \(k\) divides. Section~\ref{sec:fivelengths} classifies the effective ratios for five lengths and bounds \(\kzero(\ell)\) for every length vector with at most four rational relations.

\subsection{The boundary unitary}

We recall the three statements of \cite{bell2026four} that the section rests on. They are stated and proved there for every number \(n\) of intervals.

\begin{lemma}[{\cite[Lemma~3.1]{bell2026four}}]
\label{lem:boundaryunitary}
Let \(D\) be the self-adjoint operator on \(L^2(\Omega)\) which is diagonal in the orthonormal basis \(\{e_\lambda : \lambda \in \Lambda\}\) with eigenvalue \(\lambda\) at \(e_\lambda\). There is a unique unitary matrix \(U \in U(n)\) such that the domain of \(D\) consists of the functions \(f\) on \(\Omega\) which are absolutely continuous on each closed component with derivative in \(L^2\) and satisfy \(f(b) = Uf(a)\), where \(f(a) = (f(a_j))_j\) and \(f(b) = (f(b_j))_j\); on this domain \(Df = (2\pi i)^{-1}f'\).
\end{lemma}

The matrix \(U\) is the \emph{boundary unitary} of the pair \((\Omega,\Lambda)\). For \(\zeta \in \C\) put \(z_j(\zeta) = \ex{\zeta\ell_j}\), \(D_\ell(\zeta) = \operatorname{diag}(z_1(\zeta),\dots,z_n(\zeta))\), \(M(\zeta) = D_\ell(\zeta) - U\), and \(y_\zeta = (\ex{\zeta a_j})_{j=1}^{n}\). The \emph{secular determinant} is \(F(\zeta) = \det M(\zeta)\).

\begin{proposition}[{\cite[Proposition~3.2]{bell2026four}}]
\label{prop:secular}
The zeros of \(F\) are exactly the points of \(\Lambda\); in particular every zero is real. For \(\lambda \in \Lambda\) the kernel of \(M(\lambda)\) is the line \(\C y_\lambda\), the adjugate of \(M(\lambda)\) is the rank-one matrix \(\kappa(\lambda)\, y_\lambda\bigl(D_\ell(\lambda)y_\lambda\bigr)^*\) for a scalar \(\kappa(\lambda) \ne 0\), and \(F'(\lambda) = 2\pi i L\,\kappa(\lambda)\). Thus every zero of \(F\) is simple.
\end{proposition}

\begin{lemma}[{\cite[Lemma~3.3]{bell2026four}}]
\label{lem:freqwidth}
Let \(G(\zeta) = \sum_{\beta \in \mathcal{B}} g_\beta \ex{\beta\zeta}\) with \(\mathcal{B} \subset \R\) finite and complex coefficients \(g_\beta\). If \(\mathcal{B}\) lies in an interval of length strictly less than \(L\) and \(G\) vanishes on \(\Lambda\), then every \(g_\beta\) is zero.
\end{lemma}

For \(T \subset \nset\) write \(z_T(\zeta) = \ex{\zeta\ell_T} = \prod_{j \in T} z_j(\zeta)\), with \(z_\emptyset = 1\) and \(z_{\nset}(\zeta) = \ex{\zeta L}\). Let \(U[T]\) denote the principal minor of \(U\) on the index set \(T\), with \(U[\emptyset] = 1\), and put
\[
c_T = (-1)^{\,n - \lvert T\rvert}\, U[T^c], \qquad T \subset \nset,
\]
so that \(c_\emptyset = (-1)^n\det U\) and \(c_{\nset} = 1\). For a real number \(\alpha\) let
\[
s_\alpha = \sum_{T \subset \nset,\ \ell_T = \alpha} c_T ,
\]
with \(s_\alpha = 0\) when \(\alpha\) is not a subset sum. Let \(C_j(\zeta)\) be the \((j,j)\) cofactor of \(M(\zeta)\), the determinant of the matrix obtained by deleting row and column \(j\).

\subsection{The cofactor identity in level form}

The next theorem is \cite[Theorem~4.1]{bell2026four}, which is stated there for four intervals. Its proof uses nothing about the number four, and we give it for \(n\) intervals.

\begin{theorem}
\label{thm:cofactor}
The secular determinant and its diagonal cofactors are
\[
F(\zeta) = \sum_{T \subset \nset} c_T\, z_T(\zeta), \qquad z_j(\zeta)\,C_j(\zeta) = \sum_{T \ni j} c_T\, z_T(\zeta),
\]
and
\[
z_j(\zeta)\,C_j(\zeta) = \frac{F'(\zeta)}{2\pi i L}\quad\forall j \in \nset,\ \zeta \in \C.
\]
Equivalently, for every real \(\alpha\),
\[
\sum_{\ell_T = \alpha} c_T\, 1_T = \frac{\alpha}{L}\, s_\alpha\, \one \quad\text{in } \C^n,
\]
that is, the sum \(\sum_{\ell_T = \alpha,\ T \ni j} c_T\) does not depend on \(j\). In particular, if a nonempty proper subset \(T\) is the only subset with length sum \(\ell_T\), then \(c_T = 0\), that is \(U[T^c] = 0\).
\end{theorem}

\begin{proof}
Row \(i\) of \(M(\zeta)\) is \(z_i(\zeta)e_i^{T} - u_i\), where \(u_i\) is the \(i\)-th row of \(U\). The determinant is multilinear in the rows, so \(F(\zeta)\) is the sum, over all \(T \subset \nset\), of the determinant of the matrix whose row \(i\) is \(z_i e_i^T\) for \(i \in T\) and \(-u_i\) for \(i \notin T\). Expanding that determinant along the rows indexed by \(T\), each of which has its only nonzero entry on the diagonal, leaves \(\prod_{i \in T} z_i\) times the determinant of \(-U\) restricted to the rows and columns in \(T^c\), which is \((-1)^{\lvert T^c\rvert}U[T^c]\). The summand is therefore \(c_T z_T(\zeta)\), which proves the first identity. The cofactor \(C_j\) is the secular determinant of the principal submatrix of \(D_\ell(\zeta) - U\) on the index set \(\nset \setminus \{j\}\). The same expansion gives \(C_j = \sum_{T \not\ni j} (-1)^{n-1-\lvert T\rvert}\, U[(T \cup \{j\})^c]\, z_T\), and after multiplication by \(z_j\) the summand indexed by \(T\) is \(c_{T \cup \{j\}}\, z_{T \cup \{j\}}\). Reindexing by \(T' = T \cup \{j\}\) gives the second identity.

Let \(\lambda \in \Lambda\). By Proposition~\ref{prop:secular} the \((j,j)\) entry of the adjugate of \(M(\lambda)\) is
\[
C_j(\lambda) = \kappa(\lambda)\,(y_\lambda)_j\,\overline{z_j(\lambda)\,(y_\lambda)_j} = \kappa(\lambda)\,\overline{z_j(\lambda)},
\]
because \(\lambda\) is real and \(\lvert (y_\lambda)_j\rvert = 1\). Hence \(z_j(\lambda)C_j(\lambda) = \kappa(\lambda) = F'(\lambda)/(2\pi i L)\). The function
\[
E_j(\zeta) = z_j(\zeta)C_j(\zeta) - \frac{F'(\zeta)}{2\pi i L}
\]
is an exponential polynomial which vanishes on \(\Lambda\). Its coefficient at the frequency \(L\) is \(c_{\nset} - (2\pi i L)/(2\pi i L) = 0\), and it has no constant term, since every \(T \ni j\) is nonempty and the derivative removes the constant term of \(F\). For \(n=1\) no frequencies remain, so \(E_j\equiv0\). For \(n\ge2\), every remaining frequency is a sum \(\ell_T\) with \(T\) nonempty and proper, hence lies in the interval \([\min_i \ell_i,\, L - \min_i\ell_i]\), whose length \(L - 2\min_i\ell_i\) is strictly less than \(L\). Lemma~\ref{lem:freqwidth} gives \(E_j \equiv 0\), which is the third identity.

The coefficient of \(\ex{\alpha\zeta}\) in \(z_jC_j\) is the \(j\)-th coordinate of \(\sum_{\ell_T = \alpha} c_T 1_T\), and in \(F'/(2\pi i L)\) it is \((2\pi i\alpha/2\pi i L)\, s_\alpha = (\alpha/L)s_\alpha\). Exponentials with distinct real frequencies are linearly independent, so the two coefficients agree for every \(\alpha\) and every \(j\), which is the displayed vector identity. Finally, if \(T\) is nonempty, proper, and alone at its level, the identity reads \(c_T 1_T = (\ell_T/L)c_T\one\); the vector \(1_T\) is not a multiple of \(\one\), so \(c_T = 0\).
\end{proof}

Grouping the first identity by levels gives
\begin{equation}
\label{eq:Flevels}
F(\zeta) = \sum_{\alpha} s_\alpha\, \ex{\alpha\zeta},
\end{equation}
a finite sum over the subset sums \(\alpha\) of \(\ell\), with \(s_0 = (-1)^n\det U \ne 0\) and \(s_L = 1\). The sums \(s_\alpha\) depend on the unknown unitary \(U\). The theorem shows that an interior level \(0<\alpha<L\) can only survive in \eqref{eq:Flevels}, in the sense that \(s_\alpha \ne 0\), when the indicators of the subsets at that level span a space containing \(\one\). This is a condition on the lengths alone, and we give it a name.

\subsection{Effective levels}

The definitions are made for an arbitrary positive weight vector, because Section~\ref{sec:fivelengths} applies them to auxiliary vectors with fewer than \(n\) entries.

\begin{definition}
\label{def:effective}
Let \(E\) be a finite set and \(w \in \R^E_{>0}\) a vector of positive weights, with total weight \(L(w) = \sum_{j \in E} w_j\); for \(T \subset E\) let \(w(T) = \sum_{j \in T} w_j\) and let \(1_T \in \{0,1\}^E\) be the indicator of \(T\). A \emph{level} of \(w\) is a number \(\alpha \in (0, L(w))\) which is a subset sum \(w(T)\); its \emph{slice} is the family \(\slice{\alpha}(w) = \{T \subset E : w(T) = \alpha\}\) of nonempty proper subsets. A family \(\mathcal{T}\) of subsets of \(E\) with a common weight \(w(T) = \alpha\) is called \emph{equilevel}; it is \emph{effective} when \(1_E \in \operatorname{span}_{\Q}\{1_T : T \in \mathcal{T}\}\), and \emph{minimal effective} when it is effective and no proper subfamily is effective. The level \(\alpha\) is \emph{effective} when its slice is effective. The \emph{pinned ratio} of an effective level is \(r_\alpha = \alpha/L(w)\). We write \(R(w)\) for the set of pinned ratios of \(w\), \(D(w)\) for the set of denominators of the members of \(R(w)\) written in lowest terms, and
\[
\kzero(w) = \operatorname{lcm} D(w),
\]
with \(\kzero(w) = 1\) when \(R(w)\) is empty. For the length vector \(\ell\) of \(\Omega\) we write \(\slice{\alpha}\), \(\Rat\), \(\Den\) and \(\kzero(\ell)\).
\end{definition}

Since the indicators are rational vectors, the rational span and the complex span of a family of indicators have the same dimension, and a rational vector lies in one exactly when it lies in the other. The next lemma collects the elementary properties of minimal effective families.

\begin{lemma}
\label{lem:minimalfamily}
Let \(\mathcal{T}\) be a minimal effective equilevel family on \((E,w)\) at the level \(\alpha\). Then the indicators \(1_T\), \(T \in \mathcal{T}\), are linearly independent, so \(\lvert\mathcal{T}\rvert \le \lvert E\rvert\); the representation \(1_E = \sum_{T \in \mathcal{T}} \gamma_T 1_T\) is unique, and every \(\gamma_T\) is a nonzero rational number; every element of \(E\) lies in some member of \(\mathcal{T}\); and \(\alpha \sum_T \gamma_T = L(w)\).
\end{lemma}

\begin{proof}
Fix a representation \(1_E = \sum_T \gamma_T 1_T\). If \(\sum_T \delta_T 1_T = 0\) with some \(\delta_{T_0} \ne 0\), then \(1_E = \sum_T (\gamma_T + s\delta_T)1_T\) for every scalar \(s\), and the choice \(s = -\gamma_{T_0}/\delta_{T_0}\) expresses \(1_E\) through the family \(\mathcal{T} \setminus \{T_0\}\), against minimality. So the indicators are independent, the representation is unique, and it is rational because it is the unique solution of a linear system with rational coefficients. A zero coefficient would again express \(1_E\) through a proper subfamily. An element of \(E\) in no member of \(\mathcal{T}\) would have coordinate \(0\) in \(\sum_T \gamma_T 1_T = 1_E\). Pairing the representation with \(w\) gives \(L(w) = \sum_T \gamma_T w(T) = \alpha\sum_T\gamma_T\).
\end{proof}

\begin{proposition}
\label{prop:pinned}
Let \(w \in \R^E_{>0}\).
\begin{enumerate}
\item[\rm (i)] For the length vector \(\ell\) of \(\Omega\): if \(0 < \alpha < L\) and \(s_\alpha \ne 0\), then \(\alpha\) is an effective level of \(\ell\).
\item[\rm (ii)] If \(\alpha\) is an effective level of \(w\) and \(1_E = \sum_{T \in \slice{\alpha}(w)} \gamma_T 1_T\) with rational \(\gamma_T\), then \(\sum_T \gamma_T \ne 0\) and
\[
r_\alpha = \frac{\alpha}{L(w)} = \frac{1}{\sum_T \gamma_T} \in \Q \cap (0,1).
\]
In particular \(R(w)\) is a finite set of rational numbers.
\item[\rm (iii)] Let \(\alpha\) be an effective level of \(w\), and let \(w' \in \R^{E'}_{>0}\) be a second weight vector together with a map \(T \mapsto T'\) from \(\slice{\alpha}(w)\) to subsets of \(E'\) such that \(1_{T'}\) is the image of \(1_T\) under a fixed linear map \(\varphi\colon \Q^E \to \Q^{E'}\) with \(\varphi(1_E) = 1_{E'}\), and such that all the sets \(T'\) have one common weight \(\alpha' = w'(T')\) with \(0 < \alpha' < L(w')\). Then \(\alpha'\) is an effective level of \(w'\) and \(\alpha'/L(w') = \alpha/L(w)\).
\item[\rm (iv)] If \(\alpha\) is an effective level of \(w\), then so is \(L(w) - \alpha\), with pinned ratio \(1 - r_\alpha\).
\end{enumerate}
\end{proposition}

\begin{proof}
(i) By Theorem~\ref{thm:cofactor}, \(\one = \frac{L}{\alpha s_\alpha}\sum_{T \in \slice{\alpha}} c_T 1_T\) lies in the complex span of the indicators of the slice, hence in their rational span.

(ii) Pairing the representation with \(w\) gives \(\alpha\sum_T\gamma_T = L(w) > 0\), so \(\sum_T \gamma_T = L(w)/\alpha \ne 0\) and \(r_\alpha = 1/\sum_T\gamma_T\) is rational; it lies in \((0,1)\) because \(0 < \alpha < L(w)\). There are finitely many levels.

(iii) Let \(1_E = \sum_T \gamma_T 1_T\) be a rational representation. Applying \(\varphi\) gives \(1_{E'} = \sum_T \gamma_T 1_{T'}\), and every \(T'\) lies in \(\slice{\alpha'}(w')\), so \(\alpha'\) is effective for \(w'\). Pairing the two representations with \(w\) and with \(w'\) gives \(\alpha\sum_T\gamma_T = L(w)\) and \(\alpha'\sum_T\gamma_T = L(w')\), and the ratios agree.

(iv) The map \(T \mapsto E \setminus T\) is a bijection from the slice \(\slice{\alpha}(w)\) onto the slice \(\slice{L(w)-\alpha}(w)\). If \(1_E = \sum_T \gamma_T 1_T\), then \(\sum_T \gamma_T 1_{E \setminus T} = (\sum_T\gamma_T - 1)1_E\), and \(\sum_T \gamma_T = 1/r_\alpha \ne 1\). So \(1_E\) lies in the span of the complementary slice, with the representation \(\sum_T \gamma_T 1_{E\setminus T}/(1/r_\alpha - 1)\), whose coefficient sum is \(1/(1 - r_\alpha)\).
\end{proof}

Part (ii) says that the ratio \(\alpha/L\) of an effective level is determined by the slice as a family of subsets, independently of the particular lengths realizing it: this is why we call it pinned. Part (iii) is the form in which this is used in Section~\ref{sec:fivelengths}: when the lengths are changed in a way that keeps the members of a slice at one common level, the effective level survives with the same ratio.

\subsection{Support and period}

The next proposition is \cite[Proposition~4.3]{bell2026four}. There it is formulated for the individual coefficients \(c_T\), requiring \(\ell_T \in (L/k_0)\Z\) whenever \(c_T \ne 0\). The proof uses only the consequence that every frequency of \eqref{eq:Flevels} with a nonzero coefficient lies in \((L/k_0)\Z\), as assumed here; the converse is already formulated for the sums \(s_\alpha\). The proof therefore applies without change, and we do not repeat it.

\begin{proposition}
\label{prop:periodsupport}
Let \(k_0\) be a positive integer. If \(s_\alpha = 0\) \(\forall \alpha \notin (L/k_0)\Z\), then \(F(\zeta) = q\bigl(\ex{(L/k_0)\zeta}\bigr)\) for a monic polynomial \(q\) of degree \(k_0\) with \(q(0) \ne 0\) and with \(k_0\) distinct unimodular roots, and \(k_0/L\) is a period of \(\Lambda\). Conversely, if \(p > 0\) is a period of \(\Lambda\), then \(\ex{p\alpha} = 1\) \(\forall \alpha\) with \(s_\alpha \ne 0\); in particular \(pL \in \Z\).
\end{proposition}

\subsection{The verdict rule}

\begin{theorem}
\label{thm:verdict}
Let \(\Omega\) be a spectral union of \(n\) pairwise disjoint non-adjacent bounded intervals with length vector \(\ell\), and let \(\Lambda\), \(h\), \(k\) be as above. Then \(\kzero(\ell)/L\) is a period of \(\Lambda\), and
\[
k \mid \kzero(\ell).
\]
More precisely, \(k\) is the least positive integer \(m\) such that \(s_\alpha = 0\) \(\forall \alpha \notin (L/m)\Z\), and every such \(m\) is a multiple of \(k\). If \(\ell\) has no effective level, then \(k = 1\) and \(\Lambda = \lambda_0 + L^{-1}\Z\) for some real \(\lambda_0\).
\end{theorem}

\begin{proof}
Let \(0 < \alpha < L\) with \(s_\alpha \ne 0\). By Proposition~\ref{prop:pinned}(i), \(\alpha\) is effective, and \(\alpha = r_\alpha L\) with \(r_\alpha\) a rational number whose denominator divides \(\kzero(\ell)\); hence \(\alpha \in (L/\kzero(\ell))\Z\). The levels \(0\) and \(L\) lie there as well. Proposition~\ref{prop:periodsupport} with \(k_0 = \kzero(\ell)\) shows that \(\kzero(\ell)/L\) is a period of \(\Lambda\). By Theorem~\ref{thm:periodic} the periods of \(\Lambda\) form the group \(h^{-1}\Z\), and \(h^{-1} = k/L\) by Theorem~\ref{thm:fiber}. So \(\kzero(\ell)/L \in (k/L)\Z\), that is \(k \mid \kzero(\ell)\).

If \(m\) is a positive integer with \(s_\alpha = 0\) \(\forall \alpha \notin (L/m)\Z\), the same argument gives the period \(m/L\), so \(m/L \in (k/L)\Z\) and \(k \mid m\). Conversely \(k/L\) is a period, and the converse part of Proposition~\ref{prop:periodsupport} gives \(\ex{k\alpha/L} = 1\), that is \(\alpha \in (L/k)\Z\), \(\forall \alpha\) with \(s_\alpha \ne 0\). So \(m = k\) has the property, and it is the least such \(m\).

If \(\ell\) has no effective level, then \(s_\alpha = 0\) for \(0 < \alpha < L\) by Proposition~\ref{prop:pinned}(i), \(F(\zeta) = \ex{L\zeta} + (-1)^n\det U\), and the zero set of \(F\) is one coset \(\lambda_0 + L^{-1}\Z\). Then \(\kzero(\ell) = 1\) and \(k = 1\).
\end{proof}

The theorem is the \emph{verdict rule}: the integer \(\kzero(\ell)\) is read off from the coincidences among the subset sums of the lengths, with no knowledge of \(U\) or of \(\Lambda\), and it bounds the number \(k\) of cosets of the spectrum. The periodicity theorem enters the proof only through Lemma~\ref{lem:freqwidth} and through the identification of the minimal period; the period \(\kzero(\ell)/L\) itself is produced by Proposition~\ref{prop:periodsupport}.

\subsection{A bound in terms of the number of intervals}

The denominator of a pinned ratio is bounded by a quantity depending on \(n\) alone, because a pinned ratio is a quotient of two determinants of \(0\)-\(1\) matrices. Put
\[
d(n) = \Bigl\lfloor \frac{(n+2)^{(n+2)/2}}{2^{\,n+1}} \Bigr\rfloor, \qquad \Kn = \operatorname{lcm}(1,2,\dots,d(n)).
\]

\begin{theorem}
\label{thm:Kn}
Let \(\ell \in \R^n_{>0}\) and let \(\alpha\) be an effective level of \(\ell\). Then the denominator of \(r_\alpha\) in lowest terms divides the absolute value of a nonzero determinant of a \(0\)-\(1\) matrix of order \(n+1\), and is therefore at most \(d(n)\). Consequently \(\kzero(\ell)\) divides \(\Kn\), and for every spectral union of \(n\) pairwise disjoint non-adjacent bounded intervals \(k \le \kzero(\ell) \le \Kn\).
\end{theorem}

\begin{proof}
The slice \(\slice{\alpha}\) is finite and effective, so it contains a minimal effective subfamily \(\mathcal{T}\). By Lemma~\ref{lem:minimalfamily} the \(t = \lvert\mathcal{T}\rvert\) indicators \(1_T\), \(T \in \mathcal{T}\), are linearly independent in \(\Q^n\), so \(t \le n\), and \(\one = \sum_{T \in \mathcal{T}} \gamma_T 1_T\) with a unique rational \(\gamma\), \(r_\alpha = 1/\sum_T \gamma_T\) by Proposition~\ref{prop:pinned}(ii). Let \(N\) be the \(t \times n\) matrix with rows \(1_T\). It has rank \(t\), so some \(t\) of its columns, indexed by a set \(J \subset \nset\), form an invertible \(t \times t\) matrix \(M\) with entries in \(\{0,1\}\). The representation reads \(N^T\gamma = \one_n\); restricting to the coordinates in \(J\) gives \(M^T\gamma = \one_t\), so \(\gamma = M^{-T}\one_t\) and
\[
\sum_T \gamma_T = \one_t^T M^{-T}\one_t = \one_t^T M^{-1}\one_t .
\]
Consider the \(0\)-\(1\) matrix of order \(t+1\)
\[
B = \begin{pmatrix} M & \one_t \\ \one_t^T & 0 \end{pmatrix}.
\]
The Schur complement formula gives \(\det B = \det M \cdot (0 - \one_t^T M^{-1}\one_t) = -\det M \cdot \sum_T \gamma_T\). Since \(\sum_T\gamma_T = L/\alpha > 0\), \(\det B \ne 0\), and
\[
r_\alpha = \frac{1}{\sum_T\gamma_T} = -\frac{\det M}{\det B}
\]
is a quotient of two nonzero integers, so its denominator in lowest terms divides \(\lvert\det B\rvert\). The matrix \(B' = \operatorname{diag}(B, I_{n-t})\) is a \(0\)-\(1\) matrix of order \(n+1\) with \(\det B' = \det B\).

It remains to bound the determinant of a \(0\)-\(1\) matrix \(B'\) of order \(s = n+1\). Let \(J_s\) be the all-ones matrix and consider the matrix of order \(s+1\) with entries \(\pm 1\)
\[
\widetilde{B} = \begin{pmatrix} 1 & \one_s^T \\ \one_s & J_s - 2B' \end{pmatrix}.
\]
Subtracting the first row from each of the other rows replaces the lower block row \((\one_s,\ J_s - 2B')\) by \((0,\ -2B')\), so \(\det\widetilde{B} = (-2)^s\det B'\). Hadamard's determinant inequality \cite{horn-johnson2013} bounds the determinant of a matrix of order \(s+1\) with entries of modulus one by \((s+1)^{(s+1)/2}\). Hence
\[
\lvert\det B'\rvert = \frac{\lvert\det\widetilde{B}\rvert}{2^s} \le \frac{(s+1)^{(s+1)/2}}{2^s} = \frac{(n+2)^{(n+2)/2}}{2^{n+1}},
\]
and since \(\lvert\det B'\rvert\) is an integer, \(\lvert\det B'\rvert \le d(n)\). Thus every element of \(\Den\) is at most \(d(n)\), so each divides \(\Kn\), and so does their least common multiple \(\kzero(\ell)\). The last assertion is Theorem~\ref{thm:verdict}.
\end{proof}

The values are \(d(1), \dots, d(8) = 1, 2, 3, 6, 14, 32, 76, 195\), so \(K(1) = 1\), \(K(2) = 2\), \(K(3) = 6\), \(K(4) = 60\), \(K(5) = \operatorname{lcm}(1,\dots,14) = 360360\). The bound is far from sharp: \(\kzero(\ell) \le 6\) for four lengths \cite[Theorem~4.4]{bell2026four}, and Section~\ref{sec:fivelengths} shows \(\kzero(\ell) \le 10\) for five. Its role is to make the number of possible values of \(k\) finite for each \(n\), which is what the general reduction of Section~\ref{sec:general} needs.

\subsection{The integer form of the verdict rule}

A finite spectral pair \((A,S)\) with \(A \subset \Z\) is the same thing as a spectral union of unit intervals with a single fiber, and the verdict rule then becomes a statement about the run vector of \(A\). A \emph{run} of a finite set \(A \subset \Z\) is a maximal set of consecutive integers contained in \(A\), and the \emph{run vector} \(d = (d_1,\dots,d_m)\) of \(A\) lists the lengths of its runs, so that \(\sum_i d_i = \lvert A\rvert\). The quantities \(R(d)\), \(D(d)\) and \(\kzero(d)\) are those of Definition~\ref{def:effective} for the weight vector \(d \in \Z^m_{>0}\); its levels are the integers \(\delta\) with \(0 < \delta < \lvert A\rvert\) that are sums of some of the \(d_i\).

\begin{lemma}
\label{lem:runvector}
Let \((A,S)\) be a spectral pair with \(A \subset \Z\), \(S \subset \R/\Z\), \(\lvert A\rvert = \lvert S\rvert = k\), such that \(S\) has no nontrivial translation symmetry, and let \(d\) be the run vector of \(A\). Then \(\kzero(d) = k\). Equivalently, for every prime \(p\) dividing \(k\) there is an effective level of \(d\) not divisible by \(p\); if \(k > 1\), then \(d\) has an effective level.
\end{lemma}

\begin{proof}
Let \(R_1,\dots,R_m\) be the maximal runs of \(A\), of lengths \(d_1,\dots,d_m\), and put \(E = A + [0,1) = \bigsqcup_i [\min R_i, \max R_i + 1)\). The components of \(E\) have lengths \(d_i\); they are pairwise disjoint, and non-adjacent because two adjacent components would come from two runs whose union is a run, against maximality. So \(E\) is a union of \(m\) intervals in the normal form of Definition~\ref{def:standing} with length vector \(d\) and measure \(L = k\). Take \(h = 1\) and let \(\Gamma \subset [0,1)\) be the set of representatives of \(S\), so that \(\Lambda = \Gamma + \Z\) and \(h\Gamma = S\) in the notation of Definition~\ref{def:fiber}. For \(x \in [0,1)\) the fiber of \(E\) at \(x\) is \(\{j \in \Z : x + j \in E\} = A\), because \(x + j \in [a, a+1)\) exactly when \(j = a\) for \(a \in A\). Every fiber is a spectrum of \(S\), so by Theorem~\ref{thm:fiber} the set \(E\) is spectral with spectrum \(\Lambda\). A period \(p\) of \(\Lambda\) satisfies \(S + p = S\) in \(\R/\Z\), so \(p \in \Z\); the minimal period is \(1 = k/L\), and the number of cosets of \(\Lambda\) is \(\lvert\Gamma\rvert = k\). Theorem~\ref{thm:verdict} applied to \(E\) gives \(k \mid \kzero(d)\).

Conversely every level of \(d\) is an integer \(\delta\) with \(0 < \delta < k\), with ratio \(\delta/k\), whose denominator in lowest terms divides \(k\). So \(\kzero(d)\), the least common multiple of these denominators, divides \(k\). Hence \(\kzero(d) = k\).

For the equivalent form, let \(\Delta\) be the set of effective levels of \(d\). For a prime \(p\) and \(\delta \in \Delta\), the exponent of \(p\) in the denominator \(k/\gcd(k,\delta)\) of \(\delta/k\) is \(v_p(k) - \min(v_p(k), v_p(\delta))\), which equals \(v_p(k)\) exactly when \(p \nmid \delta\). Thus \(v_p(\kzero(d)) = v_p(k)\) \(\forall p \mid k\) exactly when every prime divisor of \(k\) fails to divide some member of \(\Delta\); when \(k > 1\) this requires \(\Delta \ne \emptyset\).
\end{proof}

Thus an aperiodic spectrum can occur only when \(\kzero(d)=\sum_i d_i\): the run vector alone can rule out aperiodicity. This restriction is used in Section~\ref{sec:general}.

\section{Five lengths and the specialization lemma}
\label{sec:fivelengths}

By the verdict rule, the number \(k\) of cosets of the spectrum divides \(\kzero(\ell)\), and \(\kzero(\ell)\) is determined by the coincidences among the subset sums of the lengths. We determine the effective ratios for five lengths and then extend the resulting bound by specialization. For four lengths the corresponding classification is given in \cite[Theorem~4.4]{bell2026four}. For five lengths it is carried out here in full, with the following organizing principle: a minimal effective family on at most five elements at a level \(\alpha\) is, after the elements of equal incidence type have been merged, a square invertible \(0\)-\(1\) matrix whose inverse has positive row sums, and such matrices of size at most five can be listed. The list gives seventeen possible ratios and, for the large denominators \(7\), \(8\), \(9\), the length vectors themselves; the incompatibilities among the ratios then give \(\kzero(\ell) \le 10\). The second half of the section shows that a length vector with \(\dimQ\Krel \le 4\) can be specialized to a vector with at most five positive entries without losing a coincidence of subset sums, so that \(\kzero(\ell) \le 10\) holds for every \(n\).

We use the notation of Definition~\ref{def:effective} for an arbitrary finite set \(E\) with positive weights \(w\). Subsets of a set with at most nine elements are written as strings of their elements, so that \(134\) stands for \(\{1,3,4\}\). For a family \(\mathcal{T}\) of subsets of \(E\) and \(j \in E\), the \emph{type} of \(j\) is the subfamily \(S_j = \{T \in \mathcal{T} : j \in T\}\).

\subsection{Contraction and configurations}

\begin{definition}
\label{def:configuration}
Let \(\mathcal{T}\) be a minimal effective equilevel family on \((E,w)\) at the level \(\alpha\). Its \emph{contraction} is obtained by merging the elements of equal type: the new element set is the set \(\overline{E}\) of types that occur, with weights \(\overline{w}(S) = \sum_{j : S_j = S} w_j\), and each member \(T \in \mathcal{T}\) is replaced by \(\overline{T} = \{S \in \overline{E} : T \in S\}\). An \emph{\(m\)-configuration} is an invertible \(m \times m\) matrix \(M\) with entries in \(\{0,1\}\) such that the vector \(w = M^{-1}\one\) has all coordinates positive. Its \emph{rows} are the subsets \(T_1,\dots,T_m\) of \([m]\) whose indicators are the rows of \(M\), its \emph{weight vector} is \(w\), its \emph{mass} is \(L(M) = \one^T M^{-1}\one = \sum_j w_j\), and its \emph{ratio} is \(r(M) = 1/L(M)\). Two configurations are \emph{isomorphic} when they differ by a permutation of the rows and a permutation of the columns; isomorphic configurations have the same ratio and weight vectors differing by a permutation.
\end{definition}

A configuration is thus an equilevel family on \(([m], w)\), with common weight \(1\): \(Mw = \one\) says that every row has weight \(1\). Since \(M\) is invertible, its rows are distinct and nonempty, the indicators are linearly independent, and \(\one = \sum_i \gamma_i 1_{T_i}\) has the unique solution \(\gamma = M^{-T}\one\); so a configuration is an effective family, and it is minimal effective exactly when no coordinate of \(M^{-T}\one\) vanishes. For \(m \ge 2\), its ratio is the pinned ratio \(1/L(w)\) of the level \(1\). The weight vector is determined by the \(0\)-\(1\) pattern alone: the weights of a configuration are pinned. Two further facts are used throughout. The rows of a configuration form an antichain, since \(T_i \subsetneq T_k\) would give \(w(T_k \setminus T_i) = 0\). For \(m \ge 2\) no row equals \([m]\), because it would contain every other row, and therefore \(L(M) = w([m]) > w(T_1) = 1\) and \(r(M) < 1\).

\begin{lemma}
\label{lem:contraction}
Let \(\mathcal{T}\) be a minimal effective equilevel family on \((E,w)\) at the level \(\alpha\), with coefficients \(\gamma\) as in Lemma~\ref{lem:minimalfamily}. Its contraction \(\overline{\mathcal{T}} = \{\overline{T} : T \in \mathcal{T}\}\) is a minimal effective equilevel family on \((\overline{E},\overline{w})\) at the same level \(\alpha\), with the same coefficients \(\gamma\) and the same ratio, whose elements have pairwise distinct types, and \(\lvert\overline{E}\rvert \le \lvert E\rvert\).
\end{lemma}

\begin{proof}
Let \(\varphi\colon \Q^{\overline{E}} \to \Q^{E}\) be the linear map which repeats the coordinate \(S\) at every \(j\) with \(S_j = S\). It is injective because every type in \(\overline{E}\) occurs, \(\varphi(1_{\overline{T}})_j = [T \in S_j] = [j \in T]\), so \(\varphi(1_{\overline{T}}) = 1_T\), and \(\varphi(1_{\overline{E}}) = 1_E\). Hence the linear relations among the vectors \(1_{\overline{T}}\) and \(1_{\overline{E}}\) are exactly those among the \(1_T\) and \(1_E\): the contracted family is effective with the same coefficients, its indicators are independent (so the \(\overline{T}\) are distinct and nonempty), and every effective subfamily of \(\overline{\mathcal{T}}\) corresponds to an effective subfamily of \(\mathcal{T}\), which gives minimality. The weight of \(\overline{T}\) is \(\sum_{S \ni T}\sum_{j : S_j = S} w_j = \sum_{j \in T} w_j = \alpha\). The type of \(S \in \overline{E}\) in the contracted family is \(\{\overline{T} : S \in \overline{T}\} = \{\overline{T} : T \in S\}\), which determines \(S\); so types are pairwise distinct. The ratio \(\alpha/\overline{w}(\overline{E}) = \alpha/w(E)\) is unchanged.
\end{proof}

\begin{lemma}
\label{lem:configbasic}
Let \(M\) be an \(m\)-configuration with \(m \ge 2\), weight vector \(w\) and mass \(L\).
\begin{enumerate}
\item[\rm (a)] If a row is a singleton, \(T_i = \{j\}\), then \(w_j = 1\) and \(j\) lies in no other row, so column \(j\) of \(M\) is the unit vector \(e_i\).
\item[\rm (b)] The \(m\)-configurations having a column with a single \(1\) are, up to isomorphism, exactly the matrices
\[
M = \begin{pmatrix} 1 & 1_R^T \\ 0 & N \end{pmatrix}
\]
with \(N\) an \((m-1)\)-configuration with weight vector \(w'\) and \(R\) a subset of the column set of \(N\) with \(w'(R) < 1\). The weight vector of \(M\) is \((1 - w'(R), w')\) and \(L(M) = L(N) + 1 - w'(R)\). We call \(M\) the extension of \(N\) by \(R\).
\item[\rm (c)] Let \(J\) be the all-ones matrix. Then \(J - M\) is an \(m\)-configuration, the complement of \(M\); its rows are the complements \([m] \setminus T_i\), its weight vector is \(w/(L-1)\), its mass is \(L/(L-1)\) and its ratio is \(1 - r(M)\). A column of weight \(c\) (that is, containing \(c\) ones) in \(M\) has weight \(m - c\) in \(J - M\).
\item[\rm (d)] Every column of \(M\) has weight in \([1, m-1]\). Exactly one of the following holds: {\rm (I)} some column has weight \(1\); {\rm (II)} no column has weight \(1\) and some column has weight \(m-1\); {\rm (III)} every column weight lies in \([2, m-2]\). The complement of a configuration of class {\rm (II)} is of class {\rm (I)}, the class {\rm (III)} is closed under complementation, and in class {\rm (III)} every row has size in \([2, m-2]\).
\end{enumerate}
\end{lemma}

\begin{proof}
(a) \(w_j = w(T_i) = 1\). A row \(T_k \ne T_i\) containing \(j\) contains a second element, of positive weight, so \(w(T_k) > 1\).

(b) Let column \(j\) of \(M\) be the unit vector \(e_i\). After permuting rows and columns so that \(i = j = 1\), the matrix has the displayed block form, where \(N\) is the submatrix on the other rows and columns and \(R\) is the set of columns in which row \(1\) has further ones. Then \(\det M = \det N \ne 0\). The rows \(2,\dots,m\) have a zero in column \(1\), so \(Nw' = \one\) with \(w' = (w_2,\dots,w_m) > 0\), which makes \(N\) an \((m-1)\)-configuration, and row \(1\) gives \(w_1 = 1 - w'(R) > 0\). Conversely the displayed matrix is invertible with \(M(1 - w'(R), w')^T = \one\) and positive weights.

(c) \((J - M)w = L\one - \one = (L-1)\one\), and \(L > 1\). If \((J - M)x = 0\), then \(Mx = \sigma\one\) with \(\sigma = \sum_j x_j\); if \(\sigma = 0\) then \(x = 0\), and otherwise \(x = \sigma w\) and summing the coordinates gives \(\sigma = \sigma L\), so \(L = 1\), which is excluded. Hence \(J - M\) is invertible and \((J-M)^{-1}\one = w/(L-1) > 0\). Its mass is \(L/(L-1)\) and its ratio is \((L-1)/L\).

(d) A zero column makes \(M\) singular, and a column of weight \(m\) would be a zero column of the configuration \(J - M\). The three classes are exhaustive and disjoint by definition. If \(M\) is of class (II), some column of \(J - M\) has weight \(1\). Complementation maps weights in \([2, m-2]\) to weights in \([2, m-2]\). In class (III) a singleton row would give a column of weight \(1\) by (a), and a row of size at least \(m-1\) would give a row of size at most \(1\) in \(J - M\), which is also of class (III); an empty row is impossible in an invertible matrix.
\end{proof}

\subsection{The configurations of size at most four}

\begin{proposition}
\label{prop:smallconfigs}
Up to isomorphism, the configurations of size \(m \le 4\) are the following, listed with their rows, weight vector, mass and ratio.
\[
\setlength{\arraycolsep}{2pt}
\begin{array}{lllll}
\text{size } 1: & [1] & 1 & w = (1) & L = 1,\ r = 1 \\
\text{size } 2: & I_2 & 1,\ 2 & w = (1,1) & L = 2,\ r = 1/2 \\
\text{size } 3: & I_3 & 1,\ 2,\ 3 & w = (1,1,1) & L = 3,\ r = 1/3 \\
       & J - I_3 & 23,\ 13,\ 12 & w = \tfrac12(1,1,1) & L = 3/2,\ r = 2/3 \\
\text{size } 4: & I_4 & 1,\ 2,\ 3,\ 4 & w = (1,1,1,1) & L = 4,\ r = 1/4 \\
       & J - I_4 & 234,\ 134,\ 124,\ 123 & w = \tfrac13(1,1,1,1) & L = 4/3,\ r = 3/4 \\
       & A & 1,\ 34,\ 24,\ 23 & w = \tfrac12(2,1,1,1) & L = 5/2,\ r = 2/5 \\
       & A' = J - A & 234,\ 12,\ 13,\ 14 & w = \tfrac13(2,1,1,1) & L = 5/3,\ r = 3/5 \\
       & D & 12,\ 34,\ 24,\ 23 & w = \tfrac12(1,1,1,1) & L = 2,\ r = 1/2
\end{array}
\]
The ratios of configurations of sizes \(2\) through \(4\) are therefore \(1/2\), \(1/3\), \(2/3\), \(1/4\), \(3/4\), \(2/5\) and \(3/5\). The only configurations of size at most \(4\) with ratio \(1/3\), \(1/4\), \(2/5\) are \(I_3\), \(I_4\), \(A\), and the only ones with ratio \(2/3\), \(3/4\), \(3/5\) are their complements. The configuration \(D\) is not minimal effective, since \(\one = 1_{12} + 1_{34}\).
\end{proposition}

The proof enumerates the extensions and complements of the configurations of each size, and is given in Appendix~\ref{app:cases}.

\subsection{Minimal families on at most five elements are square}

The next proposition is the reason the five-length case can be done in full: every minimal effective family on at most five elements contracts to a configuration.

\begin{proposition}
\label{prop:fewmembers}
Let \(E\) be a set with \(\lvert E\rvert \le 5\), \(w \in \R^E_{>0}\), \(\alpha\) an effective level of \(w\), and \(\mathcal{T} \subset \slice{\alpha}(w)\) a minimal effective subfamily with \(t\) members. Then \(2 \le t \le 5\), the elements of \(E\) fall into exactly \(t\) types, which partition \(E\) into classes \(C_1,\dots,C_t\), the contraction of \(\mathcal{T}\) is a \(t\)-configuration \(M\) with \(w(C_i) = \alpha\,(M^{-1}\one)_i\), and \(r_\alpha = r(M)\). Conversely, let \(t \ge 2\), let \(M\) be a \(t\)-configuration with weight vector \(w\), and let \(v \in \R^{E'}_{>0}\) be obtained by splitting each weight \(w_i\) into one or more positive parts. Then \(1\) is an effective level of \(v\) with pinned ratio \(r(M)\).
\end{proposition}

The proof, a case analysis on the number \(t\) of members and the singleton types, is given in Appendix~\ref{app:cases}.

\subsection{The configurations of size five}

\begin{theorem}
\label{thm:fiveconfigs}
Every \(5\)-configuration is isomorphic to a configuration in the following lists, which need not be disjoint.
\begin{enumerate}
\item[\rm (I)] The extensions, by Lemma~{\rm \ref{lem:configbasic}(b)}, of the five \(4\)-configurations \(N\) on the elements \(2,3,4,5\) by the sets \(R\) with \(w'(R) < 1\), the new element \(1\) having weight \(1 - w'(R)\):
\[
\begin{array}{lllllll}
 & N & R & w & L & r \\
E_1 & I_4 & \emptyset & (1,1,1,1,1) & 5 & 1/5 \\
E_2 & J - I_4 & \emptyset & \tfrac13(3,1,1,1,1) & 7/3 & 3/7 \\
E_3 & J - I_4 & \text{one element} & \tfrac13(2,1,1,1,1) & 2 & 1/2 \\
E_4 & J - I_4 & \text{two elements} & \tfrac13(1,1,1,1,1) & 5/3 & 3/5 \\
E_5 & A & \emptyset & \tfrac12(2,2,1,1,1) & 7/2 & 2/7 \\
E_6 & A & \text{one light element} & \tfrac12(1,2,1,1,1) & 3 & 1/3 \\
E_7 & A' & \emptyset & \tfrac13(3,2,1,1,1) & 8/3 & 3/8 \\
E_8 & A' & \text{one light element} & \tfrac13(2,2,1,1,1) & 7/3 & 3/7 \\
E_9 & A' & \text{two light elements} & \tfrac13(1,2,1,1,1) & 2 & 1/2 \\
E_{10} & A' & \text{the heavy element} & \tfrac13(1,2,1,1,1) & 2 & 1/2 \\
E_{11} & D & \emptyset & \tfrac12(2,1,1,1,1) & 3 & 1/3 \\
E_{12} & D & \text{one element} & \tfrac12(1,1,1,1,1) & 5/2 & 2/5
\end{array}
\]
Here the light elements of \(A\) are the three of weight \(\frac12\), the heavy element of \(A'\) is the one of weight \(\frac23\) and its light elements are the three of weight \(\frac13\); \(E_{12}\) includes the extensions by each of the four singleton subsets of \(D\), which need not be isomorphic.
\item[\rm (II)] The complements of the configurations in {\rm (I)}. The weight vectors are proportional, and the ratios are
\[
\begin{gathered}
4/5;\quad 4/7,1/2,2/5;\quad 5/7,2/3;\\
5/8,4/7;\quad 1/2,1/2;\quad 2/3,3/5.
\end{gathered}
\]
\item[\rm (III)] Three configurations of class {\rm (III)} and their complements. With the five elements named \(a,b,c,d,e\) and the rows written as strings of elements:
\[
\begin{array}{llll}
C_1: & ae,\ ab,\ bc,\ cd,\ de & w = \tfrac12(1,1,1,1,1) & r = 2/5 \\
C_2: & abe,\ ce,\ de,\ ac,\ bd & w = \tfrac13(1,1,2,2,1) & r = 3/7 \\
C_3: & ade,\ bde,\ cd,\ ce,\ ab & w = \tfrac14(2,2,3,1,1) & r = 4/9
\end{array}
\]
and the complements have the ratios \(3/5\), \(4/7\), \(5/9\).
\end{enumerate}
The set of ratios of \(5\)-configurations is
\[
R_5 = \{1/5, 4/5, 2/7, 5/7, 3/7, 4/7, 1/3, 2/3, 3/8, 5/8, 2/5, 3/5, 4/9, 5/9, 1/2\}.
\]
The \(5\)-configurations whose ratio has denominator \(7\), \(8\) or \(9\), or equals \(1/5\) or \(4/5\), have, up to a scalar, the following weight vectors: \((3,1,1,1,1)\) for \(E_2\), \((2,2,1,1,1)\) for \(E_5\) and \(E_8\), \((1,1,2,2,1)\) for \(C_2\), and their complements, at the denominator \(7\); \((3,2,1,1,1)\) for \(E_7\) and its complement, at the denominator \(8\); \((2,2,3,1,1)\) for \(C_3\) and its complement, at the denominator \(9\); \((1,1,1,1,1)\) for \(E_1\) and its complement. The \(5\)-configurations with ratio \(1/3\) are \(E_6\) and \(E_{11}\), with weight vector proportional to \((1,2,1,1,1)\); none has ratio \(1/4\) or \(3/4\); every one with ratio \(2/5\) has all weights equal.
\end{theorem}

The proof, which lists the extensions of the configurations of size four and the configurations of class {\rm (III)} by the number of columns of weight three, is given in Appendix~\ref{app:cases}.

\subsection{Pinned ratios and length vectors}

Let
\[
\mathcal{R} = \bigl\{\tfrac12,\ \tfrac13,\ \tfrac23,\ \tfrac14,\ \tfrac34,\ \tfrac15,\ \tfrac25,\ \tfrac35,\ \tfrac45,\ \tfrac27,\ \tfrac37,\ \tfrac47,\ \tfrac57,\ \tfrac38,\ \tfrac58,\ \tfrac49,\ \tfrac59\bigr\}
\]
be the union of the ratios in \((0,1)\) listed in Proposition~\ref{prop:smallconfigs} and Theorem~\ref{thm:fiveconfigs}.

\begin{theorem}
\label{thm:fiveratios}
Let \(v \in \R^5_{>0}\) with \(L = L(v)\).
\begin{enumerate}
\item[\rm (a)] Every pinned ratio of \(v\) lies in \(\mathcal{R}\), so \(D(v) \subset \{2,3,4,5,7,8,9\}\). Every member of \(\mathcal{R}\) is a pinned ratio of some \(v \in \R^5_{>0}\).
\item[\rm (b)] If some pinned ratio of \(v\) has denominator \(7\), then \(v\) is proportional to a permutation of \((1,1,1,1,3)\) or of \((1,1,1,2,2)\); if one has denominator \(8\), to a permutation of \((1,1,1,2,3)\); if one has denominator \(9\), to a permutation of \((1,1,2,2,3)\); if \(1/5\) or \(4/5\) is a pinned ratio, \(v\) is proportional to \((1,1,1,1,1)\).
\item[\rm (c)] If \(1/3\) or \(2/3\) is a pinned ratio of \(v\), then \([5]\) is a disjoint union of three blocks each of \(v\)-sum \(L/3\). If \(1/4\) or \(3/4\) is a pinned ratio, \([5]\) is a disjoint union of four blocks each of sum \(L/4\); so \(v\) is a permutation of \((c, d, u, u, u)\) with \(c + d = u\). If \(2/5\) or \(3/5\) is a pinned ratio, then \(v\) is a permutation of \((a, b, s, s, s)\) with \(a + b = 2s\) or of \((2s, s, s, a, b)\) with \(a + b = s\).
\end{enumerate}
\end{theorem}

\begin{proof}
(a) Let \(\alpha\) be an effective level of \(v\) and \(\mathcal{T} \subset \slice{\alpha}(v)\) a minimal effective subfamily. By Proposition~\ref{prop:fewmembers} its contraction is a \(t\)-configuration \(M\) with \(t \le 5\) and \(r_\alpha = r(M)\), which lies in \(\mathcal{R}\) by Proposition~\ref{prop:smallconfigs} and Theorem~\ref{thm:fiveconfigs}. Conversely every configuration of size between two and five is realized by the converse of Proposition~\ref{prop:fewmembers}, splitting weights until five lengths are obtained.

(b) A ratio with denominator \(7\), \(8\), \(9\), or equal to \(1/5\) or \(4/5\), is not the ratio of a configuration of size at most four. So in Proposition~\ref{prop:fewmembers} \(t = 5\), the five classes are singletons, and \(v = \alpha M^{-1}\one\) is proportional to the weight vector of a \(5\)-configuration with that ratio; Theorem~\ref{thm:fiveconfigs} lists these weight vectors, and they are the stated ones up to permutation.

(c) By Proposition~\ref{prop:pinned}(iv) it suffices to treat the ratios \(1/3\), \(1/4\), \(2/5\). For \(1/3\), the configuration \(M\) of Proposition~\ref{prop:fewmembers} is \(I_3\) or has size five. If \(M = I_3\), the three classes \(C_1, C_2, C_3\) partition \([5]\) and each has sum \(\alpha = L/3\). If \(M\) has size five, then \(v\) is proportional to \((1,2,1,1,1)\) by Theorem~\ref{thm:fiveconfigs}, which has the blocks \(2\), \(13\), \(45\). For \(1/4\), the only configuration of ratio \(1/4\) is \(I_4\), so the four classes partition \([5]\) into blocks of sum \(L/4\); one block has two elements, the other three are singletons of weight \(u = L/4\), and the two elements of the fourth block have sum \(u\). For \(2/5\), \(M\) is \(A\) or a \(5\)-configuration of ratio \(2/5\). In the second case all weights are equal, so \(v\) is proportional to \((1,1,1,1,1)\), which is \((a,b,s,s,s)\) with \(a = b = s\). If \(M = A\), the classes have sums \(\alpha M^{-1}\one = \alpha(1,\frac12,\frac12,\frac12) = (2s, s, s, s)\) with \(s = \alpha/2\), and exactly one class has two elements; if it is the heavy class, \(v\) is a permutation of \((a, b, s, s, s)\) with \(a + b = 2s\), and if it is a light class, of \((2s, s, s, a, b)\) with \(a + b = s\).
\end{proof}

\begin{lemma}
\label{lem:incompatible}
Let \(v \in \R^5_{>0}\).
\begin{enumerate}
\item[\rm (a)] The ratios \(1/3\) and \(1/4\) are not both pinned ratios of \(v\).
\item[\rm (b)] If \(2/5\) is a pinned ratio of \(v\), then neither \(1/3\) nor \(1/4\) is.
\item[\rm (c)] If \(2/5\) and \(1/2\) are both pinned ratios of \(v\), then \(v\) is proportional to a permutation of \((1,2,2,2,3)\) or of \((1,1,2,2,4)\).
\end{enumerate}
The same statements hold with \(2/3\), \(3/4\), \(3/5\) in place of \(1/3\), \(1/4\), \(2/5\), by Proposition~{\rm \ref{prop:pinned}(iv)}.
\end{lemma}

\begin{proof}
Only Theorem~\ref{thm:fiveratios}(c) is used, together with the fact that a pinned ratio \(r\) requires \(rL\) to be a subset sum of \(v\), and \(1/3\) requires three pairwise disjoint subsets of sum \(L/3\).

(a) By the ratio \(1/4\), \(v\) is a permutation of \((c, d, u, u, u)\) with \(c + d = u\), \(L = 4u\). The subset sums equal to \(4u/3\) are the pairs \(c + u\) if \(c = u/3\) and \(d + u\) if \(d = u/3\): the singletons are at most \(u\), the other pairs are \(c + d = u\) and \(2u\), and the triples are at least \(c + d + u = 2u\). Since \(c = d = u/3\) contradicts \(c + d = u\), all subsets of sum \(4u/3\) contain one fixed element, and three pairwise disjoint ones do not exist.

(b) By the ratio \(2/5\), \(L = 5s\) and \(v\) is of the first form \((a,b,s,s,s)\) with \(a + b = 2s\), where we may take \(a \le s \le b\), or of the second form \((2s,s,s,a,b)\) with \(a + b = s\), \(a \le b < s\). For the ratio \(1/3\) we list the subset sums equal to \(5s/3\). In the first form: \(b\) if \(a = s/3\); \(a + s\) if \(a = 2s/3\); the other sums are \(a \le s\), \(s\), \(a + b = 2s\), \(b + s \ge 2s\), \(2s\), and the triples, which are at least \(a + 2s\). For a given \(a\) at most one of the two kinds occurs, and neither kind supplies three pairwise disjoint sets. In the second form: \(s + b\) if \(b = 2s/3\); the other sums are \(2s\), \(s\), \(a\), \(b\), \(a + b = s\), \(s + a \le 3s/2\), the pairs containing \(2s\), which exceed \(2s\), \(2s\), and the triples, which are at least \(s + a + b = 2s\). Two sets of the form \(\{s_i, b\}\) are not disjoint. For the ratio \(1/4\), \(v\) would be a permutation of \((c,d,u,u,u)\) with \(c + d = u\), so some entry would occur three times and the other two entries would sum to it. In the first form the entry occurring three times is \(s\) (if \(a = s\) or \(b = s\) then all entries equal \(s\)), and the other two sum to \(2s \ne s\) or, when all are equal, to \(2s\) as well. In the second form no entry occurs three times, since \(a \le b < s < 2s\).

(c) In the first form we list the subset sums equal to \(L/2 = 5s/2\). Singletons are at most \(b \le 2s\). Among the pairs, \(a + b = 2s\), \(a + s \le 2s\), \(s + s = 2s\), and \(b + s = 5s/2\) exactly when \(b = 3s/2\). Among the triples, \(a + b + s = 3s\), \(b + s + s \ge 3s\), \(3s\), and \(a + s + s = 5s/2\) exactly when \(a = s/2\). A four-element set has sum at least \(5s - b \ge 3s\). So \(L/2\) is a subset sum exactly when \((a,b) = (s/2, 3s/2)\), and \(v\) is proportional to a permutation of \((1,3,2,2,2)\). In the second form, the singletons are at most \(2s\); the pairs are \(2s + s = 3s\), \(s + s = 2s\), \(s + a < 2s\), \(s + b < 2s\), \(a + b = s\), and \(2s + a\), \(2s + b\), which equal \(5s/2\) exactly when \(a = s/2\), respectively \(b = s/2\); the triples are \(4s\), \(2s + s + a > 3s\), \(2s + s + b > 3s\), \(2s + a + b = 3s\), \(s + a + b = 2s\), and \(s + s + a\), \(s + s + b\), which equal \(5s/2\) exactly when \(a = s/2\), respectively \(b = s/2\); a four-element set has sum at least \(5s - 2s = 3s\). Since \(a + b = s\), \(L/2\) is a subset sum exactly when \(a = b = s/2\), and \(v\) is proportional to a permutation of \((4,2,2,1,1)\).
\end{proof}

\begin{theorem}
\label{thm:fivekzero}
For every \(v \in \R^5_{>0}\), \(\kzero(v) \le 10\). More precisely, if \(7 \in D(v)\) then \(\kzero(v) = 7\); if \(8 \in D(v)\) then \(\kzero(v) = 8\); if \(9 \in D(v)\) then \(\kzero(v) = 9\); otherwise \(D(v)\) is one of \(\emptyset\), \(\{2\}\), \(\{3\}\), \(\{4\}\), \(\{5\}\), \(\{2,3\}\), \(\{2,4\}\), \(\{2,5\}\), and \(\kzero(v)\) is, respectively, \(1\), \(2\), \(3\), \(4\), \(5\), \(6\), \(4\), \(10\).
\end{theorem}

\begin{proof}
By Theorem~\ref{thm:fiveratios}(a), \(D(v) \subset \{2,3,4,5,7,8,9\}\). Pinned ratios are invariant under scaling \(v\). If \(7 \in D(v)\), Theorem~\ref{thm:fiveratios}(b) lets us scale \(v\) to a permutation of \((1,1,1,1,3)\) or \((1,1,1,2,2)\), with \(L = 7\); every level is an integer \(\alpha \in (0,7)\), so every pinned ratio \(\alpha/7\) has denominator \(7\), \(D(v) = \{7\}\) and \(\kzero(v) = 7\). If \(8 \in D(v)\), we scale to \((1,1,1,2,3)\), \(L = 8\), every pinned ratio has the form \(\alpha/8\) with denominator dividing \(8\), and since \(8 \in D(v)\), \(\kzero(v) = 8\). If \(9 \in D(v)\), we scale to \((1,1,2,2,3)\), \(L = 9\), every denominator divides \(9\), and \(\kzero(v) = 9\). The three cases are exclusive, the vectors being pairwise non-proportional.

Otherwise \(D(v) \subset \{2,3,4,5\}\). A denominator \(3\) means that \(1/3\) or \(2/3\) is a pinned ratio, a denominator \(4\) that \(1/4\) or \(3/4\) is, and a denominator \(5\) that one of \(1/5, 2/5, 3/5, 4/5\) is. By Lemma~\ref{lem:incompatible}(a), \(3\) and \(4\) do not both lie in \(D(v)\). If \(1/5\) or \(4/5\) is a pinned ratio, then \(v\) is proportional to \((1,1,1,1,1)\) by Theorem~\ref{thm:fiveratios}(b), where every pinned ratio has the form \(\alpha/5\) with \(\alpha \in \{1,2,3,4\}\), so \(D(v) = \{5\}\). If \(2/5\) or \(3/5\) is a pinned ratio, Lemma~\ref{lem:incompatible}(b) excludes \(3\) and \(4\) from \(D(v)\). Hence \(D(v)\) is one of the eight listed sets, with the least common multiples \(1, 2, 3, 4, 5, 6, 4, 10\).
\end{proof}

\begin{remark}
\label{rem:fivecells}
The arguments above also identify the length vectors with \(\kzero(v) \ge 6\), although this is not needed in the sequel. By Lemma~\ref{lem:incompatible}(c), \(D(v) = \{2,5\}\) forces \(v\) to be proportional to a permutation of \((1,2,2,2,3)\) or of \((1,1,2,2,4)\), and both have \(\kzero = 10\): for \((1,2,2,2,3)\) the slice at the level \(4\) contains \(15, 23, 24, 34\) with \(2\cdot\one = 2\cdot 1_{15} + 1_{23} + 1_{24} + 1_{34}\), and \(v(\{2,5\}) = 5 = L/2\); for \((1,1,2,2,4)\) the slice at the level \(4\) contains \(5, 34, 123, 124\) with \(2\cdot\one = 2\cdot 1_{5} + 1_{34} + 1_{123} + 1_{124}\), and \(v(\{1,5\}) = 5\). By Theorem~\ref{thm:fiveratios}(c), \(D(v) = \{2,3\}\), with \(\kzero(v) = 6\), occurs exactly when \([5]\) has a partition into three blocks of sum \(L/3\) and a subset of sum \(L/2\); the simplest instance is \((1,1,1,1,2)\). The vectors with \(\kzero(v) = 7, 8, 9\) are those of Theorem~\ref{thm:fiveratios}(b): in \((1,1,1,1,3)\) the slice at the level \(3\) contains \(5, 123, 124, 134, 234\) with \(3\cdot\one = 3\cdot 1_5 + 1_{123} + 1_{124} + 1_{134} + 1_{234}\), of ratio \(3/7\); in \((1,1,1,2,2)\) the slice at the level \(2\) contains \(4, 5, 12, 13, 23\) with \(2\cdot\one = 2\cdot 1_4 + 2\cdot 1_5 + 1_{12} + 1_{13} + 1_{23}\), of ratio \(2/7\); in \((1,1,1,2,3)\) the slice at the level \(3\) contains \(5, 14, 24, 34, 123\) with \(3\cdot\one = 3\cdot 1_5 + 1_{14} + 1_{24} + 1_{34} + 2\cdot 1_{123}\), of ratio \(3/8\); in \((1,1,2,2,3)\) the slice at the level \(4\) contains \(15, 25, 34, 123, 124\) with \(4\cdot\one = 2\cdot 1_{15} + 2\cdot 1_{25} + 3\cdot 1_{34} + 1_{123} + 1_{124}\), of ratio \(4/9\). For every other \(v \in \R^5_{>0}\), \(\kzero(v) \le 5\).
\end{remark}

\subsection{The specialization lemma}

We now return to \(n\) lengths \(\ell \in \R^n_{>0}\) with \(\Krel = \{r \in \Q^n : \sum_j r_j\ell_j = 0\}\) and \(\rho = \dimQ \operatorname{span}_\Q(\ell_1,\dots,\ell_n)\) as in Definition~\ref{def:standing}; the map \(r \mapsto \sum_j r_j\ell_j\) from \(\Q^n\) onto that span has kernel \(\Krel\), so \(\dimQ\Krel = n - \rho\). The idea is that the coincidences \(\ell_T = \ell_{T'}\) among subset sums are rational linear relations among the lengths, and that the lengths can be moved along the space of all vectors satisfying the same rational relations until some of them vanish, without destroying any coincidence.

\begin{lemma}
\label{lem:specialization}
Let \(\ell \in \R^n_{>0}\), let \(\tau_1,\dots,\tau_\rho\) be a \(\Q\)-basis of \(\operatorname{span}_\Q(\ell_1,\dots,\ell_n)\), and let \(A \in \Q^{n\times\rho}\) be the matrix with \(\ell_j = \sum_i A_{ji}\tau_i\). Then \(A\) has rank \(\rho\). For \(t \in \R^\rho\) put \(\ell'(t) = At \in \R^n\).
\begin{enumerate}
\item[\rm (a)] If \(\ell_T = \ell_{T'}\) for two subsets \(T, T' \subset \nset\), then \(\ell'(t)_T = \ell'(t)_{T'}\) \(\forall t\).
\item[\rm (b)] Suppose \(\ell'(t) \ge 0\) coordinatewise and \(\ell'(t) \ne 0\), let \(N = \{j : \ell'(t)_j > 0\}\), and let \(v \in \R^N_{>0}\) be the restriction of \(\ell'(t)\) to \(N\). Then every pinned ratio of \(\ell\) is a pinned ratio of \(v\): \(\Rat \subset R(v)\), \(\Den \subset D(v)\), and \(\kzero(\ell)\) divides \(\kzero(v)\).
\end{enumerate}
\end{lemma}

\begin{proof}
The rows \(A_j\) of \(A\) are the coordinate vectors of the \(\ell_j\) in the basis \(\tau\), and the \(\ell_j\) span the space, so the rows span \(\Q^\rho\) and \(A\) has rank \(\rho\).

(a) \(\ell_T = \sum_i (\sum_{j \in T}A_{ji})\tau_i\), and the \(\tau_i\) are linearly independent over \(\Q\), so \(\ell_T = \ell_{T'}\) exactly when the rational row sums \(A_T = \sum_{j \in T}A_j\) and \(A_{T'}\) coincide, and then \(\ell'(t)_T = A_T\cdot t = A_{T'}\cdot t = \ell'(t)_{T'}\).

(b) Write \(\ell' = \ell'(t)\) and \(L' = \sum_j\ell'_j = L(v) > 0\). Let \(\alpha\) be an effective level of \(\ell\) with a rational representation \(\one = \sum_{T \in \slice{\alpha}}\gamma_T 1_T\). By (a), all \(T \in \slice{\alpha}\) have one common value \(\ell'_T =: \alpha'\). Pairing the representation with \(\ell\) and with \(\ell'\) gives \(\alpha\sum_T\gamma_T = L\) and \(\alpha'\sum_T\gamma_T = L'\), so \(\alpha'/L' = \alpha/L \in (0,1)\) and \(0 < \alpha' < L'\). Let \(\varphi\colon \Q^n \to \Q^N\) be the restriction of coordinates to \(N\); then \(\varphi(1_T) = 1_{T \cap N}\), \(\varphi(\one) = 1_N\), and \(v(T \cap N) = \ell'_T = \alpha'\) \(\forall T \in \slice{\alpha}\), because the coordinates outside \(N\) carry the length \(0\). Proposition~\ref{prop:pinned}(iii) shows that \(\alpha'\) is an effective level of \(v\) with pinned ratio \(\alpha/L\). Hence \(\Rat \subset R(v)\), \(\Den \subset D(v)\), and the least common multiple \(\kzero(\ell)\) of \(\Den\) divides \(\kzero(v)\).
\end{proof}

\begin{lemma}
\label{lem:extremeray}
In the notation of Lemma~{\rm \ref{lem:specialization}} there is \(t^* \in \R^\rho\) with \(\ell'(t^*) \ge 0\), \(\ell'(t^*) \ne 0\), and at least \(\rho - 1\) vanishing coordinates. Consequently the vector \(v\) of Lemma~{\rm \ref{lem:specialization}(b)} can be chosen with at most \(n - \rho + 1 = \dimQ\Krel + 1\) entries.
\end{lemma}

\begin{proof}
For \(t \in \R^\rho\) with \(At \ge 0\) let \(Z(t) = \{j : A_j\cdot t = 0\}\) be the set of tight rows and let \(\operatorname{rk}(t)\) be the rank of \(\{A_j : j \in Z(t)\}\). Since \(A\) is injective, a nonzero \(t\) has \(\operatorname{rk}(t) \le \rho - 1\), because the tight rows annihilate \(t\). Start with \(t_0 = (\tau_1,\dots,\tau_\rho)\), for which \(\ell'(t_0) = \ell > 0\), \(Z(t_0) = \emptyset\) and \(\operatorname{rk}(t_0) = 0\). Suppose \(t \ne 0\) satisfies \(At \ge 0\) and \(\operatorname{rk}(t) \le \rho - 2\). The space of vectors annihilated by the tight rows has dimension at least \(2\), so it contains a vector \(u\) not proportional to \(t\). Then \(Au \ne 0\), and since \(A_j\cdot u = 0\) for \(j \in Z(t)\), some \(j \notin Z(t)\) has \(A_j\cdot u \ne 0\); replacing \(u\) by \(-u\) if necessary, some \(A_j\cdot u < 0\). Put
\[
s^* = \min\Bigl\{\frac{A_j\cdot t}{-A_j\cdot u} : A_j\cdot u < 0\Bigr\} > 0, \qquad t_1 = t + s^*u .
\]
For \(j\) with \(A_j\cdot u \ge 0\), \(A_j\cdot t_1 \ge A_j\cdot t \ge 0\); for \(j\) with \(A_j\cdot u < 0\), \(A_j\cdot t_1 = A_j\cdot t + s^*A_j\cdot u \ge 0\), with equality for a minimizing index \(j^*\). Every \(j \in Z(t)\) remains tight, and \(j^*\) becomes tight, with \(A_{j^*}\cdot u \ne 0\) while \(A_j\cdot u = 0\) on \(Z(t)\), so \(A_{j^*}\) is not in the span of the old tight rows and \(\operatorname{rk}(t_1) > \operatorname{rk}(t)\). Since \(u\) is not proportional to \(t\), \(t_1 \ne 0\), so \(\operatorname{rk}(t_1) \le \rho - 1\). Repeating the step at most \(\rho - 1\) times produces \(t^* \ne 0\) with \(At^* \ge 0\) and \(\operatorname{rk}(t^*) = \rho - 1\); then \(\ell'(t^*) = At^* \ne 0\) and at least \(\rho - 1\) of its coordinates vanish. The vector \(v\) of Lemma~\ref{lem:specialization}(b) then has at most \(n - \rho + 1\) entries.
\end{proof}

\begin{theorem}
\label{thm:kzeroten}
Let \(\ell \in \R^n_{>0}\) with \(\dimQ\Krel \le 4\). Then there is a vector \(v\) with at most five positive entries such that \(\Rat \subset R(v)\) and \(\kzero(\ell)\) divides \(\kzero(v)\). Consequently \(\Den \subset \{2,3,4,5,7,8,9\}\) and \(\kzero(\ell) \le 10\). In particular, for every spectral union of \(n\) pairwise disjoint non-adjacent bounded intervals with \(\dimQ\Krel \le 4\), the number \(k\) of cosets of the spectrum divides \(\kzero(\ell) \le 10\).
\end{theorem}

\begin{proof}
Lemmas~\ref{lem:specialization} and~\ref{lem:extremeray} give a vector \(v\) with \(m \le \dimQ\Krel + 1 \le 5\) positive entries, \(\Rat \subset R(v)\) and \(\kzero(\ell) \mid \kzero(v)\). If \(m < 5\), split one entry of \(v\) into \(6 - m\) positive parts to obtain \(v' \in \R^5_{>0}\). Every effective level \(\alpha\) of \(v\) is an effective level of \(v'\) with the same ratio: the map \(\varphi\) which repeats the split coordinate on its parts sends \(1_T\) to the indicator of the set \(T'\) obtained by replacing the split element, if present, by all its parts, sends the all-ones vector to the all-ones vector, and \(v'(T') = v(T) = \alpha\) with \(0 < \alpha < L(v) = L(v')\); so Proposition~\ref{prop:pinned}(iii) applies. Hence \(R(v) \subset R(v')\) and \(\kzero(v) \mid \kzero(v')\); if \(m = 5\) put \(v' = v\). Theorems~\ref{thm:fiveratios}(a) and~\ref{thm:fivekzero} give \(D(v') \subset \{2,3,4,5,7,8,9\}\) and \(\kzero(v') \le 10\). So \(\Den \subset D(v) \subset D(v')\) and \(\kzero(\ell)\) divides \(\kzero(v') \le 10\). The last sentence is Theorem~\ref{thm:verdict}.
\end{proof}

The bound \(\kzero(\ell) \le 10\) is attained: \(\ell = (1,2,2,2,3)\) has \(\dimQ\Krel = 4\) and \(\kzero(\ell) = 10\) by Remark~\ref{rem:fivecells}. For four lengths \cite[Theorem~4.4]{bell2026four} gives the sharper bounds \(D(v) \subset \{2,3,4,5\}\) and \(\kzero(v) \le 6\); they are not needed here, since a vector with fewer than five positive entries is handled above by splitting an entry.

\begin{figure}[ht]
\centering
\input{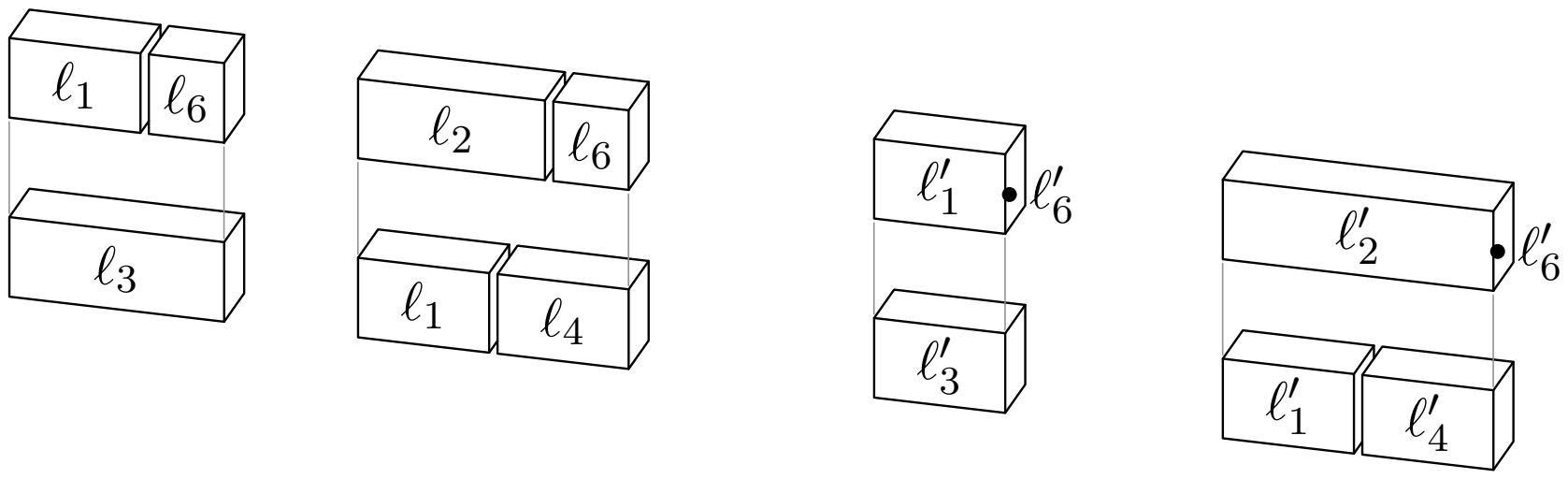}
\caption{The lengths \(\ell = (1, 1 + \theta, 2 - \theta, 1, 2 + \theta, 1 - \theta)\), with \(\dimQ\Krel = 4\) when \(\theta\) is irrational, drawn schematically as boxes at \(\theta \approx 0.4\) on the left, and their specialization \(\ell' = (1, 2, 1, 1, 3, 0)\) at \(\theta = 1\) on the right, in which \(\ell'_6\) vanishes and is drawn as a dot. Each pair of rows shows one of the coincidences \(\ell_3 = \ell_1 + \ell_6\) and \(\ell_1 + \ell_4 = \ell_2 + \ell_6\): the boxes of the upper row cover exactly the extent of the lower row, as the gray lines at the two ends show, and the coincidence survives the specialization.}
\label{fig:specialization}
\end{figure}

\begin{remark}
\label{rem:specialization}
The specialization has a concrete meaning in the fiber picture of Section~\ref{sec:fibers}. With the endpoint classes \(c_1 < \dots < c_r\) of Definition~\ref{def:classes}, every length is \(\ell_j = h d_j + c_{v(j)} - c_{u(j)}\) with an integer \(d_j\). When \(\rho = r\), the numbers \(h, c_2 - c_1, \dots, c_r - c_1\) form a basis of the rational span of the lengths (Lemma~\ref{lem:maximalrank}), and the specializations \(\ell'(t)\) are exactly the vectors \(h'd_j + \varphi_{v(j)} - \varphi_{u(j)}\) obtained by varying \(h\) and replacing the class positions by an arbitrary real potential \(\varphi\) on the classes; a coincidence \(\ell_T = \ell_{T'}\) then says that \(\sum_{j \in T}d_j = \sum_{j \in T'}d_j\) and that \(T\) and \(T'\) have the same net incidence at every class, and it survives every choice of \(h'\) and \(\varphi\). Lemma~\ref{lem:extremeray} chooses \(h'\) and the potential so that the lengths stay nonnegative and \(r - 1\) of them vanish. Equality of the two rational spans is not enough for this description. For \(\Omega = [0, \tfrac12) \cup [\tfrac32, 2)\) with the spectrum \(\Z\), one has \(h = 1\), the classes \(0\) and \(\tfrac12\), and both spans equal to \(\Q\), but moving the class \(\tfrac12\) to \(\tfrac15\) changes the equal lengths \(\tfrac12, \tfrac12\) into \(\tfrac15, \tfrac45\). Figure~\ref{fig:specialization} shows an instance with six lengths.
\end{remark}

\section{Spectral pairs with a consecutive pair}
\label{sec:consecutive}

This section proves the consecutive-pair statement up to ten points, which is the range the proof of Theorem~\ref{thm:A} needs. A fiber of a union of \(n\) intervals with more than \(n\) points contains two consecutive integers, and the consecutive-pair statement \(\CP{k}\) of Definition~\ref{def:namedstatements} says that a finite spectral pair whose spectrum contains two consecutive integers has a spectral set with a translation symmetry. The proof uses a common scheme: the consecutive pair makes the \(k\) points of \(S\) a vanishing sum of \(k\) distinct roots of unity, the classification of such sums by Poonen and Rubinstein lists the possible shapes, and residue-class bounds exclude the shapes without a symmetry. Most cases have a zero set disjoint from \(m\Z\) for some \(m < k\), so a spectrum has at most \(m\) points; two ten-point cases use the refinement that each occupied residue class contains at most two points. The same method gives eleven points, and at twelve points the statement is false; both are treated in Section~\ref{sec:twelve}, after the proof of Theorem~\ref{thm:A}.

\subsection{The reduction to vanishing sums}

Throughout this section \(N \ge 2\), \(\zroot{N} = e^{2\pi i/N}\), and for \(S' \subset \ZN{N}\) the character sum, the zero set \(\zeroset{S'}\), the spectra of \(S'\) in \(\ZN{N}\), periodicity and primitivity are those of Definition~\ref{def:zeroset}. For \(S \subset \R/\Z\) finite and \(t \in \Z\) we write \(\widehat{\muS}(t) = \sum_{s \in S} \ex{st}\) as in Definition~\ref{def:spectralpair}, and
\[
\zeroset{S} = \{ t \in \Z \setminus \{0\} : \widehat{\muS}(t) = 0 \}.
\]
A spectral pair \((A, S)\) with \(A \subset \Z\) thus has \((A - A) \setminus \{0\} \subset \zeroset{S}\). Translating \(S\) by \(s_0\) multiplies \(\widehat{\muS}(t)\) by \(\ex{s_0 t}\), so \(\zeroset{S}\) and the translation symmetries of \(S\) are unchanged by translation; translating \(A\) changes nothing. The two zero sets are compatible: if \(S = s_0 + S'/N\) with \(S' \subset \ZN{N}\), then \(\widehat{\muS}(t) = \ex{s_0 t}\,\widehat{\mu_{S'}}(t \bmod N)\), so \(\zeroset{S}\) is the preimage of \(\zeroset{S'}\) under reduction modulo \(N\).

\begin{lemma}
\label{lem:cpreduction}
Let \(k \ge 2\).
\begin{enumerate}
\item[(a)] Let \((A, S)\) be a spectral pair with \(A \subset \Z\), \(S \subset \R/\Z\) and \(\lvert A\rvert = \lvert S\rvert = k\). Let \(N\), \(s_0\) and \(S' = N(S - s_0) \subset \ZN{N}\) be as in Lemma~\ref{lem:reduction}, and let \(A' \subset \ZN{N}\) be the reduction of \(A\) modulo \(N\). Then \(S'\) is primitive, \(A'\) is a spectrum of \(S'\) in \(\ZN{N}\), and for an integer \(m \ge 2\) one has \(S + 1/m = S\) iff \(m \mid N\) and \(S' + N/m = S'\). In particular \(S\) has a translation symmetry iff \(S'\) is periodic. If \(A\) contains two consecutive integers, then \(1 \in A' - A'\), hence \(1 \in \zeroset{S'}\), that is, \(\sum_{s \in S'} \zroot{N}^{s} = 0\).
\item[(b)] The statement \(\CP{k}\) holds iff, for every \(N\), every primitive \(S' \subset \ZN{N}\) with \(\lvert S'\rvert = k\) which has a spectrum \(A' \subset \ZN{N}\) such that \(A' - A'\) contains a unit of \(\ZN{N}\) is periodic.
\item[(c)] The statement \(\CP{k}\) holds as soon as, for every \(N\), every primitive \(S' \subset \ZN{N}\) with \(\lvert S'\rvert = k\) and \(\sum_{s \in S'} \zroot{N}^{s} = 0\) which has a spectrum in \(\ZN{N}\) is periodic.
\end{enumerate}
\end{lemma}

\begin{proof}
(a) By Theorem~\ref{thm:rational}, \(S - S \subset \Q/\Z\), so Lemma~\ref{lem:reduction} applies: \(0 \in S'\), \(\gcd(N, S') = 1\), and \(A'\) is a spectrum of \(S'\), the reduction being injective on \(A\). Since \(0 \in S'\), the subgroup generated by \(S' - S'\) contains \(S'\) and is \(\ZN{N}\), so \(S'\) is primitive. If \(S + 1/m = S\), then \(1/m \in S - S \subset (1/N)\Z/\Z\), so \(m \mid N\), and \(S' = N(S - s_0)\) satisfies \(S' + N/m = S'\); the converse is the same identity read backward. If \(S'\) is periodic, \(S' + t = S'\) for some \(t \ne 0\), then \(S'\) is a union of cosets of the subgroup generated by \(t\), which is \((N/m)\ZN{N}\) for the order \(m \ge 2\) of \(t\), so \(S' + N/m = S'\) and \(S + 1/m = S\). Finally, two consecutive integers in \(A\) reduce to two elements of \(A'\) with difference \(1\), and \((A' - A') \setminus \{0\} \subset \zeroset{S'}\) gives \(1 \in \zeroset{S'}\).

(b) Suppose the condition holds and let \((A, S)\) be as in \(\CP{k}\). By (a), \(S'\) is primitive with the spectrum \(A'\), and \(1 \in A' - A'\) is a unit, so \(S'\) is periodic and \(S\) has a translation symmetry. Conversely suppose \(\CP{k}\) holds, and let \(S' \subset \ZN{N}\) be primitive with a spectrum \(A'\) and a unit \(u \in A' - A'\). The zero set \(\zeroset{S'}\) is stable under multiplication by units: for a unit \(v\) the automorphism \(\sigma_v\) of \(\Q(\zroot{N})\) with \(\sigma_v(\zroot{N}) = \zroot{N}^{v}\) satisfies \(\sigma_v(\widehat{\mu_{S'}}(t)) = \widehat{\mu_{S'}}(vt)\). Hence \(u^{-1}A'\) is a spectrum of \(S'\) with \(1 \in u^{-1}A' - u^{-1}A'\). Choose \(b \in u^{-1}A'\) with \(b + 1 \in u^{-1}A'\), and let \(A \subset \Z\) be a set of representatives of the elements of \(u^{-1}A'\) in which the representatives of \(b\) and \(b + 1\) are consecutive integers. Put \(S = S'/N \subset \R/\Z\). The matrix \((\ex{as/N})_{s \in S', a \in A}\) depends only on the residues of the elements of \(A\) modulo \(N\), so \((A, S)\) is a spectral pair with \(\lvert A\rvert = \lvert S\rvert = k\) and a consecutive pair in \(A\). By \(\CP{k}\), \(S + 1/m = S\) for some \(m \ge 2\), and by (a) applied to this pair, with the integer \(N\) of Lemma~\ref{lem:reduction} for \(S\) equal to the given \(N\) because \(S'\) is primitive, \(S' + N/m = S'\).

(c) If \(A' - A'\) contains a unit \(u\), then \(u^{-1}A'\) is a spectrum of \(S'\) with \(1\) among its differences, so \(1 \in \zeroset{S'}\) and the hypothesis of (c) applies. Thus the condition of (c) implies the condition of (b).
\end{proof}

The condition in (c) is stronger than the one in (b): it concerns every spectrum of a vanishing \(S'\), and not only the spectra with a unit difference. It is the condition that the proofs below establish for \(k \le 11\).

\subsection{Vanishing sums of distinct roots of unity}
\label{subsec:vanishing}

A \emph{relation} is an identity \(\sum_{i=1}^{r} a_i \eta_i = 0\) with positive integers \(a_i\) and distinct roots of unity \(\eta_i\); its \emph{weight} is \(\sum_i a_i\). It is \emph{minimal} when the only relations \(\sum_i b_i \eta_i = 0\) with integers \(0 \le b_i \le a_i\) are the trivial one and itself. Rotating a relation means multiplying every \(\eta_i\) by one root of unity. Poonen and Rubinstein \cite{poonen-rubinstein1998} observe that every relation is a sum of minimal relations, and their Theorem~3 states that their Table~1 is the complete list, up to rotation, of the minimal relations of weight at most twelve. For every prime \(p\) the sum \(\Rp{p}\) of all \(p\)th roots of unity is minimal; the notation \((R : T_1, \dots, T_j)\) denotes a relation obtained from a relation \(R\) by rotating each \(T_i\) so that it shares exactly one root with \(R\), these shared roots being distinct, subtracting the rotated \(T_i\) from \(R\), and absorbing the minus signs into the roots. Up to weight eleven the table reads: weight \(2\), \(\Rp{2}\); weight \(3\), \(\Rp{3}\); weight \(5\), \(\Rp{5}\); weight \(6\), \((\Rp{5} : \Rp{3})\); weight \(7\), \((\Rp{5} : 2\Rp{3})\) and \(\Rp{7}\); weight \(8\), \((\Rp{5} : 3\Rp{3})\) and \((\Rp{7} : \Rp{3})\); weight \(9\), \((\Rp{5} : 4\Rp{3})\) and \((\Rp{7} : 2\Rp{3})\); weight \(10\), \((\Rp{7} : 3\Rp{3})\) and \((\Rp{7} : \Rp{5})\); weight \(11\), \((\Rp{7} : 4\Rp{3})\), \((\Rp{7} : \Rp{5}, \Rp{3})\), \((\Rp{7} : (\Rp{5} : \Rp{3}))\) and \(\Rp{11}\). There is no minimal relation of weight \(1\) or \(4\).

We write \(\mu_p\) for the set of \(p\)th roots of unity. For \(p \in \{5, 7\}\) and a set \(X \subset \ZN{p}\) with \(1 \le \lvert X\rvert \le p - 1\), put
\[
T(p, X) = \{ \zroot{p}^{x} : x \notin X \} \cup \bigcup_{x \in X} \zroot{p}^{x}\{\zroot{6}, \zroot{6}^{-1}\}, \qquad T_p = T(p, \{0\}).
\]
Subtracting the rotated triangle \(\zroot{p}^{x}\{1, \zroot{3}, \zroot{3}^{2}\}\) from \(\Rp{p}\) replaces the root \(\zroot{p}^{x}\) by \(-\zroot{p}^{x}\zroot{3} = \zroot{p}^{x}\zroot{6}^{-1}\) and \(-\zroot{p}^{x}\zroot{3}^{2} = \zroot{p}^{x}\zroot{6}\), so the roots of a relation of type \((\Rp{p} : j\Rp{3})\) form a set \(T(p, X)\) with \(\lvert X\rvert = j\), of weight \(p + j\), and \(T(p, \{x\}) = \zroot{p}^{x} T_p\). The set \(T_5\) is the configuration of \cite[Proposition~6.2]{bell2026four}. The relation \((\Rp{7} : \Rp{5})\) has the root set
\[
T_{7,5} = (\mu_7 \setminus \{1\}) \cup \bigl(-(\mu_5 \setminus \{1\})\bigr),
\]
and the two remaining relations of weight eleven have the root sets
\begin{align*}
T_{7,5,3}(y) &= (\mu_7 \setminus \{1, \zroot{7}^{y}\}) \cup \bigl(-(\mu_5 \setminus \{1\})\bigr) \cup \zroot{7}^{y}\{\zroot{6}, \zroot{6}^{-1}\}, \qquad y \in \ZN{7} \setminus \{0\}, \\
T_{7,(5,3)}(W) &= (\mu_7 \setminus \{1\}) \cup \bigl(-(W \setminus \{1\})\bigr), \qquad W = u^{-1}T_5 \text{ for some } u \in T_5,
\end{align*}
after a rotation that places the root shared with \(\Rp{5}\), respectively with \((\Rp{5} : \Rp{3})\), at \(1\); here \(W\) is a rotation of \(T_5\) containing \(1\). In each of these sets the roots are distinct. Within one constituent this holds by construction: two roots \(\zroot{p}^{x}\zroot{6}^{\varepsilon}\) and \(\zroot{p}^{x'}\zroot{6}^{\varepsilon'}\) coincide only when \(\zroot{p}^{x - x'} \in \{1, \zroot{3}, \zroot{3}^{2}\}\), that is, when \(x = x'\) and \(\varepsilon = \varepsilon'\). Roots from different constituents cannot coincide: a root of \(\mu_p\) has order \(1\) or \(p\), whereas \(\zroot{p}^{x}\zroot{6}^{\pm 1}\) has order \(6\) when \(x = 0\) and order \(6p\) otherwise. A root \(-\zroot{5}^{i}\) with \(i \ne 0\) has order \(10\). In \(T_{7,(5,3)}(W)\), the roots in \(\mu_7 \setminus \{1\}\) have order \(7\), whereas the roots \(-u'/u\) with \(u' \in T_5 \setminus \{u\}\) have orders dividing \(30\). Every minimal relation of weight at most eleven therefore has all its coefficients equal to one, and its root set is a rotation of one of the sets \(\mu_p\), \(T(p, X)\), \(T_{7,5}\), \(T_{7,5,3}(y)\), \(T_{7,(5,3)}(W)\).

The consequence used below is the following. Let \(S\) be a set of \(k \le 11\) distinct roots of unity with \(\sum_{\omega \in S} \omega = 0\). Then \(S\) is a disjoint union of \emph{blocks}, each block being a rotation of one of the root sets above. Indeed, the relation \(\sum_{\omega \in S} \omega = 0\) contains a minimal subrelation \(\sum b_\omega \omega = 0\) with \(0 \le b_\omega \le 1\); its support is a block, and the remainder \(\sum (1 - b_\omega)\omega = 0\) is again a relation with coefficients one, of smaller weight, to which the argument applies again. A decomposition into blocks is not unique in general; the arguments below take any one decomposition. The blocks are listed by their weights, so for a given \(k\) the possible shapes of \(S\) are indexed by the partitions of \(k\) into parts from \(\{2, 3, 5, 6, 7, 8, 9, 10, 11\}\), with each part carrying every relation of that weight, repeated parts allowed, and arbitrary relative rotations of the blocks.

We also record how a translation symmetry of \(S\) reads on the roots. If \(S \subset \R/\Z\) has \(k\) elements and \(S + 1/m = S\) for some integer \(m \ge 2\), then rotation by \(\ex{1/m}\) permutes the \(k\) roots \(\ex{s}\) in orbits of size \(m\), so \(m \mid k\). For \(k\) prime a set with a translation symmetry is therefore a coset \(s_0 + (1/k)\Z/\Z\), that is, its roots form a rotation of \(\mu_k\); for \(k = 8\) the symmetric sets are those with \(S + 1/2 = S\), the unions of four antipodal pairs; for \(k = 9\) those with \(S + 1/3 = S\), the unions of three rotated triangles; for \(k = 10\) those with \(S + 1/2 = S\) or \(S + 1/5 = S\).

\subsection{Three lemmas}

For \(t \in \Z\) put \(c(t) = 2\cos(\pi t/3) - 1\), so that \(c(t) = 1, 0, -2, -3\) according as \(t \equiv 0, \pm 1, \pm 2, 3 \pmod 6\). We write \([p \mid t]\) for \(1\) when \(p \mid t\) and for \(0\) otherwise, and for a finite set \(T\) of roots of unity \(P_T(t) = \sum_{\omega \in T} \omega^{t}\) and \(Z(T)=\{t\in\Z\setminus\{0\}:P_T(t)=0\}\).

\begin{lemma}
\label{lem:blockpowersums}
Let \(t \in \Z\) and let \(\rho\) be a root of unity.
\begin{enumerate}
\item[(i)] \(P_{\rho\mu_p}(t) = p\rho^{t}[p \mid t]\).
\item[(ii)] For \(p \in \{5, 7\}\) and \(X \subset \ZN{p}\) with \(1 \le \lvert X\rvert = j \le p - 1\),
\[
P_{T(p, X)}(t) = p[p \mid t] + c(t)E_X(t), \qquad E_X(t) = \sum_{x \in X}\zroot{p}^{xt}.
\]
Here \(E_X(t) = j\) when \(p \mid t\), and \(E_X(t) \ne 0\) when \(p \nmid t\); moreover \(\lvert E_X(t)\rvert = 1\) when \(p \nmid t\) and \(j \in \{1, p - 1\}\). In particular \(P_{T_p}(t) = p[p \mid t] + c(t)\).
\item[(iii)] \(P_{T_{7,5}}(t) = (7[7 \mid t] - 1) + (-1)^{t}(5[5 \mid t] - 1)\).
\end{enumerate}
Consequently \(\zeroset{\rho T(p, X)} = \{ t : p \nmid t,\ t \equiv \pm 1 \pmod 6 \}\) \(\forall p \in \{5, 7\}\) and every admissible \(X\), and \(\zeroset{\rho T_{7,5}} = \{ t \text{ odd} : 5 \nmid t,\ 7 \nmid t \}\). Both zero sets consist of odd integers.
\end{lemma}

\begin{proof}
(i) is the geometric sum. (ii) The roots \(\zroot{p}^{x}\), \(x \notin X\), contribute \(p[p \mid t] - E_X(t)\), and each \(x \in X\) contributes \(\zroot{p}^{xt}(\zroot{6}^{t} + \zroot{6}^{-t}) = 2\cos(\pi t/3)\zroot{p}^{xt}\); adding gives the formula. If \(p \mid t\) every term of \(E_X(t)\) equals \(1\). If \(p \nmid t\), multiplication by \(t\) permutes \(\ZN{p}\), so \(E_X(t) = \sum_{y \in tX}\zroot{p}^{y}\) with \(\lvert tX\rvert = j\); since \(\cyc{p}\) is the minimal polynomial of \(\zroot{p}\), an integer relation \(\sum_{y \in \ZN{p}} n_y\zroot{p}^{y} = 0\) holds only when \(y \mapsto n_y\) is constant, and the indicator of \(tX\) is not constant because \(0 < j < p\). For \(j = 1\), \(E_X(t)\) is a root of unity; for \(j = p - 1\) and \(p \nmid t\), \(E_X(t) = -\zroot{p}^{yt}\) for the element \(y \notin X\). (iii) \(\sum_{j \ne 0}\zroot{7}^{jt} = 7[7 \mid t] - 1\) and \(\sum_{i \ne 0}(-\zroot{5}^{i})^{t} = (-1)^{t}(5[5 \mid t] - 1)\).

Zero sets. The factor \(\rho^{t}\) does not affect vanishing. In (ii), if \(p \mid t\) then \(P = p + jc(t) \in \{p + j, p, p - 2j, p - 3j\}\), and \(p = 2j\) or \(p = 3j\) is impossible for \(p \in \{5, 7\}\); if \(p \nmid t\) then \(P = c(t)E_X(t)\) with \(E_X(t) \ne 0\), which vanishes iff \(c(t) = 0\), that is, \(t \equiv \pm 1 \pmod 6\). In (iii), if \(5 \nmid t\) and \(7 \nmid t\) the sum is \(-1 - (-1)^{t}\), zero exactly for odd \(t\); if \(7 \mid t\) and \(5 \nmid t\) it is \(6 - (-1)^{t}\); if \(5 \mid t\) and \(7 \nmid t\) it is \(-1 + 4(-1)^{t}\); if \(35 \mid t\) it is \(6 + 4(-1)^{t}\); the last three are never zero.
\end{proof}

\begin{lemma}
\label{lem:clique}
Let \(Z \subset \Z \setminus \{0\}\) and let \(A \subset \Z\) be finite with \((A - A) \setminus \{0\} \subset Z\).
\begin{enumerate}
\item[(i)] If \(m \ge 2\) and \(Z \cap m\Z = \emptyset\), then \(\lvert A\rvert \le m\).
\item[(ii)] If \(m \ge 1\) and there is a root of unity \(\alpha\) with \(Z \cap m\Z \subset \{ t \in m\Z : \alpha^{t/m} = -1 \}\), then every residue class modulo \(m\) contains at most two elements of \(A\), so \(\lvert A\rvert \le 2m\).
\end{enumerate}
\end{lemma}

\begin{proof}
(i) Two elements of \(A\) congruent modulo \(m\) have a nonzero difference in \(Z \cap m\Z\). (ii) If \(x < y < z\) lie in one class, then \(t_1 = y - x\), \(t_2 = z - y\) and \(t_1 + t_2\) lie in \(Z \cap m\Z\), so \(\alpha^{t_1/m} = \alpha^{t_2/m} = -1\), and then \(\alpha^{(t_1 + t_2)/m} = 1 \ne -1\).
\end{proof}

\begin{lemma}
\label{lem:rationalmoduli}
The following facts about algebraic integers are used repeatedly.
\begin{enumerate}
\item[(a)] Let \(x\) be an algebraic integer in a cyclotomic field. Then \(\lvert x\rvert^{2}\) is an algebraic integer, and if \(\lvert x\rvert^{2}\) is rational it is a nonnegative integer. In particular, if the modulus of a sum of roots of unity is rational, it is an integer.
\item[(b)] Let \(p \ge 5\) be prime and \(p \nmid v\). Then \(\lvert 1 + \zroot{p}^{v}\rvert^{2} = 2 + 2\cos(2\pi v/p)\) is irrational. Hence \(\lvert \zroot{p}^{a} + \zroot{p}^{b}\rvert^{2}\) is irrational whenever \(a \not\equiv b \pmod p\).
\item[(c)] If \(a\zroot{7}^{u} + b\zroot{7}^{w} = 0\) with \(a, b \in \Q\) and \(u \not\equiv w \pmod 7\), then \(a = b = 0\).
\end{enumerate}
\end{lemma}

\begin{proof}
(a) Complex conjugation restricts to an automorphism of the cyclotomic field, so \(\bar x\) is an algebraic integer and so is \(x\bar x = \lvert x\rvert^{2}\); a rational algebraic integer is an integer. A rational number whose square is an integer is itself an integer. (b) If \(2\cos(2\pi v/p) = r\) were rational, then \(\zroot{p}^{v}\) would be a root of \(z^{2} - rz + 1\), while its minimal polynomial \(\cyc{p}\) has degree \(p - 1 \ge 4\). The second claim follows from \(\lvert \zroot{p}^{a} + \zroot{p}^{b}\rvert = \lvert 1 + \zroot{p}^{b - a}\rvert\). (c) If \(a \ne 0\) then \(\zroot{7}^{u - w} = -b/a\) is rational, while a primitive seventh root of unity is not real.
\end{proof}

\subsection{At most six points}

\begin{proposition}
\label{prop:cpsmall}
Let \((A, S)\) be a spectral pair with \(A \subset \Z\), \(S \subset \R/\Z\), \(\lvert A\rvert = \lvert S\rvert = k \le 6\), and suppose \(\widehat{\muS}(1) = \sum_{s \in S}\ex{s} = 0\); this holds when \(A\) contains two consecutive integers. Then \(S + 1/m = S\) for some \(m \ge 2\). More precisely, for \(k \in \{2, 3, 5\}\) the set \(S\) is a coset \(s_0 + (1/k)\Z/\Z\); for \(k = 4\), \(S + 1/2 = S\); for \(k = 6\), \(S + 1/2 = S\) or \(S + 1/3 = S\). The case \(k = 1\) does not occur.
\end{proposition}

\begin{proof}
By Theorem~\ref{thm:rational}, \(S - S \subset \Q/\Z\); translating \(S\), which changes neither \(\zeroset{S}\) nor the symmetries of \(S\) nor the vanishing of \(\widehat{\muS}(1)\), we may assume that the \(k\) numbers \(\omega_s = \ex{s}\) are distinct roots of unity, with \(\sum_s\omega_s = 0\). A single root of unity does not vanish, so \(k \ge 2\). Two distinct roots with sum zero are antipodal, so for \(k = 2\) one has \(S = \{s_0, s_0 + 1/2\}\). For \(k = 3\), after rotating so that one root is \(1\), the other two roots \(u, v\) satisfy \(\lvert 1 + u\rvert = \lvert v\rvert = 1\), so \(\operatorname{Re}u = -1/2\), \(u \in \{\zroot{3}, \zroot{3}^{2}\}\) and \(v = -1 - u\) is the other cube root of unity; the roots form a rotation of \(\mu_3\), and \(S = s_0 + (1/3)\Z/\Z\). For \(k = 4\), the table of Section~\ref{subsec:vanishing} has no minimal relation of weight \(4\) or \(1\), so the four roots split into two blocks of weight \(2\), two antipodal pairs, and \(S + 1/2 = S\). For \(k = 5\) the blocks have weights \(5\) or \(3 + 2\). In the first case the roots form a rotation of \(\mu_5\) and \(S = s_0 + (1/5)\Z/\Z\). In the second case \(S\) corresponds to \(\eta\mu_3 \cup \{\rho, -\rho\}\) for roots of unity \(\eta, \rho\), and \(\widehat{\muS}(t) = 3\eta^{t}[3 \mid t] + \rho^{t}(1 + (-1)^{t})\) by Lemma~\ref{lem:blockpowersums}(i). For \(3 \mid t\) this equals \(3\eta^{t}\) when \(t\) is odd and \(3\eta^{t} + 2\rho^{t}\), of modulus at least \(1\), when \(t\) is even. Thus \(\zeroset{S} \cap 3\Z = \emptyset\), and Lemma~\ref{lem:clique}(i) gives \(\lvert A\rvert \le 3 < 5\), so this case does not occur. For \(k = 6\), \cite[Proposition~6.2]{bell2026four} says that the six roots form three antipodal pairs, or two rotated triangles, or a rotation of \(T_5\). In the first case \(S + 1/2 = S\), in the second \(S + 1/3 = S\). In the third, \cite[Lemma~6.3]{bell2026four} (which is the case \(p = 5\), \(X = \{0\}\) of Lemma~\ref{lem:blockpowersums}) gives \(\zeroset{S} = \{t : 5 \nmid t,\ t \equiv \pm 1 \pmod 6\}\), a set of odd integers, so \(\zeroset{S} \cap 2\Z = \emptyset\) and \(\lvert A\rvert \le 2 < 6\) by Lemma~\ref{lem:clique}(i); this case does not occur either.
\end{proof}

The case \(k = 6\) with a consecutive pair in \(A\) is \cite[Lemma~6.4]{bell2026four}. The argument for \(7 \le k \le 11\) has the same three steps: list the shapes, restrict the zero set of each shape without a symmetry using Lemma~\ref{lem:blockpowersums}, and bound the size of a spectrum by Lemma~\ref{lem:clique}. In every case the bound obtained is strictly smaller than \(k\), and never exceeds \(7\).

\subsection{Seven to ten points}
\label{subsec:seventen}

In the four theorems of this subsection \((A, S)\) is a spectral pair with \(A \subset \Z\), \(S \subset \R/\Z\), \(\lvert A\rvert = \lvert S\rvert = k\), and \(\widehat{\muS}(1) = 0\), which holds when \(A\) contains two consecutive integers. As in the proof of Proposition~\ref{prop:cpsmall}, \(S\) is translated so that its elements correspond to \(k\) distinct roots of unity with sum zero, and we identify \(S\) with this set of roots; \(Z = \zeroset{S}\), and \(\rho, \eta, \xi, \rho_i, \eta_i, \xi_i\) denote roots of unity. The shapes of \(S\) are read off the partitions of \(k\) listed in Section~\ref{subsec:vanishing}; in each shape the power sum \(\widehat{\muS}(t)\) is the sum of the power sums of the blocks, given by Lemma~\ref{lem:blockpowersums}. A block of type \((\Rp{p} : \Rp{3})\) is a rotation of \(T_p\) and is written \(\rho T_p\), with \(P_{T_p}(t) = p[p \mid t] + c(t)\); explicitly \(P_{T_5}(t) = 1, 0, -2, -3\) for \(5 \nmid t\) and \(t \equiv 0, \pm 1, \pm 2, 3 \pmod 6\), and \(P_{T_5}(t) = 6, 5, 3, 2\) in the same four cases when \(5 \mid t\).

\begin{theorem}
\label{thm:cpseven}
If \(k = 7\), then \(S = s_0 + \tfrac17\Z/\Z\) for some \(s_0\); in particular \(S + 1/7 = S\).
\end{theorem}

\begin{proof}
The partitions of \(7\) are \(7\), \(5 + 2\) and \(3 + 2 + 2\), so \(S\) has one of the shapes (a) \(\rho\mu_7\), (b) \(\rho T(5, X)\) with \(\lvert X\rvert = 2\), (c) \(\rho\mu_5 \cup \{\eta, -\eta\}\), (d) \(\eta\mu_3 \cup \{\pm\rho_1\} \cup \{\pm\rho_2\}\). Shape (a) is the conclusion. In shape (b), \(Z\) consists of odd integers by Lemma~\ref{lem:blockpowersums}, so \(Z \cap 2\Z = \emptyset\) and \(\lvert A\rvert \le 2\) by Lemma~\ref{lem:clique}(i). In shape (c), \(\widehat{\muS}(t) = 5\rho^{t}[5 \mid t] + \eta^{t}(1 + (-1)^{t})\); for even \(t\) this is \(2\eta^{t}\) or a sum of two terms of moduli \(5\) and \(2\), so it is nonzero, \(Z \cap 2\Z = \emptyset\) and \(\lvert A\rvert \le 2\). In shape (d), put \(\theta = \rho_1/\rho_2\), a root of unity with \(\theta \ne \pm 1\) since the two antipodal pairs are disjoint; then \(\widehat{\muS}(t) = 3\eta^{t}[3 \mid t] + (1 + (-1)^{t})\rho_2^{t}(1 + \theta^{t})\). For odd \(t\) with \(3 \mid t\) this is \(3\eta^{t} \ne 0\); for \(6 \mid t\) it is \(3\eta^{t} + 2\rho_2^{t}(1 + \theta^{t})\), and its vanishing would force \(\lvert 1 + \theta^{t}\rvert = 3/2\), excluded by Lemma~\ref{lem:rationalmoduli}(a). So \(Z \cap 3\Z = \emptyset\) and \(\lvert A\rvert \le 3\). Since \(\lvert A\rvert = 7\), shapes (b), (c), (d) do not occur.
\end{proof}

\begin{theorem}
\label{thm:cpeight}
If \(k = 8\), then \(S + 1/2 = S\); the roots of \(S\) form four antipodal pairs.
\end{theorem}

\begin{proof}
The partitions of \(8\) are \(8\), \(6 + 2\), \(5 + 3\), \(3 + 3 + 2\) and \(2 + 2 + 2 + 2\), so \(S\) has one of the shapes (e) four antipodal pairs, (f) \(\rho T(5, X)\) with \(\lvert X\rvert = 3\), (g) \(\rho T(7, X)\) with \(\lvert X\rvert = 1\), (h) \(\{\eta, -\eta\} \cup \rho T_5\), (i) \(\rho\mu_5 \cup \eta\mu_3\), (j) \(\rho_1\mu_3 \cup \rho_2\mu_3 \cup \{\eta, -\eta\}\) with \(\theta = \rho_1/\rho_2 \notin \mu_3\). Shape (e) is the conclusion. In shapes (f) and (g), \(Z\) consists of odd integers, so \(\lvert A\rvert \le 2\). In shape (h), \(\widehat{\muS}(t) = \eta^{t}(1 + (-1)^{t}) + \rho^{t}P_{T_5}(t)\). For odd \(t\) with \(3 \mid t\) this is \(\rho^{t}P_{T_5}(t)\) with \(P_{T_5}(t) \in \{-3, 2\}\), nonzero; for \(6 \mid t\) it is \(2\eta^{t} + \rho^{t}P_{T_5}(t)\) with \(P_{T_5}(t) \in \{1, 6\}\), and the two moduli differ. So \(Z \cap 3\Z = \emptyset\) and \(\lvert A\rvert \le 3\). In shape (i), \(\widehat{\muS}(t) = 5\rho^{t}[5 \mid t] + 3\eta^{t}[3 \mid t]\), which for \(3 \mid t\) equals \(3\eta^{t}\) or \(5\rho^{t} + 3\eta^{t}\), nonzero in both cases; \(\lvert A\rvert \le 3\). In shape (j), \(\widehat{\muS}(t) = 3\rho_2^{t}(1 + \theta^{t})[3 \mid t] + \eta^{t}(1 + (-1)^{t})\). For even \(t\) with \(3 \nmid t\) this is \(2\eta^{t} \ne 0\); for \(6 \mid t\) it is \(3\rho_2^{t}(1 + \theta^{t}) + 2\eta^{t}\), whose vanishing would force \(\lvert 1 + \theta^{t}\rvert = 2/3\), excluded by Lemma~\ref{lem:rationalmoduli}(a). So \(Z \cap 2\Z = \emptyset\) and \(\lvert A\rvert \le 2\), whether or not \(\theta = -1\). Since \(\lvert A\rvert = 8\), only shape (e) occurs.
\end{proof}

\begin{theorem}
\label{thm:cpnine}
If \(k = 9\), then \(S + 1/3 = S\); the roots of \(S\) form three rotated triangles.
\end{theorem}

\begin{proof}
The partitions of \(9\) are \(9\), \(7 + 2\), \(6 + 3\), \(5 + 2 + 2\), \(3 + 3 + 3\) and \(3 + 2 + 2 + 2\). The shape \(\rho_1\mu_3 \cup \rho_2\mu_3 \cup \rho_3\mu_3\) is the conclusion. The others are excluded as follows, with the bound of Lemma~\ref{lem:clique}(i) in each case.

\(\rho T(5, X)\) with \(\lvert X\rvert = 4\), and \(\rho T(7, X)\) with \(\lvert X\rvert = 2\): \(Z\) consists of odd integers, \(\lvert A\rvert \le 2\).

\(\rho\mu_7 \cup \{\pm\eta\}\): for even \(t\), \(\widehat{\muS}(t) = 7\rho^{t}[7 \mid t] + 2\eta^{t}\) has modulus \(2\) or at least \(5\); \(Z \cap 2\Z = \emptyset\), \(\lvert A\rvert \le 2\).

\(\rho T(5, X) \cup \{\pm\eta\}\) with \(\lvert X\rvert = 2\): \(\widehat{\muS}(t) = \rho^{t}P(t) + \eta^{t}(1 + (-1)^{t})\) with \(P = P_{T(5, X)}\), and for even \(t\) a zero needs \(\lvert P(t)\rvert = 2\). For even \(t\), \(c(t) \in \{1, -2\}\). If \(5 \mid t\), then \(P(t) = 5 + 2c(t) \in \{7, 1\}\). If \(5 \nmid t\), then \(P(t) = c(t)E_X(t)\) and \(\lvert P(t)\rvert^{2} = c(t)^{2}\lvert E_X(t)\rvert^{2}\) is irrational by Lemma~\ref{lem:rationalmoduli}(b), so it is not \(4\). Thus \(Z \cap 2\Z = \emptyset\) and \(\lvert A\rvert \le 2\).

\(\rho T_5 \cup \eta\mu_3\): \(\widehat{\muS}(t) = \rho^{t}P_{T_5}(t) + 3\eta^{t}[3 \mid t]\). For even \(t\) with \(3 \nmid t\) this is \(\rho^{t}P_{T_5}(t)\) with \(P_{T_5}(t) \in \{-2, 3\}\); for \(6 \mid t\), \(P_{T_5}(t) \in \{1, 6\}\) and a zero would need \(\lvert P_{T_5}(t)\rvert = 3\). Thus \(Z \cap 2\Z = \emptyset\), \(\lvert A\rvert \le 2\).

\(\rho\mu_5 \cup \{\pm\eta_1\} \cup \{\pm\eta_2\}\): for \(5 \mid t\), \(\widehat{\muS}(t) = 5\rho^{t} + (1 + (-1)^{t})(\eta_1^{t} + \eta_2^{t})\), and the second term has modulus at most \(4\). Thus \(Z \cap 5\Z = \emptyset\), \(\lvert A\rvert \le 5\).

\(\rho\mu_3 \cup \{\pm\eta_1\} \cup \{\pm\eta_2\} \cup \{\pm\eta_3\}\): \(\widehat{\muS}(t) = 3\rho^{t}[3 \mid t] + (1 + (-1)^{t})(\eta_1^{t} + \eta_2^{t} + \eta_3^{t})\). For odd \(t\) with \(3 \mid t\) this is \(3\rho^{t}\); for \(6 \mid t\) a zero would need \(\lvert \eta_1^{t} + \eta_2^{t} + \eta_3^{t}\rvert = 3/2\), excluded by Lemma~\ref{lem:rationalmoduli}(a). Thus \(Z \cap 3\Z = \emptyset\), \(\lvert A\rvert \le 3\).

Since \(\lvert A\rvert = 9 > 5\), only three rotated triangles remain.
\end{proof}

\begin{theorem}
\label{thm:cpten}
If \(k = 10\), then \(S + 1/2 = S\) or \(S + 1/5 = S\); the roots of \(S\) form five antipodal pairs or two rotated pentagons.
\end{theorem}

\begin{proof}
The partitions of \(10\) are \(10\), \(8 + 2\), \(7 + 3\), \(6 + 2 + 2\), \(5 + 5\), \(5 + 3 + 2\), \(3 + 3 + 2 + 2\) and \(2 + 2 + 2 + 2 + 2\). The shapes \(\rho_1\mu_5 \cup \rho_2\mu_5\) and five antipodal pairs are the conclusion. The others are excluded as follows.

\(\rho T(7, X)\) with \(\lvert X\rvert = 3\), and \(\rho T_{7,5}\): \(Z\) consists of odd integers, \(\lvert A\rvert \le 2\).

\(\rho T(5, X) \cup \{\pm\eta\}\) with \(\lvert X\rvert = 3\): as for \(\lvert X\rvert = 2\) in the proof of Theorem~\ref{thm:cpnine}. For even \(t\) a zero needs \(\lvert P_{T(5, X)}(t)\rvert = 2\); if \(5 \mid t\) then \(P_{T(5, X)}(t) = 5 + 3c(t) \in \{8, -1\}\), and if \(5 \nmid t\) then \(E_X(t) = -(\zroot{5}^{at} + \zroot{5}^{bt})\) for the two elements \(a, b \notin X\), so \(\lvert P_{T(5, X)}(t)\rvert^{2} = c(t)^{2}\lvert E_X(t)\rvert^{2}\) is irrational by Lemma~\ref{lem:rationalmoduli}(b). Thus \(Z \cap 2\Z = \emptyset\), \(\lvert A\rvert \le 2\).

\(\rho T_7 \cup \{\pm\eta\}\): \(\widehat{\muS}(t) = \rho^{t}P_{T_7}(t) + \eta^{t}(1 + (-1)^{t})\) with \(P_{T_7}(t) = 7[7 \mid t] + c(t)\). For odd \(t\) with \(3 \mid t\), \(P_{T_7}(t) \in \{-3, 4\}\); for \(6 \mid t\), \(P_{T_7}(t) \in \{1, 8\}\), of modulus different from \(2\). Thus \(Z \cap 3\Z = \emptyset\), \(\lvert A\rvert \le 3\).

\(\rho\mu_7 \cup \eta\mu_3\): for \(3 \mid t\), \(\widehat{\muS}(t) = 7\rho^{t}[7 \mid t] + 3\eta^{t} \ne 0\); \(\lvert A\rvert \le 3\).

\(\rho T(5, X) \cup \eta\mu_3\) with \(\lvert X\rvert = 2\): \(\widehat{\muS}(t) = \rho^{t}P_{T(5, X)}(t) + 3\eta^{t}[3 \mid t]\). For even \(t\) with \(3 \nmid t\), \(c(t) = -2\) and \(P_{T(5, X)}(t)\) equals \(5 - 4 = 1\) when \(5 \mid t\) and \(-2E_X(t) \ne 0\) when \(5 \nmid t\). For \(6 \mid t\), \(c(t) = 1\) and a zero needs \(\lvert P_{T(5, X)}(t)\rvert = 3\), while \(P_{T(5, X)}(t) = 7\) when \(5 \mid t\) and \(\lvert P_{T(5, X)}(t)\rvert^{2} = \lvert E_X(t)\rvert^{2}\) is irrational when \(5 \nmid t\). Thus \(Z \cap 2\Z = \emptyset\), \(\lvert A\rvert \le 2\).

\(\rho T_5 \cup \{\pm\eta_1\} \cup \{\pm\eta_2\}\): \(\widehat{\muS}(t) = \rho^{t}P_{T_5}(t) + (1 + (-1)^{t})(\eta_1^{t} + \eta_2^{t})\). For odd \(t\) with \(3 \mid t\) this is \(\rho^{t}P_{T_5}(t)\) with \(P_{T_5}(t) \in \{-3, 2\}\); for \(6 \mid t\), \(P_{T_5}(t) \in \{1, 6\}\) and a zero would need \(\lvert \eta_1^{t} + \eta_2^{t}\rvert \in \{1/2, 3\}\), excluded by Lemma~\ref{lem:rationalmoduli}(a) and by \(\lvert \eta_1^{t} + \eta_2^{t}\rvert \le 2\). Thus \(Z \cap 3\Z = \emptyset\), \(\lvert A\rvert \le 3\).

\(\rho\mu_5 \cup \eta\mu_3 \cup \{\pm\xi\}\): \(\widehat{\muS}(t) = 5\rho^{t}[5 \mid t] + 3\eta^{t}[3 \mid t] + \xi^{t}(1 + (-1)^{t})\). If exactly one or exactly two of the primes \(2, 3, 5\) divide \(t\), the terms present have distinct moduli taken from \(\{2, 3, 5\}\) and cannot cancel. Hence \(Z \subset \{t : \gcd(t, 30) = 1\} \cup 30\Z\), and for \(30 \mid t\) a zero forces \(\lvert 3\eta^{t} + 2\xi^{t}\rvert = 5\), that is \(\eta^{t} = \xi^{t}\), and then \(\rho^{t} = -\eta^{t}\), so \((\rho/\eta)^{t} = -1\). Two elements of \(A\) in different residue classes modulo \(30\) have a difference coprime to \(30\), hence odd; three elements of \(A\) pairwise in different classes would have three pairwise odd differences, which is impossible since \(a - c = (a - b) + (b - c)\). So \(A\) meets at most two classes modulo \(30\), and Lemma~\ref{lem:clique}(ii) with \(m = 30\) and \(\alpha = (\rho/\eta)^{30}\) puts at most two elements of \(A\) in each class: \(\lvert A\rvert \le 4\).

\(\rho_1\mu_3 \cup \rho_2\mu_3 \cup \{\pm\eta_1\} \cup \{\pm\eta_2\}\) with \(\theta = \rho_1/\rho_2 \notin \mu_3\) and \(\psi = \eta_1/\eta_2 \ne \pm 1\): \(\widehat{\muS}(t) = 3\rho_2^{t}(1 + \theta^{t})[3 \mid t] + \eta_2^{t}(1 + \psi^{t})(1 + (-1)^{t})\). For even \(t\) with \(3 \nmid t\) a zero means \(\psi^{t} = -1\). For \(6 \mid t\) a zero gives \(3\lvert 1 + \theta^{t}\rvert = 2\lvert 1 + \psi^{t}\rvert\); squaring, with \(x = \theta^{t} + \theta^{-t}\) and \(y = \psi^{t} + \psi^{-t}\), this reads \(9(2 + x) = 4(2 + y)\), that is \(9x - 4y = -10\). The numbers \(x\) and \(y\) are algebraic integers, and all their conjugates are of the form \(\zeta + \zeta^{-1}\) with \(\zeta\) a root of unity, hence real and in \([-2, 2]\). From \(9x = 2(2y - 5)\) and \(8x = 2 \cdot 4x\) we get \(x = 9x - 8x \in 2\mathcal{O}\), where \(\mathcal{O}\) is the ring of algebraic integers; write \(x = 2w\). Every conjugate of \(x\) equals \((4y' - 10)/9\) for a conjugate \(y' \in [-2, 2]\) of \(y\), hence lies in \([-2, -2/9]\), so every conjugate of \(w\) lies in \([-1, -1/9]\). The norm of the nonzero algebraic integer \(w\) has absolute value at least \(1\), so every conjugate of \(w\) has absolute value \(1\), and \(w = -1\). Then \(x = -2\), so \(\theta^{t} = -1\), and \(y = (9x + 10)/4 = -2\), so \(\psi^{t} = -1\). Thus every even \(t \in Z\) satisfies \(\psi^{t} = -1\), and Lemma~\ref{lem:clique}(ii) with \(m = 2\) and \(\alpha = \psi^{2}\) gives \(\lvert A\rvert \le 4\).

Since \(\lvert A\rvert = 10 > 4\), only the two symmetric shapes remain.
\end{proof}

\begin{corollary}
\label{cor:cpten}
Let \((A, S)\) be a spectral pair with \(A \subset \Z\), \(S \subset \R/\Z\) and \(\lvert A\rvert = \lvert S\rvert = k \le 10\), and suppose \(\sum_{s \in S}\ex{s} = 0\). Then \(S + 1/m = S\) for some \(m \ge 2\). In particular \(\CP{k}\) holds for every \(k \le 10\).
\end{corollary}

\begin{proof}
This is Proposition~\ref{prop:cpsmall} for \(k \le 6\) and Theorems~\ref{thm:cpseven} to~\ref{thm:cpten} for \(7 \le k \le 10\). A consecutive pair in \(A\) gives \(\widehat{\muS}(1) = 0\) by Definition~\ref{def:spectralpair}.
\end{proof}

This is the range that the proof of Theorem~\ref{thm:A} uses. Section~\ref{sec:twelve} treats eleven points, where the statement still holds, and twelve points, where it fails.

\section{The class graph and the reduction to eight intervals}
\label{sec:classgraph}

Throughout this section \(\Omega = \bigsqcup_{j=1}^{n} [a_j, b_j)\) is a spectral union of \(n\) intervals in the standing reduction of Definition~\ref{def:standing}, \(\Lambda\) is a spectrum of \(\Omega\) with minimal period \(h^{-1}\), and \(k\), \(S\) and the fibers \(A_x\) are as in Definition~\ref{def:fiber}. By Theorem~\ref{thm:fiber} and Lemma~\ref{lem:everyfiber}, \(L = kh\) and every fiber is a spectrum of \(S\); in particular \(\lvert A_x\rvert = k\) \(\forall x \in \R\). The endpoint classes \(c_1 < \dots < c_r\) in \([0, h)\), the multiplicities \(m_c\), the class indices \(u(j), v(j) \in [r]\) of \(a_j\) and \(b_j\), the integers \(p_j, q_j\) with \(a_j = c_{u(j)} + hp_j\) and \(b_j = c_{v(j)} + hq_j\), the integer parts \(d_j = q_j - p_j\), and the sets \(\Out{c}\) and \(\In{c}\) of leaving and entering heights are those of Definition~\ref{def:classes}. The counting law of Definition~\ref{def:classes} says that at every class \(c\) the number of intervals ending at \(c\) equals the number \(m_c\) of intervals starting at \(c\).

This section records what the positions of the endpoints, read modulo \(h\), say about the lengths and about the fibers. The information is encoded in a directed multigraph on the classes. Its first use is the bound \(\dimQ \Krel \ge n - r\), which together with the merging along rigid pairs (Proposition~\ref{prop:mergechain}) confines the proof of Theorem~A to fully mixing configurations with at most eight intervals. Its second use is the description of the fibers of a configuration in which no fiber contains two consecutive integers: every integer part is then forced, each interval carries at most one point of each fiber, and the points carried by a loop or a backward edge take two values whose locations are determined by the endpoint classes.

\subsection{The class graph}

\begin{definition}
\label{def:classgraph}
The \emph{class graph} \(\clsgraph\) of \(\Omega\) (relative to \(h\)) is the directed multigraph with vertex set \([r]\) and one edge \(u(j) \to v(j)\) for each interval \(j \in [n]\). The edge of \(j\) is a \emph{loop} when \(u(j) = v(j)\), \emph{forward} when \(c_{u(j)} < c_{v(j)}\), and \emph{backward} when \(c_{u(j)} > c_{v(j)}\). The numbers of loops, forward edges and backward edges are \(\nloop\), \(\nfwd\) and \(\nback\), so that \(n = \nloop + \nfwd + \nback\).
\end{definition}

By the counting law, the in-degree and the out-degree of vertex \(i\) of \(\clsgraph\) both equal \(m_{c_i}\): the class graph is balanced. A balanced directed multigraph decomposes into edge-disjoint directed cycles, so every edge of \(\clsgraph\) lies on a directed cycle; a loop is a cycle by itself.

\begin{lemma}
\label{lem:arcs}
For every \(j \in [n]\),
\[
\ell_j = c_{v(j)} - c_{u(j)} + h d_j,
\]
and for every \(x \in [0, h)\) the run of the \(j\)-th interval in the fiber \(A_x\) has length
\[
\lvert I_j(x)\rvert = d_j + [x < c_{v(j)}] - [x < c_{u(j)}],
\]
where \([\,\cdot\,]\) is \(1\) when the condition holds and \(0\) otherwise. Consequently, if the edge of \(j\) is a loop, then \(\ell_j = hd_j \in h\Z\) and \(\lvert I_j(x)\rvert = d_j\) \(\forall x\). If it is forward, then \(d_j = \lfloor \ell_j/h\rfloor\), and \(\lvert I_j(x)\rvert = d_j + 1\) for \(x \in [c_{u(j)}, c_{v(j)})\) and \(\lvert I_j(x)\rvert = d_j\) for the other \(x \in [0, h)\). If it is backward, then \(d_j = \lceil \ell_j/h\rceil\), and \(\lvert I_j(x)\rvert = d_j - 1\) for \(x \in [c_{v(j)}, c_{u(j)})\) and \(\lvert I_j(x)\rvert = d_j\) for the other \(x \in [0, h)\). In every case \(\ell_j/h \equiv (c_{v(j)} - c_{u(j)})/h \pmod 1\). If the edge of \(j\) is not a loop, the run of \(j\) has length \(\lceil \ell_j/h\rceil = \lfloor \ell_j/h\rfloor + 1\) exactly when \(x\) lies on the arc of \(\R/h\Z\) that runs from \(c_{u(j)}\) (included) upward to \(c_{v(j)}\) (excluded), and length \(\lfloor \ell_j/h\rfloor\) otherwise; for a loop the length is \(\ell_j/h\) at every \(x\). The number \(e = k - \sum_j \lfloor \ell_j/h\rfloor = \sum_j \{\ell_j/h\}\) of edges that are not loops and whose arc contains \(x\) is the same for every \(x\).
\end{lemma}

\begin{proof}
The first identity is \(b_j - a_j\) written with the chosen representatives. For the second, Proposition~\ref{prop:intervalfibers} gives \(I_j(x) = \Z \cap [\lceil (a_j - x)/h\rceil, \lceil (b_j - x)/h\rceil)\), and \((a_j - x)/h = p_j + (c_{u(j)} - x)/h\) with \((c_{u(j)} - x)/h \in (-1, 1)\) because \(c_{u(j)}\) and \(x\) both lie in \([0, h)\). The ceiling of a number in \((-1, 1)\) is \(1\) when the number is positive and \(0\) otherwise, so \(\lceil (a_j - x)/h\rceil = p_j + [x < c_{u(j)}]\), and likewise \(\lceil (b_j - x)/h\rceil = q_j + [x < c_{v(j)}]\). Subtracting gives the formula for \(\lvert I_j(x)\rvert\).

If \(u(j) = v(j)\), the two brackets cancel, and \(\ell_j = hd_j\). If \(c_{u(j)} < c_{v(j)}\), the bracket difference equals \(1\) on \([c_{u(j)}, c_{v(j)})\) and \(0\) elsewhere in \([0, h)\); since \(0 < (c_{v(j)} - c_{u(j)})/h < 1\), the first identity gives \(d_j = \lfloor \ell_j/h\rfloor\). If \(c_{u(j)} > c_{v(j)}\), the bracket difference equals \(-1\) on \([c_{v(j)}, c_{u(j)})\) and \(0\) elsewhere, and \(\ell_j/h = d_j - (c_{u(j)} - c_{v(j)})/h\) with the subtracted term in \((0, 1)\), so \(d_j = \lceil \ell_j/h\rceil\); the set \([0, c_{v(j)}) \cup [c_{u(j)}, h)\) on which the run is long is the arc from \(c_{u(j)}\) upward to \(c_{v(j)}\) on the circle \(\R/h\Z\). The congruence modulo \(1\) is the first identity divided by \(h\). Finally, summing the run lengths over \(j\) and using \(\lvert A_x\rvert = k\) shows that the number of edges that are not loops and whose arc contains \(x\) is \(k - \sum_j \lfloor \ell_j/h\rfloor\), independent of \(x\), and \(\sum_j \{\ell_j/h\} = L/h - \sum_j \lfloor \ell_j/h\rfloor\) equals the same number because \(L = kh\).
\end{proof}

The lengths of the integer blocks in a fiber can thus be read off the class graph drawn on the circle \(\R/h\Z\): a block whose edge is not a loop is long exactly on the arc of its edge, and these arcs cover every point of the circle exactly \(e\) times. On \([c_1,c_1+h)\), the fiber is constant on the \(r\) intervals between consecutive endpoint classes; after one full turn it is shifted by \(-1\).

\begin{lemma}
\label{lem:cyclesums}
At each class choose a bijection from the set of intervals ending there onto the set of intervals starting there, and let \(\pi\) be the resulting permutation of \([n]\), so that \(b_j \equiv a_{\pi(j)} \pmod{h\Z}\) \(\forall j\). Then for every cycle \(C\) of \(\pi\) the sum \(\ell_C = \sum_{j \in C} \ell_j\) equals \(k_C h\) for a positive integer \(k_C\), and \(\sum_C k_C = k\). If \(\pi\) has \(t\) cycles, then \(\dimQ \Krel \ge t - 1\).
\end{lemma}

\begin{proof}
The bijections exist by the counting law, and \(\pi\) is a bijection of \([n]\) because every interval ends in exactly one class and starts in exactly one class. For a cycle \(C\) one has \(\pi(C) = C\), so \(\ell_C = \sum_{j \in C} (b_j - a_j) \equiv \sum_{j \in C} (a_{\pi(j)} - a_j) = 0 \pmod{h\Z}\), and \(\ell_C > 0\); thus \(\ell_C = k_C h\) with \(k_C \ge 1\). Summing over the cycles gives \(L = h\sum_C k_C\), and \(L = kh\). For each cycle put \(w_C = k 1_C - k_C \one \in \Z^n\); then \(\sum_j w_C(j)\ell_j = k\ell_C - k_C L = 0\), so \(w_C \in \Krel\). Let the cycles be \(C_1, \dots, C_t\) and suppose \(\sum_{i < t} \lambda_i w_{C_i} = 0\) with rational \(\lambda_i\). A coordinate \(j \in C_t\) gives \(\sum_{i < t} \lambda_i k_{C_i} = 0\), and then a coordinate \(j \in C_m\) with \(m < t\) gives \(k\lambda_m - \sum_{i<t} \lambda_i k_{C_i} = k\lambda_m = 0\). So \(w_{C_1}, \dots, w_{C_{t-1}}\) are linearly independent in \(\Krel\).
\end{proof}

The permutation \(\pi\) depends on the chosen bijections, and so do its cycles; the lemma holds for every choice.

\begin{lemma}
\label{lem:rankbound}
Let \(\rho = \dimQ \operatorname{span}_{\Q}(\ell_1, \dots, \ell_n)\). Then
\[
\rho \le \dimQ \operatorname{span}_{\Q}(h, c_2 - c_1, \dots, c_r - c_1) \le r,
\]
so \(\dimQ \Krel = n - \rho \ge n - r\). If the class graph is connected as an undirected graph, the first inequality is an equality.
\end{lemma}

\begin{proof}
By Lemma~\ref{lem:arcs}, \(\ell_j = (c_{v(j)} - c_1) - (c_{u(j)} - c_1) + d_j h\) is an integer combination of \(h\) and the differences \(c_i - c_1\), which gives the inclusion of spans and the first inequality; the second is the count of generators. Suppose \(\clsgraph\) is connected. Then \(h = L/k\) lies in \(\operatorname{span}_{\Q}(\ell)\), each edge \(j\) gives \(c_{v(j)} - c_{u(j)} = \ell_j - d_j h \in \operatorname{span}_{\Q}(\ell)\), and along an undirected path from the vertex \(1\) to the vertex \(i\) the signed sum of these differences is \(c_i - c_1\). So the two spans coincide.
\end{proof}

\subsection{Fully mixing configurations and the reduction}

\begin{corollary}
\label{cor:fullymixing}
Let the configuration be fully mixing, in the sense of Definition~\ref{def:rigid}. Then \(m_c \ge 2\) for every class, \(r \le n/2\), and \(\dimQ \Krel \ge n - r \ge n/2\). Consequently, when \(\dimQ \Krel \le 4\), then \(n \le 8\), and for \(n \ge 5\) the multiset \((m_{c_i})_{i=1}^r\) of class sizes, together with \(\rho = n - \dimQ \Krel\), is one of the following:
\begin{enumerate}
\item[\rm (i)] \(n = 5\): the pattern \((5)\) with \(\rho = 1\), or the pattern \((3, 2)\) with \(\rho \in \{1, 2\}\);
\item[\rm (ii)] \(n = 6\): the patterns \((4, 2)\) and \((3, 3)\) with \(\rho = 2\), or the pattern \((2, 2, 2)\) with \(\rho \in \{2, 3\}\);
\item[\rm (iii)] \(n = 7\): the pattern \((3, 2, 2)\) with \(\rho = 3\);
\item[\rm (iv)] \(n = 8\): the pattern \((2, 2, 2, 2)\) with \(\rho = 4\).
\end{enumerate}
\end{corollary}

\begin{proof}
A class with \(m_c = 1\) is rigid by Lemma~\ref{lem:rigidpair}, so a fully mixing configuration has \(m_c \ge 2\) at every class, whence \(n = \sum_c m_c \ge 2r\), and Lemma~\ref{lem:rankbound} gives \(\dimQ \Krel \ge n - r \ge n/2\). Now suppose \(\dimQ \Krel \le 4\). Then \(n \le 8\), and \(\rho = n - \dimQ \Krel \ge n - 4\), while \(\rho \le r\) by Lemma~\ref{lem:rankbound}. The patterns are the partitions of \(n\) into \(r\) parts, each at least \(2\), with \(r \ge n - 4\). For \(n = 5\) they are \((5)\) and \((3, 2)\); in the pattern \((5)\) all endpoints lie in one class, so every \(\ell_j\) lies in \(h\Z\) by Lemma~\ref{lem:arcs} and \(\rho = 1\). For \(n = 6\) one needs \(r \ge 2\), which excludes \((6)\) and leaves \((4, 2)\), \((3, 3)\) with \(r = 2\), hence \(\rho = 2\), and \((2, 2, 2)\) with \(\rho \in \{2, 3\}\). For \(n = 7\) one needs \(r \ge 3\), which excludes \((7)\), \((5, 2)\), \((4, 3)\) and leaves \((3, 2, 2)\) with \(\rho = r = 3\). For \(n = 8\) one needs \(r \ge 4\), which leaves only \((2, 2, 2, 2)\) with \(\rho = r = 4\). For \(n \ge 9\) the requirement \(r \ge n - 4 > n/2\) contradicts \(r \le n/2\).
\end{proof}

In the single-class pattern \((n)\), every length lies in \(h\Z\), and \(\Omega = c_1 + h(A + [0, 1))\) with \(A = A_{c_1}\). The fiber equals \(A\) for \(x \in [c_1, c_1+h)\); elsewhere it is obtained from these fibers by \(A_{x+h}=A_x-1\). In this case \(\rho = 1\) and \(\dimQ \Krel = n - 1\). The bound \(\dimQ \Krel \le 4\) therefore forces \(n \le 5\), where it is covered by Theorem~\ref{thm:four} and by the five-interval case of the proof of Theorem~A.

\begin{lemma}
\label{lem:mergerank}
Let \(\Omega'\) be obtained from \(\Omega\) by sliding an interval along a rigid pair as in Lemma~\ref{lem:slide}, followed by the coalescence of any adjacencies that arise, and let \(\ell'\) be its vector of lengths. Then \(\dimQ K_{\ell'} \le \dimQ \Krel\), where \(K_{\ell'}\) denotes the space of rational relations among the entries of \(\ell'\).
\end{lemma}

\begin{proof}
Let \((i, j)\) be the rigid pair. By Lemma~\ref{lem:slide} the components of \(\Omega'\), before coalescence, are the intervals \([a_m, b_m)\) for \(m \ne i, j\) and the interval \([a_i, b_i + \ell_j)\), so the new length vector consists of the \(\ell_m\) with \(m \ne i, j\) and of \(\ell_i + \ell_j\), and has \(n - 1\) entries. Its rational span \(V'\) is contained in \(V = \operatorname{span}_{\Q}(\ell)\), and \(V \subset V' + \Q\ell_i\) because \(\ell_j = (\ell_i + \ell_j) - \ell_i\). Hence \(\rho - 1 \le \rho' \le \rho\) for the dimensions, and \(\dimQ K_{\ell'} = (n - 1) - \rho' \in \{\dimQ \Krel - 1, \dimQ \Krel\}\). A coalescence of two adjacent components replaces two lengths by their sum and is the same computation, so each coalescence also leaves the dimension unchanged or lowers it by one.
\end{proof}

\begin{proposition}
\label{prop:reduction}
Let \(\Omega\) be a spectral union of \(n\) intervals with \(\dimQ \Krel \le 4\), with spectrum \(\Lambda\) and the data \(h\), \(k\), \(S\). Then either \(k = 1\), or there is a fully mixing spectral union \(\Omega_r\) of \(n_r \le \min(n, 8)\) intervals, with the same spectrum \(\Lambda\), the same \(h\), \(k\) and \(S\), and with \(\dimQ K_{\ell^{(r)}} \le 4\) for its length vector \(\ell^{(r)}\); when \(n_r \ge 5\), its class pattern is one of those listed in Corollary~\ref{cor:fullymixing}. Consequently, if \(k \le 5\) holds for every fully mixing spectral union of at most eight intervals with \(\dimQ \Krel \le 4\), then \(k \le 5\) holds for every spectral union of intervals with \(\dimQ \Krel \le 4\), whatever the number of intervals.
\end{proposition}

\begin{proof}
By Proposition~\ref{prop:mergechain}, iterated sliding along rigid pairs ends either at a single interval, in which case \(k = 1\), or at a fully mixing union \(\Omega_r\) of \(n_r \le n\) components with the same spectrum \(\Lambda\), hence the same minimal period \(h^{-1}\), the same \(k\) and the same \(S\). Each step preserves spectrality and the spectrum (Lemma~\ref{lem:slide}) and does not increase the dimension of the space of rational relations (Lemma~\ref{lem:mergerank}), so \(\dimQ K_{\ell^{(r)}} \le 4\). Corollary~\ref{cor:fullymixing} applied to \(\Omega_r\) gives \(n_r \le 8\) and the list of patterns. The last sentence follows because \(k\) is the same for \(\Omega\) and \(\Omega_r\).
\end{proof}

Proposition~\ref{prop:reduction} is used in the proof of Theorem~A. The configurations with \(n_r \le 5\) are handled there by the bound \(k\le 10\), the pigeonhole argument below and the consecutive-pair statements (Proposition~\ref{prop:cpsmall} and Theorems~\ref{thm:cpseven} to~\ref{thm:cpten}), and the patterns with \(6 \le n_r \le 8\) by the fiber-side mechanisms (Theorem~\ref{thm:fullymixing}).

\subsection{The pigeonhole and the forced integer parts}

\begin{lemma}
\label{lem:pigeonhole}
\begin{enumerate}
\item[\rm (i)] If \(k > n\), then every fiber \(A_x\) contains two consecutive integers.
\item[\rm (ii)] Suppose that the class graph has an edge which is not a loop; equivalently, some \(\ell_j\) does not lie in \(h\Z\); this is the case whenever \(\rho \ge 2\). If \(k \ge n\), then some fiber \(A_x\) contains two consecutive integers.
\end{enumerate}
\end{lemma}

\begin{proof}
(i) By Proposition~\ref{prop:intervalfibers}, \(A_x\) is the disjoint union of the \(n\) integer intervals \(I_j(x)\). If \(k > n\), one of them has at least two elements, and these are consecutive integers.

(ii) An edge \(j\) is a loop exactly when \(\ell_j \in h\Z\), by Lemma~\ref{lem:arcs}; if every edge is a loop then every length is a rational multiple of \(h\) and \(\rho = 1\). Let \(j\) be an edge which is not a loop. Then \(\ell_j/h \notin \Z\), and by Lemma~\ref{lem:arcs} the run \(I_j(x)\) has \(\lceil \ell_j/h\rceil\) elements for \(x\) on the arc of \(j\) and \(\lfloor \ell_j/h\rfloor\) elements for \(x\) off it; both sets of \(x\) are nonempty since \(c_{u(j)} \ne c_{v(j)}\). If \(\lfloor \ell_j/h\rfloor \ge 1\), the run has at least two elements on the arc. If \(\lfloor \ell_j/h\rfloor = 0\), the run is empty off the arc, and there the \(k \ge n\) elements of \(A_x\) lie in the remaining \(n - 1\) integer intervals, one of which has two elements.
\end{proof}

\begin{theorem}
\label{thm:forced}
Suppose that no fiber \(A_x\) contains two consecutive integers. Then \(d_j = 1\) for every loop and every backward edge, \(d_j = 0\) for every forward edge, every run \(I_j(x)\) has at most one element, and
\[
k = \nloop + \nback = n - \nfwd .
\]
Every simple directed cycle of the class graph through at least two classes contains a forward edge and a backward edge. Hence if the class graph has an edge which is not a loop, then \(\nfwd \ge 1\) and \(k \le n - 1\).
\end{theorem}

\begin{proof}
A run with two or more elements contains two consecutive integers, so every run has at most one element at every \(x\). For a loop, Lemma~\ref{lem:arcs} gives \(\lvert I_j(x)\rvert = d_j\) \(\forall x\), and \(d_j = \ell_j/h \ge 1\); hence \(d_j = 1\). For a forward edge the run has \(d_j + 1\) elements on the arc of \(j\), so \(d_j + 1 \le 1\), and \(d_j = \lfloor \ell_j/h\rfloor \ge 0\); hence \(d_j = 0\). For a backward edge the run has \(d_j\) elements on the arc, so \(d_j \le 1\), and \(d_j = \lceil \ell_j/h\rceil \ge 1\); hence \(d_j = 1\).

For the count, fix \(x \in [0, h)\) and sum the formula of Lemma~\ref{lem:arcs} over \(j\):
\[
k = \lvert A_x\rvert = \sum_j \lvert I_j(x)\rvert = \sum_j d_j + \sum_{c} m_c [x < c] - \sum_{c} m_c [x < c] = \sum_j d_j,
\]
where the two middle sums count the edges ending at \(c\) and the edges starting at \(c\), both equal to \(m_c\) by the counting law. With the forced values, \(\sum_j d_j = \nloop + \nback\), and \(n = \nloop + \nfwd + \nback\).

Along a simple directed cycle through at least two classes let \(c^{*}\) be the class with the largest representative in \([0, h)\). The cycle visits \(c^{*}\) once and contains no loop, so the edge of the cycle entering \(c^{*}\) comes from a class with a smaller representative and is forward, and the edge leaving \(c^{*}\) goes to a class with a smaller representative and is backward. Finally, the balanced graph \(\clsgraph\) is an edge-disjoint union of simple directed cycles, so an edge which is not a loop lies on a simple directed cycle, which passes through at least two classes; so \(\nfwd \ge 1\) and \(k = n - \nfwd \le n - 1\).
\end{proof}

When no fiber contains a consecutive pair, the lengths are therefore of three kinds: a loop is an interval of length exactly \(h\), a forward edge has length \(c_{v(j)} - c_{u(j)} \in (0, h)\), and a backward edge has length \(h - (c_{u(j)} - c_{v(j)}) \in (0, h)\). The next lemma says which integer each interval contributes to the fiber.

\begin{lemma}
\label{lem:avatars}
Suppose that \(d_j = 1\) for every loop and every backward edge and \(d_j = 0\) for every forward edge, as in Theorem~\ref{thm:forced}. Then for \(x \in [0, h)\):
\begin{enumerate}
\item[\rm (i)] if the edge of \(j\) is a loop at the class \(c\), then \(I_j(x) = \{p_j + 1\}\) for \(x < c\) and \(I_j(x) = \{p_j\}\) for \(x \ge c\);
\item[\rm (ii)] if the edge of \(j\) is forward, from \(u = u(j)\) to \(v = v(j)\), then \(I_j(x) = \{p_j\}\) for \(x \in [c_u, c_v)\) and \(I_j(x) = \emptyset\) otherwise;
\item[\rm (iii)] if the edge of \(j\) is backward, from \(u = u(j)\) to \(v = v(j)\) with \(c_v < c_u\), so that \(q_j = p_j + 1\), then \(I_j(x) = \{q_j\}\) for \(x \in [0, c_v)\), \(I_j(x) = \emptyset\) for \(x \in [c_v, c_u)\), and \(I_j(x) = \{p_j\}\) for \(x \in [c_u, h)\).
\end{enumerate}
Consequently, at every class \(c\) the sets \(\Out{c} = \{q_i : c_{v(i)} = c\}\) and \(\In{c} = \{p_j : c_{u(j)} = c\}\) have \(m_c\) elements each, and if \(C\) denotes the set of points of the fiber carried by the intervals with neither endpoint in the class \(c\), then the fiber is \(C \sqcup \Out{c}\) immediately before \(c\) and \(C \sqcup \In{c}\) at \(c\) and immediately after it. The sets \(\Out{c}\) and \(\In{c}\) are disjoint, since a common height would give \(b_i = a_j\), contrary to the standing reduction.
\end{lemma}

\begin{proof}
As in the proof of Lemma~\ref{lem:arcs}, for \(x \in [0, h)\) one has \(I_j(x) = \Z \cap [p_j + [x < c_{u(j)}],\ q_j + [x < c_{v(j)}])\). For a loop at \(c\), \(q_j = p_j + 1\) and both brackets equal \([x < c]\), so \(I_j(x) = \{p_j + [x < c]\}\). For a forward edge, \(q_j = p_j\), and the integer interval \([p_j + [x < c_u], p_j + [x < c_v])\) is nonempty exactly when \([x < c_u] = 0\) and \([x < c_v] = 1\), that is for \(c_u \le x < c_v\), where it is \(\{p_j\}\). For a backward edge, \(q_j = p_j + 1\): for \(x < c_v\) the integer run is \(\Z\cap[p_j + 1, p_j + 2) = \{q_j\}\); for \(c_v \le x < c_u\) it is \(\Z\cap[p_j + 1, p_j + 1) = \emptyset\); for \(x \ge c_u\) it is \(\Z\cap[p_j, p_j + 1) = \{p_j\}\).

For the consequence, consider a class \(c\) and the intervals with an endpoint in it. At \(c=0\), the left-hand fiber is computed using \(A_{x+h}=A_x-1\). By (i) to (iii), an interval \(i\) ending at \(c\) carries the point \(q_i\) immediately before \(c\) (for a loop \(q_i = p_i + 1\); for a forward edge \(q_i = p_i\), carried on \([c_{u(i)}, c)\); for a backward edge \(q_i\) is carried on \([0, c)\) when \(c>0\), and immediately to the left of \(0\) when \(c=0\)) and no longer carries it from \(c\) on; an interval \(j\) starting at \(c\) carries no point immediately before \(c\) unless it also ends there, and carries \(p_j\) from \(c\) on. Intervals with no endpoint in the class carry the same point, or no point, on both sides of \(c\), since by (i) to (iii) their runs change only at their own classes. The points carried by different intervals at one \(x\) are distinct because the runs \(I_j(x)\) are pairwise disjoint (Proposition~\ref{prop:intervalfibers}); applied immediately before \(c\) this shows that the \(m_c\) heights \(q_i\) with \(c_{v(i)} = c\) are distinct, and applied at \(c\) that the \(m_c\) heights \(p_j\) with \(c_{u(j)} = c\) are distinct. This is the description of the two fibers.
\end{proof}

For fibers without consecutive pairs, Lemma~\ref{lem:avatars} gives the form of the fibers used in the arguments leading to Theorem~\ref{thm:fullymixing}. The two integer heights \(p_j\) and \(q_j\) associated with a loop or a backward edge are called its two avatars; they differ by one, and the lemma specifies where each occurs. The fibers \(A_x\) with \(x \in [0, h)\) are all spectra of the same set \(S\), and a point and its other avatar never belong to the same fiber, since the fibers contain no consecutive pair. What the maximal rank \(\rho = r\) of Corollary~\ref{cor:fullymixing} says about the coincidences of subset sums is recorded where it is used, in Lemma~\ref{lem:maximalrank}.

\section{Fiber-side mechanisms and the fully mixing case}
\label{sec:mechanisms}

Section~\ref{sec:classgraph} reduced Theorem~\ref{thm:A} to the fully mixing configurations with at most eight intervals, and it described, when no fiber contains two consecutive integers, how a fiber changes as the point \(x\) crosses an endpoint class: the integer parts are forced (Theorem~\ref{thm:forced}), every nonempty run is a single integer, and at a class \(c\) exactly the \(m_c\) points carried by the edges ending at \(c\) leave the fiber while the \(m_c\) points carried by the edges starting at \(c\) enter (Lemma~\ref{lem:avatars}). This section exploits the fact that every one of these fibers is a spectrum of the same \(k\)-point set \(S\). Each crossing therefore relates two orthogonal bases of \(\C^S\) that differ in at most \(m_c\) vectors, and the shape of the exponential vectors turns this relation into arithmetic constraints on \(S\). Four mechanisms result. A loop at a class of size two forces a spectrum with a consecutive pair (loop exchange). For \(k\in\{3,5\}\) the fibers are cosets of one subgroup and a class that moves fewer than \(k\) points cannot change the coset (cosets). For \(k=4\) the fibers are pairs of antipodal pairs and a class that moves at most three points keeps one pair and translates the other inside a fixed subgroup (lineage). A closed walk of classes that carries whole transition spaces around gives an invariant subspace of a diagonal unitary (monodromy). Together with the class graphs of the two patterns at \(n=7\) and \(n=8\), these mechanisms prove that a fully mixing spectral union of six, seven or eight intervals with \(\dimQ\Krel\le 4\) has \(k\le 4\), \(3\), \(2\) respectively (Theorem~\ref{thm:fullymixing}).

Throughout this section \(\Omega\) is a spectral union of \(n\) intervals \([a_j,b_j)\) in the standing reduction of Definition~\ref{def:standing}, \(\Lambda\) is a spectrum of \(\Omega\) with minimal period \(h^{-1}\), and \(\Gamma\), \(k\) and \(S=h\Gamma\subset\R/\Z\) are as in Theorem~\ref{thm:periodic}, so that \(\lvert S\rvert=k\) and \(S\) has no nontrivial translation symmetry. By Lemma~\ref{lem:everyfiber} every fiber \(A_x\), \(x\in\R\), is a spectrum of \(S\) in the sense of Definition~\ref{def:spectralpair}. The classes \(c_1<\dots<c_r\) in \([0,h)\), their sizes \(m_c\), the integers \(u(j)\), \(v(j)\), \(p_j\), \(q_j\), \(d_j=q_j-p_j\), the sets \(\Out{c}\), \(\In{c}\) and the class graph \(\clsgraph\) are those of Definitions~\ref{def:classes} and~\ref{def:classgraph}. By Proposition~\ref{prop:intervalfibers} the fiber is constant on each arc between two consecutive classes, and \(A_{x+h}=A_x-1\). For a class \(c\) we write \(A_{c-}\) for the fiber on the arc of \([0,h)\) that ends at \(c\) and \(A_c\) for the fiber on the arc that starts at \(c\); when \(c=c_1=0\), the arc ending at \(c\) is \([c_r-h,0)\), and the relation \(A_{x+h}=A_x-1\) transports every statement about the crossing of \(c\) to every point of \(c+h\Z\).

We assume \(k\ge 2\). Then \(S-S\subset\Q/\Z\) by Theorem~\ref{thm:rational}, and we fix \(s_0\in S\) and the integer \(N\) of Lemma~\ref{lem:reduction}, so that \(S'=N(S-s_0)\subset\ZN{N}\) has \(k\) elements, contains \(0\), satisfies \(\gcd(N,S')=1\), and every spectrum \(A\subset\Z\) of \(S\) reduces injectively modulo \(N\) onto a spectrum of \(S'\) in \(\ZN{N}\). Every \(s\in S\) is written \(s=s_0+s'/N\) with \(s'\in S'\), and we put \(z_s=\ex{s}\); the \(k\) numbers \(z_s\) are distinct points of the unit circle.

For the exchange arguments, assume that
\begin{equation}
\label{eq:nocp}
\text{no fiber } A_x \text{ of } \Omega \text{ contains two consecutive integers.}
\end{equation}
Under \eqref{eq:nocp}, Theorem~\ref{thm:forced} gives \(d_j=1\) on every loop and every backward edge, \(d_j=0\) on every forward edge, every run of every fiber has at most one element, and \(k=\nloop+\nback=n-\nfwd\); Lemma~\ref{lem:avatars} describes the runs: a loop \(j\) at \(c\) has the run \(\{p_j+1\}\) before \(c\) and \(\{p_j\}\) after \(c\); a forward edge \(j\) has the run \(\{p_j\}\) on the arc from \(c_{u(j)}\) to \(c_{v(j)}\) and is empty elsewhere; a backward edge \(j\) has the run \(\{q_j\}\) before \(c_{v(j)}\), is empty between \(c_{v(j)}\) and \(c_{u(j)}\), and has the run \(\{p_j\}=\{q_j-1\}\) after \(c_{u(j)}\). In every case the point that leaves the fiber at the class \(c_{v(j)}\) is \(q_j\) and the point that enters at \(c_{u(j)}\) is \(p_j\), so that \(\Out{c}=\{q_j: c_{v(j)}=c\}\) and \(\In{c}=\{p_j: c_{u(j)}=c\}\), and
\begin{equation}
\label{eq:crossing}
A_{c-}=C_c\sqcup\Out{c},\quad A_c=C_c\sqcup\In{c},\quad \lvert\Out{c}\rvert=\lvert\In{c}\rvert=m_c,\quad \lvert C_c\rvert=k-m_c,
\end{equation}
where \(C_c\) is the set of points untouched at \(c\). The unions are disjoint because the runs of one fiber are pairwise disjoint, so that the \(m_c\) leaving points are distinct, and because both fibers have \(k\) elements. Moreover, \(\Out{c}\cap\In{c}=\emptyset\): an equality \(q_i=p_j\) at the same class would give \(b_i=a_j\), contrary to the standing reduction.

The consecutive-pair statement gives \eqref{eq:nocp}: if some fiber \(A_x\) contains two consecutive integers, then \((A_x,S)\) is a spectral pair with a consecutive pair in \(A_x\), and \(\CP{k}\) (Definition~\ref{def:namedstatements}) gives \(S+1/m=S\) for some \(m\ge 2\), against the choice of \(S\). Since \(\CP{k}\) holds for \(k\le 10\) by Corollary~\ref{cor:cpten}, \eqref{eq:nocp} holds whenever \(k\le 10\), and by Theorem~D(i) also for \(k=11\).

\subsection{The transition space at a class}

\begin{definition}
\label{def:transition}
Let \(V=\C^S\) with the inner product \(\langle u,v\rangle=\sum_{s\in S}u(s)\overline{v(s)}\). For \(a\in\Z\) the exponential vector \(\ev{a}\in V\) is \(\ev{a}(s)=\ex{as}\), and \(\Dop\) is the diagonal unitary operator on \(V\) with entries \(z_s=\ex{s}\), so that \(\Dop\ev{a}=\ev{a+1}\) \(\forall a\in\Z\). For a class \(c\) the \emph{transition space} is \(W_c=\operatorname{span}\{\ev{o}: o\in\Out{c}\}\).
\end{definition}

Every coordinate of \(\ev{a}\) has modulus one, and \(\langle\ev{a},\ev{b}\rangle=\sum_{s\in S}\ex{(a-b)s}=\widehat{\mu_S}(a-b)\). By Definition~\ref{def:spectralpair}, a set \(A\subset\Z\) of \(k\) integers is a spectrum of \(S\) exactly when the matrix \(\bigl(\ex{sa}\bigr)_{s\in S,\,a\in A}\) is a complex Hadamard matrix, and a square matrix with unimodular entries and pairwise orthogonal rows has pairwise orthogonal columns; the columns are the vectors \(\ev{a}\). Hence
\begin{equation}
\label{eq:orthobasis}
A \text{ is a spectrum of } S \iff \{\ev{a}: a\in A\} \text{ is an orthogonal basis of } V.
\end{equation}

\begin{lemma}
\label{lem:transition}
Assume \eqref{eq:nocp} and let \(c\) be a class. Then
\[
W_c=\operatorname{span}\{\ev{o}: o\in\Out{c}\}=\operatorname{span}\{\ev{i}: i\in\In{c}\}=\bigl(\operatorname{span}\{\ev{f}: f\in C_c\}\bigr)^{\perp},
\]
and \(\dim W_c=m_c\). A forward edge \(j\) contributes the same vector \(\ev{p_j}\) to \(W_{c_{u(j)}}\) and to \(W_{c_{v(j)}}\). A backward edge \(j\) contributes \(\ev{q_j}\) to \(W_{c_{v(j)}}\) and \(\ev{p_j}=\Dop^{-1}\ev{q_j}\) to \(W_{c_{u(j)}}\). A loop \(j\) at \(c\) contributes both \(\ev{p_j+1}\) and \(\ev{p_j}=\Dop^{-1}\ev{p_j+1}\) to \(W_c\).
\end{lemma}

\begin{proof}
By \eqref{eq:crossing} and \eqref{eq:orthobasis}, both \(\{\ev{a}: a\in A_{c-}\}\) and \(\{\ev{a}: a\in A_c\}\) are orthogonal bases of \(V\) containing the family \(\{\ev{f}: f\in C_c\}\). The span of the remaining vectors of an orthogonal basis is the orthogonal complement of the span of the common part, in both cases, which gives the two equalities; the dimension is \(k-\lvert C_c\rvert=m_c\). The edge \(j\) contributes \(\ev{p_j}\) to \(W_{c_{u(j)}}\) through \(\In{c_{u(j)}}\) and \(\ev{q_j}\) to \(W_{c_{v(j)}}\) through \(\Out{c_{v(j)}}\), and \(q_j=p_j+d_j\) with the forced values \(d_j=0\) on forward edges and \(d_j=1\) on backward edges and loops; finally \(\ev{a}=\Dop^{-1}\ev{a+1}\).
\end{proof}

\subsection{Loop exchange}

\begin{lemma}
\label{lem:loopexchange}
Let \(\lvert S\rvert=k\ge 3\), let \(C\subset\Z\) have \(k-2\) elements, and let \(\lambda,a,b\in\Z\) be such that the two \(k\)-element sets \(C\cup\{\lambda+1,a\}\) and \(C\cup\{\lambda,b\}\) are spectra of \(S\). Then \(\ev{\lambda+1}\) is a unimodular multiple of \(\ev{b}\), \(\ev{\lambda}\) is a unimodular multiple of \(\ev{a}\), and \(C\cup\{\lambda,\lambda+1\}\) is a spectrum of \(S\). In particular \(\sum_{s\in S}\ex{s}=0\).
\end{lemma}

\begin{proof}
Put \(W=\{\ev{c}: c\in C\}^{\perp}\). By \eqref{eq:orthobasis}, \(W\) has dimension \(2\), and both \(\{\ev{\lambda+1},\ev{a}\}\) and \(\{\ev{\lambda},\ev{b}\}\) are orthogonal bases of \(W\). Write \(\ev{\lambda+1}=\alpha\ev{\lambda}+\beta\ev{b}\) with \(\alpha,\beta\in\C\). Evaluating at \(s\in S\) and dividing by \(\ev{\lambda}(s)\), which has modulus one,
\[
z_s=\alpha+\beta\,\ex{(b-\lambda)s}\qquad(s\in S).
\]
If \(\beta=0\), the \(k\) distinct numbers \(z_s\) all equal \(\alpha\), so \(k\le 1\). If \(\beta\ne 0\), every \(z_s\) lies on the circle \(\{\alpha+\beta w: \lvert w\rvert=1\}\) of center \(\alpha\) and radius \(\lvert\beta\rvert\) as well as on the unit circle; two distinct circles meet in at most two points, so \(k\le 2\) unless the two circles coincide, that is unless \(\alpha=0\) and \(\lvert\beta\rvert=1\). Since \(k\ge 3\), we conclude \(\ev{\lambda+1}=\beta\ev{b}\) with \(\lvert\beta\rvert=1\). Consequently \(\ev{\lambda+1}\) is orthogonal to \(\ev{\lambda}\), because \(\ev{b}\) is, and to every \(\ev{c}\), \(c\in C\). The \(k\) vectors \(\ev{c}\) (\(c\in C\)), \(\ev{\lambda}\), \(\ev{\lambda+1}\) are pairwise orthogonal and nonzero, hence an orthogonal basis of \(V\), and their indices are \(k\) distinct integers (\(\lambda\notin C\) and \(\lambda+1\notin C\) because the two given sets have \(k\) elements). By \eqref{eq:orthobasis}, \(C\cup\{\lambda,\lambda+1\}\) is a spectrum of \(S\), and \(0=\langle\ev{\lambda+1},\ev{\lambda}\rangle=\sum_{s\in S}\ex{s}\). Finally, \(\ev{\lambda}\) and \(\ev{a}\) span the same one-dimensional orthogonal complement of \(\ev{\lambda+1}\) in \(W\); since both have norm \(\sqrt{k}\), they differ by a unimodular factor.
\end{proof}

\begin{theorem}
\label{thm:noloop}
Assume \eqref{eq:nocp}, \(k\ge 3\), and that \(\CP{k}\) holds, as it does for \(k\le 10\) by Corollary~\ref{cor:cpten} and for \(k=11\) by Theorem~D(i). Then no class of size two carries a loop.
\end{theorem}

\begin{proof}
Suppose the class \(c\) has \(m_c=2\) and carries the loop \(j\). Under \eqref{eq:nocp} the loop has \(d_j=1\), and by the description of the runs its point is \(p_j+1\) before \(c\) and \(p_j\) after \(c\); thus \(p_j+1\in\Out{c}\) and \(p_j\in\In{c}\). By \eqref{eq:crossing} there are integers \(a\) and \(b\) with \(\Out{c}=\{p_j+1,a\}\), \(\In{c}=\{p_j,b\}\), and
\[
A_{c-}=C_c\cup\{p_j+1,a\},\qquad A_c=C_c\cup\{p_j,b\},\qquad \lvert C_c\rvert=k-2,
\]
both unions disjoint. Both fibers are spectra of \(S\). Lemma~\ref{lem:loopexchange}, applied with \(\lambda=p_j\), shows that \(C_c\cup\{p_j,p_j+1\}\) is a spectrum of \(S\). It contains two consecutive integers, so \(\CP{k}\) gives \(S+1/m=S\) for some integer \(m\ge 2\), a translation symmetry of \(S\), which contradicts Theorem~\ref{thm:periodic}.
\end{proof}

\begin{remark}
The argument uses that the transition space has dimension two. At a class of size \(m_c\ge 3\) carrying a loop, the vector \(\ev{p_j+1}\) lies in \(W_c\) together with \(\ev{p_j}=\Dop^{-1}\ev{p_j+1}\), but the orthogonal complement of \(\ev{p_j+1}\) in \(W_c\) has dimension \(m_c-1\ge 2\), so the two-circle argument no longer applies. Larger classes can carry loops: for any \(k\ge 3\), the set \(\Omega=2\{0,\dots,k-1\}+[0,1)\) has spectrum \(\Lambda=\frac1{2k}\{0,\dots,k-1\}+\Z\), since its fibers give the Fourier matrix of order \(k\). Here the minimal period is \(1\), the unique endpoint class has size \(k\), and all its edges are loops. The proof of Theorem~\ref{thm:fullymixing} uses loop exchange only at classes of size two.
\end{remark}

\subsection{Cosets: three and five points}

\begin{lemma}
\label{lem:cosets}
Let \(p\in\{2,3,5\}\), let \(S\subset\R/\Z\) have \(p\) elements, and suppose that some \(A\subset\Z\) is a spectrum of \(S\). Then \(S-S\subset\Q/\Z\), the prime \(p\) divides the integer \(N\) of Lemma~\ref{lem:reduction}, and, with \(K=N/p\), every spectrum \(A\subset\Z\) of \(S\) reduces modulo \(N\) onto a coset of the subgroup \(K\ZN{N}\) of order \(p\). This is \cite[Theorem~5.5]{bell2026four}.
\end{lemma}

\begin{theorem}
\label{thm:cosetclass}
Assume \eqref{eq:nocp} and \(k\in\{3,5\}\). Then some class has size \(m_c\ge k\).
\end{theorem}

\begin{proof}
Suppose that every class has \(m_c<k\). Put \(K=N/k\), an integer by Lemma~\ref{lem:cosets}. Every fiber \(A_x\) is a spectrum of \(S\), so by Lemma~\ref{lem:cosets} its reduction modulo \(N\) is a coset of \(K\ZN{N}\); as a set of integers, \(A_x\) is therefore contained in a single residue class modulo \(K\). As \(x\) increases through \(\R\), the fiber changes only at the points of \(c+h\Z\), \(c\) a class, and at such a point the fibers before and after share the \(k-m_c\ge 1\) untouched points of \eqref{eq:crossing}, so they lie in the same residue class modulo \(K\). By induction along the finitely many crossings in \([0,h]\), the fibers \(A_0\) and \(A_h\) lie in the same residue class modulo \(K\). But \(A_h=A_0-1\) and \(A_0\ne\emptyset\), so \(1\equiv 0\pmod K\), that is \(K=1\) and \(N=k\). Then \(S'\) is a \(k\)-element subset of \(\ZN{k}\), so \(S'=\ZN{k}\) and \(S=s_0+(1/k)\Z/\Z\), which has the translation symmetry \(1/k\), against Theorem~\ref{thm:periodic}.
\end{proof}

The proof uses only that a class moves fewer than \(k\) points and that a full turn shifts the fiber by one; for \(k=2\), \(m_c<2\) at every class makes every class rigid (Lemma~\ref{lem:rigidpair}), so there is nothing to prove.

\subsection{Lineage: four points}

\begin{lemma}
\label{lem:fourpairs}
Let \(S\subset\R/\Z\) have four elements and a spectrum \(A\subset\Z\), and let \(N\) and \(S'\) be as in Lemma~\ref{lem:reduction}. Then \(N\) is even; put \(M=N/2\). Exactly one nonzero element of \(S'\) is even, and every spectrum \(A\) of \(S'\) in \(\ZN{N}\) satisfies \(A+M=A\), so that \(A\) is the union of two cosets of the subgroup \(\{0,M\}\). This is \cite[Lemma~5.3 and Theorem~5.6]{bell2026four}.
\end{lemma}

In this subsection \(k=4\), and we write \(S'=\{0,e,o_1,o_2\}\) with \(e\) even and \(o_1,o_2\) odd, as Lemma~\ref{lem:fourpairs} allows. Let
\[
\mathcal{L}=\{t\in\ZN{N}: te\equiv 0 \text{ and } t(o_1-o_2)\equiv 0 \pmod N\},
\]
a subgroup of \(\ZN{N}\). It contains \(M\), because \(e\) and \(o_1-o_2\) are even. Let \(E\) be the diagonal unitary on \(V\) with \((Ev)(s)=(-1)^{s'}v(s)\). An \emph{antipodal pair} is a coset \(P=\{x,x+M\}\) of \(\{0,M\}\) in \(\ZN{N}\). For an integer \(y\) whose residue lies in \(P\) we have
\begin{equation}
\label{eq:planes}
\ev{y+N}=\ex{Ns_0}\,\ev{y},\qquad \ev{y+M}=\ex{Ms_0}\,E\ev{y},
\end{equation}
because \(\ex{Ns}=\ex{Ns_0}\ex{s'}=\ex{Ns_0}\) and \(\ex{Ms}=\ex{Ms_0}\ex{s'/2}=\ex{Ms_0}(-1)^{s'}\). Since \(E^2\) is the identity, the plane
\[
\Pi(P)=\operatorname{span}\{\ev{y},E\ev{y}\}
\]
is the same for every integer \(y\) with residue in \(P\): changing \(y\) by \(N\) multiplies both spanning vectors by a constant, and changing \(y\) by \(M\) exchanges them up to constants. It is a plane, because \(\langle\ev{y},E\ev{y}\rangle=\sum_{s'\in S'}(-1)^{s'}=0\). If \(A\subset\Z\) is a spectrum of \(S\) whose reduction is \(P_1\cup P_2\), the two elements of \(A\) over \(P_1\) give two orthogonal vectors of \(\Pi(P_1)\), which therefore span it, and likewise for \(P_2\); as the four vectors form an orthogonal basis of \(V\),
\begin{equation}
\label{eq:planecomplement}
\Pi(P_2)=\Pi(P_1)^{\perp}.
\end{equation}

\begin{lemma}
\label{lem:onestep}
Assume \eqref{eq:nocp} and \(k=4\), and let \(c\) be a class with \(m_c\le 3\). Write the reductions modulo \(N\) of \(A_{c-}\) and \(A_c\) as \(P_1\cup P_2\) and \(P_1'\cup P_2'\), unions of two antipodal pairs. Then, after renumbering, \(P_1'=P_1\) and \(P_2'=P_2+t\) for some \(t\in\mathcal{L}\).
\end{lemma}

\begin{proof}
Both fibers are spectra of \(S\) with four elements, so their reductions are unions of two antipodal pairs by Lemma~\ref{lem:fourpairs}. Since \(m_c\le 3\), the set \(C_c\) of \eqref{eq:crossing} is nonempty; let \(f\in C_c\). The antipodal pair containing the residue of \(f\) is a pair of both reductions; call it \(P_1=P_1'\). By \eqref{eq:planecomplement}, applied to both fibers, \(\Pi(P_2)=\Pi(P_1)^{\perp}=\Pi(P_2')\). Choose \(x\in A_{c-}\) with residue in \(P_2\) and \(y\in A_c\) with residue in \(P_2'\). Then \(\ev{y}\in\Pi(P_2')=\Pi(P_2)=\operatorname{span}\{\ev{x},E\ev{x}\}\), say \(\ev{y}=\alpha\ev{x}+\beta E\ev{x}\). Dividing coordinatewise by \(\ev{x}(s)\),
\[
\ex{(y-x)s}=\alpha+\beta(-1)^{s'}\qquad(s\in S).
\]
The right side takes the same value at the two even elements \(s'=0\) and \(s'=e\), so \(\ex{(y-x)s_0}=\ex{(y-x)(s_0+e/N)}\), that is \((y-x)e\equiv 0\pmod N\); and it takes the same value at \(s'=o_1\) and \(s'=o_2\), so \((y-x)(o_1-o_2)\equiv 0\pmod N\). Hence \(t=y-x\bmod N\) lies in \(\mathcal{L}\), and \(P_2'\), the antipodal pair of the residue of \(y\), equals \(P_2+t\).
\end{proof}

\begin{theorem}
\label{thm:lineage}
Assume \eqref{eq:nocp} and \(k=4\). Then some class has size \(m_c=4\).
\end{theorem}

\begin{proof}
Suppose that every class has \(m_c\le 3\). Choose \(x_0\in\R\) that is not congruent to a class modulo \(h\), and write the reduction of \(A_{x_0}\) as \(P_1\cup P_2\). As \(x\) increases from \(x_0\) to \(x_0+h\), the fiber changes exactly at the \(r\) points congruent to the classes, and Lemma~\ref{lem:onestep} applies at each of them (transported by \(A_{x+h}=A_x-1\) when the crossing point lies outside \([0,h)\)). We claim that after each crossing the reduction of the fiber is \((P_1+\tau_1)\cup(P_2+\tau_2)\) with \(\tau_1,\tau_2\in\mathcal{L}\). This holds before the first crossing with \(\tau_1=\tau_2=0\). At a crossing, the lemma keeps one of the two current pairs and replaces the other, \(P_i+\tau_i\) say, by \(P_i+\tau_i+t\) with \(t\in\mathcal{L}\); since \(\mathcal{L}\) is a subgroup, the claim persists, each pair keeping its lineage. After the last crossing the fiber is \(A_{x_0+h}=A_{x_0}-1\), whose reduction is \((P_1-1)\cup(P_2-1)\). Thus
\[
(P_1+\tau_1)\cup(P_2+\tau_2)=(P_1-1)\cup(P_2-1)
\]
as unions of two distinct antipodal pairs, and there are two cases. If \(P_1+\tau_1=P_1-1\), then \(\tau_1+1\in\{0,M\}\), so \(-1\in\tau_1+\{0,M\}\subset\mathcal{L}\) (recall \(M\in\mathcal{L}\)), and the condition \(te\equiv 0\) at \(t=-1\) gives \(e\equiv 0\pmod N\), which is false because \(e\) is a nonzero element of \(S'\). Otherwise \(P_1+\tau_1=P_2-1\) and \(P_2+\tau_2=P_1-1\). Writing \(P_i=x_i+\{0,M\}\), these read \(x_1+\tau_1-x_2+1\in\{0,M\}\) and \(x_2+\tau_2-x_1+1\in\{0,M\}\); adding, \(\tau_1+\tau_2+2\in\{0,M\}\), so \(-2\in\mathcal{L}\): \(2e\equiv 0\) and \(2(o_1-o_2)\equiv 0\pmod N\). As \(e\ne 0\) and \(o_1\ne o_2\) in \(\ZN{N}\), this forces \(e=M\) and \(o_1-o_2=M\), so \(S'=\{0,M,o_1,o_1+M\}=S'+M\) and \(S+1/2=S\), against Theorem~\ref{thm:periodic}.
\end{proof}

\begin{corollary}
\label{cor:smallk}
Assume \eqref{eq:nocp} and \(k\in\{3,4,5\}\). Then \(k\le\max_c m_c\).
\end{corollary}

\begin{proof}
Theorem~\ref{thm:cosetclass} for \(k\in\{3,5\}\) and Theorem~\ref{thm:lineage} for \(k=4\).
\end{proof}

\subsection{Monodromy of the transition spaces}

\begin{proposition}
\label{prop:monodromy}
\textup{(a)} Let \(b\ge 1\) be an integer and let \(W\subset V\) be a subspace of dimension \(m\) with \(\Dop^{b}W=W\) that contains a vector with no zero coordinate. Then the map \(s\mapsto z_s^{\,b}\) takes at most \(m\) values on \(S\), so that \(k\le bm\); and if \(k=bm\), then \(S\) is a union of \(m\) cosets of \((1/b)\Z/\Z\), so that \(S+1/b=S\).

\textup{(b)} Assume \eqref{eq:nocp}. Let \(c^{(0)},c^{(1)},\dots,c^{(\ell)}=c^{(0)}\) be classes, \(\ell\ge 1\), such that for each \(i\) every edge starting at \(c^{(i-1)}\) ends at \(c^{(i)}\). Then all these classes have the same size \(m\), and if \(b\) denotes the number of indices \(i\) for which the edges from \(c^{(i-1)}\) to \(c^{(i)}\) are backward edges or loops, then \(\Dop^{b}W_{c^{(0)}}=W_{c^{(0)}}\). Consequently \(k\le bm\) when \(b\ge 1\), and \(k=bm\) forces \(S+1/b=S\).
\end{proposition}

\begin{proof}
(a) Let \(w_1,\dots,w_t\) be the distinct values of \(z_s^{\,b}\), \(s\in S\), let \(S_i=\{s\in S: z_s^{\,b}=w_i\}\), and let \(Q_i\) be the orthogonal projection of \(V\) onto the coordinate subspace \(\C^{S_i}\). The operator \(\Dop^{b}\) is diagonal with entry \(z_s^{\,b}\) at \(s\), so \(Q_i=\prod_{j\ne i}(\Dop^{b}-w_j)/(w_i-w_j)\) is a polynomial in \(\Dop^{b}\), and \(Q_iW\subset W\). Let \(v\in W\) have no zero coordinate. The vectors \(Q_1v,\dots,Q_tv\) are nonzero, have pairwise disjoint supports, hence are linearly independent, and lie in \(W\); so \(t\le m\). The equation \(z^b=w_i\) has exactly \(b\) roots, and the \(z_s\) are distinct, so \(\lvert S_i\rvert\le b\) and \(k=\sum_i\lvert S_i\rvert\le bt\le bm\). If \(k=bm\), then \(t=m\) and every \(S_i\) has exactly \(b\) elements, so \(\{z_s: s\in S_i\}\) is the full set of \(b\)-th roots of \(w_i\), and \(S_i\) is a coset of \((1/b)\Z/\Z\) in \(\R/\Z\).

(b) The edges starting at \(c^{(i-1)}\) number \(m_{c^{(i-1)}}\) and all end at \(c^{(i)}\), where \(m_{c^{(i)}}\) edges end; so \(m_{c^{(i-1)}}\le m_{c^{(i)}}\) \(\forall i\), and around the closed sequence all sizes are equal, to \(m\) say. Then the edges ending at \(c^{(i)}\) are exactly the edges starting at \(c^{(i-1)}\). Since all of them join the same two classes (or are loops at the same class when \(c^{(i)}=c^{(i-1)}\)), they are all forward, with \(d_j=0\), or all backward or loops, with \(d_j=1\); write \(b_i\in\{0,1\}\) for this common value. By Lemma~\ref{lem:transition},
\[
W_{c^{(i)}}=\operatorname{span}\{\ev{q_j}: c_{v(j)}=c^{(i)}\}=\operatorname{span}\{\Dop^{b_i}\ev{p_j}: c_{u(j)}=c^{(i-1)}\}=\Dop^{b_i}W_{c^{(i-1)}}.
\]
Composing, \(W_{c^{(0)}}=W_{c^{(\ell)}}=\Dop^{b}W_{c^{(0)}}\) with \(b=\sum_i b_i\). The space \(W_{c^{(0)}}\) has dimension \(m\) and contains the vector \(\ev{p_j}\) of any edge \(j\) starting at \(c^{(0)}\), whose coordinates are unimodular. Part (a) gives the rest.
\end{proof}

\subsection{The class graphs at seven and eight intervals}

In this subsection \(\Omega\) is fully mixing with \(\dimQ\Krel\le 4\) and \(n\in\{7,8\}\). By Corollary~\ref{cor:fullymixing} the pattern is \((3,2,2)\) at \(n=7\) and \((2,2,2,2)\) at \(n=8\), and in both cases \(\rho=r\), where \(\rho=\dimQ\operatorname{span}_{\Q}(\ell_1,\dots,\ell_n)\). We first record what the maximal rank \(\rho=r\) means for the class graph and for the levels. For \(T\subset\nset\) and a class \(c\) let \(\varphi_T(c)\) be the number of edges of \(T\) ending at \(c\) minus the number starting at \(c\), the \emph{flow} of \(T\) at \(c\), and let \(d(T)=\sum_{j\in T}d_j\) be the \emph{weight} of \(T\). The set \(T\) is \emph{balanced} when \(\varphi_T=0\); a nonempty balanced set is an edge-disjoint union of directed cycles (a longest trail in it must close up, and removing the cycle so found leaves a balanced set).

\begin{lemma}
\label{lem:maximalrank}
Suppose \(\rho=r\). Then \(\clsgraph\) is connected, the numbers \(h,c_2-c_1,\dots,c_r-c_1\) are linearly independent over \(\Q\), and for \(T,T'\subset\nset\) one has \(\ell_T=\ell_{T'}\) iff \(d(T)=d(T')\) and \(\varphi_T=\varphi_{T'}\). In particular, \(\ell_T\in h\Z\) iff \(T\) is balanced, and every effective level is of the form \(\alpha=wh\), where the slice \(\slice{\alpha}\) consists of balanced sets of weight \(w\).
\end{lemma}

\begin{proof}
By Lemma~\ref{lem:arcs}, \(\ell_j=c_{v(j)}-c_{u(j)}+hd_j\). If \(\clsgraph\) had \(\kappa\ge 2\) connected components, every \(\ell_j\) would lie in the \(\Q\)-span of \(h\) together with the differences \(c_v-c_u\) between classes in the same component, a space of dimension at most \(1+(r-\kappa)<r\), against \(\rho=r\). So \(\clsgraph\) is connected, and Lemma~\ref{lem:rankbound} gives \(\operatorname{span}_{\Q}(\ell_1,\dots,\ell_n)=\operatorname{span}_{\Q}(h,c_2-c_1,\dots,c_r-c_1)\); since the left side has dimension \(r\), the \(r\) generators on the right are linearly independent. Summing the expression of \(\ell_j\) over \(j\in T\) gives
\[
\ell_T=\sum_{c}\varphi_T(c)\,c+h\,d(T)=\sum_{i=2}^{r}\varphi_T(c_i)(c_i-c_1)+h\,d(T),
\]
because \(\sum_c\varphi_T(c)=0\). By linear independence, \(\ell_T=\ell_{T'}\) exactly when the coefficients agree, and \(\ell_T\in h\Z\) exactly when \(\varphi_T=0\). Finally let \(\alpha\) be an effective level and \(\one=\sum_{T\in\slice{\alpha}}\gamma_T1_T\) with rational \(\gamma_T\). All members of the slice have the same flow \(\varphi\) and the same weight \(w\). The flow is linear in the indicator, and \(\varphi_{\nset}=0\) because \(\clsgraph\) is balanced; so \(0=\bigl(\sum_T\gamma_T\bigr)\varphi\), and \(\sum_T\gamma_T=L/\alpha\ne 0\) by Proposition~\ref{prop:pinned}. Hence \(\varphi=0\) and \(\alpha=\ell_T=wh\) \(\forall T\) in the slice.
\end{proof}

The classes are numbered by their positions, \(c_1<\dots<c_r\), and an edge from class \(u\) to class \(v\) is written \(u\to v\); it is forward when \(u<v\) and backward when \(u>v\). Under \eqref{eq:nocp}, which holds here since \(k\le\kzero(\ell)\le 10\) by Theorem~\ref{thm:kzeroten}, the integer parts are forced and \(k=n-\nfwd\). By Theorem~\ref{thm:forced}, every simple directed cycle through at least two classes contains a forward edge, namely the edge entering its highest class, and a backward edge, the edge leaving it.

\begin{lemma}
\label{lem:sevensix}
Let \(n=7\), so that the pattern is \((3,2,2)\). Then \(k\ne 6\).
\end{lemma}

\begin{proof}
Suppose \(k=6\), so \(\nfwd=1\). The non-loop edges of \(\clsgraph\) form a balanced set (the loops are balanced, and so is \(\clsgraph\)), hence an edge-disjoint union of simple directed cycles, each of which passes through at least two classes and therefore contains a forward edge. As there is only one forward edge, the non-loop edges form a single simple directed cycle. By Lemma~\ref{lem:maximalrank}, \(\clsgraph\) is connected, so every class is incident with a non-loop edge and the cycle passes through all three classes. A directed cycle through the classes \(1,2,3\) is either \(1\to 2\to 3\to 1\), with two forward edges, or \(3\to 2\to 1\to 3\), with the single forward edge \(1\to 3\) and the backward edges \(3\to 2\) and \(2\to 1\). So the latter occurs. The cycle contributes in-degree one at every class, and the remaining four edges are loops, \(m_c-1\) of them at the class \(c\): two at the class of size three and one at each class of size two. Thus a class of size two carries a loop, which Theorem~\ref{thm:noloop} excludes.
\end{proof}

\begin{remark}
The verdict rule does not exclude these configurations. Write \(Z\) for the set of the three edges of the cycle, of weight \(2\), and \(\lambda_1,\dots,\lambda_4\) for the loops. The balanced sets \(Z\cup\{\lambda_i\}\) and the triples of loops have weight \(3\), and
\[
\one=\tfrac14\sum_{i=1}^{4}1_{Z\cup\{\lambda_i\}}+\tfrac14\sum_{i<j<l}1_{\{\lambda_i,\lambda_j,\lambda_l\}},
\]
so the level \(3h=L/2\) is effective; the balanced sets \(Z\) and \(\{\lambda_i,\lambda_j\}\) have weight \(2\), and \(\one=1_Z+\tfrac13\sum_{i<j}1_{\{\lambda_i,\lambda_j\}}\), so the level \(2h=L/3\) is effective. Hence \(6\) divides \(\kzero(\ell)\), and the exclusion of \(k=6\) at \(n=7\) needs the fibers.
\end{remark}

\begin{lemma}
\label{lem:descent}
Let \(n=8\), so that the pattern is \((2,2,2,2)\). Then \(k\le 6\), and \(k\ne 6\). Hence \(k\notin\{6,7\}\).
\end{lemma}

\begin{proof}
No fiber contains a consecutive pair, \(\clsgraph\) is connected (Lemma~\ref{lem:maximalrank}) and therefore has a non-loop edge, and Theorem~\ref{thm:forced} gives \(k\le 7\). If \(k\le 2\) there is nothing to prove, so assume \(k\ge 3\). Every class has size two, so Theorem~\ref{thm:noloop} excludes loops, and \(\clsgraph\) is connected, loopless and balanced; it therefore has an Eulerian circuit, a closed trail using every edge once. Reading the circuit cyclically, a maximal run of consecutive backward arcs strictly decreases the position at each arc, so it contains at most \(r-1=3\) arcs, and two such runs are separated by at least one forward arc. Hence the number of backward runs is at most \(\nfwd\), and \(\nback\le 3\nfwd\). Since there are no loops, \(k=\nback=8-\nfwd\), and \(8-\nfwd\le 3\nfwd\) gives \(\nfwd\ge 2\) and \(k\le 6\).

Suppose \(k=6\). Then \(\nfwd=2\) and \(\nback=6=3\nfwd\), so equality holds throughout: there are exactly two backward runs, each of exactly three arcs, hence each the full descent \(4\to 3\to 2\to 1\), and each is followed by exactly one forward arc, which must lead from position \(1\) back to position \(4\). Thus \(\clsgraph\) consists of two copies of the cycle \(4\to 3\to 2\to 1\to 4\): it is the graph called (D4) in Appendix~\ref{app:cases}, with the backward arcs \(4\to 3\), \(3\to 2\), \(2\to 1\) doubled, each of weight \(1\), and the forward arc \(1\to 4\) doubled, of weight \(0\).

We determine the effective levels. In a balanced subset of this doubled cycle, the numbers of selected edges on the four successive links are equal. Every nonempty proper balanced subset therefore selects one edge from each parallel pair and has weight \(3\). By Lemma~\ref{lem:maximalrank}, \(3h=L/2\) is consequently the only possible effective level in \((0,L)\). It is effective: choose one edge from each parallel pair to form a cycle \(C\); its complement is another cycle of weight \(3\), and \(1_C+1_{C^c}=\one\). Thus \(\Rat=\{1/2\}\), \(\kzero(\ell)=2\), and Theorem~\ref{thm:verdict} gives \(k\mid 2\), against \(k=6\). Hence \(k\ne 6\), and \(k\notin\{6,7\}\).
\end{proof}

The monodromy proposition excludes the same configuration without the verdict rule. In the graph (D4) placed as above, the two edges starting at position \(4\) end at \(3\), those starting at \(3\) end at \(2\), those starting at \(2\) end at \(1\), and those starting at \(1\) end at \(4\). Proposition~\ref{prop:monodromy}(b) applies to the closed sequence of classes \(4,3,2,1,4\) with \(m=2\) and \(b=3\) backward steps: \(\Dop^{3}W_{c_4}=W_{c_4}\), so \(k\le 6\), and \(k=6\) forces \(S+1/3=S\), against Theorem~\ref{thm:periodic}. In this case \(z_s^{\,3}\) takes at most two values on the six points of \(S\), which must then be two full cosets of \((1/3)\Z/\Z\).

\subsection{The fully mixing theorem}

\begin{theorem}
\label{thm:fullymixing}
Let \(\Omega\) be a fully mixing spectral union of \(n\in\{6,7,8\}\) intervals with \(\dimQ\Krel\le 4\), with spectrum \(\Lambda\) of minimal period \(h^{-1}\) and \(k=L/h\). Then:
\begin{enumerate}
\item[\rm (i)] for \(n=8\), in the pattern \((2,2,2,2)\), \(k\le 2\);
\item[\rm (ii)] for \(n=7\), in the pattern \((3,2,2)\), \(k\le 3\);
\item[\rm (iii)] for \(n=6\), \(k\le 4\); moreover \(k=4\) occurs only in the pattern \((4,2)\), \(k=3\) only in the patterns \((4,2)\) and \((3,3)\), and \(k\le 2\) in the pattern \((2,2,2)\).
\end{enumerate}
In every case \(k\le\max_c m_c\).
\end{theorem}

\begin{proof}
By Theorem~\ref{thm:kzeroten}, \(k\) divides \(\kzero(\ell)\le 10\), so \(k\le 10\) and no fiber contains a consecutive pair. By Corollary~\ref{cor:fullymixing}, \(\rho\ge 2\) in every fully mixing pattern with \(n\ge 6\); if every edge of \(\clsgraph\) were a loop, every \(\ell_j\) would lie in \(h\Z\) and \(\rho\) would be \(1\); so \(\clsgraph\) has a non-loop edge, and Theorem~\ref{thm:forced} gives \(k\le n-1\). The possible patterns are those of Corollary~\ref{cor:fullymixing}.

(i) \(n=8\): \(k\le 7\), and \(k\notin\{6,7\}\) by Lemma~\ref{lem:descent}. Every class has size two, so \(k\notin\{3,5\}\) by Theorem~\ref{thm:cosetclass} and \(k\ne 4\) by Theorem~\ref{thm:lineage}. Hence \(k\le 2\).

(ii) \(n=7\): \(k\le 6\), and \(k\ne 6\) by Lemma~\ref{lem:sevensix}. Every class has size at most three, so \(k\ne 5\) by Theorem~\ref{thm:cosetclass} and \(k\ne 4\) by Theorem~\ref{thm:lineage}. Hence \(k\le 3\).

(iii) \(n=6\): \(k\le 5\). Every class has size at most four, so \(k\ne 5\) by Theorem~\ref{thm:cosetclass}. If \(k=4\), Theorem~\ref{thm:lineage} gives a class of size four, which exists only in the pattern \((4,2)\). If \(k=3\), Theorem~\ref{thm:cosetclass} gives a class of size at least three, which exists only in the patterns \((4,2)\) and \((3,3)\). In the pattern \((2,2,2)\), therefore, \(k\le 2\). The pattern \((2,2,2)\) admits \(\rho\in\{2,3\}\); nothing above depends on the rank, since Theorem~\ref{thm:kzeroten} is a statement about \(\ell\) alone and Theorems~\ref{thm:cosetclass} and~\ref{thm:lineage} use only the crossing structure \eqref{eq:crossing}.

The final assertion follows from \(\max_c m_c=2\) in (i), \(\max_c m_c=3\) in (ii), and the restrictions on the patterns in (iii).
\end{proof}

\begin{remark}
\label{rem:openmechanisms}
The bound \(k\le\max_c m_c\) does not extend to all fully mixing configurations without consecutive pairs; the example below has \(k=6\) and \(\max_c m_c=3\). The obstruction is visible in Theorems~\ref{thm:cosetclass} and~\ref{thm:lineage}: their proofs use the relation between spectra that share a point. For three and five points their reductions are the same coset; for four points the shared antipodal pair \(P_1\) determines the plane \(\Pi(P_1)^{\perp}\) of the other pair. The six-point dichotomy of Theorem~\ref{thm:sixdichotomy} describes each spectrum separately as three cosets of the subgroup of order two or two cosets of the subgroup of order three, and both shapes can occur for one set (Remark~\ref{rem:digitbox}). If two spectra with three antipodal pairs share one pair, their remaining pairs span the same four-dimensional complement in \(\C^S\), but that complement does not determine the pairs; a crossing may also pass from one shape to the other.
\end{remark}

\begin{remark}
\label{rem:omegatwelve}
Let \(S'=\{0,1,8,9,17,28\}\subset\ZN{36}\), and put
\[
\begin{aligned}
A_0&=\{2,12,14,24,26,36\},& A_1&=\{2,8,14,20,26,32\},\\
A_2&=\{2,5,8,20,23,26\},& A_3&=\{2,5,17,20,23,35\},\\
A_4&=\{5,11,17,23,29,35\}.&&
\end{aligned}
\]
Each \(A_i\) reduces to a spectrum of \(S'\): every nonzero difference has order in \(\{2,3,6,12,18\}\), and the character sum of \(S'\) vanishes at these orders. For orders \(2\), \(3\) and \(6\), this follows from its equal residue counts; modulo \(12\) it is two cosets of \(\{0,4,8\}\), and modulo \(18\) it is three antipodal pairs. The set \(S'\) has no period, since it contains \(0\) but neither \(12\) nor \(18\). Define
\[
\Omega_{12}=\bigcup_{i=0}^{4}\bigl(A_i+[i/5,(i+1)/5)\bigr).
\]
Its fibers for \(h=1\) on the five successive fifths of \([0,1)\) are \(A_0,\dots,A_4\), so \(\Lambda=S'/36+\Z\) is a spectrum of minimal period \(1\) and \(k=6\). At the four interior crossings the numbers of entering points are \(3,2,2,2\); at \(0\), comparing \(A_4+1\) with \(A_0\), the number is \(3\). Thus \(\Omega_{12}\) has twelve components and class pattern \((3,3,2,2,2)\). No fiber contains consecutive integers. Since \(1/36\in(S'-S')/36\), a rigid slide would be a nonzero multiple of \(36\), whereas every difference between an ending and a starting height at one class has absolute value less than \(36\). The configuration is therefore fully mixing, with \(\max_c m_c=3<6=k\). Its lengths are commensurable, giving \(\dimQ K_{\ell}=11\), so Theorem~\ref{thm:A} does not apply.
\end{remark}

\section{Proof of Theorem A and Corollary B}
\label{sec:assembly}

The pieces are now in place: the reduction to fully mixing configurations with at most eight components (Section~\ref{sec:fibers} and Proposition~\ref{prop:reduction}), the bound \(k \mid \kzero(\ell) \le 10\) from Theorem~\ref{thm:kzeroten}, the consecutive-pair statements \(\CP{k}\) for \(k \le 10\) (Proposition~\ref{prop:cpsmall} and Theorems~\ref{thm:cpseven} to~\ref{thm:cpten}), the fiber-side exclusions for six, seven and eight components (Theorem~\ref{thm:fullymixing}), and the fiber lemma for at most five cosets (Theorem~\ref{thm:fiberfew}). For four components the bound on \(k\) is also proved in \cite{bell2026four}, and we record it.

\begin{corollary}[{\cite[Corollary~6.5]{bell2026four}}]
\label{cor:kfour}
Let \(\Omega\) be a union of four pairwise disjoint non-adjacent bounded intervals, and let \(\Lambda\) be a spectrum of \(\Omega\) with minimal period \(h^{-1}\). Then \(k = \lvert\Omega\rvert/h \le 4\).
\end{corollary}

\begin{proof}[Proof of Theorem~\ref{thm:A}]
Let \(\Omega\) be a bounded measurable set which agrees up to a null set with a finite union of pairwise disjoint bounded intervals with lengths \(\ell_1, \dots, \ell_n\) satisfying \(\dimQ \Krel \le 4\), and suppose that \(\Omega\) is spectral. Apply the standing reduction of Definition~\ref{def:standing}. It changes \(\Omega\) only by a null set, which affects neither spectrality nor tiling, and it replaces the length vector by one obtained from \(\ell\) by deleting the lengths of the discarded null intervals and by replacing groups of adjacent intervals by single intervals, whose lengths are the sums of the original ones; by the computation in the proof of Lemma~\ref{lem:mergerank}, each such replacement leaves the dimension of the space of rational relations unchanged or lowers it. We may therefore assume that \(\Omega = \bigsqcup_{j=1}^{n} [a_j, b_j)\) with pairwise disjoint non-adjacent components, \(\ell_j > 0\), \(L = \sum_j \ell_j\) and \(\dimQ \Krel \le 4\). Let \(\Lambda\) be a spectrum of \(\Omega\). Theorem~\ref{thm:periodic} gives the minimal period \(h^{-1}\), the integer \(k\), the set \(\Gamma\) with \(\Lambda = \Gamma + h^{-1}\Z\), and the set \(S = h\Gamma \subset \R/\Z\) of \(k\) points without translation symmetry; by Lemma~\ref{lem:everyfiber}, \(L = kh\) and every fiber of \(\Omega\) is a spectrum of \(S\).

We first show that \(k \le 5\). By Proposition~\ref{prop:reduction}, either \(k = 1\), or there is a fully mixing spectral union \(\Omega_r\) of \(n_r\) components with \(2 \le n_r \le 8\), with the same spectrum \(\Lambda\), hence the same \(h\), \(k\) and \(S\), and whose length vector \(\ell^{(r)}\) satisfies \(\dimQ K_{\ell^{(r)}} \le 4\). In the second case Theorem~\ref{thm:kzeroten}, applied to \(\Omega_r\), gives \(k \mid \kzero(\ell^{(r)}) \le 10\), so \(k \le 10\). Suppose that \(k > n_r\). Then Lemma~\ref{lem:pigeonhole}(i), applied to \(\Omega_r\), says that every fiber of \(\Omega_r\) contains two consecutive integers. Fix one fiber \(A\) of \(\Omega_r\); by Lemma~\ref{lem:everyfiber} it is a spectrum of \(S\), so \((A, S)\) is a spectral pair with \(A \subset \Z\), \(\lvert A\rvert = k\), \(3 \le k \le 10\), and \(A\) contains a consecutive pair. The statement \(\CP{k}\) holds for every such \(k\): for \(k \le 6\) by Proposition~\ref{prop:cpsmall}, and for \(k = 7, 8, 9, 10\) by Theorems~\ref{thm:cpseven}, \ref{thm:cpeight}, \ref{thm:cpnine} and~\ref{thm:cpten}. Hence \(S + 1/m = S\) for some integer \(m \ge 2\), a nontrivial translation symmetry of \(S\), which Theorem~\ref{thm:periodic} excludes. Therefore \(k \le n_r\). If \(n_r \le 5\), this is the bound \(k \le 5\); for \(n_r = 4\) it is also Corollary~\ref{cor:kfour}. If \(n_r \in \{6, 7, 8\}\), Theorem~\ref{thm:fullymixing}, applied to the fully mixing configuration \(\Omega_r\), gives \(k \le 4\). In every case \(k \le 5\).

It remains to produce the tiling, and for this we return to \(\Omega\) itself. By Lemma~\ref{lem:everyfiber}, \(\sum_{m \in \Z} 1_\Omega(x + mh) = k\) \(\forall x\), and \(\Lambda = \Gamma + h^{-1}\Z\) is a spectrum of \(\Omega\) with \(\lvert\Gamma\rvert = k \le 5\), the elements of \(\Gamma\) being distinct modulo \(h^{-1}\Z\). Theorem~\ref{thm:fiberfew}, applied to \(F = \Omega\), gives a set \(B \subset \Z\), depending only on \(S\), with \(\sum_{b \in B} 1_\Omega(t - hb) = 1\) for almost every \(t \in \R\). Thus the translates \(\Omega + hb\), \(b \in B\), are pairwise disjoint up to null sets and cover \(\R\) up to a null set: \(\Omega\) tiles \(\R\) by translations, by the translation set \(hB\). When \(k = 1\) the set \(B\) is \(\Z\) and the tiling is by the lattice \(L\Z\), as Proposition~\ref{prop:mergechain}(i) also shows directly.
\end{proof}

The proof uses the fiber lemma only for \(k \le 5\), and the consecutive-pair statement only for \(k \le 10\). The case \(k = 5\) arises only when \(n_r = 5\), that is, for configurations that merge along rigid pairs to a five-interval union; in a fully mixing configuration with six, seven or eight components the proof gives \(k \le 4\), and in fact \(k \le 2\), \(k \le 3\) and \(k \le 4\) at \(n_r = 8\), \(7\) and \(6\), by Theorem~\ref{thm:fullymixing}.

\begin{proof}[Proof of Corollary~\ref{cor:B}]
After the standing reduction, \(\Omega\) agrees up to a null set with a union of \(n \le 5\) pairwise disjoint non-adjacent intervals of positive lengths \(\ell_1, \dots, \ell_n\). The lengths span a rational vector space of dimension \(\rho \ge 1\), so \(\dimQ \Krel = n - \rho \le n - 1 \le 4\), and Theorem~\ref{thm:A} applies: if \(\Omega\) is spectral, it tiles \(\R\) by translations, and the number \(k\) of cosets of its spectrum is at most \(5\).
\end{proof}

The bound \(k \le 5\) in Theorem~\ref{thm:A} cannot be improved to \(k \le 4\), even for six or more intervals, and the examples that attain it all contain a rigid pair.

\begin{figure}[tp]
\centering
\input{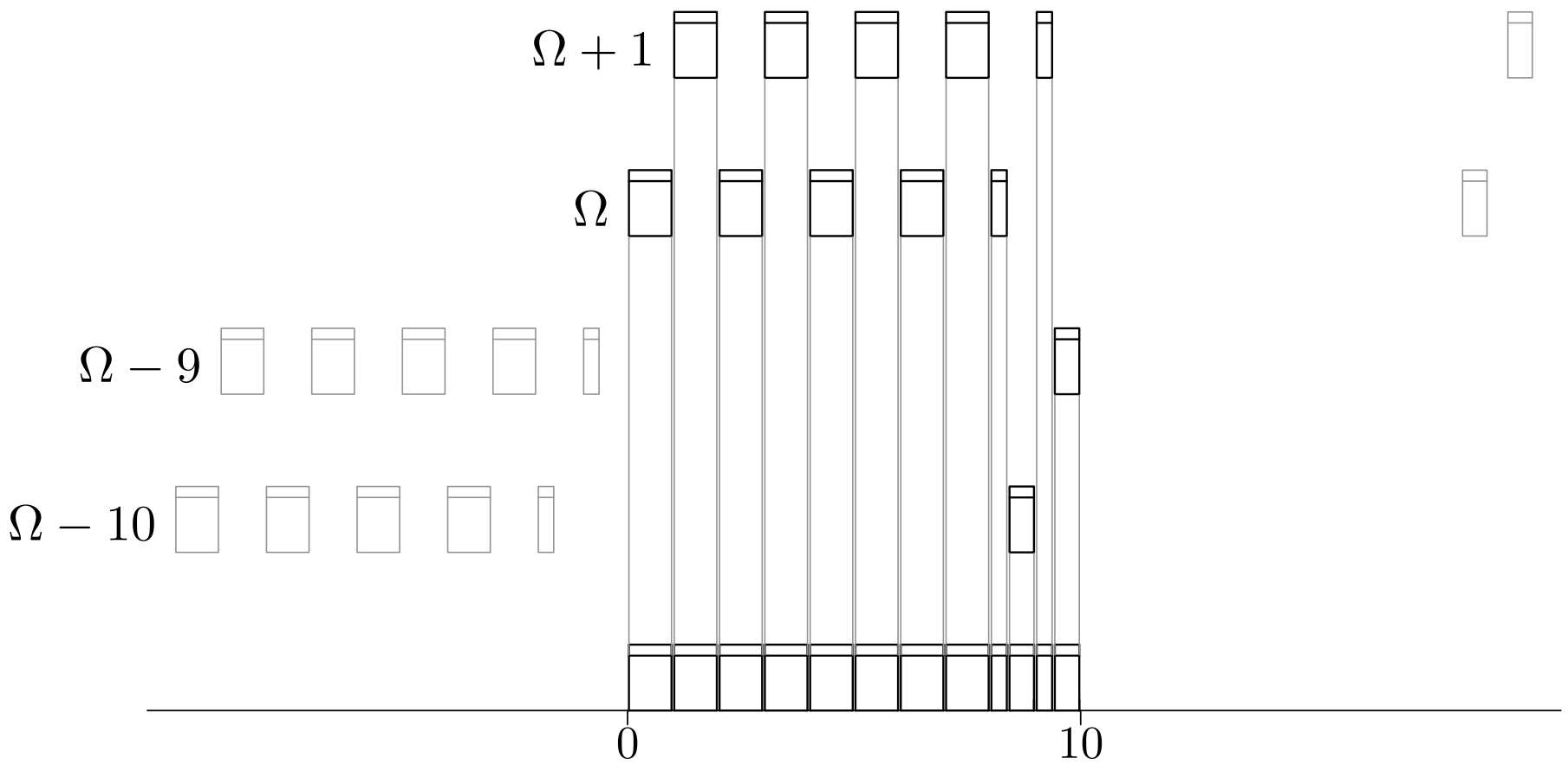}
\caption{The set \(\Omega\) of Proposition~\ref{prop:sharp} and its translates by \(1\), \(-9\) and \(-10\), raised to separate levels. The translates by \(\{0, 1\} + 10\Z\) tile the real line: the pieces lying over one period \([0, 10)\), between the marks \(0\) and \(10\) on the base line, fill it exactly, and the parts outside the period are gray.}
\label{fig:tiling}
\end{figure}

\begin{proposition}
\label{prop:sharp}
Let \(a \in (0, 1)\) be irrational and let
\[
\Omega = \bigl(\{0, 2, 4, 6\} + [0, 1)\bigr) \cup [8, 8 + a) \cup [18 + a, 19).
\]
Then \(\Omega\) is a union of six pairwise disjoint non-adjacent intervals, with lengths \((1, 1, 1, 1, a, 1 - a)\) and \(\dimQ \Krel = 4\); it is spectral, with the spectrum \(\Lambda = \tfrac{1}{10}\{0, 1, 2, 3, 4\} + \Z\) of minimal period \(1\), so that \(h = 1\) and \(k = 5\); the pair formed by \([8, 8 + a)\) and \([18 + a, 19)\) is rigid, with slide \(-10\), and the slide along it produces the five-interval set \(\{0, 2, 4, 6, 8\} + [0, 1)\), all of whose endpoints lie in one class; and \(\Omega\) tiles \(\R\) by the translation set \(\{0, 1\} + 10\Z\). For \(n \in \{7, 8\}\) there are spectral unions of \(n\) intervals with \(\dimQ \Krel = 4\) and \(k = 5\) as well. In every spectral union of \(n \ge 6\) intervals with \(\dimQ \Krel \le 4\) and \(k = 5\), some pair is rigid.
\end{proposition}

\begin{proof}
The six intervals \([0,1)\), \([2,3)\), \([4,5)\), \([6,7)\), \([8, 8+a)\) and \([18+a, 19)\) are pairwise disjoint and non-adjacent, since \(8 + a < 18 + a\) and \(a\) is not an integer. The rational span of the lengths is \(\Q + \Q a\), of dimension \(2\) because \(a\) is irrational, so \(\dimQ \Krel = 6 - 2 = 4\); the relations are the three equalities among the unit lengths and \(a + (1 - a) = 1\).

Take \(h = 1\), \(\Gamma = \tfrac{1}{10}\{0, 1, 2, 3, 4\}\) and \(S = \Gamma \bmod 1\), whose elements are distinct modulo \(\Z\). For \(x \in [0, 1)\) the fiber \(A_x = \{m \in \Z : x + m \in \Omega\}\) contains \(0, 2, 4, 6\), contains \(8\) exactly when \(x < a\), and contains \(18\) exactly when \(x \ge a\); so \(A_x = \{0, 2, 4, 6, 8\}\) for \(x < a\) and \(A_x = \{0, 2, 4, 6, 18\}\) for \(x \ge a\). Both sets are congruent to \(\{0, 2, 4, 6, 8\}\) modulo \(10\), and for \(\sigma = j/10\) with \(1 \le \lvert j\rvert \le 4\), which are the nonzero elements of \(S - S\),
\[
\sum_{t=0}^{4} \ex{(j/10)\, 2t} = \sum_{t=0}^{4} \ex{jt/5} = 0,
\]
since \(5 \nmid j\); the exponentials \(\ex{\sigma m}\) depend only on \(m \bmod 10\). Thus every fiber is a spectrum of \(S\), and Theorem~\ref{thm:fiber} says that \(\Lambda = \Gamma + \Z\) is a spectrum of \(\Omega\), with \(L = 5 = kh\). If \(1/m\) were a period of \(\Lambda\) with \(m \ge 2\), then \(S + 1/m = S\); but \(S + t = S\) forces \(t \in S - S = \tfrac{1}{10}\{-4, \dots, 4\}\), and \(\{0, \dots, 4\} + j \equiv \{0, \dots, 4\} \pmod{10}\) only for \(j = 0\). So the minimal period is \(1\), \(h = 1\) and \(k = 5\).

The endpoints \(8 + a\) and \(18 + a\) are congruent modulo \(\Z\), and the pair \((i, j)\) with \([a_i, b_i) = [8, 8 + a)\) and \([a_j, b_j) = [18 + a, 19)\) has slide \(d = (b_i - a_j)/h = -10\), which satisfies \(d(S - S) \subset \Z\); the pair is rigid. The slide of Lemma~\ref{lem:slide} moves \([18 + a, 19)\) to \([8 + a, 9)\), which completes \([8, 8 + a)\) to \([8, 9)\), and the result is \(\{0, 2, 4, 6, 8\} + [0, 1)\), whose ten endpoints are integers and lie in the single class \(0\). The other pair of congruent endpoints, \(b_j = 19\) and \(a_i = 8\), has \((19 - 8)/h = 11\), and \(11 \cdot \tfrac{1}{10} \notin \Z\): sliding \([8, 8 + a)\) by \(11\) would not preserve the spectrum, and indeed the pair \(([8, 8+a), [18+a, 19))\) is rigid only in the one order.

For the tiling, \(\{0, 2, 4, 6, 8\} \oplus \{0, 1\} = \{0, 1, \dots, 9\}\), so \(A_x \oplus (\{0, 1\} + 10\Z) = \Z\) \(\forall x\), because both fibers are congruent to \(\{0, 2, 4, 6, 8\}\) modulo \(10\). As in the proof of Lemma~\ref{lem:ucfl}, this gives \(\sum_{b \in \{0,1\} + 10\Z} 1_\Omega(t - b) = 1\) \(\forall t \in \R\), which is the tiling by \(\{0, 1\} + 10\Z\). Figure~\ref{fig:tiling} shows it.

For \(n = 7\), let \(b \in (0, 1)\) be such that \(1, a, b\) are linearly independent over \(\Q\), and replace the component \([6, 7)\) by \([6, 6 + b) \cup [16 + b, 17)\); for \(n = 8\), let \(c \in (0,1)\) be such that \(1, a, b, c\) are linearly independent over \(\Q\), and replace in addition \([4, 5)\) by \([4, 4 + c) \cup [14 + c, 15)\). The resulting intervals are pairwise disjoint and non-adjacent, since \(8 + a < 14 + c < 15 < 16 + b < 17 < 18 + a\), the rational span of the lengths has dimension \(3\), respectively \(4\), so \(\dimQ \Krel = 7 - 3 = 4\), respectively \(8 - 4 = 4\), and every fiber over \([0, 1)\) is again congruent to \(\{0, 2, 4, 6, 8\}\) modulo \(10\), by the same computation: the point \(6\) of the fiber is replaced by \(16\) when \(x \ge b\), and the point \(4\) by \(14\) when \(x \ge c\). The same \(\Lambda\) is a spectrum, with \(k = 5\).

Finally, a fully mixing spectral union of \(n \in \{6, 7, 8\}\) intervals with \(\dimQ \Krel \le 4\) has \(k \le 4\) by Theorem~\ref{thm:fullymixing}, and by Corollary~\ref{cor:fullymixing} a fully mixing configuration with \(\dimQ \Krel \le 4\) has at most eight components. So a spectral union of \(n \ge 6\) intervals with \(\dimQ \Krel \le 4\) and \(k = 5\) is not fully mixing: some pair is rigid.
\end{proof}

\section{Eleven points, and the failure at twelve}
\label{sec:twelve}

This section proves parts (i)--(iii) of Theorem~D. Section~\ref{sec:consecutive} proved the consecutive-pair statement up to ten points, which is all that the proof of Theorem~\ref{thm:A} uses. Here the same method settles eleven points, and at twelve points the statement fails. Every failure with a periodic spectrum arises by gluing smaller spectral pairs, and we determine when the glued set is aperiodic.

\subsection{Eleven points}

\begin{theorem}
\label{thm:cpeleven}
Let \((A, S)\) be a spectral pair with \(A \subset \Z\), \(S \subset \R/\Z\), \(\lvert A\rvert = \lvert S\rvert = 11\) and \(\widehat{\muS}(1) = 0\), which holds when \(A\) contains two consecutive integers. Then \(S = s_0 + \tfrac1{11}\Z/\Z\) for some \(s_0\).
\end{theorem}

\begin{proof}
With the conventions of Section~\ref{subsec:seventen}, the partitions of \(11\) into parts from \(\{2, 3, 5, 6, 7, 8, 9, 10, 11\}\) are
\(11\), \(9 + 2\), \(8 + 3\), \(7 + 2 + 2\), \(6 + 5\), \(6 + 3 + 2\), \(5 + 3 + 3\), \(5 + 2 + 2 + 2\), \(3 + 3 + 3 + 2\) and \(3 + 2 + 2 + 2 + 2\); the part \(10\) does not occur because \(1\) is not a weight. Since \(11\) is prime, \(S\) has a translation symmetry iff its roots form a rotation of \(\mu_{11}\), which is the conclusion. The remaining shapes are the following fifteen, and for each of them a residue class \(m\Z\) with \(m \le 7\) disjoint from \(Z\) is exhibited, so that \(\lvert A\rvert \le 7 < 11\) by Lemma~\ref{lem:clique}(i).

(1) \(\rho T(7, X)\), \(\lvert X\rvert = 4\): \(Z\) consists of odd integers; \(Z \cap 2\Z = \emptyset\).

(2) \(\rho T_{7,5,3}(y)\), \(y \ne 0\): by Lemma~\ref{lem:blockpowersums}, the power sum of \(T_{7,5,3}(y)\) is
\[
7[7 \mid t] - 1 - \zroot{7}^{yt} + (-1)^{t}(5[5 \mid t] - 1) + \zroot{7}^{yt}(\zroot{6}^{t} + \zroot{6}^{-t}).
\]
For even \(t\) this equals \(7[7 \mid t] + (5[5 \mid t] - 2) + c(t)\zroot{7}^{yt}\), with \(5[5 \mid t] - 2 \in \{3, -2\}\) and \(c(t) \in \{1, -2\}\). If \(7 \nmid t\), the roots \(1\) and \(\zroot{7}^{yt}\) are distinct and carry nonzero rational coefficients, so the sum is nonzero by Lemma~\ref{lem:rationalmoduli}(c). If \(7 \mid t\), the sum is at least \(7 - 2 - 2 > 0\). Thus \(Z \cap 2\Z = \emptyset\).

(3) \(\rho T_{7,(5,3)}(W)\), \(W = u^{-1}T_5 \ni 1\): the power sum of \(T_{7,(5,3)}(W)\) is
\[
7[7 \mid t] - 1 + (-1)^{t}\bigl(P_W(t) - 1\bigr), \qquad P_W(t) = u^{-t}P_{T_5}(t).
\]
For \(6 \mid t\) this is \(7[7 \mid t] - 2 + u^{-t}P_{T_5}(t)\) with \(P_{T_5}(t) \in \{1, 6\}\). A zero would need \(\lvert P_{T_5}(t)\rvert = 2\) when \(7 \nmid t\) and \(\lvert P_{T_5}(t)\rvert = 5\) when \(7 \mid t\), neither of which holds. Thus \(Z \cap 6\Z = \emptyset\) and \(\lvert A\rvert \le 6\).

(4) \(\rho T(5, X) \cup \{\pm\xi\}\), \(\lvert X\rvert = 4\): \(\widehat{\muS}(t) = \rho^{t}P(t) + \xi^{t}(1 + (-1)^{t})\) with \(P(t) = 5[5 \mid t] + c(t)E_X(t)\), where \(E_X(t) = 4\) when \(5 \mid t\) and \(\lvert E_X(t)\rvert = 1\) when \(5 \nmid t\). For odd \(t\) with \(3 \mid t\), \(c(t) = -3\) and \(P(t) = -7\) or \(\lvert P(t)\rvert = 3\), nonzero. For \(6 \mid t\), \(c(t) = 1\) and \(P(t) = 9\) or \(\lvert P(t)\rvert = 1\), both different from the modulus \(2\) of the other term. Thus \(Z \cap 3\Z = \emptyset\).

(5) \(\rho T(7, X) \cup \{\pm\xi\}\), \(\lvert X\rvert = 2\): now \(P(t) = 7[7 \mid t] + c(t)E_X(t)\) with \(E_X(t) = 2\) when \(7 \mid t\) and \(\lvert E_X(t)\rvert^{2}\) irrational when \(7 \nmid t\), by Lemma~\ref{lem:rationalmoduli}(b). For odd \(t\) with \(3 \mid t\), \(P(t) = 7 - 6 = 1\) or \(P(t) = -3E_X(t) \ne 0\). For \(6 \mid t\), \(P(t) = 9\) or \(\lvert P(t)\rvert^{2} = \lvert E_X(t)\rvert^{2}\) is irrational, so \(\lvert P(t)\rvert \ne 2\). Thus \(Z \cap 3\Z = \emptyset\).

(6) \(\rho T(5, X) \cup \eta\mu_3\), \(\lvert X\rvert = 3\): \(\widehat{\muS}(t) = \rho^{t}P(t) + 3\eta^{t}[3 \mid t]\), with \(E_X(t) = 3\) when \(5 \mid t\) and \(E_X(t) = -(\zroot{5}^{at} + \zroot{5}^{bt})\) for the two elements \(a, b \notin X\) when \(5 \nmid t\), of irrational squared modulus. For even \(t\) with \(3 \nmid t\), \(c(t) = -2\) and \(P(t) = 5 - 6 = -1\) or \(P(t) = -2E_X(t) \ne 0\). For \(6 \mid t\), \(P(t) = 8\) or \(\lvert P(t)\rvert^{2}\) is irrational, so \(\lvert P(t)\rvert \ne 3\). Thus \(Z \cap 2\Z = \emptyset\).

(7) \(\rho T_7 \cup \eta\mu_3\): \(\widehat{\muS}(t) = \rho^{t}(7[7 \mid t] + c(t)) + 3\eta^{t}[3 \mid t]\). For even \(t\) with \(3 \nmid t\) the coefficient of \(\rho^t\) is \(7[7 \mid t] - 2 \in \{5, -2\}\); for \(6 \mid t\) it is \(7[7 \mid t] + 1 \in \{8, 1\}\), of modulus different from \(3\). Thus \(Z \cap 2\Z = \emptyset\).

(8) \(\rho T(5, X) \cup \{\pm\eta_1\} \cup \{\pm\eta_2\}\), \(\lvert X\rvert = 2\): \(\widehat{\muS}(t) = \rho^{t}P(t) + (1 + (-1)^{t})(\eta_1^{t} + \eta_2^{t})\), with \(E_X(t) = 2\) when \(5 \mid t\) and \(\lvert E_X(t)\rvert^{2} \in \{(3 + \sqrt5)/2, (3 - \sqrt5)/2\}\) when \(5 \nmid t\), the two values of \(2 + 2\cos(2\pi v/5)\) for \(5 \nmid v\). For odd \(t\) with \(3 \mid t\), \(P(t) = 5 - 6 = -1\) or \(P(t) = -3E_X(t) \ne 0\). For \(6 \mid t\), \(c(t) = 1\); if \(5 \mid t\), then \(P(t) = 7 > 4 \ge 2\lvert \eta_1^{t} + \eta_2^{t}\rvert\); if \(5 \nmid t\), a zero would give \(\lvert E_X(t)\rvert^{2} = 4\lvert \eta_1^{t} + \eta_2^{t}\rvert^{2} = 4(2 + \psi^{t} + \psi^{-t})\) with \(\psi = \eta_1/\eta_2\), so \(\psi^{t} + \psi^{-t} = (-13 \pm \sqrt5)/8\), a number whose two conjugates sum to \(-13/4\); an algebraic integer has an integer trace, so this is impossible. Thus \(Z \cap 3\Z = \emptyset\).

(9) \(\rho\mu_7 \cup \{\pm\eta_1\} \cup \{\pm\eta_2\}\): for \(7 \mid t\), \(\widehat{\muS}(t) = 7\rho^{t} + (1 + (-1)^{t})(\eta_1^{t} + \eta_2^{t})\), and the second term has modulus at most \(4\). Thus \(Z \cap 7\Z = \emptyset\) and \(\lvert A\rvert \le 7\).

(10) \(\rho T_5 \cup \eta\mu_5\): \(\widehat{\muS}(t) = \rho^{t}P_{T_5}(t) + 5\eta^{t}[5 \mid t]\). For even \(t\) with \(5 \nmid t\) this is \(\rho^{t}c(t)\) with \(c(t) \in \{1, -2\}\); for even \(t\) with \(5 \mid t\), \(P_{T_5}(t) \in \{6, 3\}\) has modulus different from \(5\). Thus \(Z \cap 2\Z = \emptyset\).

(11) \(\rho T_5 \cup \eta\mu_3 \cup \{\pm\xi\}\): for \(5 \mid t\), \(P_{T_5}(t) = 6, 5, 3, 2\) according as \(t \equiv 0, \pm 1, \pm 2, 3 \pmod 6\), and \(\widehat{\muS}(t) = \rho^{t}P_{T_5}(t) + 3\eta^{t}[3 \mid t] + \xi^{t}(1 + (-1)^{t})\) equals, in the four cases, \(6\rho^{t} + 3\eta^{t} + 2\xi^{t}\), \(5\rho^{t}\), \(3\rho^{t} + 2\xi^{t}\), \(2\rho^{t} + 3\eta^{t}\). None vanishes: the first has modulus at least \(6 - 5 = 1\), and the others are a single term or two terms of distinct moduli. Thus \(Z \cap 5\Z = \emptyset\) and \(\lvert A\rvert \le 5\).

(12) \(\rho\mu_5 \cup \eta_1\mu_3 \cup \eta_2\mu_3\): for \(5 \mid t\), \(\widehat{\muS}(t) = 5\rho^{t} + 3[3 \mid t](\eta_1^{t} + \eta_2^{t})\), and a zero would need \(\lvert \eta_1^{t} + \eta_2^{t}\rvert = 5/3\), excluded by Lemma~\ref{lem:rationalmoduli}(a). Thus \(Z \cap 5\Z = \emptyset\).

(13) \(\rho\mu_5 \cup \{\pm\eta_1\} \cup \{\pm\eta_2\} \cup \{\pm\eta_3\}\): for \(5 \mid t\), \(\widehat{\muS}(t) = 5\rho^{t} + (1 + (-1)^{t})(\eta_1^{t} + \eta_2^{t} + \eta_3^{t})\), and a zero would need \(\lvert \eta_1^{t} + \eta_2^{t} + \eta_3^{t}\rvert = 5/2\), excluded. Thus \(Z \cap 5\Z = \emptyset\).

(14) \(\eta_1\mu_3 \cup \eta_2\mu_3 \cup \eta_3\mu_3 \cup \{\pm\xi\}\): for even \(t\) with \(3 \nmid t\), \(\widehat{\muS}(t) = 2\xi^{t}\); for \(6 \mid t\) a zero would need \(\lvert \eta_1^{t} + \eta_2^{t} + \eta_3^{t}\rvert = 2/3\), excluded. Thus \(Z \cap 2\Z = \emptyset\).

(15) \(\eta\mu_3 \cup \{\pm\xi_1\} \cup \{\pm\xi_2\} \cup \{\pm\xi_3\} \cup \{\pm\xi_4\}\): for odd \(t\) with \(3 \mid t\), \(\widehat{\muS}(t) = 3\eta^{t}\); for \(6 \mid t\) a zero would need \(\lvert \xi_1^{t} + \xi_2^{t} + \xi_3^{t} + \xi_4^{t}\rvert = 3/2\), excluded. Thus \(Z \cap 3\Z = \emptyset\).

Together with a rotation of \(\mu_{11}\), these fifteen shapes exhaust the possibilities. The three other minimal relations of weight \(11\) give (1)--(3); the parts \(9\), \(8\), \(7\) and \(6\) carry respectively two, two, two and one minimal relations, which gives (4) and (5), (6) and (7), (8) and (9), and (10) and (11), and the partitions into parts \(2\), \(3\) and \(5\) give (12) to (15). Since \(\lvert A\rvert = 11 > 7\), the roots of \(S\) form a rotation of \(\mu_{11}\).
\end{proof}

\begin{corollary}
\label{cor:cpvanishing}
Let \((A, S)\) be a spectral pair with \(A \subset \Z\), \(S \subset \R/\Z\) and \(\lvert A\rvert = \lvert S\rvert = k \le 11\), and suppose \(\sum_{s \in S}\ex{s} = 0\). Then \(S + 1/m = S\) for some \(m \ge 2\). In particular \(\CP{k}\) holds for every \(k \le 11\), and for \(k \in \{2, 3, 5, 7, 11\}\) the set \(S\) in \(\CP{k}\) is a coset of \(\tfrac1k\Z/\Z\).
\end{corollary}

\begin{proof}
This is Corollary~\ref{cor:cpten} for \(k \le 10\) and Theorem~\ref{thm:cpeleven} for \(k = 11\). A consecutive pair in \(A\) gives \(\widehat{\muS}(1) = 0\) by Definition~\ref{def:spectralpair}.
\end{proof}

The exclusions only require \(A\) to be an orthogonal set of integer frequencies of sufficient size: for \(k=11\), eight frequencies suffice, and smaller bounds suffice for smaller \(k\). In the language of Lemma~\ref{lem:cpreduction}(c), every primitive \(S' \subset \ZN{N}\) with at most eleven elements, \(\sum_{s \in S'}\zroot{N}^{s} = 0\) and a spectrum in \(\ZN{N}\) is periodic. The method stops here, and not for lack of a classification: Christie, Dykema and Klep classify the minimal relations up to weight sixteen \cite{cdk2020}. The obstruction is that the next subsection exhibits, at \(k = 12\), a shape without a symmetry that carries spectra with a consecutive pair, so no argument about the zero set of \(S\) alone can exclude it.

\subsection{Gluing}

\begin{definition}
\label{def:gluing}
Let \(m \mid N\) and let \(\pi_m : \ZN{N} \to \ZN{N/m}\) be reduction modulo \(N/m\), whose kernel is the subgroup \((N/m)\ZN{N}\) of order \(m\). We represent each \(j \in \ZN{m}\) by the unique integer \(0 \le j < m\). For \(S' \subset \ZN{N}\), its \emph{pieces} are
\[
S'_j = S' \cap (j + m\ZN{N}) = j + m\widetilde S_j, \qquad \widetilde S_j = \{ (s - j)/m : s \in S'_j \} \subset \ZN{N/m},
\]
where \((s - j)/m\) means \((\tilde s - j)/m\) modulo \(N/m\), with \(\tilde s\) any integer representative of \(s\). Since \(s \in S'_j\), the numerator is divisible by \(m\); changing \(\tilde s\) by \(N\) does not change the resulting residue. For \(\bar B \subset \ZN{N/m}\), its full inverse image is \(\pi_m^{-1}(\bar B) = B + (N/m)\ZN{N}\), with \(B\) any set of lifts of \(\bar B\).
\end{definition}

\begin{lemma}
\label{lem:gluing}
Let \(m \mid N\), let \(S' \subset \ZN{N}\) have the pieces \(\widetilde S_j \subset \ZN{N/m}\), and let \(A \subset \ZN{N}\) satisfy \(A + N/m = A\), with \(\bar B = \pi_m(A)\), so that \(A = \pi_m^{-1}(\bar B)\) and \(\lvert A\rvert = m\lvert \bar B\rvert\).
\begin{enumerate}
\item[(a)] \(A\) is a spectrum of \(S'\) in \(\ZN{N}\) iff \(\lvert \widetilde S_j\rvert = \lvert \bar B\rvert\) \(\forall j \in \ZN{m}\) and \(\bar B\) is a spectrum of \(\widetilde S_j\) in \(\ZN{N/m}\) \(\forall j\).
\item[(b)] Suppose \(N/m \ge 2\). Then \(A - A\) contains a unit of \(\ZN{N}\) iff \(\bar B - \bar B\) contains a unit of \(\ZN{N/m}\). In that case, when \(A\) is a spectrum of \(S'\), every piece satisfies \(\sum_{\tilde s \in \widetilde S_j}\zroot{N/m}^{\tilde s} = 0\).
\item[(c)] Let \(p\) be a prime divisor of \(N\). If \(p \mid N/m\), then \(S' + N/p = S'\) iff \(\widetilde S_j + (N/m)/p = \widetilde S_j\) \(\forall j\). If \(p \nmid N/m\), then \(S' + N/p = S'\) iff \(\widetilde S_{j'} = \widetilde S_j + c_j\) \(\forall j\), where \(j' \in \ZN{m}\) is the residue of \(j + N/p\) and \(c_j = (j + N/p - j')/m \in \ZN{N/m}\), using the fixed representatives of \(j\) and \(j'\); in this case \(j \mapsto j'\) permutes \(\ZN{m}\) in orbits of length \(p\), and the pieces along each orbit are translates of one another.
\end{enumerate}
\end{lemma}

\begin{proof}
(a) For \(s = j + m\tilde s \in S'_j\) and \(t \in \ZN{N}\), \(\zroot{N}^{st} = \zroot{N}^{jt}\zroot{N/m}^{\tilde s\,\pi_m(t)}\), so
\[
\widehat{\mu_{S'}}(t) = \sum_{j \in \ZN{m}} \zroot{N}^{jt}\,\widehat{\mu_{\widetilde S_j}}(\pi_m(t)).
\]
Let \(B\) be a set of lifts of \(\bar B\), so that \(A - A = (B - B) + (N/m)\ZN{N}\). For \(b, b' \in B\) and \(i \in \ZN{m}\), with \(t = b - b' + iN/m\), we have \(\zroot{N}^{jt} = \zroot{N}^{j(b - b')}\zroot{m}^{ij}\) and \(\pi_m(t) = \pi_m(b - b')\), so
\[
\widehat{\mu_{S'}}(b - b' + iN/m) = \sum_{j \in \ZN{m}} \zroot{m}^{ij}\,c_j, \qquad c_j = \zroot{N}^{j(b - b')}\,\widehat{\mu_{\widetilde S_j}}(\pi_m(b - b')).
\]
This is the Fourier transform on \(\ZN{m}\) of the vector \((c_j)\), evaluated at \(i\). It vanishes for every \(i \in \ZN{m}\) iff every \(c_j\) vanishes, and it vanishes for every \(i \ne 0\) iff \(j \mapsto c_j\) is constant. For \(b \ne b'\) in \(B\), the elements \(b - b' + iN/m\) are the nonzero elements of \(A - A\) outside \((N/m)\ZN{N}\), and \(c_j = 0\) means \(\pi_m(b - b') \in \zeroset{\widetilde S_j}\); for \(b = b'\), \(c_j = \lvert \widetilde S_j\rvert\), and the elements \(iN/m\), \(i \ne 0\), are the nonzero elements of \(A - A\) inside \((N/m)\ZN{N}\). Hence \((A - A) \setminus \{0\} \subset \zeroset{S'}\) iff all \(\lvert \widetilde S_j\rvert\) are equal and every nonzero difference of \(\bar B\) lies in every \(\zeroset{\widetilde S_j}\). Under these conditions \(\lvert S'\rvert = m\lvert \widetilde S_j\rvert\), and \(\lvert A\rvert = m\lvert \bar B\rvert\); so \(\lvert A\rvert = \lvert S'\rvert\) holds exactly when \(\lvert \widetilde S_j\rvert = \lvert \bar B\rvert\), which is the statement.

(b) Since \(A = \pi_m^{-1}(\bar B)\), we have \(A - A = \pi_m^{-1}(\bar B - \bar B)\). A unit of \(\ZN{N}\) reduces to a unit of \(\ZN{N/m}\), which gives one direction. Conversely let \(v \in \bar B - \bar B\) be a unit of \(\ZN{N/m}\) and let \(t_0 \in A - A\) be one lift of \(v\); every lift \(t_0 + iN/m\), \(i \in \ZN{m}\), lies in \(A - A\). For a prime \(p \mid N\) with \(p \mid N/m\), every lift is prime to \(p\) because \(v\) is. For a prime \(p \mid N\) with \(p \nmid N/m\), the integer \(N/m\) is a unit modulo \(p\) and \(p \mid m\), so the residues modulo \(p\) of the lifts \(t_0 + iN/m\) depend only on \(i\) modulo \(p\) and run through all of \(\ZN{p}\); choose for each such \(p\) a residue class of \(i\) modulo \(p\) in which the lift is prime to \(p\). By the Chinese remainder theorem some \(i \in \ZN{m}\) satisfies all these finitely many congruences at once, and the corresponding lift is prime to every prime divisor of \(N\), hence a unit of \(\ZN{N}\) in \(A - A\). For the last claim, if \(v \in \bar B - \bar B\) is a unit, then by (a) \(v \in \zeroset{\widetilde S_j}\), that is \(\sum_{\tilde s}\zroot{N/m}^{\tilde s v} = 0\), and applying the automorphism of \(\Q(\zroot{N/m})\) that sends \(\zroot{N/m}^{v}\) to \(\zroot{N/m}\) gives \(\sum_{\tilde s}\zroot{N/m}^{\tilde s} = 0\).

(c) If \(p \mid N/m\), then \(N/p = m \cdot (N/m)/p\), so translation by \(N/p\) preserves each residue class modulo \(m\) and acts on the piece \(\widetilde S_j\) as translation by \((N/m)/p\); the claim follows. If \(p \nmid N/m\), then \(S'_j + N/p = j' + m(\widetilde S_j + c_j)\), and \(S' + N/p = S'\) iff \(S'_j + N/p = S'_{j'}\) \(\forall j\), which is the stated condition. The residue of \(N/p\) modulo \(m\) has order \(p\) in \(\ZN{m}\): its order is \(m/\gcd(m, N/p)\), and \(\gcd(m, N/p) = m/p\) because \(p\) divides \(m\) exactly as often as it divides \(N\), while every other prime divides \(N/p\) at least as often as it divides \(m\).
\end{proof}

The lemma says that a spectral pair whose spectrum has the period \(N/m\) is the same thing as a common spectrum \(\bar B\) of the \(m\) pieces of \(S'\), each piece having \(\lvert \bar B\rvert\) elements, and, when \(N/m \ge 2\), that a consecutive pair in \(A\) gives the difference \(1\) in the common spectrum \(\bar B\). Conversely, any \(m\) spectral pairs \((\bar B, \widetilde S_j)\) in \(\ZN{N/m}\) with a common spectrum \(\bar B\) glue to the spectral pair \((\pi_m^{-1}(\bar B), \bigsqcup_j (j + m\widetilde S_j))\) in \(\ZN{N}\). By (c), the glued set \(S'\) is periodic exactly when, for some prime \(p\) dividing \(N\), either \(p \mid N/m\) and every piece has the period \((N/m)/p\), or \(p \nmid N/m\) and \(\widetilde S_{j'} = \widetilde S_j + c_j\) \(\forall j\), so that translation by \(N/p\) carries the piece of each class \(j\) onto the piece of the class \(j'\). In particular, when the pieces have no common period \((N/m)/p\) and no translation \(N/p\) with \(p \nmid N/m\) permutes them in this way, the glued set has no period. This is how the consecutive-pair statement fails.

\begin{theorem}
\label{thm:cptwelve}
Let \(N = 72\) and
\begin{align*}
\Ssev &= \{5, 6, 10, 11, 30, 34, 35, 41, 47, 54, 58, 71\} \subset \ZN{72}, \\
A' &= \{0, 1, 8, 9, 17, 28\} + \{0, 36\} \subset \ZN{72}.
\end{align*}
Then \(\Ssev\) is primitive, \(A'\) is a spectrum of \(\Ssev\) in \(\ZN{72}\), \(A'\) contains \(0\) and \(1\), and \(\Ssev\) is not periodic. Consequently \(\CP{12}\) fails: with \(S = \Ssev/72 \subset \R/\Z\) and \(A = \{0, 1, 8, 9, 17, 28, 36, 37, 44, 45, 53, 64\} \subset \Z\), the pair \((A, S)\) is a spectral pair with \(\lvert A\rvert = \lvert S\rvert = 12\), \(A\) contains two consecutive integers, and \(S\) has no translation symmetry. The set \(A\) is a union of eight runs.
\end{theorem}

\begin{proof}
We apply Lemma~\ref{lem:gluing} with \(m = 2\), so \(N/m = 36\). The pieces of \(\Ssev\) are
\begin{align*}
\widetilde S_0 &= \{3, 5, 15, 17, 27, 29\} = \{3, 5\} + 12\ZN{36}, \\
\widetilde S_1 &= \{2, 5, 17, 20, 23, 35\} = \{2, 5, 17\} + \{0, 18\},
\end{align*}
since the even elements of \(\Ssev\) are \(2 \cdot \{3, 5, 15, 17, 27, 29\}\) and the odd ones are \(1 + 2 \cdot \{2, 5, 17, 20, 23, 35\}\). Thus \(\widetilde S_0\) is a union of two cosets of the subgroup of order \(3\), with \(\widetilde S_0 + 12 = \widetilde S_0\), and \(\widetilde S_1\) is a union of three antipodal pairs, with \(\widetilde S_1 + 18 = \widetilde S_1\). Put \(\bar B = \{0, 1, 8, 9, 17, 28\} \subset \ZN{36}\), so that \(A' = \pi_2^{-1}(\bar B)\) and \(\lvert \bar B\rvert = 6 = \lvert \widetilde S_0\rvert = \lvert \widetilde S_1\rvert\).

The zero sets of the pieces. By Lemma~\ref{lem:blockpowersums}(i), for \(t \in \ZN{36}\),
\[
\widehat{\mu_{\widetilde S_0}}(t) = 3[3 \mid t]\bigl(\zroot{36}^{3t} + \zroot{36}^{5t}\bigr) = 3[3 \mid t]\,\zroot{12}^{t}\bigl(1 + \zroot{18}^{t}\bigr),
\]
which vanishes exactly when \(3 \nmid t\) or \(t \equiv 9 \pmod{18}\). Hence \(\zeroset{\widetilde S_0} = \{t \ne 0 : 3 \nmid t\} \cup \{9, 27\}\). Likewise
\[
\widehat{\mu_{\widetilde S_1}}(t) = (1 + (-1)^{t})\bigl(\zroot{36}^{2t} + \zroot{36}^{5t} + \zroot{36}^{17t}\bigr) = (1 + (-1)^{t})\,\zroot{18}^{t}\bigl(1 + \zroot{12}^{t} + \zroot{12}^{5t}\bigr).
\]
For even \(t = 2u\) the last factor is \(1 + \zroot{6}^{u} + \zroot{6}^{-u} = 1 + 2\cos(\pi u/3)\), which vanishes exactly when \(u \equiv \pm 2 \pmod 6\), that is, \(t \equiv \pm 4 \pmod{12}\). Hence \(\zeroset{\widetilde S_1} = \{t \text{ odd}\} \cup \{4, 8, 16, 20, 28, 32\}\).

The differences of \(\bar B\). The positive differences of \(0, 1, 8, 9, 17, 28\) are
\[
1, 8, 9, 17, 28, 7, 8, 16, 27, 1, 9, 20, 8, 19, 11,
\]
so the set of nonzero differences of \(\bar B\) in \(\ZN{36}\) is
\begin{align*}
(\bar B - \bar B) \setminus \{0\} &= \pm\{1, 7, 8, 9, 11, 16, 17, 19, 20, 27, 28\} \\
&= \{1, 7, 8, 9, 11, 16, 17, 19, 20, 25, 27, 28, 29, 35\}.
\end{align*}
Every element of this set is prime to \(3\) or lies in \(\{9, 27\}\), so it lies in \(\zeroset{\widetilde S_0}\); and every element is odd or lies in \(\{8, 16, 20, 28\} \subset \{4, 8, 16, 20, 28, 32\}\), so it lies in \(\zeroset{\widetilde S_1}\). Thus \(\bar B\) is a common spectrum of \(\widetilde S_0\) and \(\widetilde S_1\), and Lemma~\ref{lem:gluing}(a) shows that \(A' = \bar B \cup (\bar B + 36)\) is a spectrum of \(\Ssev\) in \(\ZN{72}\). It contains \(0\) and \(1\).

Aperiodicity and primitivity. If \(\Ssev\) were periodic, then \(36\) or \(24\) would be a period, since the only prime divisors of \(72\) are \(2\) and \(3\); but \(6 \in \Ssev\) and \(6 + 36 = 42 \notin \Ssev\), and \(5 \in \Ssev\) and \(5 + 24 = 29 \notin \Ssev\). (In the terms of Lemma~\ref{lem:gluing}(c): the piece \(\widetilde S_0\) does not have the period \(18\), and the piece \(\widetilde S_1\) does not have the period \(12\).) The set \(\Ssev\) is primitive because \(6 - 5 = 1 \in \Ssev - \Ssev\).

The pair over \(\Z\). The matrix \((\ex{as/72})_{s \in \Ssev, a \in A}\) depends only on the residues of the elements of \(A\) modulo \(72\), and these form \(A'\); so \((A, S)\) is a spectral pair by Definition~\ref{def:spectralpair}. A translation symmetry \(S + 1/m = S\) would give, by Lemma~\ref{lem:cpreduction}(a) applied with the primitive set \(\Ssev\), a period of \(\Ssev\). Finally \(A\) consists of the runs \(\{0, 1\}\), \(\{8, 9\}\), \(\{17\}\), \(\{28\}\), \(\{36, 37\}\), \(\{44, 45\}\), \(\{53\}\), \(\{64\}\).
\end{proof}

\begin{remark}
\label{rem:strongform}
The zero set of \(\Ssev\) can be described exactly: it is the set of nonzero \(t \in \ZN{72}\) with \(\gcd(t, 72) \in \{1, 4, 8, 9, 36\}\), that is, of the elements of order \(72\), \(18\), \(9\), \(8\) or \(2\). Indeed, by the identity in the proof of Lemma~\ref{lem:gluing}(a), \(\widehat{\mu_{\Ssev}}(t) = \widehat{\mu_{\widetilde S_0}}(\bar t) + \zroot{72}^{t}\,\widehat{\mu_{\widetilde S_1}}(\bar t)\) with \(\bar t = \pi_2(t) \in \ZN{36}\). Both terms vanish when \(\bar t \in \zeroset{\widetilde S_0} \cap \zeroset{\widetilde S_1}\), which is the set of \(\bar t\) with \(\gcd(\bar t, 36) \in \{1, 4, 9\}\), and the lifts of these \(\bar t\) are the \(t\) with \(\gcd(t, 72) \in \{1, 4, 8, 9\}\). For \(\bar t = 0\) and \(t = 36\) the two terms are \(6\) and \(-6\), so \(36 \in \zeroset{\Ssev}\). For the remaining \(\bar t \in \{6, 12, 18, 24, 30\}\) the two terms have the moduli \(3\lvert 1 + \zroot{3}^{\bar t/6}\rvert \in \{3, 6\}\) and \(2\lvert 1 + 2\cos(\pi\bar t/6)\rvert \in \{2, 6\}\), equal to \(6\) exactly when \(\bar t = 18\), respectively \(\bar t \in \{12, 24\}\), so they never cancel; and for the other \(\bar t\) exactly one term vanishes. The nonzero differences of \(A'\) are the lifts of the differences of \(\bar B\) together with \(36\), and their orders are \(72\) (lifts of the differences prime to \(6\)), \(8\) (lifts of \(9\) and \(27\)), \(9\) and \(18\) (lifts of \(8, 16, 20, 28\)), and \(2\): every difference has its order in \(\{2, 8, 9, 18, 72\}\), as the theorem requires.

The description shows why the method of Section~\ref{sec:consecutive} and of Theorem~\ref{thm:cpeleven} cannot see this example: \(\zeroset{\Ssev}\) meets \(m\ZN{72} \setminus \{0\}\) \(\forall 2 \le m < 12\). This subgroup contains \(36\) unless \(m = 8\), because \(\gcd(m,72)\) then divides \(36\); for \(m = 8\) it contains \(8\). \end{remark}

\subsection{The failures with a periodic spectrum}

The example of Theorem~\ref{thm:cptwelve} has a periodic spectrum and an aperiodic spectral set. We show that every failure of the consecutive-pair statement of this kind is a gluing, and that twelve is the first size at which such a gluing can be aperiodic. The following terms are used in this subsection only. A \emph{failure of size \(k\)} is a spectral pair \((A, S')\) in \(\ZN{N}\), for some \(N\ge2\), with \(\lvert A\rvert = \lvert S'\rvert = k\), such that \(A - A\) contains a unit of \(\ZN{N}\) and \(S'\) is not periodic; it is a \emph{double failure} when \(A\) is not periodic either. The pair \((A', \Ssev)\) is a failure of size \(12\) which is not a double failure.

\begin{lemma}
\label{lem:nofailure}
There is no failure of size \(k \le 11\).
\end{lemma}

\begin{proof}
Let \((A, S')\) be a spectral pair in \(\ZN{N}\) with \(\lvert A\rvert = \lvert S'\rvert = k \le 11\) and a unit \(u \in A - A\). As in the proof of Lemma~\ref{lem:cpreduction}(b), \(u^{-1}A\) is a spectrum of \(S'\) with \(1 \in u^{-1}A - u^{-1}A\), and lifting \(u^{-1}A\) to a set \(A_0 \subset \Z\) of representatives containing two consecutive integers gives a spectral pair \((A_0, S'/N)\) with \(\widehat{\mu_{S'/N}}(1) = 0\). Corollary~\ref{cor:cpvanishing} gives \(S'/N + 1/m = S'/N\) for some \(m \ge 2\), and then \(1/m \in (1/N)\Z/\Z\), \(m \mid N\) and \(S' + N/m = S'\): \(S'\) is periodic.
\end{proof}

\begin{theorem}
\label{thm:periodicfailures}
Let \((A, S')\) be a failure of size \(k\) in \(\ZN{N}\) whose spectrum \(A\) is periodic, and let \((N/m)\ZN{N}\) be the stabilizer \(\{t \in \ZN{N} : A + t = A\}\) of \(A\), so that \(2 \le m \le N/2\). Then \(\bar B = \pi_m(A) \subset \ZN{N/m}\) is not periodic, \(\bar B - \bar B\) contains a unit of \(\ZN{N/m}\), and the pieces \(\widetilde S_j\) of \(S'\), \(j \in \ZN{m}\), are \(m\) subsets of \(\ZN{N/m}\) of size \(k/m\), each with \(\sum_{\tilde s \in \widetilde S_j}\zroot{N/m}^{\tilde s} = 0\), each having \(\bar B\) as a spectrum. Each piece is periodic or is, together with \(\bar B\), a double failure of size \(k/m\). There is no prime \(p\) dividing \(N/m\) such that every piece has the period \((N/m)/p\).
\end{theorem}

\begin{proof}
The stabilizer is a subgroup of \(\ZN{N}\), hence of the form \((N/m)\ZN{N}\) with \(m \mid N\). It is not \(\{0\}\) since \(A\) is periodic, so \(m \ge 2\); and it is not \(\ZN{N}\), since then \(A = \ZN{N}\), \(k = N\), \(S' = \ZN{N}\) would be periodic. Since \(A + N/m = A\), Lemma~\ref{lem:gluing}(a) applies: every piece has \(\lvert \bar B\rvert = k/m\) elements and \(\bar B\) as a spectrum, Lemma~\ref{lem:gluing}(b) gives a unit in \(\bar B - \bar B\) and the vanishing of \(\sum_{\tilde s}\zroot{N/m}^{\tilde s}\). If \(\bar B\) were periodic, \(\bar B + (N/m)/r = \bar B\) for a prime \(r \mid N/m\), then \(A = \pi_m^{-1}(\bar B)\) would satisfy \(A + N/(mr) = A\), and \(N/(mr)\), of order \(mr > m\), would lie in the stabilizer, which has order \(m\); so \(\bar B\) is not periodic. A piece which is not periodic forms with \(\bar B\) a spectral pair with a unit difference in which both members are aperiodic, a double failure of size \(k/m\). Finally, if every piece had the period \((N/m)/p\) for one prime \(p \mid N/m\), Lemma~\ref{lem:gluing}(c) would give \(S' + N/p = S'\).
\end{proof}

\begin{proposition}
\label{prop:firstsize}
Let \((A, S')\) be a failure of size \(k \le 23\) whose spectrum \(A\) is periodic, with stabilizer of order \(m\) as in Theorem~\ref{thm:periodicfailures}, and put \(k' = k/m\). Then every piece \(\widetilde S_j\) is periodic, \(k'\) has at least two distinct prime divisors, and
\[
(m, k') \in \{(2, 6), (3, 6), (2, 10)\}, \qquad k \in \{12, 18, 20\}.
\]
In particular there is no failure of size \(k < 12\) whose spectrum is periodic, and a failure of size \(12\) with a periodic spectrum has \(m = 2\) and six-point pieces.
\end{proposition}

\begin{proof}
Since \(m \ge 2\) and \(k \le 23\), \(k' \le 11\), so by Lemma~\ref{lem:nofailure} no piece is a double failure, and Theorem~\ref{thm:periodicfailures} says that every piece is periodic. For each \(j\) let \(\Pi_j\) be the set of primes \(r\) with \(\widetilde S_j + (N/m)/r = \widetilde S_j\); it is nonempty, and each \(r \in \Pi_j\) divides \(k'\), because \(\widetilde S_j\) is then a union of cosets of the subgroup of order \(r\). By the last assertion of Theorem~\ref{thm:periodicfailures}, \(\bigcap_j \Pi_j = \emptyset\). Take \(r \in \Pi_0\); some \(j\) has \(r \notin \Pi_j\), and any \(r' \in \Pi_j\) is a second prime dividing \(k'\). Hence \(k' \ge 6\), and \(k = mk' \ge 12\). The values \(k' \le 11\) with two distinct prime divisors are \(6\) and \(10\), and \(mk' \le 23\) leaves \((m, k') \in \{(2, 6), (3, 6), (2, 10)\}\).
\end{proof}

Size \(18\) is also attained: Corollary~\ref{cor:alllargern} uses three of the six-point pieces above to construct such a failure. The size \(20\) case is not addressed here. At size \(12\) the structure is exact.

\begin{proposition}
\label{prop:twelvestructure}
Let \((A, S')\) be a failure of size \(12\) in \(\ZN{N}\) with \(A\) periodic. Then \(N = 2M\) with \(6 \mid M\), \(A = \pi_2^{-1}(\bar B)\) for a six-element set \(\bar B \subset \ZN{M}\) which is not periodic and whose difference set contains a unit, and the two pieces of \(S'\) are, in some order,
\[
C + \{0, M/2\} \quad\text{and}\quad D + (M/3)\ZN{M},
\]
where \(C \subset \ZN{M}\) has three elements pairwise distinct modulo \(M/2\) and is not a coset of \((M/6)\ZN{M}\) modulo \(M/2\), and \(D \subset \ZN{M}\) has two elements distinct modulo \(M/3\) whose difference is not congruent to \(M/6\) modulo \(M/3\); and \(\bar B\) is a common spectrum of the two pieces. Conversely, for every such \(M\), \(C\), \(D\) and \(\bar B\), the set \(S' = 2(D + (M/3)\ZN{M}) \cup (1 + 2(C + \{0, M/2\}))\), or the set with the roles of the two pieces exchanged, is not periodic, and \((\pi_2^{-1}(\bar B), S')\) is a failure of size \(12\) in \(\ZN{2M}\).
\end{proposition}

\begin{proof}
By Proposition~\ref{prop:firstsize}, \(m = 2\), \(N = 2M\), and the two pieces \(\widetilde S_0, \widetilde S_1 \subset \ZN{M}\) are six-point sets with \(\sum_{\tilde s}\zroot{M}^{\tilde s} = 0\), both having the aperiodic spectrum \(\bar B = \pi_2(A)\) with a unit difference. Lifting a spectrum with a unit difference to \(\Z\) as in the proof of Lemma~\ref{lem:nofailure}, Proposition~\ref{prop:cpsmall} applies to each piece and gives \(\widetilde S_j + M/2 = \widetilde S_j\) or \(\widetilde S_j + M/3 = \widetilde S_j\). In the first case \(M\) is even and \(\widetilde S_j\) is a union of three cosets of \(\{0, M/2\}\), that is, \(\widetilde S_j = C + \{0, M/2\}\) with \(C\) three elements pairwise distinct modulo \(M/2\); in the second case \(3 \mid M\) and \(\widetilde S_j = D + (M/3)\ZN{M}\) with \(D\) two elements distinct modulo \(M/3\). A piece of both forms is a coset of the subgroup \((M/6)\ZN{M}\) of order six, and a set of the first form is a coset of that subgroup exactly when \(C\) is a coset of \((M/6)\ZN{M}\) modulo \(M/2\), a set of the second form exactly when the difference of the two elements of \(D\) is congruent to \(M/6\) modulo \(M/3\). By the last assertion of Theorem~\ref{thm:periodicfailures}, the two pieces do not share the period \(M/2\) and do not share the period \(M/3\); so neither piece is a coset of order six, exactly one piece is of each form, and \(6 \mid M\).

Conversely, let \(M\), \(C\), \(D\), \(\bar B\) be as stated. By Lemma~\ref{lem:gluing}(a), \(\pi_2^{-1}(\bar B)\) is a spectrum of \(S'\), and by (b) its difference set contains a unit. It remains to see that \(S'\) is not periodic, by Lemma~\ref{lem:gluing}(c). For \(p = 2\) and \(p = 3\), both dividing \(M = N/m\), the piece \(D + (M/3)\ZN{M}\) does not have the period \(M/2\) (it is a union of cosets of \((M/3)\ZN{M}\), and if it were also a union of cosets of \((M/2)\ZN{M}\) it would be a union of cosets of \((M/6)\ZN{M}\), hence one coset of order six, excluded), and the piece \(C + \{0, M/2\}\) does not have the period \(M/3\) for the same reason. For a prime \(p \ge 5\) dividing \(M\), a six-point set with the period \(M/p\) would have \(p \mid 6\). No prime \(p\) divides \(N = 2M\) without dividing \(M\). So \(S'\) has no period, and \((\pi_2^{-1}(\bar B), S')\) is a failure of size \(12\).
\end{proof}

\begin{figure}[ht]
\centering
\input{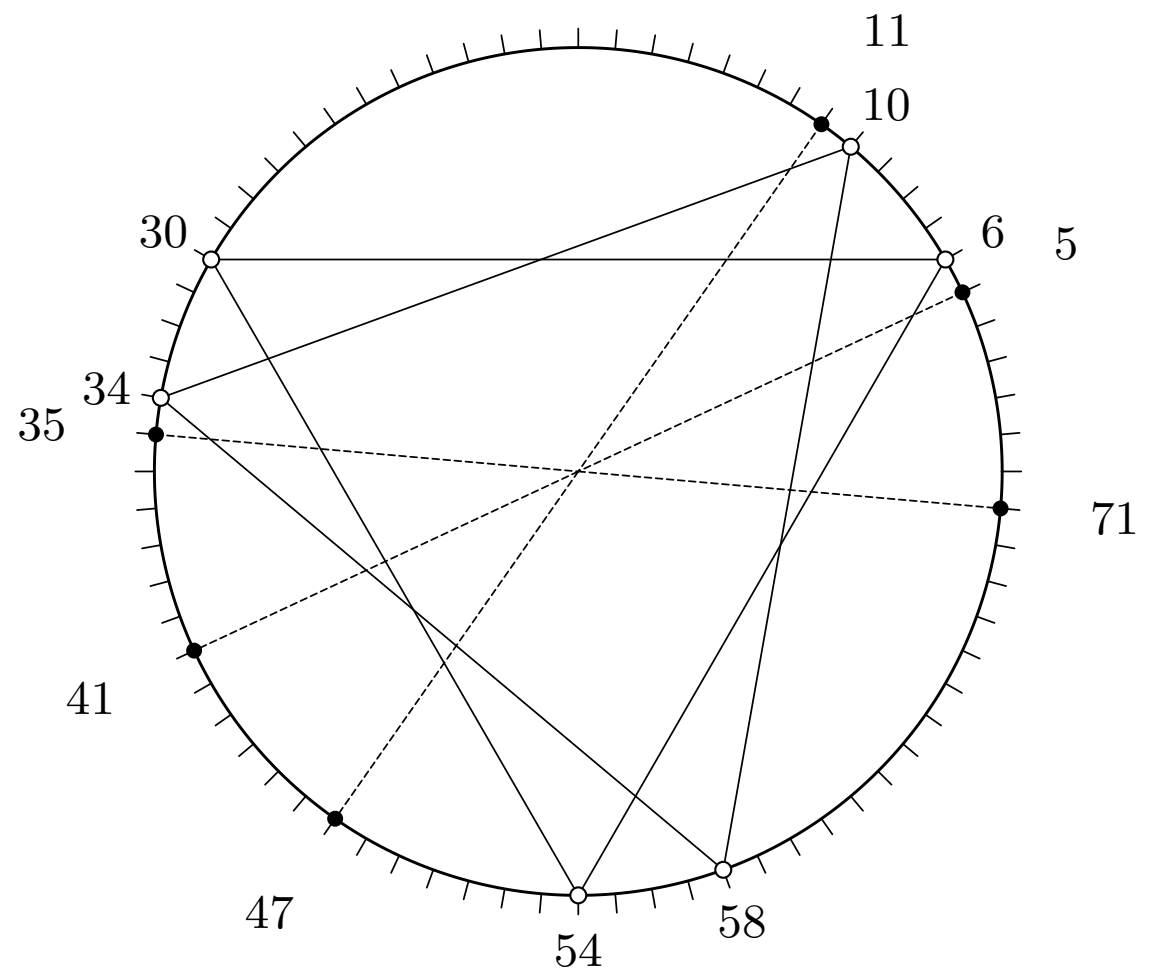}
\caption{The twelve roots of unity \(\zroot{72}^{s}\), \(s \in \Ssev\), labeled by \(s\), with the \(72\) residues marked on the circle. The six roots with \(s\) odd, filled, form three antipodal pairs, joined by dashed diameters; the six roots with \(s\) even, open, form two equilateral triangles. The whole configuration is stable neither under a half turn nor under a third of a turn.}
\label{fig:twelvepoints}
\end{figure}

The pair of Theorem~\ref{thm:cptwelve} is the instance \(M = 36\), \(C = \{2, 5, 17\}\), \(D = \{3, 5\}\), \(\bar B = \{0, 1, 8, 9, 17, 28\}\). Figure~\ref{fig:twelvepoints} shows the twelve roots \(\zroot{72}^{s}\), \(s \in \Ssev\), on the unit circle, the three antipodal pairs and the two triangles being distinguished.

\begin{remark}
\label{rem:dichotomy}
The conclusion that one member of a primitive spectral pair must be periodic is false in general. The factorization of Lagarias and Szab\'o \cite[Theorems~2.1 and~4.1]{lagarias-szabo2001} gives a double failure of size \(30\). In \(\ZN{900}\), put
\[
A = \{0,36,72,108,144\}+\{0,100,200\}+\{0,225\},
\]
and obtain \(S'\) from \(H=30\ZN{900}\) by replacing the three disjoint cosets
\[
C_1=90+180\ZN{900},\qquad C_2=120+300\ZN{900},\qquad C_3=150+450\ZN{900}
\]
by their translates by \(36\), \(100\), and \(225\), respectively. Thus \(\lvert A\rvert=\lvert S'\rvert=30\). The set \(A\) is a complete residue system modulo \(30\), so \(A\oplus H=\ZN{900}\). Moreover, \(A+C_1\), \(A+C_2\), and \(A+C_3\) are invariant under the respective shifts: the corresponding factor of \(A\), added to the subgroup in \(C_i\), fills the subgroup generated by that shift. The replacements therefore preserve the factorization, giving \(A\oplus S'=\ZN{900}\).

The character sum of \(A\) is
\[
\widehat{\mu_A}(t)=\Bigl(\sum_{i=0}^{4}\zroot{25}^{it}\Bigr)
\Bigl(\sum_{j=0}^{2}\zroot{9}^{jt}\Bigr)(1+\zroot{4}^{t});
\]
it vanishes exactly when \(p\mid t\) but \(p^2\nmid t\) for some \(p\in\{2,3,5\}\). If \(t\in A-A\), the bounds on the three digit differences show instead that \(p\mid t\) implies \(p^2\mid t\) for each such prime. Hence \(\widehat{\mu_A}(t)\ne0\) on \(A-A\), and the Fourier transform of the tiling identity gives \(\widehat{\mu_{S'}}(t)=0\) on \((A-A)\setminus\{0\}\). Thus \((A,S')\) is a spectral pair.

Both members are primitive: \(361=36+100+225\in A-A\) is a unit modulo \(900\), while \(0,126,220,375\in S'\) and \(\gcd(126,220,375)=1\). Neither member is periodic. For \(A\), the point \(0\) excludes the possible prime-order periods \(450\), \(300\), and \(180\). For \(S'\), these periods are excluded respectively by \(0+450\notin S'\), \(126+300\notin S'\), and \(220+180\notin S'\). Theorem~\ref{thm:periodicfailures} requires a periodic spectrum, so it does not cover this pair, although both members tile.
\end{remark}

\begin{proof}[Proof of Theorem~\ref{thm:D}, parts (i) to (iii)]
Part (i) is Corollary~\ref{cor:cpvanishing}: for \(k \le 11\) a spectral pair \((A, S)\) with \(A \subset \Z\), \(\lvert A\rvert = \lvert S\rvert = k\) and two consecutive integers in \(A\) has \(S + 1/m = S\) for some \(m \ge 2\), and for \(k \in \{7, 11\}\) the set \(S\) is a coset of \(\tfrac1k\Z/\Z\) by Theorems~\ref{thm:cpseven} and~\ref{thm:cpeleven}. Part (ii) is Theorem~\ref{thm:cptwelve}. Part (iii) follows from Lemma~\ref{lem:gluing}, with \(m\) the order of the stabilizer of \(A\): the \(m\) pieces have the common spectrum \(\bar B = \pi_m(A)\), and the stated periodicity criterion is precisely part (c) of that lemma. Part (iv) is proved in Section~\ref{sec:general}.
\end{proof}

\section{Fiber lemmas for six and seven points}
\label{sec:fiberlemmas}

This section proves Theorem~\ref{thm:C}. Nothing in it refers to intervals: the objects are finite spectral pairs in cyclic groups, and the letter \(n\) denotes the order of such a group. The cases \(k\le5\) follow from the theorem of Dutkay and Haussermann recalled as Theorem~\ref{thm:fiberfew}; by Lemma~\ref{lem:ucfl} it suffices to prove the common complement statements \(\UC{6}\) and \(\UC{7}\). For six points the spectra of a given set are no longer cosets of a fixed subgroup, and the mechanism is different: every spectrum of a primitive six-point set in \(\ZN{N}\) has the period \(N/2\) or the period \(N/3\), and the two shapes share a tiling complement built from cyclotomic polynomials. For seven points the rigidity is complete: every spectrum is a coset of the subgroup of order seven.

\subsection{Reduced pairs, and the kinds of a row}

Throughout, for a finite subset \(X\) of a cyclic group \(\ZN{n}\) and \(t\in\ZN{n}\), the character sum is \(\widehat{\mu_X}(t)=\sum_{x\in X}\zroot{n}^{tx}\), as in Definition~\ref{def:zeroset}, and \(tX=\{tx:x\in X\}\) is a multiset with \(\lvert X\rvert\) elements counted with multiplicity. For a prime \(p\) dividing \(n\) with \(c=v_p(n)\), the reduction \(\rho_p\colon\ZN{n}\to\ZN{p^c}\) is the natural surjection, and for \(x\in\ZN{n}\) the valuation \(v_p(x)\) is \(\min(v_p(\tilde x),c)\) for any representative \(\tilde x\in\Z\); it is well defined, and \(v_p(x)=c\) exactly when \(\rho_p(x)=0\). For \(m\mid n\), \(\Tsub{m}=(n/m)\ZN{n}\) is the subgroup of order \(m\), and ``\(X\) has the period \(n/m\)'' means \(X+n/m=X\). The multiset \(\rho_p(X)\), that is \(X\) reduced modulo \(p^c\), is called \(X\) modulo \(p^c\); its pattern is the multiset of its multiplicities.

\begin{definition}
\label{def:reducedpair}
Let \(k\ge 2\). A \emph{reduced pair} \((D,S)\) in \(\ZN{n}\) is a pair of \(k\)-point subsets of \(\ZN{n}\) with \(0\in D\cap S\), both generating \(\ZN{n}\), such that \(\widehat{\mu_D}(l)=0\) \(\forall l\in L:=(S-S)\setminus\{0\}\). The multisets \(lD\) with \(l\in L\) are the \emph{rows} of the pair, and the multisets \(\delta S\) with \(\delta\in\Delta:=(D-D)\setminus\{0\}\) are its \emph{columns}.
\end{definition}

The condition is symmetric. The \(k\times k\) matrix \(H=(\zroot{n}^{sd})_{s\in S,\,d\in D}\) has pairwise orthogonal rows, since the inner product of the rows \(s\) and \(s'\) is \(\widehat{\mu_D}(s-s')\); so \(HH^*=kI\), hence \(H^*H=kI\), and the columns are orthogonal as well: \(\widehat{\mu_S}(\delta)=0\) \(\forall \delta\in\Delta\). Thus every row and every column of a reduced pair is a vanishing sum of \(k\) roots of unity of order dividing \(n\), counted with multiplicity, and every statement proved about rows holds for columns with the roles of \(D\) and \(S\) exchanged. Two consequences of generation are used constantly.
\begin{itemize}
\item[(G)] If \(p\) is a prime dividing \(n\) and \(D\) is constant modulo \(p^j\) for some \(j\ge 1\), then, since \(0\in D\), \(D\subset p^j\ZN{n}\), which contradicts the generation of \(\ZN{n}\) by \(D\). The same holds for \(S\). In particular the multiset \(D\) modulo \(p^c\) never has the pattern \((k)\), and more generally \(D\) modulo \(p^{j}\) never has a single value when \(1\le j\le c\).
\item[(R)] If \(q\) is a prime dividing \(n\), the set \(L_{(q)}\) of elements of \(L\) prime to \(q\) generates \(\ZN{n}\). Indeed \(S\) meets at least two classes modulo \(q\), because \(0\in S\) and \(S\not\subset q\ZN{n}\). If \(s_1,s_1'\) lie in one class and \(s_2\) in another, then \(s_1-s_1'=(s_1-s_2)-(s_1'-s_2)\) is a difference of two elements of \(L_{(q)}\). So the subgroup generated by \(L_{(q)}\) contains \(S-S\), which generates \(\ZN{n}\) because \(0\in S\). Consequently, if \(t\in\ZN{n}\) satisfies \(lt=0\) \(\forall l\in L_{(q)}\), then \(t=0\); and if \(U\) is a subgroup of \(\ZN{n}\) and \(lt\in U\) \(\forall l\in L_{(q)}\), then \(t\in U\). For fixed \(t\), the set \(\{l\in\ZN{n}:lt\in U\}\) is a subgroup containing \(L_{(q)}\), hence is all of \(\ZN{n}\); taking \(l=1\) gives \(t\in U\). The first assertion is the case \(U=\{0\}\). The same holds with \(\Delta\) in place of \(L\).
\end{itemize}

Every spectrum of a primitive set reduces to a reduced pair. Recall from Definition~\ref{def:zeroset} that \(S'\subset\ZN{N}\) is primitive when \(S'-S'\) generates \(\ZN{N}\), and that a spectrum of \(S'\) is a set \(A\subset\ZN{N}\) with \(\lvert A\rvert=\lvert S'\rvert\) and \((A-A)\setminus\{0\}\subset\zeroset{S'}\).

\begin{lemma}
\label{lem:reducedpair}
Let \(S'\subset\ZN{N}\) be primitive with \(\lvert S'\rvert=k\ge 2\), let \(A\) be a spectrum of \(S'\), and fix \(a_0\in A\), \(s_0\in S'\). Let \(g\) be the greatest common divisor of \(N\) and the elements of \(A-a_0\), let \(n=N/g\), let \(D\subset\ZN{n}\) be the image of \(A-a_0\subset g\ZN{N}\) under the isomorphism \(g\ZN{N}\to\ZN{n}\), \(gy\mapsto y\), and let \(S\subset\ZN{n}\) be the reduction of \(S'-s_0\) modulo \(n\). Then \((D,S)\) is a reduced pair in \(\ZN{n}\), the reduction modulo \(n\) is injective on \(S'\), and for every \(m\mid n\), if \(D+n/m=D\) then \(A+N/m=A\).
\end{lemma}

\begin{proof}
The matrix \((\zroot{N}^{sa})_{s\in S',\,a\in A}\) has orthogonal columns by the definition of a spectrum, hence orthogonal rows: \(\sum_{a\in A}\zroot{N}^{(s-s')a}=0\) for \(s\ne s'\) in \(S'\). Writing \(a=a_0+gd\) with \(d\in D\), this sum equals \(\zroot{N}^{(s-s')a_0}\sum_{d\in D}\zroot{n}^{(s-s')d}\), because \(\zroot{N}^{g}=\zroot{n}\). The inner sum depends only on \(s-s'\) modulo \(n\); it equals \(k\ne 0\) when \(s\equiv s'\pmod n\). So distinct elements of \(S'\) stay distinct modulo \(n\), \(\lvert S\rvert=k\), and \(\widehat{\mu_D}\) vanishes on \((S-S)\setminus\{0\}\). Both sets contain \(0\). If a prime \(p\) divided \(n\) and every element of \(D\), then \(gp\) would divide \(N\) and every element of \(A-a_0\), against the choice of \(g\); so \(D\) generates \(\ZN{n}\). The image of \(S'-S'\) modulo \(n\) generates \(\ZN{n}\) and equals \(S-S\), and \(0\in S\), so \(S\) generates \(\ZN{n}\). Finally \(D+n/m=D\) gives \(A-a_0+g\cdot n/m=A-a_0\), and \(g\cdot n/m=N/m\).
\end{proof}

The kinds of a row are read from the list of minimal vanishing sums of roots of unity of Poonen and Rubinstein~\cite[Table~1 and Theorem~3]{poonen-rubinstein1998}, which this section uses up to weight seven. In the form needed: every vanishing sum of roots of unity with positive integer coefficients is a sum of minimal relations, and the minimal relations of weight at most seven are, up to rotation, \(\Rp{2}\), \(\Rp{3}\), \(\Rp{5}\), \(\Rp{7}\), the relation \((\Rp{5}:\Rp{3})\) of weight six, and the two relations \((\Rp{5}:2\Rp{3})\) of weight seven, where \(\Rp{p}\) is the sum of the \(p\)th roots of unity and \((\Rp{5}:j\Rp{3})\) is obtained from \(\Rp{5}\) by subtracting \(j\) rotated copies of \(\Rp{3}\) attached to \(j\) distinct fifth roots of unity and absorbing the signs into the roots. There is no minimal relation of weight one or four. The case of weight six is \cite[Proposition~6.2]{bell2026four}, used already in Section~\ref{sec:consecutive}, where the minimal relations of a decomposition were called blocks and the list of their weights a shape; for the six- and seven-element multisets of this section the decomposition itself is called the kind, and the kinds are few enough to be named. In exponent form inside \(\ZN{n}\): a rotated \(\Rp{p}\) is a coset \(y+\Tsub{p}\), which requires \(p\mid n\); since \(-\zroot{3}=\zroot{6}^{5}\) and \(-\zroot{3}^2=\zroot{6}\), the relation \((\Rp{5}:\Rp{3})\) is \(x+\{n/6,\,5n/6,\,n/5,\,2n/5,\,3n/5,\,4n/5\}\), which requires \(30\mid n\), and \((\Rp{5}:2\Rp{3})\) is
\[
x+\bigl(\{jn/5:j\in\ZN{5}\setminus X\}\cup\{x_in/5+n/6,\ x_in/5-n/6:x_i\in X\}\bigr),
\]
with \(X=\{x_1,x_2\}\) a two-element subset of \(\ZN{5}\), seven distinct elements, again with \(30\mid n\).

\begin{lemma}
\label{lem:rowkinds}
Let \(M\) be a multiset of \(k\) elements of \(\ZN{n}\) whose roots of unity \(\zroot{n}^{x}\), \(x\in M\), sum to zero.

{\rm (a)} If \(k=6\), then \(M\) is of exactly one of the following kinds. Kind {\rm (I)}, three antipodal pairs that are not of kind {\rm (6)} below: \(M=\{x_1,x_1+n/2,x_2,x_2+n/2,x_3,x_3+n/2\}\), with \(2\mid n\); the \(x_i\) may coincide, so the pattern of \(M\) is \((1,1,1,1,1,1)\), \((2,2,1,1)\) or \((3,3)\). Kind {\rm (II)}, two triangles that are not of kind {\rm (6)} below: \(M=\{y_1,y_2\}+\Tsub{3}\), with \(3\mid n\) and pattern \((1,1,1,1,1,1)\) or \((2,2,2)\). Kind {\rm (6)}: \(M=x+\Tsub{6}\), with \(6\mid n\); such an \(M\) is both three antipodal pairs and two triangles, and it is the only multiset with both properties. Kind {\rm (III)}: \(M=x+\{n/6,5n/6,n/5,2n/5,3n/5,4n/5\}\), six distinct elements, with \(30\mid n\); such an \(M\) is neither three antipodal pairs nor two triangles. A multiset of kind {\rm (I)} or {\rm (6)} is said to be three antipodal pairs, one of kind {\rm (II)} or {\rm (6)} two triangles, and one of kind {\rm (I)} or {\rm (II)} is called pure.

{\rm (b)} If \(k=7\), then \(M\) is of exactly one of the following kinds: {\rm (a)} a coset \(x+\Tsub{7}\), with \(7\mid n\); {\rm (b)} a rotation of \((\Rp{5}:2\Rp{3})\) as displayed above, with \(30\mid n\); {\rm (c)} a coset \(y+\Tsub{5}\) together with an antipodal pair \(\{z,z+n/2\}\), with \(10\mid n\); {\rm (d)} a coset \(y+\Tsub{3}\) together with two antipodal pairs \(\{z_1,z_1+n/2\}\) and \(\{z_2,z_2+n/2\}\), with \(6\mid n\). Kinds {\rm (b)}, {\rm (c)}, {\rm (d)} are called bad; a bad multiset requires \(2\mid n\).
\end{lemma}

\begin{proof}
(a) The relation has weight six, so it is a sum of minimal relations with weights in \(\{2,3,5,6\}\) adding to six, and \(5+1\) is excluded; the partitions are \(2+2+2\), \(3+3\) and \(6\), which give three antipodal pairs, two triangles, or kind (III). Two antipodal pairs, and two cosets of \(\Tsub{3}\), are equal or disjoint, which gives the patterns. Suppose \(M\) is both three antipodal pairs and two triangles, \(M=\{y_1,y_2\}+\Tsub{3}\) with \(M+n/2=M\). If \(y_1+\Tsub{3}\ne y_2+\Tsub{3}\) and \(y_2\notin y_1+n/2+\Tsub{3}\), then \(y_1+n/2\in y_1+\Tsub{3}\), so \(n/2\in\Tsub{3}\), which is impossible since \(n/2\) has order two; if \(y_1+\Tsub{3}=y_2+\Tsub{3}\), the doubled coset is stable under \(+n/2\) for the same reason only if \(n/2\in\Tsub{3}\). Hence \(y_2\in y_1+n/2+\Tsub{3}\) and \(M=y_1+\Tsub{6}\). A multiset of kind (III) has six distinct elements; in units of \(n/30\) they are \(5,25,6,12,18,24\), no two differ by \(15\), and no element \(u\) has both \(u+10\) and \(u+20\) in the set, so it is neither three antipodal pairs nor two triangles.

(b) The partitions of seven into minimal weights from \(\{2,3,5,6,7\}\) are \(7\), \(5+2\) and \(3+2+2\), and the minimal relations of weight seven are \(\Rp{7}\) and the two \((\Rp{5}:2\Rp{3})\). This gives the four kinds. They exclude one another. A coset of \(\Tsub{7}\) contains no two elements differing by \(n/2\), \(n/3\) or \(n/5\), since all its differences have order seven; so (a) excludes the rest. A multiset of kind (b) contains no antipodal pair: in units of \(n/30\) its elements are \(6j\) for the three \(j\notin X\) and \(6x_i\pm 5\), and every difference is congruent to \(0\), \(\pm 1\), \(\pm 2\) or \(\pm 4\) modulo \(6\), never to \(3\); so (b) excludes (c) and (d). Finally the five elements of a coset of \(\Tsub{5}\) pairwise differ by elements of order five, so at most one of them lies in a coset of \(\Tsub{3}\) and at most one in each antipodal pair; three blocks of kind (d) hold at most three of them, and (c) excludes (d).
\end{proof}

For a multiset \(M\) of kind (I), (II) or (6) in \(\ZN{n}\), the reductions are read off at once, and they are used throughout. If \(M\) is three antipodal pairs and \(p\ne 2\), then \(\rho_p(n/2)=0\) and \(\rho_p(M)\) has all multiplicities even; if \(M\) is two triangles and \(p\ne 3\), then \(\rho_p(n/3)=0\) and \(\rho_p(M)\) has all multiplicities divisible by three; if \(M\) is three antipodal pairs, \(\rho_2(M)\) is stable under translation by \(\rho_2(n/2)=2^{a-1}\), where \(a=v_2(n)\); if \(M\) is two triangles, \(\rho_3(M)\) is stable under translation by \(\rho_3(n/3)\), an element of order three. The following principle converts these facts into statements about \(D\).

\begin{lemma}[collapse principle]
\label{lem:collapseprinciple}
Let \(D\subset\ZN{n}\), let \(p\) be a prime with \(c=v_p(n)\ge 1\), and let \(l\in\ZN{n}\) with \(i=v_p(l)<c\). Then the multiset \(\rho_p(lD)\) is the image of the multiset \(D\) modulo \(p^{c-i}\) under an injective homomorphism \(\ZN{p^{c-i}}\to\ZN{p^c}\) which carries the subgroup of order \(p\) onto the subgroup of order \(p\). In particular \(\rho_p(lD)\) and \(D\) modulo \(p^{c-i}\) have the same pattern, the one is stable under a translation of order \(p\) exactly when the other is, and if \(\rho_p(lD)\) lies in one coset of the subgroup of order \(p\) then \(D\) is constant modulo \(p^{c-i-1}\). If \(i\le i'\le c\), then \(D\) modulo \(p^{c-i'}\) is the image of \(D\) modulo \(p^{c-i}\) under reduction; the number of distinct values does not increase and multiplicities add.
\end{lemma}

\begin{proof}
Write \(\rho_p(l)=p^iw\) with \(w\) a unit of \(\ZN{p^c}\). The map \(\ZN{p^{c-i}}\to\ZN{p^c}\), \(y\mapsto p^iwy\), is a well-defined injective homomorphism with image \(p^i\ZN{p^c}\), it sends the subgroup \(p^{c-i-1}\ZN{p^{c-i}}\) of order \(p\) onto \(p^{c-1}\ZN{p^c}\), and \(\rho_p(ld)=p^iw\,\rho_p(d)\) depends on \(d\) only through \(d\) modulo \(p^{c-i}\). Everything follows, the assertion about constancy because a coset of \(p^{c-i-1}\ZN{p^{c-i}}\) is a single class modulo \(p^{c-i-1}\).
\end{proof}

The common complement is organized by the cyclotomic data of a finite set, which Coven and Meyerowitz introduced~\cite{coven-meyerowitz1999}.

\begin{definition}
\label{def:cyclotomicdata}
For a finite set \(A\subset\ZN{N}\) let \(\mask{A}=\sum_{a\in A}X^{\tilde a}\) with representatives \(\tilde a\in[0,N)\), and let \(\cycset{A}\) be the set of divisors \(e>1\) of \(N\) such that the cyclotomic polynomial \(\cyc{e}\) divides \(\mask{A}\). Since \(\zroot{N}^{t}\) is a primitive root of unity of order \(\operatorname{ord}(t)\) and \(\cyc{e}\) is irreducible over \(\Q\), one has \(\widehat{\mu_A}(t)=A(\zroot{N}^t)=0\) exactly when \(\cyc{\operatorname{ord}(t)}\mid\mask{A}\); so \(\cycset{A}=\{\operatorname{ord}(t):t\in\zeroset{A}\}\), and \(\cycset{A}\) is unchanged by translation of \(A\). The set \(A\) satisfies the condition {\rm (T1)} when \(\lvert A\rvert=\prod p\), the product over the prime powers \(p^c\in\cycset{A}\), and {\rm (T2)} when \(p_1^{c_1}\cdots p_t^{c_t}\in\cycset{A}\) whenever \(p_1^{c_1},\dots,p_t^{c_t}\in\cycset{A}\) are powers of distinct primes. Coven and Meyerowitz prove that a finite set of integers satisfying {\rm (T1)} and {\rm (T2)} tiles \(\Z\), with an explicit complement built from the cyclotomic polynomials of the prime powers outside \(\cycset{A}\); Theorem~\ref{thm:sixcomplement} below carries out their construction in \(\ZN{N}\), with the proof, for the spectra of a six-point set.
\end{definition}

\begin{lemma}
\label{lem:primepowers}
Let \(D\subset\ZN{n}\) with \(\lvert D\rvert=6\). {\rm (a)} The prime powers in \(\cycset{D}\) consist of at most one power of \(2\), at most one power of \(3\), and nothing else. {\rm (b)} Let \(l\in\zeroset{D}\). If \(lD\) is three antipodal pairs, then \(v_2(l)<a:=v_2(n)\) and \(2^{a-v_2(l)}\in\cycset{D}\). If \(lD\) is two triangles, then \(v_3(l)<b:=v_3(n)\) and \(3^{b-v_3(l)}\in\cycset{D}\).
\end{lemma}

\begin{proof}
(a) The polynomials \(\cyc{e}\), \(e\in\cycset{D}\), are distinct monic irreducible divisors of \(\mask{D}\) in \(\Z[X]\), so their product divides \(\mask{D}\) in \(\Z[X]\) by Gauss's lemma, and \(\prod_{e\in\cycset{D}}\cyc{e}(1)\) divides \(D(1)=6\). Since \(\cyc{p^j}(1)=p\) for a prime power and \(\cyc{e}(1)=1\) when \(e>1\) is not a prime power, the prime powers in \(\cycset{D}\) contribute a product of primes dividing \(6\) with the stated multiplicities.

(b) If \(lD\) is three antipodal pairs then \(n/2\in lD-lD\), which lies in the subgroup generated by \(l\); so \(2\) divides \(e=\operatorname{ord}(l)=n/\gcd(n,\tilde l)\), that is \(v_2(l)<a\) and \(v_2(e)=a-v_2(l)\). Put \(m=e/2^{v_2(e)}\), an odd integer. Then \(ml\) has order \(2^{v_2(e)}\), and \(mlD=\{mx_i,\,mx_i+mn/2\}_{i\le 3}=\{mx_i,\,mx_i+n/2\}_{i\le 3}\) is again three antipodal pairs, so \(ml\in\zeroset{D}\) and \(2^{a-v_2(l)}\in\cycset{D}\). If \(lD=\{y_1,y_2\}+\Tsub{3}\), then \(n/3\in lD-lD\), so \(3\mid e\), \(v_3(l)<b\), \(v_3(e)=b-v_3(l)\); with \(m=e/3^{v_3(e)}\), prime to \(3\), \(m\Tsub{3}=\Tsub{3}\) and \(mlD=\{my_1,my_2\}+\Tsub{3}\), so \(ml\in\zeroset{D}\) has order \(3^{b-v_3(l)}\).
\end{proof}

\subsection{No row of a six-point pair is the exotic relation}

\begin{theorem}
\label{thm:noexotic}
Let \((D,S)\) be a reduced six-point pair in \(\ZN{n}\). No row and no column is of kind {\rm (III)}.
\end{theorem}

\begin{proof}
By the symmetry of Definition~\ref{def:reducedpair} it suffices to exclude a row \(l_0D\), \(l_0\in L\), of kind (III). Suppose there is one; then \(30\mid n\). Put \(c=v_5(n)\) and \(D_5=\rho_5(D)\), a multiset in \(\ZN{5^c}\) containing \(0\). The possible patterns of \(\rho_5(M)\) for a vanishing six-multiset \(M\) are read from Lemma~\ref{lem:rowkinds}: all multiplicities even for kind (I), all divisible by three for kind (II), the pattern \((6)\) for kind (6), and for kind (III), since \(\rho_5(n/6)=\rho_5(5n/6)=0\) and \(u=\rho_5(n/5)\) has order five, \(\rho_5(M)=\rho_5(x)+\{0,0,u,2u,3u,4u\}\), the pattern \((2,1,1,1,1)\) with its five values in one coset of the subgroup of order five. Only kind (III) has five distinct values.

(a) The multiset \(\rho_5(l_0D)=\rho_5(l_0)D_5\) has five distinct values, so \(D_5\) has at least five. Since \(S\) generates \(\ZN{n}\), some \(l_5\in S\) is prime to \(5\); it is nonzero, so \(l_5\in L\), and \(\rho_5(l_5D)=\rho_5(l_5)D_5\) has the pattern of \(D_5\), with at least five distinct values. Hence \(l_5D\) is of kind (III), \(D_5\) has the pattern \((2,1,1,1,1)\), and \(\rho_5(l_5D)\) lies in one coset of the subgroup of order five, so \(D_5\) does as well, and since \(0\in D_5\) the coset is the subgroup \(5^{c-1}\ZN{5^c}\) itself. If \(c\ge 2\), every element of \(D\) is divisible by \(5\), against (G). So \(c=1\): \(D\) modulo \(5\) has the pattern \((2,1,1,1,1)\). Write \(D=\{d_1,d_2,d_3,d_4,e,e'\}\) with \(e\equiv e'\pmod 5\) and \(d_1,\dots,d_4,e\) in distinct classes modulo \(5\).

(b) A row \(lD\) is of kind (III) iff \(5\nmid l\): for \(5\nmid l\) the multiset \(lD\) modulo \(5\) has the pattern \((2,1,1,1,1)\), which only kind (III) admits, and for \(5\mid l\) it is constant. For a row \(lD=x+\{n/6,5n/6,n/5,\dots,4n/5\}\) of kind (III), the two elements of \(lD\) in the doubled class modulo \(5\) are \(x+n/6\) and \(x+5n/6\), and they are \(le\) and \(le'\) in some order, since \(le\equiv le'\) and the \(ld_i\), \(le\) are pairwise incongruent modulo \(5\). Hence \(l\delta_0=\pm(n/6-5n/6)=\mp 2n/3\in\Tsub{3}\setminus\{0\}\) for \(\delta_0:=e-e'\) and every \(l\in L\) prime to \(5\).

(c) Partition \(S\) into its classes modulo \(5\); there are at least two, since \(0\in S\) and some element of \(S\) is prime to \(5\). If \(s,s'\in S\) lie in distinct classes, then \(l=s-s'\in L\) is prime to \(5\) and \(\delta_0s-\delta_0s'=l\delta_0\in\Tsub{3}\setminus\{0\}\). Four elements of \(S\) in four distinct classes would therefore have \(\delta_0\)-multiples pairwise differing by nonzero elements of \(\Tsub{3}\), that is four distinct elements in one coset of \(\Tsub{3}\), which has three elements. So \(S\) meets two or three classes modulo \(5\).

(d) For \(i\ne j\) in \(\{1,2,3,4\}\) the column \(\delta=d_i-d_j\in\Delta\) is prime to \(5\), so \(\delta S\) modulo \(5\) has the pattern of \(S\) modulo \(5\). The column is a vanishing six-multiset, so this pattern is one of the admissible ones listed at the start: all even, all divisible by three, or \((2,1,1,1,1)\). With two or three classes the pattern of \(S\) modulo \(5\) is \((4,2)\), \((2,2,2)\) or \((3,3)\); the patterns \((5,1)\), \((4,1,1)\) and \((3,2,1)\) are excluded.

(e) Suppose the pattern is \((4,2)\) or \((2,2,2)\). Then every column \(\delta\) prime to \(5\) is of kind (I), three antipodal pairs, and since \(\rho_5(n/2)=0\) the two members of each antipodal pair of \(\delta S\) are congruent modulo \(5\): every class of \(\delta S\) modulo \(5\) is a union of antipodal pairs. A class \(\{s,s'\}\) of \(S\) of size two gives the class \(\{\delta s,\delta s'\}\) of \(\delta S\), which is therefore one antipodal pair: \(\delta(s-s')=n/2\). Put \(l=s-s'\in L\) and apply this to the three columns \(d_1-d_2\), \(d_2-d_3\), \(d_1-d_3\): the first two give \(l(d_1-d_3)=n/2+n/2=0\), and the third gives \(l(d_1-d_3)=n/2\). Contradiction.

(f) Suppose the pattern is \((3,3)\). Then every column \(\delta\) prime to \(5\) is of kind (II), two triangles, and since \(\rho_5(\Tsub{3})=0\) each triangle of \(\delta S\) lies in one class modulo \(5\); each class of \(\delta S\), of size three, is one triangle. For a class \(\{s,s',s''\}\) of \(S\) and \(\delta=d_i-d_j\), the elements \(\delta s,\delta s',\delta s''\) form a coset of \(\Tsub{3}\), so \(l(d_i-d_j)\in\Tsub{3}\setminus\{0\}\) for \(l=s-s'\) and all \(1\le i<j\le 4\). Then \(ld_1,\dots,ld_4\) are four distinct elements of one coset of \(\Tsub{3}\). Contradiction.

All cases are contradictory, so no row is of kind (III).
\end{proof}

From now on every row and every column of a reduced six-point pair is three antipodal pairs or two triangles, and a row that is not three antipodal pairs is pure (II), a row that is not two triangles is pure (I).

\subsection{The six-point dichotomy}

In this subsection \((D,S)\) is a reduced six-point pair in \(\ZN{n}\), and \(n=2^a3^bm\) with \(\gcd(m,6)=1\). Lemma~\ref{lem:collapseprinciple} and the reductions of the kinds give the following facts, for a row \(lD\) with \(l\in L\). If \(lD\) is three antipodal pairs, then \(i=v_2(l)<a\) and \(D\) modulo \(2^{a-i}\) is stable under translation by \(2^{a-i-1}\); indeed \(\rho_2(lD)\) is stable under translation by \(2^{a-1}\), which excludes \(i=a\), where \(\rho_2(lD)\) would be constant. If \(lD\) is two triangles and \(i=v_2(l)<a\), then \(D\) modulo \(2^{a-i}\) has the pattern \((3,3)\), since \(\rho_2(lD)=\{\rho_2(y_1)^3,\rho_2(y_2)^3\}\) and the pattern \((6)\) is excluded by (G). Symmetrically, if \(lD\) is two triangles then \(j=v_3(l)<b\) and \(D\) modulo \(3^{b-j}\) is stable under a translation of order three; and if \(lD\) is three antipodal pairs and \(j=v_3(l)<b\), then \(D\) modulo \(3^{b-j}\) has all multiplicities even. For a prime \(q\mid m\) and \(i=v_q(l)<v_q(n)\), the multiset \(D\) modulo \(q^{v_q(n)-i}\) has all multiplicities even if \(lD\) is three antipodal pairs and all divisible by three if \(lD\) is two triangles.

\begin{lemma}
\label{lem:bothkinds}
Some row is pure {\rm (I)} and some row is pure {\rm (II)}. Consequently \(6\mid n\), the set \(\cycset{D}\) contains exactly one power of \(2\), written \(2^\alpha\), and exactly one power of \(3\), written \(3^\beta\), with \(\alpha,\beta\ge 1\); every row that is three antipodal pairs has \(v_2(l)=a-\alpha\), and every row that is two triangles has \(v_3(l)=b-\beta\).
\end{lemma}

\begin{proof}
Suppose every row is three antipodal pairs. Then \(2\mid n\), and since \(S\) generates \(\ZN{n}\), some \(l_1\in L\) is odd. Lemma~\ref{lem:primepowers}(b) gives \(2^a\in\cycset{D}\). The same lemma puts \(2^{a-v_2(l)}\) in \(\cycset{D}\) \(\forall l\in L\), and part (a) allows only one power of \(2\); hence every \(l\in L\) is odd. This is impossible, since two of the six elements of \(S\) have the same parity.

Similarly, if every row is two triangles, then \(3\mid n\), and generation gives \(l_1\in L\) prime to \(3\), so \(3^b\in\cycset{D}\). Uniqueness of the power of \(3\) forces every \(l\in L\) to be prime to \(3\), whereas two elements of \(S\) are congruent modulo \(3\).

So a row \(l\) exists that is not three antipodal pairs, hence pure (II) by Theorem~\ref{thm:noexotic}, and a row \(l'\) exists that is pure (I). Then \(2\mid n\) and \(3\mid n\). Lemma~\ref{lem:primepowers}(b) puts \(2^{a-v_2(l')}\) and \(3^{b-v_3(l)}\) in \(\cycset{D}\), with positive exponents, and Lemma~\ref{lem:primepowers}(a) says that these are the only prime powers of \(2\) and \(3\) there; this defines \(\alpha\) and \(\beta\), and applied to any row that is three antipodal pairs, respectively two triangles, it gives \(a-v_2(l)=\alpha\), respectively \(b-v_3(l)=\beta\).
\end{proof}

\begin{lemma}
\label{lem:unitrow}
If some \(l\in L\) is a unit of \(\ZN{n}\), then \(D+n/2=D\) when \(lD\) is three antipodal pairs and \(D+n/3=D\) when \(lD\) is two triangles.
\end{lemma}

\begin{proof}
\(D=l^{-1}(lD)\), and multiplication by the unit \(l^{-1}\) carries \(\{x,x+n/2\}\) to \(\{l^{-1}x,l^{-1}x+n/2\}\), because \(l^{-1}\) is odd when \(2\mid n\), and carries cosets of \(\Tsub{3}\) to cosets of \(\Tsub{3}\).
\end{proof}

\begin{lemma}[the collapsing prime]
\label{lem:collapse}
Let \(q\) be a prime dividing \(n\), and let \(L_{(q)}\) be the set of elements of \(L\) prime to \(q\).

{\rm (a)} If \(q\ne 3\) and every row \(lD\), \(l\in L_{(q)}\), is two triangles, then \(D+n/3=D\).

{\rm (b)} If \(q=3\) and every row \(lD\), \(l\in L_{(3)}\), is three antipodal pairs, then \(D+n/2=D\).
\end{lemma}

\begin{proof}
Put \(c=v_q(n)\). The set \(L_{(q)}\) is nonempty by (R).

(a) Fix \(l\in L_{(q)}\). Since \(q\ne 3\), \(\rho_q(n/3)=0\) and \(\rho_q(lD)=\{\rho_q(y_1)^3,\rho_q(y_2)^3\}\); as \(\rho_q(l)\) is a unit, \(D\) modulo \(q^c\) has the pattern \((3,3)\) or \((6)\), and \((6)\) is excluded by (G). So \(D=D_0\sqcup D_1\), where \(D_0\) is the set of elements of \(D\) divisible by \(q^c\), which contains \(0\), and \(D_1\) the set of those with the other value \(W\ne 0\) of \(\rho_q\); both have three elements. For every \(l\in L_{(q)}\), each coset of \(\Tsub{3}\) in \(lD\) has constant \(\rho_q\), and the two values \(0\) and \(\rho_q(l)W\) are distinct; so the two cosets are \(lD_0\) and \(lD_1\), and \(lD_0\) is the coset through \(l\cdot 0=0\), that is \(lD_0=\Tsub{3}\). Writing \(D_0=\{0,d,d'\}\), this says \(\{ld,ld'\}=\{n/3,2n/3\}\), so \(3ld=0\) \(\forall l\in L_{(q)}\), and (R) gives \(3d=0\); likewise \(3d'=0\). So \(d,d'\) are the two nonzero elements of \(\Tsub{3}\) and \(D_0=\Tsub{3}\). Writing \(D_1=\{d_1,d_1',d_1''\}\), the multiset \(lD_1\) is a coset of \(\Tsub{3}\) \(\forall l\in L_{(q)}\), so \(l(d_1-d_1')\in\Tsub{3}\) for all such \(l\), and (R) gives \(d_1-d_1'\in\Tsub{3}\); likewise \(d_1-d_1''\in\Tsub{3}\). Hence \(D_1=d_1+\Tsub{3}\), \(D=\Tsub{3}\sqcup(d_1+\Tsub{3})\), and \(D+n/3=D\).

(b) Fix \(l\in L_{(3)}\). Since \(\rho_3(n/2)=0\), \(\rho_3(lD)\) has all multiplicities even, and so has \(D\) modulo \(3^b\); by (G) its pattern is \((2,2,2)\) or \((4,2)\). The pattern \((4,2)\) is excluded as follows. By Lemma~\ref{lem:bothkinds} some row \(l'D\) is two triangles; then \(j=v_3(l')<b\) and \(D\) modulo \(3^{b-j}\), which is the image of \(D\) modulo \(3^b\) under reduction, is stable under a translation of order three. If \(D\) modulo \(3^b\) had the pattern \((4,2)\), its image modulo \(3^{b-j}\) would have a value of multiplicity at least four, and so would the two other points of its orbit under the translation, giving total multiplicity at least twelve. So \(D\) modulo \(3^b\) has the pattern \((2,2,2)\): \(D=D_{X_1}\sqcup D_{X_2}\sqcup D_{X_3}\) according to the three values \(X_1,X_2,X_3\) of \(\rho_3\), two elements in each part. For \(l\in L_{(3)}\), the multiset \(\rho_3(lD)=\rho_3(l)\rho_3(D)\) has three distinct values of multiplicity two, each antipodal pair of \(lD\) has constant \(\rho_3\), and so each of the three classes \(lD_{X}\) is one antipodal pair: for \(D_X=\{e,e'\}\), \(l(e-e')=n/2\), hence \(2l(e-e')=0\), \(\forall l\in L_{(3)}\). By (R), \(2(e-e')=0\), and \(e\ne e'\) gives \(e'=e+n/2\). So \(D+n/2=D\).
\end{proof}

\begin{theorem}[the six-point dichotomy]
\label{thm:sixdichotomy}
Let \((D,S)\) be a reduced six-point pair in \(\ZN{n}\). Then \(D+n/2=D\) or \(D+n/3=D\), and likewise \(S+n/2=S\) or \(S+n/3=S\). Consequently, if \(S'\subset\ZN{N}\) is a primitive six-point set, every spectrum \(A\) of \(S'\) satisfies \(A+N/2=A\) or \(A+N/3=A\).
\end{theorem}

\begin{proof}
If every odd \(l\in L\) is two triangles, Lemma~\ref{lem:collapse}(a) with \(q=2\) gives \(D+n/3=D\). If every \(l\in L\) prime to \(3\) is three antipodal pairs, Lemma~\ref{lem:collapse}(b) gives \(D+n/2=D\). Otherwise some odd \(l_1\in L\) is not two triangles, hence pure (I), and some \(l_2\in L\) prime to \(3\) is pure (II). Lemma~\ref{lem:bothkinds} gives \(a-\alpha=v_2(l_1)=0\) and \(b-\beta=v_3(l_2)=0\): every row that is three antipodal pairs is odd, and every row that is two triangles is prime to \(3\). A row divisible by \(6\) would have no kind, which Theorem~\ref{thm:noexotic} forbids; so no two elements of \(S\) are congruent modulo \(6\), and \(S\) is a complete residue system modulo \(6\), with three even and three odd elements.

Let \(q\) be a prime dividing \(m\), and \(c=v_q(n)\). For \(l\in L_{(q)}\) the pattern of \(\rho_q(lD)\) is that of \(D\) modulo \(q^c\), and it is all even when \(lD\) is three antipodal pairs and all divisible by three when \(lD\) is two triangles, while \((6)\) is excluded by (G). As no pattern of six elements other than \((6)\) is both all even and all divisible by three, all rows in \(L_{(q)}\) have the same kind. If for some \(q\mid m\) they are two triangles, Lemma~\ref{lem:collapse}(a) gives \(D+n/3=D\). Suppose that for every prime \(q\mid m\) the rows in \(L_{(q)}\) are three antipodal pairs; this is vacuous when \(m=1\). Then \(D\) modulo \(q^c\) has the pattern \((2,2,2)\) or \((4,2)\). Let \(l\in L\) be even; it is two triangles. If \(i=v_q(l)<c\), the multiset \(D\) modulo \(q^{c-i}\) would have all multiplicities divisible by three while being an image of \(D\) modulo \(q^c\) under reduction; the images of the patterns \((2,2,2)\) and \((4,2)\) are \((2,2,2)\), \((4,2)\) and \((6)\), of which only \((6)\) is divisible by three, and \((6)\) is excluded by (G). Hence every even element of \(L\) is divisible by \(q^c\) \(\forall q\mid m\), that is by \(m\). The even elements of \(L\) are the differences inside the two parity classes of \(S\), so each parity class lies in a single class modulo \(m\). The two classes modulo \(m\) are distinct modulo every prime \(q\mid m\): otherwise \(S\) would lie in one class modulo \(q\), and \(0\in S\) would put \(S\) in \(q\ZN{n}\). Hence every odd element of \(L\), being a difference across the parity classes, is prime to \(m\). Now take \(s\in S\) with \(s\equiv 1\pmod 6\). Then \(l=s\in L\) is odd, prime to \(3\), and prime to \(m\): a unit of \(\ZN{n}\). Lemma~\ref{lem:unitrow} concludes.

The statement for \(S\) follows by the symmetry of Definition~\ref{def:reducedpair}. For the last statement, Lemma~\ref{lem:reducedpair} reduces \((A,S')\) to a reduced six-point pair in \(\ZN{n}\), and a period \(n/2\) or \(n/3\) of \(D\) lifts to the period \(N/2\) or \(N/3\) of \(A\).
\end{proof}

\begin{figure}[ht]
\centering
\begin{overpic}[width=0.8\textwidth]{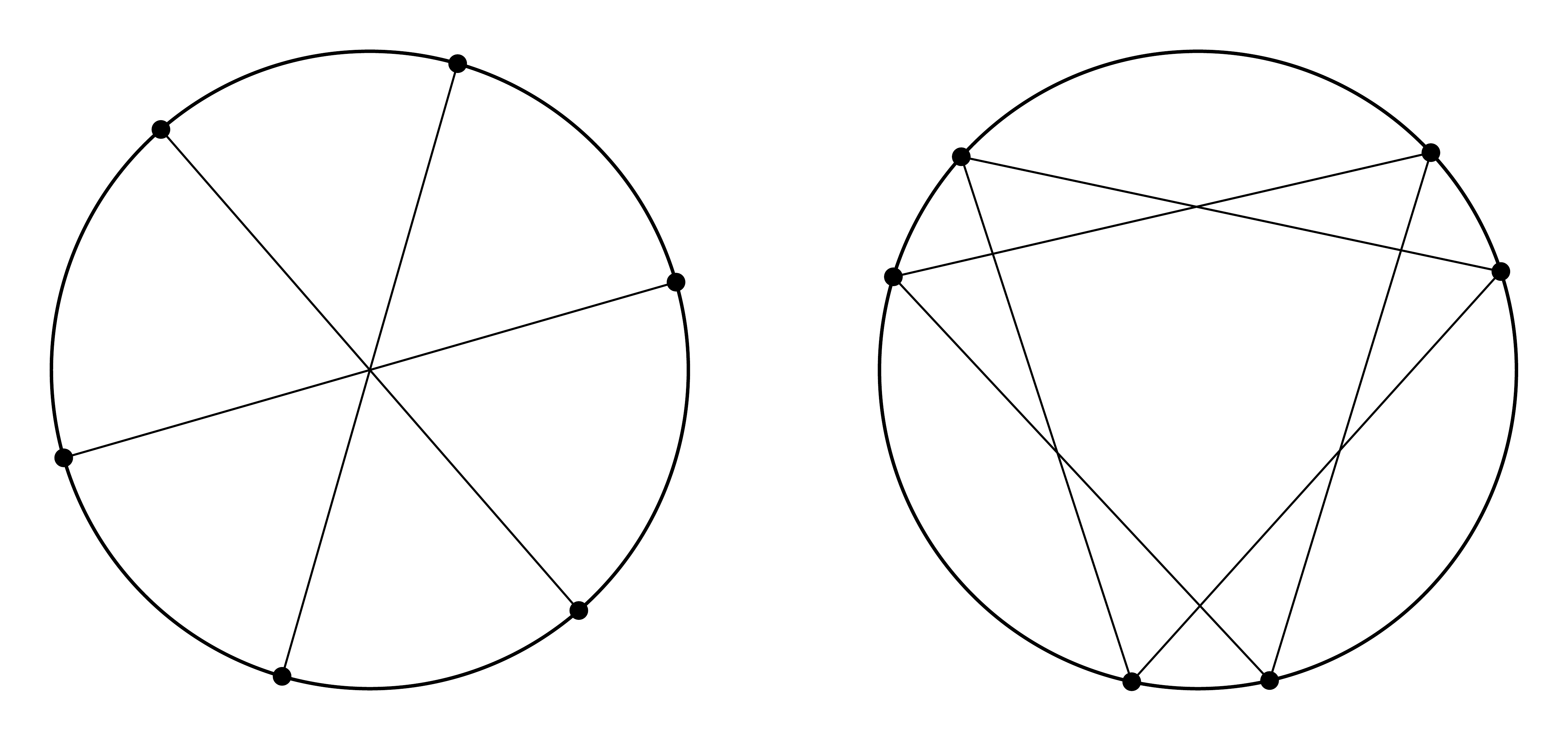}
\end{overpic}
\caption{The six-point dichotomy. Left, the roots of unity \(\zroot{n}^{d}\), \(d\in D\), for a reduced spectrum with the period \(n/2\): three antipodal pairs. Right, the same for the period \(n/3\): two rotated equilateral triangles. A reduced pair with \(n\ne 6\) has one shape on the side of \(D\) and the other on the side of \(S\).}
\label{fig:sixdichotomy}
\end{figure}

Figure~\ref{fig:sixdichotomy} shows the two shapes that a reduced six-point spectrum can take on the circle of roots of unity.

The set \(S'\) itself need not be periodic: \(S'=\{0,1,2,3,4,5\}\subset\ZN{12}\) is primitive, \(A=2\ZN{12}\) is a spectrum of it, and neither of the possible prime periods \(6\) and \(4\) fixes \(S'\). The dichotomy describes the spectra, and the reduced set \(S\) modulo \(n\), not \(S'\). In the reduced pair the two periods are moreover distributed between the two sides.

\begin{corollary}
\label{cor:sixexclusive}
Let \((D,S)\) be a reduced six-point pair in \(\ZN{n}\). Either \(n=6\) and \(D=S=\ZN{6}\), or \(n\ne 6\) and exactly one of \(D\), \(S\) has the period \(n/2\) while the other has the period \(n/3\), neither having both.
\end{corollary}

\begin{proof}
Suppose \(D+n/2=D\) and \(S+n/2=S\). Then \(n/2\in L\) and \(\widehat{\mu_D}(n/2)=\sum_{d\in D}(-1)^{d}=0\), so \(D\) has three even and three odd elements. If \(a\ge 2\), then \(n/2\) is even and \(D=D_0+\{0,n/2\}\) has an even number of even elements. So \(a=1\) and \(n/2\) is odd; the even elements \(D_e\) of \(D\) form a transversal of the antipodal pairs, \(D=D_e\sqcup(D_e+n/2)\), and likewise \(S=S_e\sqcup(S_e+n/2)\). For even \(l\in L\), \(\widehat{\mu_D}(l)=2\sum_{d\in D_e}\zroot{n}^{ld}\). The even elements of \(L\) include \((S_e-S_e)\setminus\{0\}\). Under the isomorphism \(2\ZN{n}\to\ZN{n/2}\), \(2y\mapsto y\), write \(D_e=2\widehat D\) and \(S_e=2\widehat S\); then \(\zroot{n}^{ld}=\zroot{n/2}^{2\hat l\hat d}\), so \((\widehat D,2\widehat S)\) is a three-point spectral pair in \(\ZN{n/2}\), where \(2\) is a unit. Both sets contain \(0\); since \(\langle D\rangle=\langle D_e\rangle+\{0,n/2\}=\ZN{n}\) and \(\langle D_e\rangle\subset 2\ZN{n}\), the subgroup \(\langle D_e\rangle\) has at least \(n/2\) elements, so \(\langle D_e\rangle=2\ZN{n}\) and \(\widehat D\) generates \(\ZN{n/2}\); likewise \(2\widehat S\). For a reduced three-point pair \((P,T)\) in \(\ZN{r}\) and \(\tau\in T\setminus\{0\}\), the multiset \(\tau P\) is a vanishing sum of three roots of unity containing \(0\), hence the subgroup of order three of \(\ZN{r}\); so \(3\tau P=\{0\}\), and since \(P\) generates, \(3\tau=0\). Thus \(T\) lies in the subgroup of order three, and generation gives \(r=3\). Hence \(n/2=3\) and \(n=6\).

Suppose \(D+n/3=D\) and \(S+n/3=S\). Then \(n/3\in L\) and \(\widehat{\mu_D}(n/3)=\sum_d\zroot{3}^{d}=0\) forces \(D\) to have two elements in each class modulo \(3\), because a vanishing sum \(c_0+c_1\zroot{3}+c_2\zroot{3}^2\) with nonnegative integers adding to six has \(c_0=c_1=c_2\). If \(b\ge 2\), then \(n/3\equiv 0\pmod 3\) and \(D=D_0+\Tsub{3}\) has class multiplicities divisible by three. So \(b=1\), \(\Tsub{3}\) is a transversal of the classes modulo \(3\), and \(D=D_0\sqcup(D_0+n/3)\sqcup(D_0+2n/3)\) with \(D_0=D\cap 3\ZN{n}\) of two elements containing \(0\); likewise \(S\) with \(S_0\). For \(l\in L\cap 3\ZN{n}\), \(\widehat{\mu_D}(l)=3\sum_{d\in D_0}\zroot{n}^{ld}\), and the same transfer to \(\ZN{n/3}\), where \(3\) is a unit, gives a reduced two-point pair \((\widehat D_0,3\widehat S_0)\) in \(\ZN{n/3}\). For a reduced two-point pair \((P,T)\) in \(\ZN{r}\), \(\tau P=\{0,r/2\}\) for \(\tau\in T\setminus\{0\}\), so \(2\tau=0\), \(T\subset\{0,r/2\}\), and \(r=2\). Hence \(n=6\).

For \(n=6\) the only six-point subset of \(\ZN{6}\) is \(\ZN{6}\). For \(n\ne 6\), Theorem~\ref{thm:sixdichotomy} gives each side one of the two periods, and the two cases just treated show that the sides do not share one; a side with both periods would share one with the other side.
\end{proof}

\subsection{The common complement}

In this subsection \(S'\subset\ZN{N}\) is a primitive six-point set with at least one spectrum. Translating \(S'\) changes neither \(\zeroset{S'}\) nor the spectra, so assume \(0\in S'\). Put \(a=v_2(N)\), \(b=v_3(N)\), \(L=(S'-S')\setminus\{0\}\), \(L_2=L\cap 2\ZN{N}\), \(L_3=L\cap 3\ZN{N}\). By Lemma~\ref{lem:reducedpair} and Lemma~\ref{lem:bothkinds}, \(6\) divides the modulus \(n\) of the reduced pair, which divides \(N\); so \(a,b\ge 1\), and \(L_2\), \(L_3\) are well defined. For a spectrum \(A\) the symmetry of Definition~\ref{def:reducedpair} gives \(\widehat{\mu_A}(l)=0\) \(\forall l\in L\). Since \(S'\) has six elements, two of them are congruent modulo \(3\) and two are congruent modulo \(2\), so both \(L_3\) and \(L_2\) are nonempty.

\begin{theorem}
\label{thm:sixdata}
There are integers \(\alpha\in[1,a]\) and \(\beta\in[1,b]\), depending only on \(S'\), such that for every spectrum \(A\) of \(S'\) the prime powers in \(\cycset{A}\) are exactly \(2^\alpha\) and \(3^\beta\), and \(2^\alpha3^\beta\in\cycset{A}\). Thus every spectrum satisfies {\rm (T1)} and {\rm (T2)} with the same prime powers. Explicitly: if some spectrum has the period \(N/3\), then all elements of \(L_3\) have one \(2\)-adic valuation \(i_3\), and \(\alpha=a-i_3\); if some spectrum has the period \(N/2\), then all elements of \(L_2\) have one \(3\)-adic valuation \(j_2\), and \(\beta=b-j_2\); if spectra of only one period occur, the other exponent is \(b\), respectively \(a\); if spectra of both periods occur, \((\alpha,\beta)=(a,b)\).
\end{theorem}

\begin{proof}
Let \(A\) be a spectrum; \(\cycset{A}\) is translation invariant, so assume \(0\in A\). By Theorem~\ref{thm:sixdichotomy}, \(A+N/3=A\) or \(A+N/2=A\), and by Lemma~\ref{lem:primepowers}(a) the prime powers in \(\cycset{A}\) are at most one power of \(2\) and one of \(3\). Write \(\omega=\zroot{N}\).

Case \(A+N/3=A\). Then \(A\) is a union of two cosets of \(\Tsub{3}=(N/3)\ZN{N}\), one of them \(\Tsub{3}\), so \(A=\{0,d\}+\Tsub{3}\) with \(d\notin\Tsub{3}\), and
\[
\widehat{\mu_A}(l)=(1+\omega^{ld})(1+\omega^{lN/3}+\omega^{2lN/3}).
\]
The second factor is \(0\) when \(3\nmid l\) and \(3\) when \(3\mid l\). So \(\widehat{\mu_A}(l)=0\) exactly when \(3\nmid l\) or \(ld\equiv N/2\pmod N\). An element of order \(2^j\), \(1\le j\le a\), is \(l=(N/2^j)u\) with \(u\) odd, and \(3\mid l\); the congruence \((N/2^j)ud\equiv N/2\) says \(ud\equiv 2^{j-1}\pmod{2^j}\), that is \(v_2(d)=j-1\), since \(u\) is a unit modulo \(2^j\) which fixes the class of \(2^{j-1}\). Hence \(2^j\in\cycset{A}\) iff \(v_2(d)=j-1\). For \(l\in L_3\) the vanishing of \(\widehat{\mu_A}(l)\) requires \(ld\equiv N/2\), whose \(2\)-adic valuation is \(a-1<a\); so \(v_2(l)+v_2(d)=a-1\) \(\forall l\in L_3\), all elements of \(L_3\) have the valuation \(i_3=a-1-v_2(d)\), and the power of \(2\) in \(\cycset{A}\) is \(2^{v_2(d)+1}=2^{a-i_3}\), with \(1\le a-i_3\le a\). An element of order \(3^j\) is \(l=(N/3^j)u\) with \(3\nmid u\). For \(j=b\) it is prime to \(3\), so \(3^b\in\cycset{A}\). For \(j<b\) it is divisible by \(3\), and \(ld=(N/3^j)ud\) has \(2\)-adic valuation at least \(a\) while \(N/2\) has \(a-1\); so \(3^j\notin\cycset{A}\). Elements of order \(2^{a-i_3}3^b\) are prime to \(3\), so \(2^{a-i_3}3^b\in\cycset{A}\). In this case the prime powers are \((2^{a-i_3},3^b)\).

Case \(A+N/2=A\). Then \(A=\{0,d_1,d_2\}+\Tsub{2}\) with \(0,d_1,d_2\) distinct modulo \(\Tsub{2}=\{0,N/2\}\), and
\[
\widehat{\mu_A}(l)=(1+(-1)^{l})(1+\omega^{ld_1}+\omega^{ld_2}).
\]
This vanishes when \(l\) is odd, and for even \(l\) exactly when \(1+\omega^{ld_1}+\omega^{ld_2}=0\), that is when \(\{\omega^{ld_1},\omega^{ld_2}\}=\{\zroot{3},\zroot{3}^2\}\), since two unimodular numbers summing to \(-1\) have real part \(-1/2\); that is, exactly when \(\{ld_1,ld_2\}=\{N/3,2N/3\}\). Elements of order \(2^a\) are odd, so \(2^a\in\cycset{A}\); for \(j<a\) an element \((N/2^j)u\) of order \(2^j\) is even and \((N/2^j)ud_1\) has \(3\)-adic valuation at least \(b\) while \(N/3\) has \(b-1\), so \(2^j\notin\cycset{A}\). An element of order \(3^j\) is \((N/3^j)u\) with \(3\nmid u\), and it is even. Multiplication by \(u\) permutes \(\{N/3,2N/3\}\), so \(\widehat{\mu_A}\) vanishes at all elements of order \(3^j\) exactly when \((N/3^j)\{d_1,d_2\}=\{N/3,2N/3\}\), that is when \(v_3(d_1)=v_3(d_2)=j-1\) and \(3^j\mid d_1+d_2\). For \(l\in L_2\) the vanishing of \(\widehat{\mu_A}(l)\) requires \(l\{d_1,d_2\}=\{N/3,2N/3\}\), so \(v_3(l)+v_3(d_1)=b-1\) and \(3^{b-v_3(l)}\mid d_1+d_2\). Hence all elements of \(L_2\) have the valuation \(j_2=b-1-v_3(d_1)\), the conditions for \(3^j\in\cycset{A}\) hold for \(j=b-j_2\), which lies in \([1,b]\), and the power of \(3\) in \(\cycset{A}\) is \(3^{b-j_2}\). Elements of order \(2^a3^{b-j_2}\) are odd, so \(2^a3^{b-j_2}\in\cycset{A}\). In this case the prime powers are \((2^a,3^{b-j_2})\).

Uniformity. The valuations \(i_3\) and \(j_2\) are properties of \(S'\), so within each case the pair is determined by \(S'\). Suppose \(S'\) has a spectrum \(A\) with the period \(N/3\) and a spectrum \(A'\) with the period \(N/2\). Since \(N/3\in A-A\), \(\widehat{\mu_{S'}}(N/3)=\sum_{s\in S'}\zroot{3}^{s}=0\), and as in the proof of Corollary~\ref{cor:sixexclusive} the set \(S'\) has two elements in each class modulo \(3\); since \(N/2\in A'-A'\), \(\sum_{s\in S'}(-1)^s=0\), and \(S'\) has three even and three odd elements. If \(i_3\ge 1\), every element of \(L_3\) is even, so the two elements of \(S'\) in each class modulo \(3\) have the same parity, and each parity class is a union of such pairs, of even size, which is absurd. If \(j_2\ge 1\), every element of \(L_2\) is divisible by \(3\), so the three elements of a parity class lie in one class modulo \(3\), which has two elements, again absurd. So \(i_3=j_2=0\) and both cases give \((2^a,3^b)\).

The conditions (T1) and (T2) follow: \(\lvert A\rvert=6=2\cdot 3\), and the only product of prime powers of distinct primes in \(\cycset{A}\) is \(2^\alpha3^\beta\), which lies in \(\cycset{A}\).
\end{proof}

\begin{theorem}[the common complement]
\label{thm:sixcomplement}
Let \(\alpha\), \(\beta\) be as in Theorem~\ref{thm:sixdata}, and for each prime \(p\mid N\) put \(t_p=N/p^{v_p(N)}\). Let
\[
B(X)=\prod_{p\mid N}\ \prod_{\substack{1\le c\le v_p(N)\\ p^c\notin\{2^\alpha,3^\beta\}}}\cyc{p^c}\bigl(X^{t_p}\bigr),
\]
and let \(B\subset\ZN{N}\) be the set of exponents of \(B(X)\) reduced modulo \(N\). Then \(B\) is a set with \(N/6\) elements, namely \(B=\{\sum_{p\mid N}t_pm_p:m_p\in M_p\}\), where \(M_p=\{0,1,\dots,p^{v_p(N)}-1\}\) for \(p\ge 5\), \(M_2\) is the set of integers \(m\in[0,2^a)\) whose binary digit of weight \(2^{\alpha-1}\) is \(0\), and \(M_3\) the set of integers \(m\in[0,3^b)\) whose ternary digit of weight \(3^{\beta-1}\) is \(0\). Moreover \(A\oplus B=\ZN{N}\) for every spectrum \(A\) of \(S'\).
\end{theorem}

\begin{proof}
Since \(\cyc{p^c}(Y)=\sum_{j=0}^{p-1}Y^{jp^{c-1}}\), the product over \(c\) for a fixed \(p\) expands as \(\sum_{m\in M_p}X^{t_pm}\): choosing the term \(j_c\) in each factor produces the exponent \(t_p\sum_cj_cp^{c-1}\), and \(m=\sum_cj_cp^{c-1}\) runs exactly over the integers in \([0,p^{v_p(N)})\) whose base-\(p\) digit of weight \(p^{c-1}\) is \(0\) for each omitted \(c\). The full product has the exponents \(\sum_pt_pm_p\), one for each tuple \((m_p)\). Two tuples give distinct exponents modulo \(N\): reduction modulo \(p^{v_p(N)}\) kills \(t_{p'}m_{p'}\) for \(p'\ne p\) and leaves \(t_pm_p\), and \(t_p\) is a unit modulo \(p^{v_p(N)}\), so \(m_p\) is recovered. Hence \(B\) is a set of \(\prod_p\lvert M_p\rvert=N/(2\cdot 3)\) elements and \(B(X)\equiv\sum_{\beta'\in B}X^{\beta'}\) modulo \(X^N-1\).

Let \(A\) be a spectrum, and let \(e>1\) divide \(N\). Suppose first that some prime power \(p^c\) exactly dividing \(e\) is not in \(\{2^\alpha,3^\beta\}\). Since \(\gcd(t_p,p)=1\), the roots of \(\cyc{p^c}(X^{t_p})\) are the roots of unity whose \(t_p\)th power has order \(p^c\), which are the primitive roots of order \(p^cd\) with \(d\mid t_p\); as both sides are monic and squarefree of degree \(\varphi(p^c)t_p\),
\[
\cyc{p^c}\bigl(X^{t_p}\bigr)=\prod_{d\mid t_p}\cyc{p^cd}(X).
\]
Now \(d=e/p^c\) divides \(N/p^{v_p(N)}=t_p\), so \(\cyc{e}\) divides \(B(X)\). Otherwise every prime power exactly dividing \(e\) lies in \(\{2^\alpha,3^\beta\}\), so \(e\in\{2^\alpha,3^\beta,2^\alpha3^\beta\}\subset\cycset{A}\) by Theorem~\ref{thm:sixdata}, and \(\cyc{e}\) divides \(\mask{A}\). Hence every \(\cyc{e}\) with \(e\mid N\), \(e>1\), divides \(\mask{A}B(X)\). Let \(P(X)\) be the remainder of \(\mask{A}B(X)\) modulo \(X^N-1\), the polynomial of degree less than \(N\) whose coefficient of \(X^r\) is the number of pairs \((a,\beta')\in A\times B\) with \(a+\beta'\equiv r\pmod N\). Each \(\cyc{e}\) with \(e\mid N\), \(e>1\), divides \(X^N-1\), hence divides \(P\); so \((X^N-1)/(X-1)\) divides \(P\), and \(P=c\,(1+X+\dots+X^{N-1})\) for a constant \(c\). Evaluating at \(1\) gives \(cN=\lvert A\rvert\lvert B\rvert=N\), so \(c=1\): every element of \(\ZN{N}\) is \(a+\beta'\) in exactly one way.
\end{proof}

\begin{remark}
\label{rem:digitbox}
A second common complement avoids cyclotomic polynomials. Let \(B(\alpha,\beta)\) be the set of \(g\in\ZN{N}\) whose reduction modulo \(2^a\) has binary digit of weight \(2^{\alpha-1}\) equal to \(0\) and whose reduction modulo \(3^b\) has ternary digit of weight \(3^{\beta-1}\) equal to \(0\). Then \(A\oplus B(\alpha,\beta)=\ZN{N}\) for every spectrum \(A\) of \(S'\), by a greedy matching of digits. Translating \(A\) does not change whether it tiles with a fixed complement, so we may assume \(0\in A\). The proof of Theorem~\ref{thm:sixdata} shows that in a spectrum \(\{0,d\}+\Tsub{3}\) the element \(d\) has binary digit \(1\) at the weight \(2^{\alpha-1}\) and \(0\) below, while \(\Tsub{3}\) reduces to \(0\) modulo \(2^a\) and runs through the ternary digits \(0,1,2\) at the weight \(3^{b-1}=3^{\beta-1}\) with \(0\) below; so a given \(g\) is brought into \(B(\alpha,\beta)\) by subtracting \(0\) or \(d\) and then an element of \(\Tsub{3}\), no carry disturbing the cleared digits, and the same digit reading gives uniqueness. A spectrum \(\{0,d_1,d_2\}+\Tsub{2}\) is treated with the roles of \(2\) and \(3\) exchanged, \(d_1,d_2\) carrying the ternary digits \(1\) and \(2\) at the weight \(3^{\beta-1}\). Both kinds of spectrum, beyond the coset \(\Tsub{6}\), do occur for one set: for \(N=36\) and \(S'=\{0,1,8,9,17,28\}\), whose classes modulo \(3\) are \(\{0,9\}\), \(\{1,28\}\), \(\{8,17\}\) and whose parity classes are \(\{0,8,28\}\), \(\{1,9,17\}\), one has \(L_3=\{9,27\}\) and \(L_2=\{8,16,20,28\}\); the criteria in the proof of Theorem~\ref{thm:sixdata} show that \(\{0,2\}+12\ZN{36}\) is a spectrum, since \(9\cdot 2\equiv 27\cdot 2\equiv 18\pmod{36}\), and that \(\{0,3,6\}+18\ZN{36}\) is a spectrum, since \(l\{3,6\}=\{12,24\}\) for each \(l\in L_2\). Here \((\alpha,\beta)=(2,2)\) and \(B(2,2)=\{0,1,9,20,28,29\}\).
\end{remark}

\begin{theorem}
\label{thm:fibersix}
The statements \(\UC{6}\) and \(\FL{6}\) hold. More precisely, for every primitive six-point \(S'\subset\ZN{N}\) with a spectrum, the set \(B\) of Theorem~\ref{thm:sixcomplement} satisfies \(A\oplus B=\ZN{N}\) for every spectrum \(A\) of \(S'\); and every bounded measurable \(E\subset\R\) with a spectrum \(\Gamma+h^{-1}\Z\), where the six elements of \(\Gamma\) are distinct modulo \(h^{-1}\Z\), tiles \(\R\) by \(h\widetilde B\), where \(\widetilde B=B+N\Z\) is built from \(S=h\Gamma\) through Lemma~\ref{lem:reduction} and Theorem~\ref{thm:sixcomplement}.
\end{theorem}

\begin{proof}
The first statement is Theorem~\ref{thm:sixcomplement}; translating \(S'\) to contain \(0\) changes neither its spectra nor \(B\), since \(\cycset{A}\) and hence \(\alpha,\beta\) are unchanged. The second is Lemma~\ref{lem:ucfl} with \(k=6\), whose proof produces the translation set \(h\{b\in\Z:b\bmod N\in B\}\).
\end{proof}

\subsection{Seven points}

In this subsection \((D,S)\) is a reduced seven-point pair in \(\ZN{n}\), with \(a=v_2(n)\), \(b=v_3(n)\), \(c=v_5(n)\), \(e=v_7(n)\). Every row and every column is of one of the kinds (a) to (d) of Lemma~\ref{lem:rowkinds}(b). The rows of kind (a) are described first; then the reductions of the bad kinds are listed; then the bad kinds are excluded one against the other.

\begin{lemma}
\label{lem:kinda}
A row \(lD\), \(l\in L\), is of kind {\rm (a)} iff \(l\) has order seven. In that case \(7\mid n\), \(lD=\Tsub{7}\), \(D\) is a complete residue system modulo \(7\), and \(v_q(l)=v_q(n)\) for every prime \(q\mid n\) with \(q\ne 7\) while \(v_7(l)=e-1\).
\end{lemma}

\begin{proof}
If \(lD=x+\Tsub{7}\), then \(0=l\cdot 0\in lD\), so \(lD=\Tsub{7}\) and \(7ld=0\) \(\forall d\in D\); since \(D\) generates \(\ZN{n}\), \(7l=0\), that is \(l\in\Tsub{7}\setminus\{0\}\). Conversely, suppose \(l\) has order seven, so \(7\mid n\) and \(l\in\Tsub{7}\setminus\{0\}\). Then \(lD\subset\Tsub{7}\) and \(\widehat{\mu_D}(l)=\sum_{j\in\ZN{7}}c_j\zroot{7}^{j}\) with nonnegative integers \(c_j\) adding to seven; since \(1,\zroot{7},\dots,\zroot{7}^{6}\) satisfy only the relations proportional to \(\cyc{7}\), the vanishing forces all \(c_j\) equal, so \(c_j=1\) and \(lD=\Tsub{7}\) with seven distinct elements. Writing \(l=jn/7\) with \(7\nmid j\), one has \(ld=ld'\) exactly when \(7\mid d-d'\), so \(D\) is a complete residue system modulo \(7\); and since \(l=jn/7\) with \(7\nmid j\), one has \(v_q(l)=v_q(n)\) for primes \(q\mid n\) with \(q\ne 7\), and \(v_7(l)=e-1\).
\end{proof}

The reductions of the bad kinds are read off as in the six-point case. All reductions are at primes dividing \(n\). For the kinds requiring the respective fractions, write \(h=2^{a-1}\in\ZN{2^a}\) and \(w=\rho_3(n/6)\), an element of order three, and \(u=\rho_5(n/5)\), of order five; recall \(\rho_2(n/6)=h\), \(\rho_2(n/5)=0\), \(\rho_3(n/5)=0\), \(\rho_5(n/6)=0\). For a row \(lD\) of kind (b): \(\rho_2(lD)=\rho_2(x)+\{0,0,0,h,h,h,h\}\), the pattern \((4,3)\) with the two values differing by \(h\); \(\rho_3(lD)=\rho_3(x)+\{0,0,0,w,w,-w,-w\}\), the pattern \((3,2,2)\) with three values in one coset of the subgroup of order three; \(\rho_5(lD)\) is \(\rho_5(x)\) plus the three values \(ju\), \(j\notin X\), once each and the two values \(x_iu\) twice each, the pattern \((2,2,1,1,1)\) with five values in one coset of the subgroup of order five; and for primes \(q\mid n\) with \(q\notin\{2,3,5\}\), \(\rho_q(lD)\) is constant. For kind (c), \(lD=(y+\Tsub{5})\cup\{z,z+n/2\}\): \(\rho_2(lD)=\{Y^5,Z,Z+h\}\), the pattern \((5,1,1)\) or \((6,1)\); \(\rho_5(lD)\) has the five distinct values of a coset of the subgroup of order five together with a double value; and \(\rho_q(lD)=\{Y^5,Z^2\}\) for every prime \(q\mid n\) with \(q\notin\{2,5\}\), with at most two distinct values. For kind (d), \(lD=(y+\Tsub{3})\cup\{z_1,z_1+n/2\}\cup\{z_2,z_2+n/2\}\): \(\rho_2(lD)=\{Y^3,Z_1,Z_1+h,Z_2,Z_2+h\}\); \(\rho_3(lD)=\{Y,Y+w',Y-w',Z_1^2,Z_2^2\}\) with \(w'=\rho_3(n/3)\) of order three, so with at least three distinct values; and \(\rho_q(lD)=\{Y^3,Z_1^2,Z_2^2\}\) for primes \(q\mid n\) with \(q\notin\{2,3\}\), with at most three distinct values. For kind (a), \(\rho_q(lD)\) is constant for every prime \(q\mid n\) with \(q\ne 7\). By Lemma~\ref{lem:collapseprinciple}: a row of kind (c) has \(v_5(l)<c\) and \(D\) modulo \(5^{c-v_5(l)}\) has at least five distinct values; a row of kind (d) has \(v_3(l)<b\) and \(D\) modulo \(3^{b-v_3(l)}\) has at least three distinct values; every bad row has \(v_2(l)<a\), because the displayed reductions of all three bad kinds have at least two distinct \(\rho_2\)-values; and a row with \(v_p(l)=v_p(n)\) has constant \(\rho_p(lD)\).

\begin{theorem}
\label{thm:nokindb}
No row and no column of a reduced seven-point pair is of kind {\rm (b)}.
\end{theorem}

\begin{proof}
By symmetry it suffices to treat rows. Suppose \(lD\) is of kind (b) for some \(l\in L\); then \(30\mid n\).

(1) Let \(j=v_5(l)\); then \(j<c\), and \(D\) modulo \(5^{c-j}\) has the pattern \((2,2,1,1,1)\) with its five values in one coset of the subgroup of order five, so \(D\) is constant modulo \(5^{c-j-1}\) (Lemma~\ref{lem:collapseprinciple}), and (G) gives \(j=c-1\). Hence \(D\) modulo \(5\) has the pattern \((2,2,1,1,1)\): every class modulo \(5\) meets \(D\), two of them twice.

(2) Let \(l'\in L\) with \(j'=v_5(l')<c\). The multiset \(D\) modulo \(5^{c-j'}\) reduces to \(D\) modulo \(5\), so it has at least five distinct values and no multiplicity above two. A row of kind (a) has \(v_5=c\) (Lemma~\ref{lem:kinda}); a row of kind (d) gives at most three distinct values; a row of kind (c) gives five values congruent modulo \(5^{c-j'-1}\) together with a double value, so if \(c-j'\ge 2\) the reduction of \(D\) modulo \(5^{c-j'-1}\), hence modulo \(5\), has a value of multiplicity at least five, and if \(c-j'=1\) the five values fill \(\ZN{5}\) and the double value gives the pattern \((3,1,1,1,1)\); neither is \((2,2,1,1,1)\). So \(l'D\) is of kind (b) and \(j'=c-1\) by (1). Thus every element of \(L\) has \(5\)-adic valuation at least \(c-1\), so \(S\subset 5^{c-1}\ZN{n}\), and generation gives \(c=1\). Summary: \(5\) divides \(n\) exactly once, and a row \(l'D\) is of kind (b) iff \(5\nmid l'\).

(3) For the row \(lD\) of kind (b) and \(i=v_2(l)<a\), the multiset \(D\) modulo \(2^{a-i}\) has the pattern \((4,3)\) with its two values differing by \(2^{a-i-1}\); so \(D\) is constant modulo \(2^{a-i-1}\), and (G) gives \(i=a-1\). Likewise \(v_3(l)=b-1\), since the three \(\rho_3\)-values lie in one coset of the subgroup of order three. Thus every row prime to \(5\) has \(v_2=a-1\) and \(v_3=b-1\), and \(D\) modulo \(2\) has the pattern \((4,3)\).

(4) Some \(s_1\in S\) is prime to \(5\); it lies in \(L\), is of kind (b), and \(\rho_2(s_1)=h\). If \(s,s'\in S\) lie in distinct classes modulo \(5\), then \(\rho_2(s-s')=h\); three elements in three distinct classes would give \(h=\rho_2(s-s'')=\rho_2(s-s')+\rho_2(s'-s'')=h+h=0\). So \(S=S_0\sqcup S_1\) with \(S_0\subset 5\ZN{n}\), \(0\in S_0\), and \(S_1\) in one nonzero class \(r\) modulo \(5\). For \(s\in S_1\), \(\rho_2(s)=\rho_2(s-0)=h\); for \(s\in S_0\), \(\rho_2(s)=\rho_2(s-s_1)+\rho_2(s_1)=0\). Hence \(S\subset 2^{a-1}\ZN{n}\) and (G) gives \(a=1\): \(S_0\) is even and \(S_1\) is odd. In the same way \(\rho_3(s)\in 3^{b-1}\ZN{3^b}\) \(\forall s\in S\), so \(b=1\). Thus \(n=30m'\) with \(\gcd(m',30)=1\); \(n/2\) is odd and \(n/3\) is even.

(5) Let \(\delta\in\Delta\) with \(5\nmid\delta\). The column \(\delta S\) modulo \(5\) takes the two values \(0\) and \(\delta r\), with multiplicities \(\lvert S_0\rvert\) and \(\lvert S_1\rvert\). Kinds (a), (b), (c) have one, five, and at least five distinct values modulo \(5\); so \(\delta S\) is of kind (d), and its pattern \(\{Y^3,Z_1^2,Z_2^2\}\) modulo \(5\) collapses to two values only as \((5,2)\) or \((4,3)\). In \(\delta S=(y+\Tsub{3})\cup\{z_1,z_1+n/2\}\cup\{z_2,z_2+n/2\}\), the coset of \(\Tsub{3}\) has constant parity because \(n/3\) is even, and each antipodal pair has one even and one odd element because \(n/2\) is odd; so \(\delta S\) has five elements of one parity and two of the other. If \(\delta\) were even, every element of \(\delta S\) would be even. So every column prime to \(5\) is odd.

(6) Since \(D\) modulo \(2\) has the pattern \((4,3)\) and \(0\in D\), \(D\) has at least two nonzero even elements \(d',d''\). They are even columns, so \(5\mid d'\) and \(5\mid d''\) by (5). With \(0\), the class of \(0\) modulo \(5\) contains three elements of \(D\), while (1) allows at most two. Contradiction.
\end{proof}

\begin{lemma}
\label{lem:collapseseven}
{\rm (a)} If some row is of kind {\rm (a)} and some row is bad, a contradiction follows. {\rm (b)} If \(3\mid n\) and every row prime to \(3\) is of kind {\rm (c)}, a contradiction follows.
\end{lemma}

\begin{proof}
(a) By Lemma~\ref{lem:kinda}, \(7\mid n\) and \(D\) is a complete residue system modulo \(7\), so \(D\) modulo \(7^{e}\) has seven distinct values, and so has every reduction of it modulo a positive power of \(7\). A bad row prime to \(7\) has at most three distinct \(\rho_7\)-values, while \(\rho_7(lD)\) has the pattern of \(D\) modulo \(7^e\) for \(7\nmid l\). So every bad row is divisible by \(7\). If \(e\ge 2\), rows of kind (a) have \(v_7=e-1\ge 1\), so every row is divisible by \(7\), \(S\subset 7\ZN{n}\), against generation. If \(e=1\), the rows of kind (a), the nonzero elements of \(\Tsub{7}\), are prime to \(7\), the bad rows are divisible by \(7\), and every row is of one of these types; so two elements of \(S\) in distinct classes modulo \(7\) differ by an element of \(\Tsub{7}\). The set \(S\) meets at least two classes modulo \(7\). Fix \(s'\in S\); every element of \(S\) outside the class of \(s'\) lies in \(s'+\Tsub{7}\), and every element \(s''\) in the class of \(s'\) lies in \(s+\Tsub{7}=s'+\Tsub{7}\) for any \(s\in S\) outside that class. So \(S\subset s'+\Tsub{7}\), and \(0\in S\) gives \(S\subset\Tsub{7}\); generation gives \(n=7\), where no bad row exists since \(2\nmid 7\). Contradiction.

(b) For \(l\in L\) prime to \(3\), \(\rho_3(lD)=\{\rho_3(y)^5,\rho_3(z)^2\}=\rho_3(l)\rho_3(D)\) with \(\rho_3(l)\) a unit, so \(D\) modulo \(3^b\) has the pattern \((5,2)\), the pattern \((7)\) being excluded by (G). Let \(D=D_5\sqcup D_2\) according to the two values of \(\rho_3\), with \(\lvert D_5\rvert=5\). For every \(l\in L\) prime to \(3\), the coset of \(\Tsub{5}\) in \(lD\) has constant \(\rho_3\), since \(3\mid n/5\), and it is the class of multiplicity five of \(lD\), which is \(lD_5\). So \(l(D_5-D_5)\subset\Tsub{5}\) \(\forall l\in L_{(3)}\), and (R) gives \(D_5-D_5\subset\Tsub{5}\): five distinct elements pairwise congruent modulo \(\Tsub{5}\) form a coset, \(D_5=y'+\Tsub{5}\). Write \(D_2=\{d,d'\}\). For \(t\in\ZN{n}\),
\[
\widehat{\mu_D}(t)=\zroot{n}^{y't}\sum_{j\in\ZN{5}}\zroot{5}^{jt}+\zroot{n}^{dt}\bigl(1+\zroot{n}^{(d'-d)t}\bigr).
\]
If \(5\mid t\) the first term has modulus five and the second at most two, so \(\widehat{\mu_D}(t)\ne 0\). If \(5\nmid t\) the first term vanishes, and \(\widehat{\mu_D}(t)=0\) exactly when \((d'-d)t=n/2\). So \(L\) is contained in the set \(Z'=\{t:(d'-d)t=n/2\}\), and \(Z'\) contains no sum of two of its elements, since \((d'-d)(t+t')=0\) for \(t,t'\in Z'\). But for three distinct elements \(s,s',s''\) of \(S\), the differences \(s-s'\), \(s'-s''\) and their sum \(s-s''\) all lie in \(L\). Contradiction.
\end{proof}

\begin{theorem}[seven-point rigidity]
\label{thm:sevenrigid}
Let \((D,S)\) be a reduced seven-point pair in \(\ZN{n}\). Then \(n=7\) and \(D=S=\ZN{7}\).
\end{theorem}

\begin{proof}
If no row is bad, every row is of kind (a), so \(L\subset\Tsub{7}\) by Lemma~\ref{lem:kinda}, \(S\subset\Tsub{7}\) because \(0\in S\), and generation gives \(n=7\); then \(D=S=\ZN{7}\). Suppose some row is bad. By Theorem~\ref{thm:nokindb} and Lemma~\ref{lem:collapseseven}(a), every row is of kind (c) or (d), and \(n\) is even.

Case \(3\nmid n\). Kind (d) needs \(3\mid n\), so every row is of kind (c) and \(5\mid n\). Some \(s_1\in S\) is odd, and \(\rho_2(s_1)\) is a unit, so \(D\) modulo \(2^a\) has the form \(\{Y'^5,Z',Z'+h\}\) up to the unit, which fixes \(h\). Two of the seven elements of \(S\) have the same parity, so some \(l\in L\) is even, \(i=v_2(l)\ge 1\); and \(i<a\), because a row of kind (c) has two distinct \(\rho_2\)-values. Then \(D\) modulo \(2^{a-i}\) is the reduction of \(\{Y'^5,Z',Z'+h\}\), in which \(Z'\) and \(Z'+h\) coincide because \(h\equiv 0\) modulo \(2^{a-i}\): its pattern is \((5,2)\) or \((7)\). But \(D\) modulo \(2^{a-i}\) has the pattern of \(\rho_2(lD)\), which is \((5,1,1)\) or \((6,1)\) for a row of kind (c). Contradiction.

Case \(3\mid n\). Rows prime to \(3\) exist by (R), and by Lemma~\ref{lem:collapseseven}(b) one of them, \(l_0D\), is of kind (d). So \(D\) modulo \(3^b\) has the pattern of \(\{Y,Y+w',Y-w',Z_1^2,Z_2^2\}\), with at least three distinct values, and a row of kind (c) prime to \(3\) would give at most two. Hence every row prime to \(3\) is of kind (d).

Every row divisible by \(3\) is of kind (c). A row with \(v_3(l)=b\) has constant \(\rho_3(lD)\), which kind (d) forbids. A row of kind (d) with \(1\le j=v_3(l)<b\) would make \(D\) modulo \(3^{b-j}\) have the pattern of a multiset \(\{Y,Y+w',Y-w',Z_1^2,Z_2^2\}\), which is one of \((1,1,1,2,2)\), \((1,1,1,4)\), \((3,1,1,2)\), \((3,3,1)\), \((5,1,1)\) according to the coincidences of \(Z_1,Z_2\) with the three distinct values \(Y,Y\pm w'\) and with each other. On the other hand \(D\) modulo \(3^{b-j}\) is the reduction of \(D\) modulo \(3^b\), which is a unit times \(\{Y,Y+w',Y-w',Z_1^2,Z_2^2\}\) from the row \(l_0\); since \(w'\) has order three it lies in \(3^{b-1}\ZN{3^b}\) and reduces to \(0\) modulo \(3^{b-j}\), so the three values \(Y,Y\pm w'\) merge and the reduction is \(\{\bar Y^3,\bar Z_1^2,\bar Z_2^2\}\), of pattern \((3,2,2)\), \((5,2)\), \((4,3)\) or \((7)\). The two lists are disjoint. So rows divisible by \(3\) are of kind (c). Such rows exist, since seven elements of \(S\) fall into three classes modulo \(3\); hence \(5\mid n\), and every row of kind (c) has \(v_5(l)<c\) with \(D\) modulo \(5^{c-v_5(l)}\) having at least five distinct values.

Rows prime to \(5\) exist by (R). Suppose one of them, \(l_1D\), is of kind (c). Then \(D\) modulo \(5^c\) has at least five distinct values, so every row prime to \(5\) is of kind (c), a row of kind (d) prime to \(5\) giving at most three; hence every row prime to \(5\) is divisible by \(3\), and, since every row of kind (d) is prime to \(3\), every row of kind (d) is divisible by \(5\). The set \(S\) meets at least two classes modulo \(3\). An element \(s\in S\) not divisible by \(3\) is a row prime to \(3\), of kind (d), hence divisible by \(5\). An element \(s\in S\) divisible by \(3\) satisfies \(s=(s-s')-(0-s')\) for any \(s'\in S\) in another class modulo \(3\), and both differences are rows prime to \(3\), hence divisible by \(5\). So \(S\subset 5\ZN{n}\), against generation. Therefore every row prime to \(5\) is of kind (d), and \(D\) modulo \(5^c\) has at most three distinct values. But a row of kind (c) makes \(D\) modulo \(5^{c-v_5(l)}\), a reduction of \(D\) modulo \(5^c\), have at least five distinct values, and reduction does not increase the number of values. Contradiction.

So no row is bad, and the first paragraph gives \(n=7\) and \(D=S=\ZN{7}\).
\end{proof}

\begin{corollary}
\label{cor:sevencosets}
Let \(S'\subset\ZN{N}\) be a primitive seven-point set with a spectrum \(A\). Then \(7\mid N\), \(A=a_0+\Tsub{7}=a_0+(N/7)\ZN{N}\) \(\forall a_0\in A\), and \(S'\) is a complete residue system modulo \(7\). Consequently, if \(S\subset\R/\Z\) and \(A\subset\Z\) with \(\lvert S\rvert=\lvert A\rvert=7\) form a spectral pair, there are bijections \(\varphi\colon S\to\ZN{7}\) and \(\psi\colon A\to\ZN{7}\) with \(\ex{(s-s_0)(a-a_0)}=\zroot{7}^{\varphi(s)\psi(a)}\) for fixed \(s_0\in S\), \(a_0\in A\): every complex Hadamard matrix of order seven of the form \((\ex{sa})_{s\in S,\,a\in A}\) is equivalent to the Fourier matrix \(F_7\).
\end{corollary}

\begin{proof}
Lemma~\ref{lem:reducedpair} produces from \((A,S')\) a reduced seven-point pair \((D,S)\) in \(\ZN{n}\), \(n=N/g\), with \(D\) the image of \(A-a_0\) under \(g\ZN{N}\to\ZN{n}\) and \(S\) the reduction of \(S'-s_0\) modulo \(n\). Theorem~\ref{thm:sevenrigid} gives \(n=7\) and \(D=\ZN{7}\), so \(A-a_0=g\ZN{N}=(N/7)\ZN{N}\), and \(S'-s_0\) reduces bijectively onto \(\ZN{7}\). For the second statement, Theorem~\ref{thm:rational} gives \(S-S\subset\Q/\Z\), and Lemma~\ref{lem:reduction} gives \(N\), a primitive \(S'=N(S-s_0)\subset\ZN{N}\) (Definition~\ref{def:zeroset}), and the reduction \(A\bmod N\), a spectrum of \(S'\). With \(a-a_0\equiv\psi(a)N/7\pmod N\) and \(\varphi(s)=N(s-s_0)\bmod 7\), both bijections by the first part, \(\ex{(s-s_0)(a-a_0)}=\zroot{N}^{N(s-s_0)\psi(a)N/7}=\zroot{7}^{\varphi(s)\psi(a)}\).
\end{proof}

\begin{theorem}
\label{thm:fiberseven}
The statements \(\UC{7}\) and \(\FL{7}\) hold. For every primitive seven-point \(S'\subset\ZN{N}\) with a spectrum, \(7\mid N\) and \(B=\{0,1,\dots,N/7-1\}\) satisfies \(A\oplus B=\ZN{N}\) for every spectrum \(A\) of \(S'\); and every bounded measurable \(E\subset\R\) with a spectrum \(\Gamma+h^{-1}\Z\), where the seven elements of \(\Gamma\) are distinct modulo \(h^{-1}\Z\), tiles \(\R\) by \(h(\{0,1,\dots,N/7-1\}+N\Z)\), where \(N\) is the integer of Lemma~\ref{lem:reduction} for \(S=h\Gamma\).
\end{theorem}

\begin{proof}
Every residue modulo \(N\) is uniquely \(q\cdot N/7+r\) with \(q\in\ZN{7}\) and \(0\le r<N/7\), so \(\Tsub{7}\oplus B=\ZN{N}\), and by Corollary~\ref{cor:sevencosets} every spectrum is a translate of \(\Tsub{7}\), so \(A\oplus B=\ZN{N}\). Lemma~\ref{lem:ucfl} with \(k=7\) gives the tiling, with the translation set \(h\{b\in\Z:b\bmod N\in B\}\) read from its proof.
\end{proof}

\subsection{Proof of Theorem C, and eight points}

\begin{proof}[Proof of Theorem~\ref{thm:C}]
Let \(k,h,F\) and \(\Gamma\) be as in Theorem~C, and put \(S=h\Gamma\subset\R/\Z\), so \(\lvert S\rvert=k\). For \(k\le 5\) the tiling by a set \(hB\), \(B\subset\Z\) depending only on \(S=h\Gamma\), is Theorem~\ref{thm:fiberfew}. For \(k=6\) it is Theorem~\ref{thm:fibersix}, and for \(k=7\) Theorem~\ref{thm:fiberseven}. The dichotomy for the spectra of a primitive six-point set is Theorem~\ref{thm:sixdichotomy}, and the coset statement for seven points, with \(S'\) a complete residue system modulo \(7\), is Corollary~\ref{cor:sevencosets}.
\end{proof}

\begin{remark}
\label{rem:eightpoints}
The case of eight points is not treated here. The six-point argument identifies a difficulty in extending this approach: a residue class of size four may be split into two antipodal pairs in different ways by different rows. The collapse argument needs these matchings to agree before it can recover a period of the original set. At six points the second row type supplies the additional information used in Lemma~\ref{lem:collapse}; at eight points this argument alone does not supply an analogous compatibility condition.
\end{remark}

\section{General implications for Fuglede's conjecture}
\label{sec:general}

Theorem~A uses \(\dimQ \Krel \le 4\) twice. Through the specialization lemma it gives the bound \(\kzero(\ell) \le 10\) of Theorem~\ref{thm:kzeroten}, and through the rank bound of Lemma~\ref{lem:rankbound} it limits a fully mixing configuration to at most eight intervals, with the class patterns of Corollary~\ref{cor:fullymixing} that Sections~\ref{sec:classgraph} and~\ref{sec:mechanisms} analyze. This section develops the argument for arbitrary lengths. The result, Theorem~E, reduces the problem for \(n\) intervals to statements about finite spectral pairs of finitely many possible cardinalities, each still quantified over all ambient cyclic groups, and the examples of Corollary~\ref{cor:alllargern} show that for every \(n\ge8\) the finite statements needed include fibers with more points than intervals. The section closes with the status of six intervals and with remarks on what generalizes and what does not.

\subsection{The reduction theorem}

\begin{definition}
\label{def:cprun}
For integers \(n \ge 1\) and \(k' \ge 1\), the statement \((\mathrm{CP}_{n,k'})\) is the consecutive-pair statement \(\CP{k'}\) of Definition~\ref{def:namedstatements} restricted to spectral pairs \((A, S)\) in which \(A \subset \Z\) is a union of at most \(n\) runs: if \((A, S)\) is a spectral pair with \(A \subset \Z\), \(\lvert A\rvert = \lvert S\rvert = k'\), \(A\) a union of at most \(n\) runs and \(A\) containing two consecutive integers, then \(S + 1/m = S\) for some integer \(m \ge 2\). When \(k' > n\), a consecutive pair is automatic. Clearly \(\CP{k'}\) implies \((\mathrm{CP}_{n,k'})\) \(\forall n\).
\end{definition}

\begin{theorem}
\label{thm:generalreduction}
Let \(\Omega\) be a spectral union of \(n\) intervals in the standing reduction of Definition~\ref{def:standing}, with length vector \(\ell\), and let \(\Lambda\) be a spectrum with minimal period \(h^{-1} = k/L\).
\begin{enumerate}
\item[\rm (a)] \(k\) divides \(\kzero(\ell)\), and \(\kzero(\ell)\) divides \(\Kn\). If \(\dimQ \Krel \le 4\), then \(\kzero(\ell) \le 10\).
\item[\rm (b)] If \((\mathrm{CP}_{n,k'})\) holds for every \(k'\) with \(n < k' \le \Kn\), then \(k \le n\). The same conclusion holds if \((\mathrm{CP}_{n,k'})\) holds for every \(k'\) with \(n < k' \le \kzero(\ell)\).
\item[\rm (c)] If the fiber lemma \(\FL{k'}\) holds for every \(k' \le \Kn\), then \(\Omega\) tiles \(\R\) by translations. More precisely, \(\FL{k}\) for the actual value of \(k\) suffices, and when {\rm (b)} applies it suffices that \(\FL{k'}\) hold for every \(k' \le n\).
\end{enumerate}
\end{theorem}

\begin{proof}
(a) The divisibility \(k \mid \kzero(\ell)\) is Theorem~\ref{thm:verdict}, the divisibility \(\kzero(\ell) \mid \Kn\) is Theorem~\ref{thm:Kn}, and the bound under \(\dimQ \Krel \le 4\) is Theorem~\ref{thm:kzeroten}.

(b) Suppose \(k > n\). By (a), \(k \le \kzero(\ell) \le \Kn\). Fix any \(x \in \R\). By Lemma~\ref{lem:everyfiber} the fiber \(A_x\) is a spectrum of \(S\), so \((A_x, S)\) is a spectral pair with \(A_x \subset \Z\) and \(\lvert A_x\rvert = \lvert S\rvert = k\); by Proposition~\ref{prop:intervalfibers} the set \(A_x\) is a union of at most \(n\) runs, and by Lemma~\ref{lem:pigeonhole}(i) it contains two consecutive integers. The statement \((\mathrm{CP}_{n,k})\), which is among the hypotheses in either form, gives \(S + 1/m = S\) for some integer \(m \ge 2\), a nontrivial translation symmetry of \(S\). This contradicts Theorem~\ref{thm:periodic}. Hence \(k \le n\).

(c) By Lemma~\ref{lem:everyfiber}, \(\Omega\) is a \(k\)-fold multitile of \(\R\) by \(h\Z\) with the spectrum \(\Lambda = \Gamma + h^{-1}\Z\), where the \(k\) elements of \(\Gamma\) are pairwise distinct modulo \(h^{-1}\Z\). Thus \(\FL{k}\) gives a tiling of \(\R\) by \(\Omega\). By (a), \(k \le \Kn\); when (b) applies, \(k \le n\).
\end{proof}

\begin{proof}[Proof of Theorem~E]
The divisibilities and the bound \(\kzero(\ell) \le 10\) under \(\dimQ \Krel \le 4\) are part (a) of Theorem~\ref{thm:generalreduction}, and the tiling statement is part (c).
\end{proof}

The proof of (b) uses no induction on the number of intervals and no merging: the verdict rule does not see rigid pairs. If \((\mathrm{CP}_{n,k'})\) were available for every \(k' > n\), the bound \(\Kn\) would not be needed, since Lemma~\ref{lem:pigeonhole}(i) and the consecutive-pair statement give \(k \le n\) directly; the role of the level theory is to make the set of consecutive-pair statements to be verified finite. At \(n = 5\), where \(\kzero(\ell) \le 10\) by Theorem~\ref{thm:fivekzero} and \(\CP{k'}\) holds for \(6 \le k' \le 10\) by Corollary~\ref{cor:cpten}, the theorem gives \(k \le 5\) and the tiling by Theorem~\ref{thm:fiberfew}. For six or seven runs, Theorem~D gives the consecutive-pair statements up to eleven points; the larger cardinalities are not resolved by these results. For \(n = 8\) the conclusion \(k \le n\) is false, as the next proposition shows, and Corollary~\ref{cor:alllargern} extends this to every \(n\ge8\).

\subsection{The example \texorpdfstring{\(\OmegaEight\)}{Omega\_8}}

Recall from Theorem~\ref{thm:cptwelve} the twelve-point set \(\Ssev \subset \ZN{72}\) and the set
\[
A' = \{0, 1, 8, 9, 17, 28\} + \{0, 36\} = \{0, 1, 8, 9, 17, 28, 36, 37, 44, 45, 53, 64\} \subset \Z,
\]
whose reduction modulo \(72\) is a spectrum of \(\Ssev\) in \(\ZN{72}\), while \(\Ssev\) has no period in \(\ZN{72}\). The set \(A'\) is the union of the eight runs \(\{0, 1\}\), \(\{8, 9\}\), \(\{17\}\), \(\{28\}\), \(\{36, 37\}\), \(\{44, 45\}\), \(\{53\}\), \(\{64\}\).

\begin{proposition}
\label{prop:omegaeight}
Let \(\OmegaEight = A' + [0, 1) \subset \R\). Then \(\OmegaEight\) is a union of eight pairwise disjoint non-adjacent intervals, namely \([0, 2)\), \([8, 10)\), \([17, 18)\), \([28, 29)\), \([36, 38)\), \([44, 46)\), \([53, 54)\), \([64, 65)\), of total length \(12\). The set \(\Lambda = \Ssev/72 + \Z\) is a spectrum of \(\OmegaEight\) with minimal period \(1\), so that \(h = 1\) and \(k = 12 > 8 = n\). The set \(\OmegaEight\) tiles \(\R\) by translations:
\[
\OmegaEight \oplus \bigl(\{0, 6, 12, 18, 24, 30\} + 72\Z\bigr) = \R .
\]
\end{proposition}

\begin{proof}
The eight intervals are the runs of \(A'\) each extended by one unit; any two successive components are separated by a gap of at least one unit, so they are pairwise disjoint and non-adjacent, and the total length is \(\lvert A'\rvert = 12\).

Put \(S = \Ssev/72 \subset \R/\Z\), a set of twelve points, and let \(\Gamma \subset [0, 1)\) be its set of representatives, so that \(\Lambda = \Gamma + \Z\) and the elements of \(\Gamma\) are pairwise distinct modulo \(\Z\). With \(h = 1\), the fiber of \(\OmegaEight\) at \(x \in [0, 1)\) is \(\{m \in \Z : x + m \in A' + [0, 1)\} = A'\). Since \(\Ssev\) contains the consecutive residues \(5\) and \(6\), the differences of \(S\) have least common denominator \(N = 72\), and Lemma~\ref{lem:reduction}, applied to \(S\) with \(s_0 = 5/72\), produces the set \(S' = \Ssev - 5 \subset \ZN{72}\); a translate of \(\Ssev\) has the same zero set, so the spectra of \(S'\) are those of \(\Ssev\), and by the converse part of Lemma~\ref{lem:reduction} a set of twelve integers whose reduction modulo \(72\) is a spectrum of \(\Ssev\) is a spectrum of \(S\). So \(A'\) is a spectrum of \(S\). Every fiber over \([0, 1)\) is therefore a spectrum of \(S\), and Theorem~\ref{thm:fiber} shows that \(\Lambda\) is a spectrum of \(\OmegaEight\). The number \(1\) is a period of \(\Lambda\). If \(t \in (0, 1)\) were a period of \(\Lambda\), then \(S + t = S\) in \(\R/\Z\), so \(t \in S - S \subset (1/72)\Z/\Z\) and \(72t\) would be a nonzero period of \(\Ssev\) in \(\ZN{72}\), which does not exist. Hence the minimal period is \(1\), and \(k = L/h = 12\).

For the tiling, write \(B' = \{0, 1, 8, 9, 17, 28\}\). The residues of \(B'\) modulo \(6\) are \(0, 1, 2, 3, 5, 4\), a complete residue system, so \(B' \oplus 6\Z = \Z\). Also \(6\Z = \{0, 6, 12, 18, 24, 30\} \oplus \{0, 36\} \oplus 72\Z\), since the twelve multiples of \(6\) in \([0, 72)\) are exactly the sums of one element of the first set and one of the second. Hence
\[
\Z = B' \oplus \{0, 36\} \oplus \{0, 6, \dots, 30\} \oplus 72\Z = A' \oplus T,
\]
with \(T = \{0, 6, 12, 18, 24, 30\} + 72\Z\), all sums direct. Adding \([0, 1)\) gives \(\R = \Z + [0, 1) = \OmegaEight + T\), and the union is disjoint because \(\Z \oplus [0, 1) = \R\) is.
\end{proof}

\begin{figure}[ht]
\centering
\input{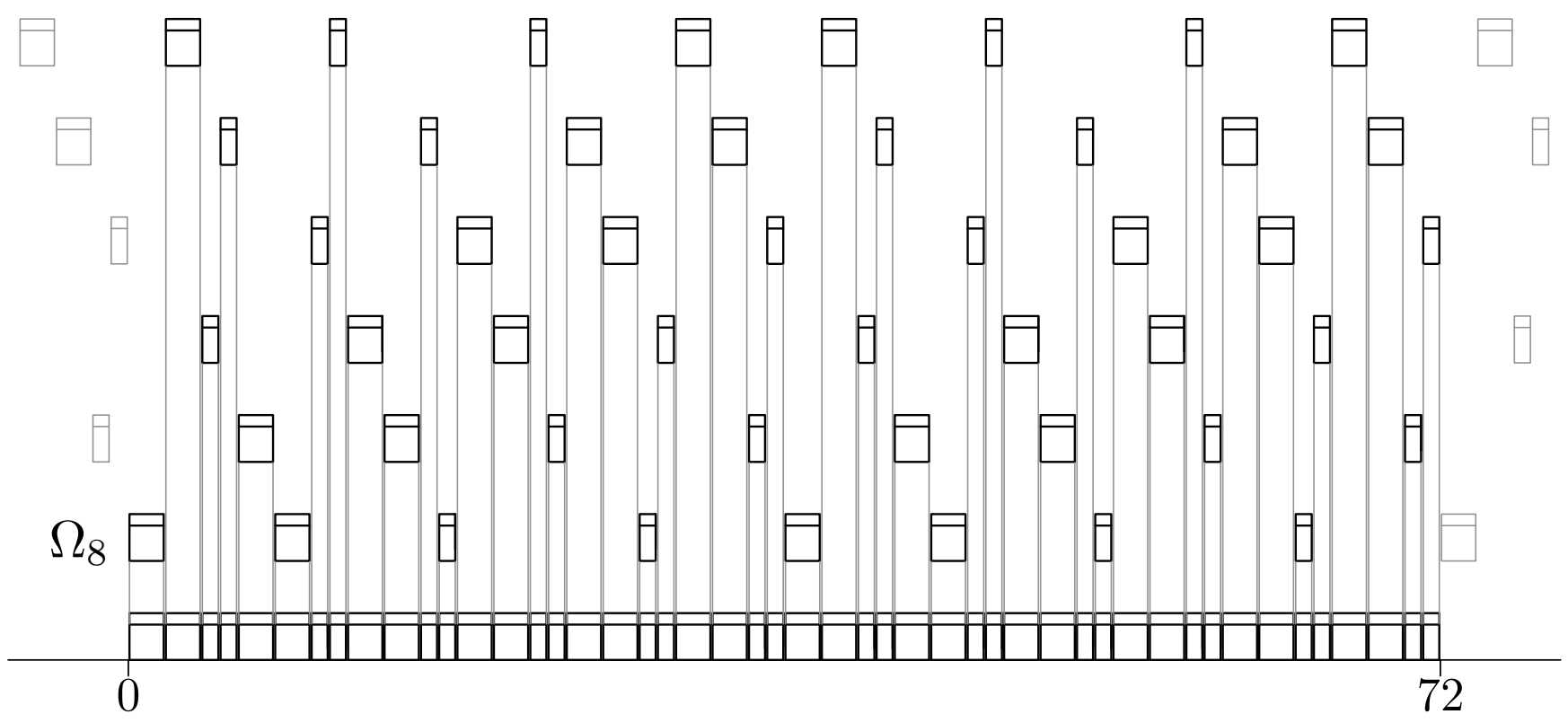}
\caption{The set \(\OmegaEight\), twelve unit intervals forming eight components, on the lowest level, and above it the translates \(\OmegaEight + t\) for \(t = 6, 12, \dots, 30\), one level for each class \(t + 72\Z\), each level also carrying the pieces of the translates by \(t - 72\) and \(t + 72\) that fall in the window. The pieces lying over the period \([0, 72)\), between the marks \(0\) and \(72\) on the base line, fill it exactly; the continuation outside the period is gray.}
\label{fig:omegaeight}
\end{figure}

Figure~\ref{fig:omegaeight} shows \(\OmegaEight\) and a stretch of the tiling, and Proposition~\ref{prop:omegaeight} proves Theorem~D(iv). Changing the integer representatives of \(A'\) modulo \(72\), and gluing copies of \(\Ssev\), gives the same phenomenon for other numbers of intervals.

\begin{proposition}
\label{prop:moreruns}
Let \(J \ge 1\) and \(N = 72J\). Let \(S_J = \bigcup_{i=0}^{J-1} (i + J\Ssev) \subset \ZN{N}\), where \(J\Ssev = \{Js : s \in \Ssev\}\) is computed with representatives of \(\Ssev\) in \([0, 72)\), and let \(A_J = A' + 72\{0, 1, \dots, J - 1\} \subset \Z\). For every \(n\) with \(8J \le n \le 12J\) there is a set \(B \subset \Z\) of \(12J\) integers, congruent to the elements of \(A_J\) modulo \(N\) one by one, which forms exactly \(n\) runs. For every such \(B\), the set \(\Omega = B + [0, 1)\) is a union of \(n\) pairwise disjoint non-adjacent intervals with the spectrum \(\Lambda_J = S_J/N + \Z\) of minimal period \(1\), so that \(h = 1\) and \(k = 12J\), and \(\Omega\) tiles \(\R\) by the translation set \(\{0, 6, 12, 18, 24, 30\} + N\Z\). In particular, for every \(n\) with \(8 \le n \le 11\) or \(n \ge 16\) there is a spectral union of \(n\) intervals that tiles the real line and whose spectrum has more than \(n\) cosets.
\end{proposition}

\begin{proof}
The elements \(i + Js\) of \(S_J\) with \(i \in \{0, \dots, J - 1\}\) and \(s \in \Ssev\) are pairwise distinct, so \(\lvert S_J\rvert = 12J\). In the notation of Definition~\ref{def:gluing} with \(m = J\), the piece of \(S_J\) in the class \(i\) modulo \(J\) is \(i + J\Ssev\), so every piece \(\widetilde S_i\) equals \(\Ssev\). Reduction \(\pi_J \colon \ZN{N} \to \ZN{72}\) maps \(A_J\) onto the reduction of \(A'\), and the reduction of \(A_J\) modulo \(N\) is \(\pi_J^{-1}(A' \bmod 72)\). Since \(A' \bmod 72\) is a spectrum of \(\Ssev\) (Theorem~\ref{thm:cptwelve}), Lemma~\ref{lem:gluing}(a) shows that the reduction of \(A_J\) modulo \(N\) is a spectrum of \(S_J\).

The set \(S_J\) has no period. Suppose \(S_J + t = S_J\) with \(t \in \{0, \dots, N - 1\}\). Translation by \(t\) maps the piece in the class \(i\) onto the piece in the class \(i'\), where \(i' \in \{0, \dots, J - 1\}\) and \(i' \equiv i + t \pmod J\). Writing \(i + t - i' = Jc\), this says \(J(\Ssev + c) = J\Ssev\) in \(\ZN{N}\), that is \(\Ssev + c = \Ssev\) in \(\ZN{72}\), and since \(\Ssev\) has no period, \(c \equiv 0 \pmod{72}\), that is \(t \equiv i' - i \pmod N\). For \(i = 0\) this gives \(t \in \{0, \dots, J - 1\}\). If \(t \ne 0\), then for \(i = J - 1\) one has \(i' = t - 1\), and \(t \equiv t - J \pmod N\), which is false since \(0 < J < N\). So \(t = 0\).

Put \(S = S_J/N \subset \R/\Z\). A period \(p\) of \(\Lambda_J\) satisfies \(S + p = S\), so \(p \in S - S \subset (1/N)\Z/\Z\) and \(Np\) is a period of \(S_J\); hence \(p \in \Z\), the minimal period is \(1\), and \(k = \lvert S\rvert = 12J\).

The integers of \(A_J\) lie in \([0, N - 8]\), and \(A' + 72i\) consists of the runs \(\{0, 1\}\), \(\{8, 9\}\), \(\{17\}\), \(\{28\}\), \(\{36, 37\}\), \(\{44, 45\}\), \(\{53\}\), \(\{64\}\) translated by \(72i\); since \(64 + 1 < 72\), the set \(A_J\) has \(8J\) runs, \(4J\) of them of length two. Choose \(s\) of the runs of length two, \(0 \le s \le 4J\), and replace the larger element of each by itself plus \(N\). The new elements lie in \([N + 1, 2N - 7]\) and differ from each other by at least \(8\), so each forms a run of length one, and each chosen run becomes a run of length one; the resulting set \(B\) has \(8J + s\) runs.

The set \(\Omega = B + [0, 1)\) is the union of the runs of \(B\), each extended by one unit, so it consists of \(n\) pairwise disjoint non-adjacent intervals of total length \(12J\), and with \(h = 1\) its fiber at every \(x \in [0, 1)\) is \(B\). The numbers \(\ex{as}\), \(s \in S\), depend only on \(a\) modulo \(N\), so \(B\) is a spectrum of \(S\), and Theorem~\ref{thm:fiber} shows that \(\Lambda_J\) is a spectrum of \(\Omega\). For the tiling, the proof of Proposition~\ref{prop:omegaeight} gives \(\Z = A' \oplus \{0, 6, \dots, 30\} \oplus 72\Z\), and \(72\Z = 72\{0, \dots, J - 1\} \oplus N\Z\), so \(\Z = A_J \oplus T_J\) with \(T_J = \{0, 6, \dots, 30\} + N\Z\). Since \(T_J + N = T_J\), the same holds with \(B\) in place of \(A_J\), and adding \([0, 1)\) gives \(\Omega \oplus T_J = \R\).

For the last statement, \(J = 1\) gives \(8 \le n \le 11\), and for \(J \ge 2\) the ranges \(8J \le n \le 12J - 1\) cover every \(n \ge 16\), since \(8(J + 1) \le 12J\).
\end{proof}

\begin{corollary}
\label{cor:alllargern}
For every \(n\ge8\), there is a spectral union of \(n\) pairwise disjoint non-adjacent bounded intervals which tiles \(\R\) and has a spectrum with more than \(n\) cosets modulo its least period.
\end{corollary}

\begin{proof}
Proposition~\ref{prop:moreruns} covers \(8\le n\le11\) and \(n\ge16\). For the remaining values, use the common spectrum and the two six-point pieces from the proof of Theorem~\ref{thm:cptwelve}:
\[
\begin{gathered}
B_0=\{0,1,8,9,17,28\},\\
T_0=\{3,5\}+\{0,12,24\},\qquad
T_1=\{2,5,17\}+\{0,18\}.
\end{gathered}
\]
in \(\ZN{36}\). Lemma~\ref{lem:gluing}(a), with \(m=3\), gives the spectral pair
\[
A_{18}=B_0+\{0,36,72\},\qquad
S_{18}=3T_0\cup(1+3T_1)\cup(2+3T_1)
\]
in \(\ZN{108}\). The set \(S_{18}\) has no nonzero period. Indeed, a nontrivial stabilizer would contain the subgroup of order \(2\) or \(3\). The period \(54\) would force \(T_0+18=T_0\), which fails since \(3\in T_0\) and \(21\notin T_0\); the period \(36\) would force \(T_1+12=T_1\), which fails since \(2\in T_1\) and \(14\notin T_1\).

As a set of integers, \(A_{18}\) has twelve runs, six of length two and six of length one. Replacing successively \(1,9,37\) by \(109,117,145\), respectively, gives sets with thirteen, fourteen and fifteen runs: each replacement splits a double run and creates an isolated singleton above the original maximum \(100\). Each resulting set \(B\), including \(A_{18}\) itself, has the same residues modulo \(108\). Thus \(B+[0,1)\) has spectrum \(S_{18}/108+\Z\), whose least period is \(1\), and \(k=18>n\).

Finally, \(B_0\) is a complete residue system modulo \(6\), so
\[
A_{18}\oplus\bigl(\{0,6,12,18,24,30\}+108\Z\bigr)=\Z.
\]
The replacements preserve this tiling, and adding \([0,1)\) gives the required tilings of \(\R\).
\end{proof}

The inequality \(k \le n\) therefore fails for every \(n\ge8\). In particular \((\mathrm{CP}_{n,12J})\) is false for every \(n \ge 8J\); when \(8J\le n<12J\), this is one of the hypotheses required by Theorem~\ref{thm:generalreduction}(b), so that theorem cannot reduce the needed fiber cardinalities to \(k'\le n\). All these sets tile, so they are no counterexamples to the conjecture that spectral unions of intervals tile; they show that a proof covering arbitrary numbers of intervals cannot impose \(k\le n\).

\subsection{The discrete direction}

\begin{proposition}
\label{prop:discrete}
Let \(A \subset \Z\) be a finite set with \(\lvert A\rvert = k\) which is spectral, in the sense that \((A, S)\) is a spectral pair for some \(S \subset \R/\Z\) with \(\lvert S\rvert = k\). If the fiber lemma \(\FL{k}\) holds, then \(A\) tiles \(\Z\) by translations: \(A \oplus T = \Z\) for some \(T \subset \Z\). The tiling can be taken periodic, so \(A\) tiles some cyclic group \(\ZN{M}\).
\end{proposition}

\begin{proof}
Let \(E = A + [0, 1)\). For every \(x \in \R\) and every \(a \in A\) there is exactly one \(m \in \Z\) with \(x + m \in [a, a + 1)\), so \(\sum_{m} 1_E(x + m) = k\) \(\forall x\): the set \(E\) is a \(k\)-fold multitile of \(\R\) by \(\Z\). Its fiber at \(x \in [0, 1)\), with \(h = 1\), is \(A\). Let \(\Gamma \subset [0, 1)\) be the set of representatives of \(S\); its elements are pairwise distinct modulo \(\Z\), and \(h\Gamma = S\). Since the fiber \(A\) is a spectrum of \(S\) at every \(x \in [0, 1)\), Theorem~\ref{thm:fiber} shows that \(\Gamma + \Z\) is a spectrum of \(E\). Thus \(\FL{k}\) gives a set \(T' \subset \R\) with \(E \oplus T' = \R\), that is, \([0, 1) \oplus (A + T') = \R\), where the sum \(A + T'\) is direct because the translates \(a + [0, 1) + t'\) are pairwise disjoint. A tiling of \(\R\) by translates of \([0,1)\) has translation set \(b_0+\Z\): disjoint interiors force consecutive left endpoints to be at least one unit apart, and coverage up to a null set permits no open gap, so each successor lies exactly one unit to the right. Hence \(A + T' = b_0 + \Z\), and with \(T = T' - b_0\) one has \(A \oplus T = \Z\) and \(T \subset \Z\). By Newman's theorem \cite{newman1977}, every tiling of \(\Z\) by a finite set is periodic, so \(T\) may be replaced by a periodic tiling set \(T_0 + M\Z\), and \(A\) tiles \(\ZN{M}\).
\end{proof}

By Theorem~\ref{thm:rational} and Lemma~\ref{lem:reduction}, the set \(A\) of Proposition~\ref{prop:discrete} reduces modulo a suitable \(N\) to a spectral subset of \(\ZN{N}\), and conversely a spectral subset of \(\ZN{N}\) with spectrum \(B\) lifts to a spectral pair \((A, B/N)\) with \(A \subset \Z\). The fiber lemma \(\FL{k}\) therefore contains the direction ``spectral implies tiling'' of Fuglede's conjecture for \(k\)-point subsets of the cyclic groups, with the proviso that the tiling obtained is one of \(\Z\), or of some \(\ZN{M}\), and not necessarily of \(\ZN{N}\) itself. That direction is known, for instance, in \(\ZN{p^a q}\) \cite{malikiosis-kolountzakis2017}, in \(\ZN{p^a q^2}\) \cite{kmsv2020}, in \(\ZN{p^a q^b}\) under conditions on the primes and the exponents \cite{malikiosis2022}, in \(\ZN{pqr}\) \cite{shi2019}, in \(\ZN{p^2 qr}\) \cite{somlai2023} and in \(\ZN{pqrs}\) \cite{kmsv2022}, and open in general; none of these results is used here. At fixed cardinality, the argument requires one complement serving all fibers, as in \(\UC{k}\), rather than separate tilings of the individual fibers. This distinction does not preclude equivalence of the unrestricted conjectures: Dutkay and Lai \cite[Theorem~1.3]{dutkay-lai2014}, together with the rational-spectrum theorem \cite{fu-song2026}, identify the spectral-to-tiling statements on \(\R\), \(\Z\), and all finite cyclic groups. What it makes visible is the shape of the problem for unions of intervals: two statements about finite spectral pairs in cyclic groups, one about a vanishing sum (a primitive spectral \(S' \subset \ZN{N}\) with \(\sum_{s \in S'} \zroot{N}^{s} = 0\) and a spectrum with few runs is periodic) and one about common complements (\(\UC{k}\)); the analytic part bounds the cardinalities needed for each \(n\), without bounding the ambient cyclic-group orders. Under \(\dimQ \Krel \le 4\), Theorem~A uses the consecutive-pair statements up to ten points and common complements up to five points, for every \(n\).

\subsection{Six intervals}

Theorem~C and Theorem~D together settle every spectral union of at most six intervals whose spectrum has at most eleven cosets, and give necessary conditions on the remaining case.

\begin{lemma}
\label{lem:threeconsecutive}
Let \((A, S)\) be a spectral pair with \(A \subset \Z\), \(\lvert A\rvert = \lvert S\rvert = k\), and \(S = s_0 + S'/N\) with \(S' \subset \ZN{N}\), and suppose that \(A\) contains three consecutive integers \(a, a+1, a+2\). Let \(\bar A \subset \ZN{N}\) be the reduction of \(A\) modulo \(N\), which has \(k\) elements and is a spectrum of \(S'\). Suppose that \(\bar A + N/m = \bar A\) for a divisor \(m \ge 2\) of \(N\). Then \(m\) divides \(k\), and if \(k/m \in \{1, 2, \dots, 11\}\), then \(S'\) is periodic in \(\ZN{N}\), so that \(S\) has a nontrivial translation symmetry.
\end{lemma}

\begin{proof}
If \(a \ne a'\) in \(A\) had \(a \equiv a' \pmod N\), then \(\sum_{s \in S} \ex{s(a - a')}\) would equal \(k\,\ex{s_0(a - a')} \ne 0\), against the spectral pair property; so the reduction is injective, and the orthogonality relations \(\sum_{s' \in S'} \zroot{N}^{s'(a - a')} = 0\) for \(a \ne a'\) descend to \(\ZN{N}\): \(\bar A\) is a spectrum of \(S'\) in the sense of Definition~\ref{def:zeroset}.

If \(N/m \le 2\), then \(\bar A\) is a union of cosets of a subgroup of \(\ZN{N}\) of index at most \(2\) and contains two consecutive residues, so \(\bar A = \ZN{N}\), \(k = N\), \(S' = \ZN{N}\), and \(S'\) is periodic. Let \(N/m \ge 3\). Then \(a, a+1, a+2\) lie in three distinct cosets of the subgroup \((N/m)\ZN{N}\), and since \(\bar A\) is a union of cosets of that subgroup we may choose a set \(B \subset \ZN{N}\) of representatives of the cosets contained in \(\bar A\) with \(a, a+1, a+2 \in B\), so that \(\bar A = B + (N/m)\ZN{N}\) with \(\lvert B\rvert = k/m\); in particular \(m \mid k\). Write \(k' = k/m\). The gluing lemma, Lemma~\ref{lem:gluing}, applied to the spectral pair \((\bar A, S')\) in \(\ZN{N}\), shows that the pieces \(S'_j = S' \cap (j + m\ZN{N})\), \(j \in \ZN{m}\), have \(k'\) elements each and that \(B\) is a common spectrum of all of them: \((B - B) \setminus \{0\} \subset \zeroset{S'_j}\) \(\forall j\). Fix \(j\), let \(T_j = \{s'/N : s' \in S'_j\} \subset \R/\Z\), a set of \(k'\) points, and let \(\widetilde B \subset \Z\) be a set of integer representatives of \(B\) containing \(a, a+1, a+2\). For \(b \ne b'\) in \(\widetilde B\) one has \(\sum_{t \in T_j} \ex{t(b - b')} = \sum_{s' \in S'_j} \zroot{N}^{s'(b - b')} = 0\), so \((\widetilde B, T_j)\) is a spectral pair with \(\widetilde B \subset \Z\), \(\lvert \widetilde B\rvert = \lvert T_j\rvert = k'\), with \(T_j - T_j \subset \Q/\Z\), and with the consecutive pair \(a, a+1\) in \(\widetilde B\); moreover \(1, 2 \in \zeroset{S'_j}\), that is, \(\sum_{t \in T_j} \ex{t} = 0\) and \(\sum_{t \in T_j} \ex{2t} = 0\).

We show that there is an integer \(g \ge 2\), depending on \(k'\) only, with \(T_j + 1/g = T_j\) \(\forall j\). The case \(k' = 1\) does not occur, since a single root of unity is not zero. For \(k' = 2\) the two points of \(T_j\) are antipodal, so \(g = 2\). For \(k' = 3\) the three unimodular numbers \(\ex{t}\), \(t \in T_j\), are distinct with zero sum, hence form a coset of the cube roots of unity, so \(g = 3\). For \(k' = 4\), Proposition~\ref{prop:cpsmall} gives \(T_j + 1/2 = T_j\), so \(g = 2\). For \(k' = 5\), Proposition~\ref{prop:cpsmall} gives that \(T_j\) is a coset of \((1/5)\Z/\Z\), so \(g = 5\). For \(k' = 6\), Proposition~\ref{prop:cpsmall} gives \(T_j + 1/3 = T_j\) or \(T_j + 1/2 = T_j\). In the second case \(T_j = \{c_1, c_2, c_3\} + \{0, 1/2\}\), and the identity \(\sum_{t \in T_j} \ex{2t} = 2\sum_{i=1}^{3} \ex{2c_i} = 0\) says that the three unimodular numbers \(\ex{2c_i}\) sum to zero; they are distinct, since a relation \(2\eta + \eta' = 0\) between unimodular numbers is impossible, so they form a coset of the cube roots of unity and \(2c_i \in 2c + \{0, 1/3, 2/3\}\) modulo \(1\) for some \(c\); hence \(c_i \in c + \{0, 1/6, 1/3\}\) modulo \(1/2\), and \(T_j = c + (1/6)\Z/\Z\), which has the period \(1/3\). So \(g = 3\) for \(k' = 6\). For \(k' = 7\), Theorem~\ref{thm:cpseven} gives that \(T_j\) is a coset of \((1/7)\Z/\Z\), so \(g = 7\). For \(k' = 8\), Theorem~\ref{thm:cpeight} gives \(T_j + 1/2 = T_j\), so \(g = 2\). For \(k' = 9\), Theorem~\ref{thm:cpnine} gives \(T_j + 1/3 = T_j\), so \(g = 3\). For \(k'=10\), Theorem~\ref{thm:cpten} gives period \(1/5\) or \(1/2\). In the latter case put \(U=2T_j\). Orthogonality of the rows indexed by \(t\) and \(t+1/2\) gives \(\sum_{b\in\widetilde B}(-1)^b=0\), so \(\widetilde B\) has five even elements. Their halves form a spectrum \(C\) of \(U\): for distinct \(c,c'\in C\),
\[
0=\sum_{t\in T_j}\ex{2t(c-c')}=2\sum_{u\in U}\ex{u(c-c')}.
\]
Also \(\sum_{u\in U}\ex{u}=\tfrac12\sum_{t\in T_j}\ex{2t}=0\). Proposition~\ref{prop:cpsmall} therefore makes \(U\) a translated regular pentagon. Its full inverse image \(T_j\) under doubling is a translated regular decagon and has period \(1/5\). Thus \(g=5\) works for \(k'=10\).

For \(k' = 11\), Theorem~\ref{thm:cpeleven} gives that \(T_j\) is a coset of \((1/11)\Z/\Z\), so \(g = 11\).

Since \(T_j \subset (1/N)\Z/\Z\) and \(1/g \in T_j - T_j\), the integer \(g\) divides \(N\), and \(T_j + 1/g = T_j\) reads \(S'_j + N/g = S'_j\) in \(\ZN{N}\). As this holds for every \(j\) with the same \(g\), \(S' + N/g = S'\): the set \(S'\) is periodic, and \(S + 1/g = S\).
\end{proof}

\begin{proposition}
\label{prop:sixintervals}
Let \(\Omega\) be a spectral union of at most six intervals in the standing reduction, with spectrum \(\Lambda\), minimal period \(h^{-1}\), and \(k\), \(S\) as in Definition~\ref{def:fiber}.
\begin{enumerate}
\item[\rm (i)] If \(k \le 11\), then \(\Omega\) tiles \(\R\) by translations.
\item[\rm (ii)] If \(k \ge 12\), then \(\Omega\) has exactly six components, and every fiber \(A_x\) is a union of exactly six runs, one of which has at least three elements. Moreover \(S - S \subset \Q/\Z\); with \(N\) and \(S' = N(S - s_0) \subset \ZN{N}\) as in Lemma~\ref{lem:reduction}, the set \(S'\) is primitive and has no period in \(\ZN{N}\); and for every \(x\) the reduction \(\bar A_x\) of \(A_x\) modulo \(N\) is a spectrum of \(S'\) with \(k\) elements which has no period \(N/m\) with \(k/m \in \{1, 2, \dots, 11\}\).
\end{enumerate}
\end{proposition}

\begin{proof}
Let \(n \le 6\) be the number of components. (i) If \(k \le 7\), then \(\Omega\) is a \(k\)-fold multitile of \(\R\) by \(h\Z\) with the spectrum \(\Gamma + h^{-1}\Z\) (Lemma~\ref{lem:everyfiber}), and Theorem~C gives the tiling. If \(8 \le k \le 11\), then \(k > n\), so by Lemma~\ref{lem:pigeonhole}(i) every fiber contains two consecutive integers; \((A_x, S)\) is a spectral pair with \(\lvert A_x\rvert = k\), and Theorem~D(i) gives \(S + 1/m = S\) for some \(m \ge 2\), against Theorem~\ref{thm:periodic}. So this case does not occur.

(ii) Let \(k \ge 12\) and fix \(x\). By Lemma~\ref{lem:everyfiber}, \((A_x, S)\) is a spectral pair and \(S\) has no nontrivial translation symmetry (Theorem~\ref{thm:periodic}). Let \(d\) be the run vector of \(A_x\), with \(m \le n \le 6\) parts and \(\lvert d\rvert = k\). By Lemma~\ref{lem:runvector}, \(k = \kzero(d)\). If \(m \le 5\), then \(\dimQ K_d \le m - 1 \le 4\) for the space \(K_d\) of rational relations among the entries of \(d\), and Theorem~\ref{thm:kzeroten} gives \(\kzero(d) \le 10 < k\), a contradiction. So \(m = 6\), hence \(n = 6\), and every fiber is a union of exactly six runs. If every run had at most two elements, then \(k \le 12\), so \(k = 12\) and \(d = (2, 2, 2, 2, 2, 2)\); but for this \(d\) every subset sum is even, the levels \(2, 4, 6, 8, 10\) are effective (the indicators of the subsets of a given size \(\beta\) sum to a positive multiple of \(\one\)), and their ratios \(\beta/6\) have denominators \(6, 3, 2, 3, 6\), so \(\kzero(d) = 6 \ne 12\), against Lemma~\ref{lem:runvector}. Hence some run has at least three elements, and in particular \(1 \in A_x - A_x\). Theorem~\ref{thm:rational} gives \(S - S \subset \Q/\Z\), and Lemma~\ref{lem:reduction} gives \(N\), \(S'\) primitive, and the injectivity of the reduction of \(A_x\) modulo \(N\) onto a spectrum \(\bar A_x\) of \(S'\). A nonzero period \(t\) of \(S'\) in \(\ZN{N}\) would give \(S + t/N = S\), a translation symmetry of \(S\); so \(S'\) has no period. Finally, if \(\bar A_x\) had a period \(N/m\) with \(k/m \in \{1, \dots, 11\}\), Lemma~\ref{lem:threeconsecutive}, applied with the three consecutive integers of the long run, would make \(S'\) periodic, which was just excluded.
\end{proof}

To prove that every spectral union of at most six intervals tiles, Proposition~\ref{prop:sixintervals} shows that it would suffice to establish the following assertion about finite spectral pairs: there is no spectral pair \((\bar A, S')\) in a cyclic group \(\ZN{N}\) with \(S'\) primitive and without period, \(\lvert \bar A\rvert = \lvert S'\rvert = k \ge 12\), \(\bar A\) the reduction of a set of \(k\) integers forming exactly six runs, one of them of length at least three, and \(\bar A\) without any period \(N/m\) with \(k/m \in \{1, \dots, 11\}\). By Theorem~\ref{thm:generalreduction}(a) only finitely many values of \(k\), the divisors of \(\kzero(\ell)\) up to \(K(6)\), need to be considered for a given length vector, and by Lemma~\ref{lem:runvector} the value \(k\) is determined by the run vector \(d\) of any fiber through \(k = \kzero(d)\), so that only finitely many integer run vectors \(d\) occur for each \(k\). The gaps between runs and the ambient cyclic-group orders remain unrestricted. The specific pair of Theorem~\ref{thm:cptwelve} has a period with quotient size six and is excluded by Proposition~\ref{prop:sixintervals}(ii). The broader gluing theorem does not by itself exclude the remaining quotient sizes. Nor can the periodicity dichotomy discussed in Remark~\ref{rem:dichotomy} settle the assertion: it is false in general, and the present reduction excludes only the specified periods of \(\bar A\), not every nonzero period. The six-run configurations remain outside the results proved here.

\subsection{Closing remarks}

\begin{remark}
\label{rem:limits}
Everything analytic in the proof of Theorem~A is free of the number of intervals: the periodicity of the spectrum, the fibers as spectra of one \(k\)-point set, Shi's boundary unitary and the cofactor identity in level form, the verdict rule \(k \mid \kzero(\ell)\) and the bound \(\kzero(\ell) \mid \Kn\), the forced integer parts in configurations without consecutive pairs, and the fiber-side mechanisms of Section~\ref{sec:mechanisms}. What is special to few intervals is the arithmetic of the lengths: for five lengths the effective ratios and the values of \(\kzero(\ell)\) are classified in full, and when \(\dimQ \Krel \le 4\), specialization reduces any number of lengths to at most five positive lengths while preserving every original subset-sum coincidence. The rank bound then limits a fully mixing configuration to at most eight intervals, and so restricts the class graphs that Sections~\ref{sec:classgraph} and~\ref{sec:mechanisms} have to treat. For unrestricted six-tuples of lengths, \(\kzero(\ell)\) is no longer bounded by ten, and Proposition~\ref{prop:sixintervals}(ii) gives necessary conditions on the configurations beyond the small-cardinality argument.

Two further obstructions arise in the general case. The consecutive-pair statement \(\CP{k}\), which excludes every \(k > n\) for \(n \le 5\), fails at \(k = 12\) (Theorem~D), and for every \(n\ge8\), the number \(k\) of cosets can exceed the number of intervals in a set that nevertheless tiles (Corollary~\ref{cor:alllargern}). So for these numbers of intervals the problem contains the discrete direction of Fuglede's conjecture for fibers with more points than intervals, together with the common complement statement for them. The bound on \(k\) alone does not provide these common complements (Proposition~\ref{prop:discrete} and the discussion following it). The fiber lemmas proved here cover \(k \le 7\) (Theorem~C), where the common complement exists because the spectra of a six-point set have the period \(N/2\) or \(N/3\) and those of a seven-point set are cosets.

The remaining questions within this argument concern common complements beyond the cardinalities treated here and the six-run configurations of Proposition~\ref{prop:sixintervals}(ii). Two tempting extensions are false: the general periodicity dichotomy (Remark~\ref{rem:dichotomy}) and the bound \(k\le\max_c m_c\) beyond the small cardinalities of Section~\ref{sec:mechanisms} (Remark~\ref{rem:openmechanisms}). The examples in those remarks indicate why a further reduction must retain more information than periodicity or the largest class size alone.
\end{remark}

\appendix
\section{Case analyses}
\label{app:cases}

This appendix contains the enumerations behind the classification of the configurations of size at most five in Section~\ref{sec:fivelengths}, and the list of the class graphs that a fully mixing configuration of eight intervals can have. The class-graph list is not needed for Theorem~\ref{thm:A}: the exclusion of \(k \in \{6, 7\}\) at eight intervals in Lemma~\ref{lem:descent} uses only that the class graph is connected, loopless and balanced.

\subsection{The configurations of size at most five}

\begin{proof}[Proof of Proposition~\ref{prop:smallconfigs}]
By Lemma~\ref{lem:configbasic}(d) it suffices to determine, at each size, the configurations of class (I), which are the extensions of the configurations of the previous size by Lemma~\ref{lem:configbasic}(b), and those of class (III); every configuration of class (II) is the complement of a configuration of class (I). For \(m = 1\) the only invertible \(0\)-\(1\) matrix is \([1]\). For \(m = 2\) and \(m = 3\), class (III) requires column weights in \([2, m-2]\), which is empty; the extensions of \([1]\), with \(w' = (1)\), need \(w'(R) < 1\), hence \(R = \emptyset\), giving \(I_2\), whose complement is \(I_2\) again; the extensions of \(I_2\), with \(w' = (1,1)\), need \(R = \emptyset\), giving \(I_3\), with complement \(J - I_3\). For \(m = 4\) the extensions of \(I_3\) need \(R = \emptyset\) and give \(I_4\); the extensions of \(J - I_3\), with \(w' = (\frac12,\frac12,\frac12)\) on the elements \(2,3,4\) and rows \(34, 24, 23\), need \(\lvert R\rvert \le 1\): \(R = \emptyset\) gives \(A\), and \(R\) a single element gives \(D\), the three choices being isomorphic under the symmetric group of \(J - I_3\). The complements are \(J - I_4\), \(A' = J - A\), and \(J - D\), whose rows \(34, 12, 13, 14\) are those of \(D\) with the roles of the elements \(1\) and \(2\) exchanged, so \(J - D \cong D\). In class (III) with \(m = 4\) every column weight and every row size equals \(2\), so the four elements are the edges of a simple graph on the four rows in which every vertex has degree two, that is, a \(4\)-cycle; the alternating sum of the four rows along the cycle vanishes, so the matrix is singular and class (III) is empty. The weight vectors are read off from \(Mw = \one\); the ones for \(J - I_3\), \(J - I_4\), \(A'\) also follow from Lemma~\ref{lem:configbasic}(c). The rest is read off from the list.
\end{proof}

\begin{proof}[Proof of Proposition~\ref{prop:fewmembers}]
Let \(\gamma\) be the coefficients of Lemma~\ref{lem:minimalfamily}. After Lemma~\ref{lem:contraction} we may assume that the elements have pairwise distinct types, and we have to show that then \(\lvert E\rvert = t\). Write \(\overline{m} = \lvert E\rvert\), so that \(t \le \overline{m} \le 5\) by Lemma~\ref{lem:minimalfamily}. Identify the types with subsets of the row index set \([t]\), where \(i\) stands for the row \(T_i\), and write \(\gamma_i = \gamma_{T_i}\). The \(j\)-th coordinate of \(\one = \sum_i\gamma_i 1_{T_i}\) says that \(\gamma(S) := \sum_{i \in S}\gamma_i = 1\) for every type \(S\) that occurs. Every \(\gamma_i\) is nonzero. We use two consequences. (N) If types \(S \subsetneq S'\) both occur, then \(\gamma(S' \setminus S) = 0\), so \(S' \setminus S\) has at least two elements. (R) The rows form an antichain by positivity (nested rows of equal weight are equal, and the rows are distinct); so for every \(i\) some type contains \(i\), and for \(i \ne k\) some type contains \(i\) and not \(k\).

If \(t = 1\), then \(T_1 = E\) and \(\alpha = w(E)\), which is not a level. So \(t \ge 2\).

Let \(t = 2\). The possible types are \(1\), \(2\), \(12\). If \(12\) occurs, then \(1\) cannot occur by (N), nor \(2\); so \(12\) is the only type, and \(T_1 = T_2 = E\), against distinctness. So the types are \(1\) and \(2\), \(\overline{m} = 2 = t\).

Let \(t = 3\), and let \(s\) be the number of singleton types. If \(s = 3\), then \(\gamma = (1,1,1)\), every other subset of \([3]\) has \(\gamma\)-sum at least \(2\), so \(\overline{m} = 3\). If \(s = 2\), say the types \(1\) and \(2\) occur, then \(\gamma_1 = \gamma_2 = 1\); the type \(12\) has sum \(2\), the types \(13\) and \(23\) would force \(\gamma_3 = 0\), so the only possible further type is \(123\), with \(\gamma_3 = -1\); then \(\overline{m} = 3\) forces it to occur, and the row \(T_3\), which consists of the element of type \(123\), is contained in \(T_1\), against (R). If \(s = 1\), say the type \(1\), the types \(12\) and \(13\) force a zero coefficient, and \(23\) and \(123\) give \(\gamma_2 + \gamma_3 = 1\) and \(\gamma_2 + \gamma_3 = 0\), so at most one of them occurs and \(\overline{m} \le 2 < t\), impossible. If \(s = 0\) and \(123\) occurs, a pair type would give a zero by (N), so \(\overline{m} = 1\), impossible; otherwise all types are pairs, \(\overline{m} \ge 3\) forces all three pairs to occur, and \(\overline{m} = 3 = t\).

Let \(t = 4\). Then \(\overline{m} \in \{4,5\}\), and we show that \(\overline{m} = 5\) is impossible. Suppose there are five distinct nonempty types, subsets of \([4]\) with \(\gamma\)-sum \(1\), and let \(s\) be the number of singletons among them.

If \(s = 4\), then \(\gamma = \one\) and a fifth type has sum at least \(2\).

If \(s = 3\), with the types \(1, 2, 3\): a further type \(S\) not containing \(4\) is a subset of \(123\) of size at least two, with sum at least \(2\); so \(4 \in S\) and \(\gamma_4 = 1 - \lvert S \cap 123\rvert\). Since \(\gamma_4 \ne 0\) and \(S \ne 4\), \(\lvert S \cap 123\rvert \in \{2,3\}\), and the two further types \(S, S'\) have the same intersection size, as both determine \(\gamma_4\). Size \(3\) gives \(S = S' = 1234\). Size \(2\) gives \(S = ij4\), \(S' = ik4\) with a common index \(i\), and then the row \(T_4\), consisting of the two elements of types \(S\) and \(S'\), is contained in \(T_i\), against (R).

If \(s = 2\), with the types \(1, 2\): by (R) some type \(S\) contains \(3\) and not \(4\); it is not a singleton, so \(\lvert S \cap 12\rvert \ge 1\), and \(\gamma_3 = 1 - \lvert S \cap 12\rvert \ne 0\) forces \(S = 123\), \(\gamma_3 = -1\). Likewise \(124\) occurs with \(\gamma_4 = -1\). A fifth type \(S''\) of size at least two meets \(34\) in \(\emptyset\), giving \(S'' = 12\) of sum \(2\), or in \(3\) or \(4\), giving \(123\) or \(124\) again, or contains \(34\), with sum \(\lvert S'' \cap 12\rvert - 2 \le 0\). None exists.

If \(s = 1\), with the type \(1\): a further type \(S\) satisfies \(\gamma(S \setminus 1) = 1 - [1 \in S]\); if \(1 \in S\) then \(\lvert S \setminus 1\rvert \ge 2\) by (N), and if \(1 \notin S\) then \(S\) is not a singleton. So \(S \setminus 1 \in \{23, 24, 34, 234\}\), and for each of these four sets \(Q\) at most one of the two versions \(Q\), \(Q \cup 1\) occurs, since they would give \(\gamma(Q) = 1\) and \(\gamma(Q) = 0\). Five types need all four, so put \(x_Q = \gamma(Q) \in \{0,1\}\). Then \(\gamma_2 + \gamma_3 = x_{23}\), \(\gamma_2 + \gamma_4 = x_{24}\), \(\gamma_3 + \gamma_4 = x_{34}\) and \(\gamma_2 + \gamma_3 + \gamma_4 = x_{234}\), so \(x_{23} + x_{24} + x_{34} = 2x_{234} \in \{0,2\}\). If all four values are \(0\), then \(\gamma_2 = \gamma_3 = \gamma_4 = 0\). Otherwise \(x_{234} = 1\) and exactly one pair value is \(0\), say \(x_{34} = 0\); then \(\gamma_2 = 1\) and \(\gamma_3 = x_{23} - \gamma_2 = 0\).

If \(s = 0\), all five types have size at least two. If \(1234\) occurs, a triple type \(ijk\) would give \(\gamma_l = 0\) by (N), so the other four types are pairs, and a pair \(ij\) gives \(\gamma_k + \gamma_l = 0\) by (N), so no two complementary pairs occur; but four of the six pairs of \([4]\) contain two complementary ones. If a triple occurs, say \(123\), and \(1234\) does not: a pair inside \(123\) gives a zero by (N), so the other four types lie among \(14, 24, 34\) (each giving \(\gamma_i + \gamma_4 = 1\)) and \(124, 134, 234\) (giving \(\gamma_4 = \gamma_3\), \(\gamma_4 = \gamma_2\), \(\gamma_4 = \gamma_1\) respectively, by comparison with \(\gamma_1 + \gamma_2 + \gamma_3 = 1\)). Three pairs and a triple give \(\gamma_1 = \gamma_2 = \gamma_3 = 1 - \gamma_4\), hence \(\gamma_4 = 2/3\), while the triple requires \(\gamma_4 = 1/3\). Two triples, say \(124\) and \(134\), give \(\gamma_2 = \gamma_3 = \gamma_4 =: g\) and \(\gamma_1 = 1 - 2g\); a pair \(24\) or \(34\) gives \(g = 1/2\) and \(\gamma_1 = 0\), and the pair \(14\) gives \(1 - g = 1\), \(g = 0\); two pairs include one of these. Three triples give \(\gamma_1 = \gamma_2 = \gamma_3 = \gamma_4 = 1/3\) and the pair gives \(2/3 = 1\). Four pairs are not available. Finally, if all five types are pairs, they are the edges of \(K_4\) minus one edge, say minus \(34\): the equations \(\gamma_1 + \gamma_2 = \gamma_1 + \gamma_3 = \gamma_1 + \gamma_4 = \gamma_2 + \gamma_3 = \gamma_2 + \gamma_4 = 1\) force \(\gamma = \frac12\one\), and the rows, with the elements named by their types, are \(T_1 = \{12, 13, 14\}\), \(T_2 = \{12, 23, 24\}\), \(T_3 = \{13, 23\}\), \(T_4 = \{14, 24\}\). The equilevel conditions \(w(T_3) = w(T_4) = \alpha\) and \(w(T_1) = w(T_2) = \alpha\) give, adding the last two, \(2w_{12} + (w_{13} + w_{23}) + (w_{14} + w_{24}) = 2\alpha\), that is \(w_{12} = 0\). Every case is contradictory, so \(\overline{m} = 4\).

Let \(t = 5\). Then \(5 = t \le \overline{m} \le 5\).

In every case the contracted family is square with linearly independent indicators, so its incidence matrix \(M\) is invertible, and the classes \(C_1,\dots,C_t\) are the sets of elements of a common type. The equilevel condition reads \(M(w(C_1),\dots,w(C_t))^T = \alpha\one\), so \((w(C_i))_i = \alpha M^{-1}\one > 0\) and \(M\) is a \(t\)-configuration. Since \(M^T\gamma = \one\), \(\sum_i\gamma_i = \one^TM^{-T}\one = L(M)\), and \(r_\alpha = 1/\sum_i\gamma_i = r(M)\) by Proposition~\ref{prop:pinned}(ii).

For the converse, let \(\varphi\colon \Q^{[t]} \to \Q^{E'}\) repeat the coordinate \(i\) on the parts of \(w_i\), and let \(T_i'\) be the set of parts of the weights \(w_j\), \(j \in T_i\). Then \(\varphi(1_{T_i}) = 1_{T_i'}\), \(\varphi(\one) = 1_{E'}\), \(v(T_i') = w(T_i) = 1\) and \(L(v) = L(M) > 1\). Proposition~\ref{prop:pinned}(iii), applied to the effective level \(1\) of \(w\) with its slice containing the rows of \(M\), shows that \(1\) is an effective level of \(v\) with ratio \(1/L(M)\).
\end{proof}

\begin{proof}[Proof of Theorem~\ref{thm:fiveconfigs}]
By Lemma~\ref{lem:configbasic}(d) every \(5\)-configuration is of class (I), (II) or (III), every configuration of class (II) is the complement of a configuration of class (I), and the complement of a class (III) configuration is of class (III).

Class (I). By Lemma~\ref{lem:configbasic}(b) and Proposition~\ref{prop:smallconfigs}, the configurations of class (I) are the extensions of \(I_4\), \(J - I_4\), \(A\), \(A'\), \(D\) by sets \(R\) of their columns with \(w'(R) < 1\), and it suffices to list every such pair \((N, R)\) up to the symmetries of \(N\). For \(I_4\), \(w' = \one\) allows only \(R = \emptyset\). For \(J - I_4\), \(w' = \frac13\one\) allows \(\lvert R\rvert \le 2\), and the symmetric group of \(J - I_4\) makes the choice of \(R\) depend on \(\lvert R\rvert\) alone. For \(A\), with weights \(1\) on the element \(2\) and \(\frac12\) on \(3,4,5\), the admissible \(R\) are \(\emptyset\) and the singletons of light elements, which are permuted by the symmetries of \(A\); the heavy element alone and two light elements give \(w'(R) = 1\). For \(A'\), whose rows \(345, 25, 24, 23\) give \(w'_3 = w'_4 = w'_5 = \frac13\) and \(w'_2 = \frac23\), the admissible \(R\) are \(\emptyset\), one light element, two light elements, and the heavy element alone, all other subsets having \(w'(R) \ge 1\). For \(D\), with all weights \(\frac12\), the admissible \(R\) are \(\emptyset\) and the singletons. The weights, masses and ratios in the table follow from Lemma~\ref{lem:configbasic}(b): the new element has weight \(1 - w'(R)\), and the mass is \(L(N) + 1 - w'(R)\). The configurations of class (II) are therefore included among the complements of this list, by Lemma~\ref{lem:configbasic}(c).

Class (III). Every column weight and every row size lies in \([2,3]\). Let \(q\) be the number of columns of weight \(3\). The sum of the row sizes equals the sum of the column weights, \(10 + q\), so exactly \(q\) rows have size \(3\). Complementation replaces \(q\) by \(5 - q\), so it suffices to treat \(q \in \{0,1,2\}\). The columns of an invertible matrix are distinct, so the elements of weight \(2\) are the edges of a simple graph \(G\) on the vertex set \([5]\) of rows, the elements of weight \(3\) are triples of rows, and \(T_i\) consists of the edges of \(G\) at \(i\) and the triples containing \(i\). We write the weight of an element as the weight of the corresponding edge or triple.

Let \(q = 0\). Then \(G\) has five edges and every vertex has degree \(\lvert T_i\rvert = 2\), so \(G\) is a disjoint union of cycles of length at least three on five vertices, hence a \(5\)-cycle with edges \(e_1,\dots,e_5\) in cyclic order. The equations \(w_{e_{i-1}} + w_{e_i} = 1\) around an odd cycle have the unique solution \(w = \frac12\one\), so the matrix is invertible with positive weights: this is \(C_1\), with \(L = 5/2\) and \(r = 2/5\).

Let \(q = 1\), with the triple \(P = 123\) of weight \(x\) and four edges. Then \(\lvert T_i\rvert = \deg(i) + [i \in P]\), so \(\deg(i) \in \{1,2\}\) for \(i \le 3\) and \(\deg(i) \in \{2,3\}\) for \(i \in \{4,5\}\), with \(\sum_i\deg(i) = 8\), one more than the sum of the minima; so exactly one vertex has degree one above its minimum. A vertex \(i \le 3\) of degree one, attached to a vertex \(v\), has \(T_i = \{iv, P\}\); if \(v \le 3\), then \(T_v\) contains \(iv\) and \(P\), against the antichain property. So every vertex \(i \le 3\) of degree one is attached to \(4\) or \(5\). If \(\deg(4) = 3\), the other degrees are minimal: \(1,2,3\) have degree one and \(5\) has degree two. The three vertices \(1,2,3\) account for three edges; the fourth edge is \(45\), since \(5\) needs a second neighbor and the vertices \(1,2,3\) are exhausted; so, up to symmetry, the edges are \(14, 24, 35, 45\). The rows give \(w_{14} = w_{24} = w_{35} = 1 - x\) and \(w_{45} = 1 - w_{35} = x\), and the row \(T_4 = \{14, 24, 45\}\) gives \(2 - x = 1\), so \(x = 1\) and \(w_{14} = 0\). The case \(\deg(5) = 3\) is symmetric. Otherwise exactly one of \(1,2,3\) has degree two, say \(1\), the vertices \(2\) and \(3\) have degree one and are attached to \(4\) or \(5\), and \(\deg(4) = \deg(5) = 2\). The two edges not at \(2\) or \(3\) lie in \(\{1,4,5\}\) and form a path of length two whose middle vertex receives no further edge; the middle vertex is \(1\), since a middle vertex \(4\) or \(5\) would leave \(1\) with a single edge. So the edges are \(14, 15\), and \(2, 3\) are attached one to \(4\) and one to \(5\), say \(24, 35\). This is \(C_2\), with \(a = 14\), \(b = 15\), \(c = 24\), \(d = 35\), \(e = P\). The rows \(T_2 = \{c, e\}\) and \(T_3 = \{d, e\}\) give \(w_c = w_d = 1 - w_e\), the rows \(T_4 = \{a, c\}\) and \(T_5 = \{b, d\}\) give \(w_a = w_b = w_e\), and \(T_1 = \{a, b, e\}\) gives \(3w_e = 1\). The solution is unique, so the matrix is invertible, with \(w = (1,1,2,2,1)/3\), \(L = 7/3\), \(r = 3/7\).

Let \(q = 2\), with triples \(P, Q\) of weights \(x, y\) and three edges, so that \(\sum_i\deg(i) = 6\) and \(\lvert T_i\rvert = \deg(i) + c_i\) with \(c_i = [i \in P] + [i \in Q]\). Two distinct triples in a five-element set share one or two elements.

If \(P = 123\) and \(Q = 345\), then \(c = (1,1,2,1,1)\), so \(\deg(3) \in \{0,1\}\) and \(\deg(i) \in \{1,2\}\) for \(i \ne 3\). Suppose \(\deg(3) = 1\). Then exactly one other vertex has degree two, and \(G\) is a path of length two together with a disjoint edge. By the symmetries exchanging \(1\) with \(2\), \(4\) with \(5\), and \(P\) with \(Q\), we may take \(1\) as the center. If \(2\) is a neighbor of \(1\), then \(T_2 = \{12, P\} \subset T_1\). If the neighbors are \(3\) and \(4\), the disjoint edge is \(25\), and the rows \(T_2 = \{25, P\}\), \(T_5 = \{25, Q\}\) give \(x = y\); \(T_4 = \{14, Q\}\) gives \(w_{14} = 1 - x\), \(T_3 = \{13, P, Q\}\) gives \(w_{13} = 1 - 2x\), and \(T_1 = \{13, 14, P\}\) gives \(2 - 2x = 1\), so \(w_{13} = 0\). If the neighbors are \(4\) and \(5\), the disjoint edge is \(23\) and \(T_2 = \{23, P\} \subset T_3\). Suppose \(\deg(3) = 0\). Then \(T_3 = \{P, Q\}\), so \(x + y = 1\), and the three edges form a path on \(\{1,2,4,5\}\) with two interior vertices of degree two. If \(12\) is an edge, both \(1\) and \(2\) are interior, since an end vertex \(i \in \{1,2\}\) would have \(T_i = \{12, P\}\) contained in the row of the other; the path is \(4\,1\,2\,5\), and the rows \(T_4, T_5, T_1\) give \(w_{14} = w_{25} = x\) and \(w_{12} = 1 - 2x\), while \(T_2\) is then automatic: the weights are not determined, so the matrix is singular (indeed \(1_{T_1} + 1_{T_5} = \one = 1_{T_2} + 1_{T_4}\)). The edge \(45\) is excluded in the same way. Otherwise the path alternates between \(\{1,2\}\) and \(\{4,5\}\), say \(1\,4\,2\,5\): the rows \(T_1 = \{14, P\}\), \(T_5 = \{25, Q\}\) give \(w_{14} = y\), \(w_{25} = x\), the row \(T_4 = \{14, 24, Q\}\) gives \(w_{24} = 1 - 2y\), and \(T_2 = \{24, 25, P\}\) gives \(1 - 2y + 2x = 1\), so \(x = y = \frac12\) and \(w_{24} = 0\).

If \(P = 123\) and \(Q = 124\), then \(c = (2,2,1,1,0)\), so \(\deg(1), \deg(2) \in \{0,1\}\), \(\deg(3), \deg(4) \in \{1,2\}\), \(\deg(5) \in \{2,3\}\). If \(\deg(1) = 0\) then \(T_1 = \{P, Q\} \subset T_2\); so \(\deg(1) = \deg(2) = 1\), and the remaining degrees sum to \(4\), which forces \((\deg 3, \deg 4, \deg 5) = (1,1,2)\). An edge \(12\) gives \(T_1 = T_2\); an edge \(1x\) with \(x \in \{3,4\}\) gives \(T_x = \{1x, P \text{ or } Q\} \subset T_1\), and likewise for \(2\). So \(1\) and \(2\) are both joined to \(5\), and the third edge is \(34\). This is \(C_3\), with \(a = 15\), \(b = 25\), \(c = 34\), \(d = P\), \(e = Q\), so that the weights of \(P, Q\) are \(w_d = x\), \(w_e = y\). The rows \(T_1, T_2\) give \(w_a = w_b = 1 - w_d - w_e\), the row \(T_5 = \{a, b\}\) gives \(w_d + w_e = \frac12\), and the rows \(T_3 = \{c, d\}\), \(T_4 = \{c, e\}\) give \(w_c = 1 - w_d = 1 - w_e\), so \(w_d = w_e = \frac14\), \(w_c = \frac34\). The solution is unique, \(w = (2,2,3,1,1)/4\), \(L = 9/4\), \(r = 4/9\).

The cases \(q = 3, 4, 5\) are the complements of \(C_3\), \(C_2\), \(C_1\), with ratios \(5/9\), \(4/7\), \(3/5\) and the same weight vectors up to a scalar.

The ratio set is the union of the class (I) ratios \(1/5\), \(3/7\), \(1/2\), \(3/5\), \(2/7\), \(1/3\), \(3/8\), \(2/5\), of their reflections \(r \mapsto 1 - r\), and of the class (III) ratios \(2/5\), \(3/7\), \(4/9\), \(5/9\), \(4/7\), \(3/5\); this union is \(R_5\). The configurations at the denominators \(7, 8, 9\) and at \(1/5\), \(4/5\) are read off from the three lists, a class (II) configuration with ratio \(r\) being the complement of a class (I) configuration with ratio \(1 - r\): \(4/7\) comes from \(E_2\), \(E_8\) and from \(C_2\); \(5/7\) from \(E_5\); \(5/8\) from \(E_7\); \(5/9\) from \(C_3\); \(4/5\) from \(E_1\). The ratio \(1/3\) occurs in class (I) for \(E_6\) and \(E_{11}\) only, in class (II) it would need a class (I) configuration of ratio \(2/3\), of which there is none, and it does not occur in class (III); the weight vectors of \(E_6\) and \(E_{11}\) are both proportional to \((1,2,1,1,1)\). The ratios \(1/4\), \(3/4\) do not occur in \(R_5\). The ratio \(2/5\) occurs for \(E_{12}\), for the complement of \(E_4\), and for \(C_1\), all with equal weights.
\end{proof}

\subsection{The class graphs at eight intervals}

\begin{lemma}
\label{lem:eightgraphs}
Let \(\Omega\) be a fully mixing spectral union of eight intervals with \(\dimQ\Krel\le 4\), so that the pattern is \((2,2,2,2)\) by Corollary~\ref{cor:fullymixing}, and suppose \(3\le k\le 10\). Then \(\clsgraph\) is a connected loopless directed multigraph on four classes with in-degree and out-degree two at every class. Up to a relabeling of the classes, there are exactly five such multigraphs:
\begin{enumerate}
\item[\rm (S2)] the bidirected complete bipartite graph on the parts \(\{1,3\}\) and \(\{2,4\}\), that is the four \(2\)-cycles \(1\leftrightarrow 2\), \(2\leftrightarrow 3\), \(3\leftrightarrow 4\), \(4\leftrightarrow 1\);
\item[\rm (S1)] the directed \(4\)-cycle \(2\to 1\to 4\to 3\to 2\) together with the two diagonal \(2\)-cycles \(1\leftrightarrow 3\) and \(2\leftrightarrow 4\);
\item[\rm (B1)] the arc \(1\to 2\) doubled, the \(2\)-cycle \(3\leftrightarrow 4\), and the single arcs \(2\to 3\), \(2\to 4\), \(3\to 1\), \(4\to 1\);
\item[\rm (C3)] the arcs \(1\to 2\) and \(3\to 4\) doubled, and the single arcs \(2\to 1\), \(2\to 3\), \(4\to 1\), \(4\to 3\);
\item[\rm (D4)] a directed \(4\)-cycle with every arc doubled.
\end{enumerate}
\end{lemma}

\begin{proof}
Every class has size two, so every loop would sit at a class of size two, and Theorem~\ref{thm:noloop} excludes loops; the graph is connected by Lemma~\ref{lem:maximalrank} and balanced with all degrees two. Since the degrees are two, every arc has multiplicity at most two; we classify by the number of doubled arcs. Call a class out-saturated when both edges starting there are known, and in-saturated when both edges ending there are known.

No doubled arc. The graph is a simple digraph with all in- and out-degrees two, and its complement in the complete loopless digraph on four vertices has all in- and out-degrees one, so the complement is a permutation without fixed points: a \(4\)-cycle or a product of two disjoint transpositions. If the complement is the \(4\)-cycle \(1\to 2\to 3\to 4\to 1\), the graph has the arcs \(1\to 3\), \(1\to 4\), \(2\to 1\), \(2\to 4\), \(3\to 1\), \(3\to 2\), \(4\to 2\), \(4\to 3\), which are the directed cycle \(2\to 1\to 4\to 3\to 2\) and the \(2\)-cycles \(1\leftrightarrow 3\), \(2\leftrightarrow 4\): this is (S1). If the complement consists of \(1\leftrightarrow 2\) and \(3\leftrightarrow 4\), the graph has every arc between \(\{1,2\}\) and \(\{3,4\}\): this is (S2) after relabeling.

Exactly one doubled arc, say \(1\to 2\). Then \(1\) is out-saturated and \(2\) is in-saturated. The two edges starting at \(3\) cannot end at \(2\) or at \(3\), and they are single, so they are \(3\to 1\) and \(3\to 4\); likewise \(4\to 1\) and \(4\to 3\). Now \(1\) is in-saturated, \(3\) and \(4\) are out-saturated, and the two edges starting at \(2\) must end at \(3\) and \(4\), one each, as \(3\) and \(4\) each still lack one incoming edge. This is (B1).

Exactly two doubled arcs. If they form a \(2\)-cycle \(1\rightrightarrows 2\rightrightarrows 1\), then \(1\) and \(2\) are saturated in both directions, and the classes \(3\) and \(4\) can only be joined to each other by single arcs, which gives them degree at most one. Two doubled arcs with a common tail or a common head give a degree four. If they form a path \(1\rightrightarrows 2\rightrightarrows 3\), then \(2\) is saturated, \(1\) is out-saturated and \(3\) is in-saturated, so both edges starting at \(4\) end at \(1\): a third doubled arc. The remaining case is two disjoint doubled arcs \(1\rightrightarrows 2\) and \(3\rightrightarrows 4\): then \(1\), \(3\) are out-saturated and \(2\), \(4\) are in-saturated, so the four remaining single arcs go from \(\{2,4\}\) to \(\{1,3\}\), each of \(2\) and \(4\) sending one arc to each of \(1\) and \(3\). This is (C3).

Exactly three doubled arcs. The doubled arcs form a loopless partial permutation of size three on four vertices: a \(3\)-cycle, which saturates its three vertices and leaves the fourth without edges; a path \(1\rightrightarrows 2\rightrightarrows 3\rightrightarrows 4\), which forces both edges starting at \(4\) to end at \(1\), a fourth doubled arc; a \(2\)-cycle \(1\rightrightarrows 2\rightrightarrows 1\) and a disjoint arc \(3\rightrightarrows 4\), which forces both edges starting at \(4\) to end at \(3\), a fourth doubled arc; or a \(2\)-cycle and an arc sharing a vertex, which gives a degree four. So three doubled arcs do not occur.

Four doubled arcs. The doubled arcs form a permutation without fixed points, doubled: a doubled \(4\)-cycle, which is (D4), or two doubled \(2\)-cycles, which is disconnected.
\end{proof}

\begin{figure}[ht]
\centering
\input{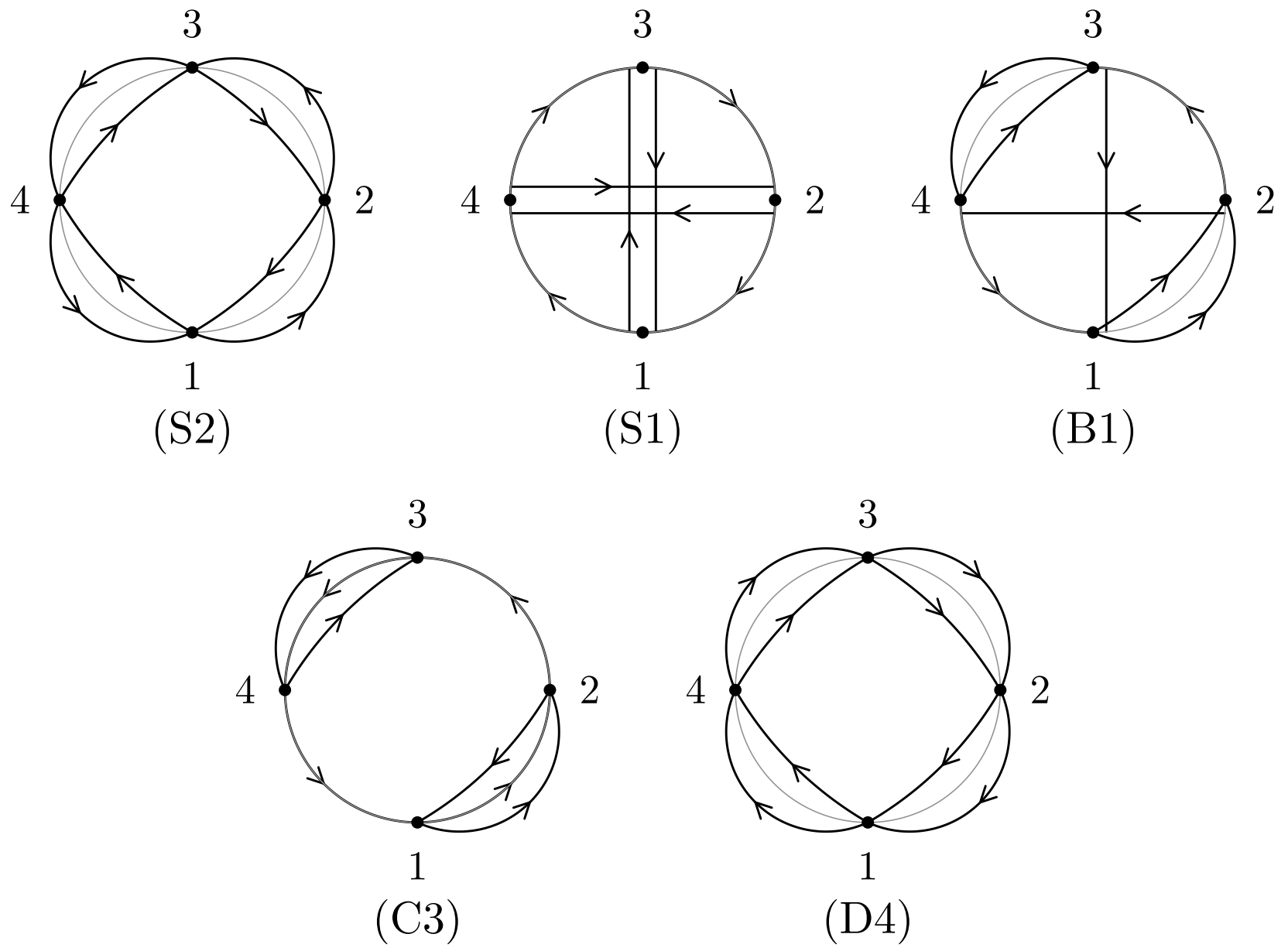}
\caption{The five connected loopless directed multigraphs on four classes with in-degree and out-degree two, drawn on the four class positions of the circle, the positions increasing counterclockwise from the bottom. The order of the classes on the circle is not part of the graph; the proof of Lemma~\ref{lem:descent} shows that, for \(k=6\), the graph is (D4) with its arcs descending through three consecutive positions and returning by a forward arc.}
\label{fig:classgraphs}
\end{figure}

Figure~\ref{fig:classgraphs} shows the five graphs.

\FloatBarrier
\section*{Acknowledgments}
I thank everyone who contributed useful discussions or took the time to read the paper and suggest edits. I am especially grateful to Zi-jian Song for helpful discussions of extensions beyond the three-interval case and of the rational-relation bound in Theorem~A.

\section*{Funding}
This research did not receive any specific grant from funding agencies in the public, commercial, or not-for-profit sectors.

\clearpage
\bibliographystyle{amsplainfull}
\makeatletter\def\@biblabel#1{\@defaultbiblabelstyle{#1}}\def\bibsetup{}\makeatother
\bibliography{refs}
\end{document}